\documentclass[12pt,uplatex]{amsart}
\usepackage{amsmath,amsfonts,amsthm,amscd,amssymb,mathrsfs,amssymb,bm,pb-diagram, here}
\usepackage{setspace}
\usepackage[all,cmtip]{xy}
\usepackage{tikz-cd}

\usepackage{multirow}
\usepackage{mathrsfs}
\usepackage{graphicx}
\usepackage[all]{xy}
\newtheorem{theorem}{Theorem}

\newtheorem{proposition}[theorem]{Proposition}
\newtheorem{lemma}[theorem]{Lemma}
\newtheorem{definition}[theorem]{Definition}

\newtheorem{convention}[theorem]{Convention}
\newtheorem{corollary}[theorem]{Corollary}

\newtheorem{remark}[theorem]{Remark}
\newtheorem{example}[theorem]{Example}
\newcommand{\aaa}{\alpha}
\newcommand{\bbb}{\beta}
\newcommand{\ccc}{\gamma}
\newcommand{\CCC}{\Gamma}
\newcommand{\ddd}{\delta}
\newcommand{\DDD}{\Delta}

\newcommand{\id}{{\rm{id}}}
\newcommand{\lmd}{\lambda}
\newcommand{\Lmd}{\Lambda}
\newcommand{\eee}{\epsilon}
\newcommand{\CP}{\mathbb{CP}}
\newcommand{\CC}{\mathbb{C}}

\newcommand{\RR}{\mathbb{R}}
\newcommand{\ZZ}{\mathbb{Z}}

\newcommand{\RP}{\mathbb{RP}}

\newcommand{\rank}{{\rm{rank}}}
\newcommand{\Sing}{{\rm{Sing}\,}}
\newcommand{\Supp}{{\rm{Supp}\,}}
\newcommand{\qdr}{\CP^1\times\CP^1}
\newcommand{\upi}{^{(i)}}
\newcommand{\upj}{^{(j)}}
\newcommand{\upk}{^{(k)}}
\newcommand{\upone}{^{(1)}}
\newcommand{\uptwo}{^{(2)}}

\newcommand{\ol}{\overline}
\newcommand{\lra}{\longrightarrow}
\newcommand{\lras}{\,\longrightarrow\,}

\newcommand{\Ra}{\Rightarrow}
\newcommand{\set}{\,|\,}
\newcommand{\proofend}{\hfill$\square$}

\newcommand{\inv}{^{-1}}

\newcommand{\Coker}{{\rm{Coker}}}

\newcommand{\Bs}{{\rm{Bs}}}
\newcommand{\Pic}{{\rm{Pic}}}

\newcommand{\ms}{\mathscr}
\newcommand{\minus}{\backslash}

\newcommand{\dpm}{_{\pm}}
\newcommand{\qandq}{\quad{\text{and}}\quad}

\newcommand{\bDh}{\bm D^{\rm h}}
\newcommand{\obDh}{\ol{\bm D}^{\rm h}}
\newcommand{\Dha}{D^{\rm h}_{\aaa}}
\newcommand{\bDha}{\ol D^{\rm h}_{\aaa}}
\newcommand{\oDha}{\ol D^{\rm h}_{\aaa}}
\newcommand{\Dh}{D^{\rm h}}
\newcommand{\oDh}{\ol D^{\rm h}}
\newcommand{\msDha}{\ms D^{\rm h}_{\aaa}}
\newcommand{\omsDha}{\ol{\ms D}^{\rm h}_{\aaa}}
\newcommand{\Ch}{\CCC^{\rm h}}
\newcommand{\oCh}{\ol\CCC^{\rm h}}
\newcommand{\Cha}{\CCC^{\rm h}_{\aaa}}
\newcommand{\ChaE}{\CCC^{\rm h}_{\aaa,E}}
\newcommand{\oChaE}{\ol\CCC^{\rm h}_{\aaa,E}}

\newcommand{\oCha}{\ol\CCC^{\rm h}_{\aaa}}

\newcommand{\ChaiE}{\big(\Cha\big)\upi_E}
\newcommand{\oChaiE}{\big(\oCha\big)\upi_E}
\newcommand{\ChaoneE}{\big(\Cha\big)\upone_E}
\newcommand{\oChaoneE}{\big(\ol\Gamma^{\rm h}_{\aaa}\big)\upone_E}
\newcommand{\ChakE}{\big(\Cha\big)\upk_E}
\newcommand{\oChakE}{\big(\oCha\big)\upk_E}

\newcommand{\ta}{\theta_{\aaa}}
\newcommand{\ota}{\ol\theta_{\aaa}}
\newcommand{\tif}{\tilde f}
\newcommand{\tiZ}{\tilde Z}
\newcommand{\Azero}{{\rm A}_0}
\newcommand{\Aone}{{\rm A}_1}
\newcommand{\Atwo}{{\rm A}_2}
\newcommand{\Athree}{{\rm A}_3}

\newcommand{\sst}{\scriptscriptstyle}

\numberwithin{equation}{section}
\numberwithin{theorem}{section}

\begin{document}
\bibliographystyle{alpha} 
\title[]
{Twistors, quartics, and del Pezzo fibrations}
\author{Nobuhiro Honda}
\address{Department of Mathematics, Tokyo Institute
of Technology}
\email{honda@math.titech.ac.jp}

\dedicatory{Dedicated to Professor Akira Fujiki on the occasion of
his 70th birthday}

\thanks{The author has been partially supported by JSPS KAKENHI Grant Number 24540061 (2012--2015) and 16H03932 (2016--2018).
Part of this work was done while I was a research member of 
MSRI in March 2016. I would like to thank the institute for hospitality.
}
\begin{abstract}
It has been known that
twistor spaces associated to self-dual metrics on 
compact 4-manifolds are source of interesting examples of 
non-projective Moishezon threefolds.
In this paper 
we investigate the structure of a variety of 
new Moishezon twistor spaces.
The anti-canonical line bundle on 
any twistor space admits a canonical half,
and we analyze the structure of 
twistor spaces by using the pluri-half-anti-canonical
map from the twistor spaces.

Specifically, each of the present 
twistor spaces is bimeromorphic to 
a double covering of a scroll of
planes over a rational normal curve,
and the branch divisor of the double cover is a cut of the scroll
by a quartic hypersurface.
In particular, 
the double covering has a pencil of Del Pezzo surfaces of degree two.
Correspondingly, the twistor spaces
have a pencil of rational surfaces
with big anti-canonical class.
The base locus of the last pencil is a 
cycle of rational curves, and 
it is an anti-canonical curve on smooth
members of the pencil.

These twistor spaces are naturally classified into
four types according to the type of singularities
of the branch divisor, or equivalently,
those of the Del Pezzo surfaces in the pencil.
We also show that the quartic hypersurface
 satisfies a strong constraint and as a result
the defining polynomial of the quartic hypersurface has to be of a specific form.

Together with our previous result in \cite{Hon_Cre1}, the present result completes
a classification of Moishezon twistor spaces
whose half-anti-canonical system is a pencil.
Twistor spaces whose half-anti-canonical system is larger than pencil
have been understood for a long time before.
In the opposite direction,
no example is known of a Moishezon twistor space whose half-anti-canonical system
is smaller than a pencil.

Twistor spaces which have a similar structure were studied in \cite{Hon_Inv} and \cite{Hon_Cre2}, and they are very special examples among 
the present twistor spaces.

\end{abstract}

\maketitle
\setcounter{tocdepth}{1}
\tableofcontents

\section{Introduction}
\subsection{Overview of the theory and the main result}
\label{ss:I}
In order to explain a motivation for the present investigation,
we first survey background materials and
what is known about Moishezon (and partially non-Moishezon)
twistor spaces.
The main theorem in this paper will be presented at the end
of this subsection.

Consider an oriented 4-dimensional Riemannian 
manifold.
By the Hodge star operator
we have a decomposition 
$\wedge^2 = \wedge^+\oplus\wedge^-$ 
of 2-forms into self-dual and anti-self-dual ones.
The Riemannian metric is called {\em self-dual} (resp.\,{\em anti-self-dual}\,)
if the anti-self-dual part (resp.\,self-dual part) of the Weyl curvature vanishes.
These conditions are invariant under conformal changes of the metric,
and self-duality and anti-self-duality make sense for a conformal class.
Also these two notions are mutually exchanged 
if we reverse the orientation of the 4-manifold.
A K\"ahler metric on a complex surface is anti-self-dual with 
respect to the complex orientation
iff it is of zero scalar curvature.
In particular, Ricci-flat K\"ahler metric on a complex surface is 
anti-self-dual,
 and for this reason
in the literature anti-self-duality  appears more often
than self-duality.
But in this paper we always consider self-dual conformal structures
because the 4-manifolds we mainly consider are the connected sum of 
copies of the complex projective planes, and they do not admit
an anti-self-dual metric
with respect to the standard orientation.

Over an oriented compact 4-manifold whose signature is non-negative,
self-dual conformal structures  attain the infimum
of the $L^2$-norm  of the Weyl curvature tensor,
determined by Gauss-Bonnet type theorem.
Thus self-dual conformal structures
may be thought of as a manifold version
of self-dual connections on a bundle
over an oriented Riemannian 4-manifold.

Another important aspect of self-dual conformal structure
is a connection with complex geometry.
Namely, the total space of the unit sphere sub-bundle of the bundle $\wedge^-$ 
of anti-self-dual 2-forms
on an oriented Riemannian 4-manifold
admits a natural almost complex structure,
and its integrability is equivalent to
self-duality of the conformal class of the metric
\cite{AHS}.
Thus for any self-dual structure,
there naturally associates a 3-dimensional
complex manifold having a sphere bundle 
structure over the 4-manifold.
This complex manifold is called the {\em twistor space}
of the self-dual structure,
and this is the main object in this paper.
One can recover a self-dual 
conformal structure from the twistor space,
and therefore these two notions are equivalent.
This is called the Penrose correspondence.
In this paper we are mainly interested in the properties
of the twistor spaces as complex manifolds or algebraic varieties.
But the self-dual conformal structures corresponding to 
the twistor spaces studied in this paper
also should be of interest from the viewpoint of differential geometry.

To be more precise,
let $M$ be a 4-manifold equipped with 
a self-dual conformal structure,
and $Z$ the associated twistor space.
Write $\pi:Z\to M$ for the projection.
Fibers of $\pi$ are complex submanifolds
of $Z$ and so isomorphic to $\CP^1$.
These are  called the {\em twistor lines}
of the twistor space.
The holomorphic normal bundle of a twistor line
in the twistor space
is non-trivial and isomorphic to
the direct sum $\ms O(1)\oplus\ms O(1)$.
This readily means that
anti-canonical bundle $-K_Z$ of $Z$ 
is of degree 4 over a twistor line.
In particular,  $\kappa(Z)=-\infty$ for the Kodaira dimension
when $Z$ is compact.
In a neighborhood of a twistor line,
there exists a 4-th square root of $-K_Z$,
and it exists globally  when the base 4-manifold $M$ admits
a spin structure.
For example, if $M$ is a K3 surface with a Ricci-flat K\"ahler metric, then $-K_Z$ admits 
a 4-th square root, and the associated
linear system induces a holomorphic submersion
$Z\to\CP^1$. 
This is so called the twistor family of K3 surfaces,
and significant for studying global structure of the moduli space of K3 surfaces.

Even if $M$ is not spin, from the construction of a twistor space,
there always exists a global and natural square root of $-K_Z$
as a holomorphic line bundle.
This is called the vertical line bundle or 
the {\em fundamental line bundle}
\cite{P92},
and is of degree two over a twistor line.
We denote it by $F\,(=K_Z^{-1/2})$ throughout this paper, and call the complete linear system $|F|$ the {\em  fundamental system}.
A member of 
the fundamental system is called a {\em fundamental divisor}.

Another important property of a twistor space
is that it is always  equipped with
an anti-holomorphic involution.
This is defined by the scalar multiplication by $(-1)$
on the bundle $\wedge^-$,
and 
is called the {\em real structure} of 
a twistor space.
We denote it by $\sigma$.
This acts on each twistor line as the anti-podal map
of $\CP^1\simeq S^2$.
The fundamental line bundle $F$
is real in the sense that $\sigma^*F\simeq \ol F$
holds, where $\ol F$ denotes the conjugate line bundle.
Hence the fundamental system $|F|$ also admits
a real structure,
and real members are parameterized by a real projective
space of the same dimension.
As we will explain later, a real fundamental divisor
plays an important role when we analyze the structure 
of twistor spaces.

Next we explain basic properties of twistor spaces
which follow from the presence of the real structure.
By Leray-Hirsch theorem on the topology of fiber bundles,
the second cohomology group of $Z$ with 
real coefficients is written as
\begin{align}\label{2ndcohom}
H^2(Z,\RR)&\simeq H^2(M,\RR)\oplus 
\langle F\rangle\\
&= H^2(M,\RR)\oplus 
\big\langle K_Z\inv\big\rangle,\notag
\end{align}
and evidently the real structure $\sigma$ acts, by pull-back,
as $\id$  and $-\id$
on the two factors
 $H^2(M,\RR)$ and $\langle F\rangle=\big\langle K_Z\inv\big\rangle$ 
 respectively.
This readily implies that,
for the algebraic dimension $a(Z)$ 
of a compact twistor space $Z$,
we have $a(Z)=\kappa(F\otimes L)$ 
for  some  real flat line bundle $L$,
where $\kappa(F\otimes L)$ is the Kodaira dimension of 
the line bundle $F\otimes L$.
If $Z$ is moreover simply connected
(which is equivalent to simply connectedness of $M$
by the homotopy exact sequence),
we always have $\Pic\,^0(Z)=1$ from gauge theory
as shown in \cite{LB92}.
Therefore the above flat line bundle $L$ is trivial, and  we have 
a simple relation $a(Z)=\kappa(F)$.
This in particular means that the twistor space of a K3 surface
always satisfies $a(Z)=1$.

As was shown in \cite{C91} by using a compactness of Chow variety
of Moishezon spaces,
simply connectedness follows  from
the assumption $a(Z)=3$.
These mean that a compact twistor space $Z$ is Moishezon 
if and only if $-K_Z$ is big.
So it might be possible to say that 
Moishezon twistor spaces are cousins of Fano threefolds.
In fact, the most basic two Fano threefolds
$\CP^3$ and the flag variety $\mathbb F$
of points and lines on $\CP^2$ are actually
twistor spaces; 
the former is the twistor space of
the 4-sphere $S^4$ equipped with the standard
 conformally flat structure, and the latter is 
the twistor space of 
$\CP^2$ equipped with the conformal class of  Fubini-Study metric.
But these are all twistor spaces which are really Fano.
More generally, Hitchin showed

\begin{proposition}
\cite{Hi81}\label{prop:Hi}
There exists no compact twistor space other than
$\CP^3$ and the flag variety $\mathbb F$ which 
admits a K\"ahler metric.
\end{proposition}

The proof of this result begins with noticing  
from \eqref{2ndcohom} and the real structure that 
the presence of a K\"ahler form
on a twistor space $Z$ implies that $F$ and $-K_Z$ are  positive line bundles,
and therefore $Z$ is really  Fano. 
Then 
analyzing the fundamental system $|F|$
by Riemann-Roch formula and
Kodaira vanishing theorem,
one finally reaches the conclusion.
Meanwhile
not only the real structure but also twistor lines
play an essential role.

The above proposition means that
besides the two spaces $\CP^3$ and $\mathbb F$
there is no compact twistor space
which is projective algebraic.
Thus, from the point of view of 
searching for twistor spaces having rich structure, 
a natural question would be whether there 
exist Moishezon twistor spaces 
which are not projective algebraic,
and nowadays a number of
such spaces are known to exist.
But before describing them,
we present the following classical result by Campana 
which means that the base 4-manifold of
a Moishezon twistor space is strongly constrained:

\begin{proposition} \cite{C91}\label{prop:C1}
If $Z$ is a twistor space on a compact 4-manifold $M$,
and if $Z$ is Moishezon,
 then $M$ is 
homeomorphic to $S^4$ or $n\CP^2$ for some 
integer $n$, where 
$n\CP^2$ denotes the connected sum of $n$ copies of 
$\CP^2$.
\end{proposition}

A key of the proof of this result is that under the Moishezon assumption,
$Z$ is covered by a {\em compact} family of rational curves
passing through a point.
Such a family forms a closed subvariety of the component of Chow variety
having twistor lines as points.
This implies simply connectedness of $Z$
(which is already mentioned),
as well as the vanishing $H^q(\ms O_Z)=0$ for 
any $q>0$.
Then by Riemann-Roch  and the formula
 $c_1c_2(Z)=12(\chi-\tau)(M)$
from \cite{Hi81},
we obtain $b_2^-(M)=0$ for the base 4-manifold $M$.
Therefore from the fundamental theorems
on closed 4-dimensional manifolds
by Freedman and Donaldson, 
$M$ has to be homeomorphic to $S^4$
or $n\CP^2$.
We also remark that the converse of Proposition \ref{prop:C1} does not hold, and twistor spaces on $n\CP^2$
are not necessarily Moishezon \cite{P92}.
This is an interesting aspect of twistor spaces
on $n\CP^2$.

From these results,
in the following, we are mainly concerned with 
twistor spaces on $n\CP^2$.
For simplicity of presentation we use the convention $0\,\CP^2=S^4$.
Let $Z$ be a twistor space on $n\CP^2$ and
$F$ the fundamental line bundle on $Z$.
The Riemann-Roch formula for the 
pluri-half-anticanonical bundle $lF$, $l\in\ZZ$,  is in effect 
computed in \cite{Hi81} and
one has
\begin{align}\label{RR1}
\chi(lF) =
 \frac13(4-n)l^3 + (4-n)l^2 + \frac13(11-2n)l
+ 1.
\end{align}
Here, as usual,
the most important in this polynomial is the sign of the 
coefficient of the top term $l^3$, which is 
positive, zero, and negative 
according as $n<4$, $n=4$ and $n>4$ respectively,
and it is natural to distinguish these three cases.
This is particularly evident when 
the corresponding self-dual conformal structure 
is of positive scalar curvature.
In the following for simplicity we
call this as self-dual structure
of {\em positive type}.
Then 
under this positivity assumption,
the Hitchin's vanishing theorem holds \cite{Hi80}:
\begin{align}\label{Hv}
H^2(lF) = 0\quad{\text{ for any $l\ge0$}}.
\end{align}
Therefore, we obtain from 
\eqref{RR1} that
if $n<4$,  $h^0(lF)=\dim H^0(lF)$ grows in cubic order,
and so we have $\kappa(F) = 3$.
Hence under the positivity assumption
the twistor space is always Moishezon if $n<4$.
When $n=0$ and $n=1$, 
the space has to be the standard ones:

\begin{itemize}
\item
Over $S^4$, any twistor space is isomorphic to $\CP^3$. This is 
a consequence of the classical result by Kuiper
on the uniqueness for conformally flat structures on the sphere.
\item
Over $\CP^2$, any twistor space
 is isomorphic to the flag variety $\mathbb F$
 if the self-dual structure is of positive type.
This was shown by Poon \cite{P86} by investigating the fundamental system $|F|$ on the twistor space.
\end{itemize}

The case $n=2$ is much more interesting:

\begin{proposition}\label{prop:Poon}
(Poon \cite{P86})
If $Z$ is a twistor space on $2\CP^2$ whose
self-dual structure is of positive type,  $\dim |F|=5$ and $\Bs\,|F|=\emptyset$ always hold. 
If $\Phi:Z\to\CP^5$ denotes the 
morphism induced by the fundamental system $|F|$, 
the image $\Phi(Z)$ is an intersection of 
 two quadrics, which has exactly four ordinary double points.
Moreover, the morphism $\Phi:Z\to\Phi(Z)$ can be identified with small resolution of 
these double points.
\end{proposition}

We mention that these twistor spaces constitute a one-dimensional moduli space.
By Proposition \ref{prop:Hi}, these twistor spaces do not admit a K\"ahler metric,
and they were the first examples of Moishezon twistor spaces
which are not projective algebraic.
A direct reason for non-projectivity is that, from the adjunction formula
we have
$F.\,\CCC=0$ for each exceptional curve $\CCC$ of the small resolution,
which cannot happen if $Z$ admits a K\"ahler metric
since if so $F$ has to be a positive line bundle as above.
We also remark that 
smooth complete intersections of two
 quadrics in $\CP^5$ are Fano threefolds.
Thus the twistor spaces on $2\CP^2$
can be obtained from such a Fano threefold
by first taking a suitably mild degeneration and
then taking small resolutions for the singularities.

Twistor spaces on $3\CP^2$ were investigated 
by Poon and Kreussler-Kurke.
They are not homogeneously described as in the case 
of $2\CP^2$, but 
can also be well understood
by the fundamental system $|F|$:

\begin{proposition}\label{prop:3}
\cite{KK92,P92}
If $Z$ is a twistor space on $3\CP^2$
whose self-dual structure 
is of positive type,
we always have $\dim|F|=3$.
Let $\Phi:Z\to\CP^3$ be the meromorphic map
induced by the fundamental system $|F|$.
Then we have:
\begin{itemize}
\item[\em (i)]
If $\Bs\,|F|\neq\emptyset$, 
the image $\Phi(Z)\subset\CP^3$ 
is a smooth quadric, and $Z$  is bimeromorphic to the total 
space of a conic bundle over the quadric.
\item[\em (ii)]
If $\Bs\,|F|=\emptyset$, 
the morphism $\Phi:Z\to\CP^3$ is surjective 
and generically two-to-one.
Moreover the branch divisor of $\Phi$ is 
a quartic surface defined by an equation of the form
\begin{align}\label{quartic1}
h_1h_2h_3h_4=Q^2,
\end{align}
where $h_i$-s are linear and $Q$ is quadratic.
\end{itemize} 
\end{proposition}

The twistor spaces in (i) are known as 
{\em LeBrun spaces},
and they exist on $n\CP^2$ for arbitrary $n$,
still having the properties $\dim |F|=3$ and $\Bs\,|F|\neq\emptyset$;
they are bimeromorphic to conic bundles over
the quadric which is the image of the meromorphic map associated to the fundamental system $|F|$.
The corresponding self-dual metrics on $n\CP^2$ are known as
LeBrun metrics  \cite{LB91}. They admit a non-trivial $S^1$-action, 
and  are explicitly constructed by hyperbolic ansatz.
Lifting the $S^1$-action and taking a complexification,
 LeBrun twistor spaces
admit a holomorphic $\CC^*$-action,
and the meromorphic map to the quadric in $\CP^3$
can be regarded as a quotient map under the $\mathbb C^*$-action.
Dependency on the number $n$ appears in
the discriminant locus of the conic bundle over the quadric.

On the other hand, for the twistor spaces in the item (ii)
of Proposition \ref{prop:3} are called of
\!{\em double solid type}.
For later use we briefly describe how the equation
\eqref{quartic1}
of the branch quartic is obtained, admitting
the fact that 
the morphism $\Phi$ is generically two-to-one and 
the branch divisor is a quartic.
First it 
is shown that the fundamental system 
has exactly four real reducible members.
Each of the corresponding hyperplane in $\CP^3$ has to be tangent
to the branch divisor along a conic
on the hyperplane.
(Such a conic is classically called a trope.)
Next it is shown that 
there exists a quadric in $\CP^3$ which contains all these four conics.
From these it follows that the branch divisor has to be 
of the form \eqref{quartic1} if we let the polynomials $h_i$ and $Q$
to be defining equations of the four hyperplanes and the quadric
respectively.
We also remark that the double covering
of $\CP^3$ branched along a smooth quartic
is a Fano threefold, and again it cannot be a twistor space. 
The quartic surface \eqref{quartic1} always has isolated singularities,
and the twistor space is obtained from 
the genuine double cover by taking small resolutions
for the singularities of the double cover.
Each of the exceptional curve $\CCC$ of the resolution
again satisfies $F.\,\CCC=0$, and 
again these are  obstructions for the twistor space
to admit a K\"ahler metric.

From these results, if $n<4$, the structure of twistor spaces 
on $n\CP^2$ is well understood by 
the fundamental system $|F|$.
If $n=4$, the first two coefficients
in the Riemann-Roch formula \eqref{RR1}
for $lF$ vanish, and 
by Hitchin's vanishing theorem \eqref{Hv}, we obtain
\begin{align}\label{RR2}
h^0(lF) \ge
\chi(lF) =
l + 1.
\end{align}
This means $\kappa(F)\ge 1$ for the line bundle $F$, and therefore
we always have $a(Z)\ge 1$.
As a matter of fact,  generic twistor spaces
on $4\CP^2$ satisfy $a(Z)=1$.
Moreover, for these spaces,
algebraic reduction is always induced by
the fundamental system, and its general members are rational surfaces.
Many examples are known of twistor spaces on $4\CP^2$
that satisfy $a(Z)=2$ \cite{CK99, HI00}. 
For Moishezon ones, we have the following
classification result which describes 
structure 
in terms of the anti-canonical map from the twistor spaces.

\begin{proposition}\label{prop:4}
\cite{Hon_JAG2}
Let $Z\to 4\CP^2$ be a Moishezon twistor space, and let $\Phi:Z\to\CP^N$
be the anti-canonical map.
Then one of the following three situations occurs.

\begin{itemize}
\item[\em (i)]
The map $\Phi$ is bimeromorphic over $\Phi(Z)$.
\item[\em (ii)] 
The map $\Phi$ is of degree two over $\Phi(Z)$.
In this case, $h^0(-K_Z)=5$ and we have 
$\Phi(Z)=p\inv(\Lmd)$ for the image, where $p:\CP^4\to \CP^2$
is a linear (meromorphic) projection and $\Lmd$ is 
a plane conic.
Further, the branch divisor of the map $\Phi$ is
a cut of the image $p\inv(\Lmd)$ by a quartic hypersurface in $\CP^4$,
whose defining equation is of the form
\begin{align}\label{quartic2}
h_1h_2h_3h_4=Q^2,
\end{align}
where $h_i$-s are linear and $Q$ is quadratic.
\item[\em (iii)] The image $\Phi(Z)$
is 2-dimensional.
\end{itemize}
\end{proposition}

Concerning the variety $p\inv(\Lmd)$ in (ii),
the indeterminacy locus of the linear projection $p:\CP^4\to\CP^2$
is a line, and fibers of $p$ are planes which contain the line.
Hence a conic in $\CP^2$ determines a 1-dimensional non-linear family of 
such planes, and 
 the variety $p\inv(\Lmd)$ is a union of such planes.
We call this variety a {\em scroll of planes over a conic},
the fixed line the {\em ridge} of the scroll,
and fibers of $p$ planes of the scroll.
The scroll has ordinary double points along the ridge.
From the equation \eqref{quartic2} of the quartic hypersurface,
the twistor spaces in (ii) may be thought of as an analogue of 
the ones over $3\CP^2$
described in the second item in Proposition \ref{prop:3}, and
we again call them as twistor spaces of {\em double solid type}.

Thus on $4\CP^2$, the structure of Moishezon twistor spaces
  is well understood by the anti-canonical system.
We mention that the result in \cite{Hon_JAG2} is much more precise
than the above form,
but for brevity we omit the detail here.
In relation with the main result in this paper,
we  mention
that the twistor spaces of double solid type
on $4\CP^2$ can be classified
into four types according to type of singularities
of  the branch divisor.
We also remark that, similarly to the case of $3\CP^2$,
the specific form of the equation \eqref{quartic2} is deduced
from a presence of real {\em reducible} members of the anti-canonical system
on the twistor space.

If $n>4$, from the Riemann-Roch formula \eqref{RR1}, we obtain $\chi(lF)<0$ for any $l>0$,
so  we do not obtain any effective lower bound 
for the algebraic dimension.
Actually generic twistor spaces on $n\CP^2$,
$n>4$, satisfy $a(Z)=0$ (\cite{P92}).
It is easy to find examples 
which satisfy $\dim |F|=a(Z)=1$.
The algebraic reduction of these examples is 
induced by the pencil $|F|$, and its general fibers
are rational surfaces
as in the case $n=4$.
But there also exist examples of twistor spaces 
on $n\CP^2$, $n>4$, satisfying $a(Z)=1$,
 whose algebraic 
reduction is induced by not the fundamental system $|F|$ but the anti-canonical system $|2F|$ (\cite{Hon_CM}).
These twistor spaces have a property
that 
 generic fibers  of the algebraic reduction are elliptic ruled surfaces.
In this direction, for twistor spaces of algebraic dimension one,
it is not completely understood as to 
which kind of surfaces really appear as
a general fiber of the algebraic reduction
(even in the case of $n\CP^2$).
Also we should  mention that no example of a twistor
space on $n\CP^2$ is known of that satisfies $a(Z)=2$ when $n>4$.
In \cite{HK17} it was shown that 
such a twistor space has to satisfy 
$h^0(F)\le 1$.
Thus basic questions on non-Moishezon twistor spaces
on $n\CP^2$ are still open when $n>4$.

Next we discuss Moishezon twistor spaces
on $n\CP^2$ when $n>4$, which 
are the main object in this paper.
LeBrun spaces described right after
Proposition \ref{prop:3} were the first
examples of such spaces.
Recall that they satisfy the property
$\dim |F|=3$, and  they are characterized 
by this property \cite{Kr99}.
There exists no twistor space
on $n\CP^2$, $n>4$, which satisfies
$\dim |F|>3$.
Other interesting examples of Moishezon twistor spaces
on $n\CP^2$, $n\ge 4$, were obtained by 
Campana and Kreussler \cite{CK98}.
Their twistor spaces satisfy
$\dim |F|=2$, and by
the associated meromorphic map $Z\to\CP^2$ 
they
are bimeromorphic to  conic bundles
over $\CP^2$.
Dependency on $n$ appears in the discriminant
curves as in the case of LeBrun spaces.
If $n>4$, these twistor spaces are characterized
by the property $\dim |F|=2$ \cite{Kr99}.

Also there are a series of distinguished families of
Moishezon twistor spaces:
they are the twistor spaces
associated to self-dual metrics on $n\CP^2$
constructed by Joyce \cite{J95}.
These metrics are toric in the sense that  they are invariant under an effective $T^2$-action on $n\CP^2$,
and some of them coincide with LeBrun metrics
whose $S^1$-action extends to an effective $T^2$-action
\cite{HV13}.
Except these, 
the twistor spaces of Joyce metrics
satisfy $\dim|F|=1$, and
the meromorphic map $Z\to\CP^1$ induced by
the pencil $|F|$
can be regarded as a quotient map
with respect to the $(\CC^*\times\CC^*)$-action
associated to the $T^2$-action for the metric.
The twistor spaces of Joyce metrics were
intensively studied by Fujiki \cite{F00},
and it was shown that 
the self-dual metrics and the twistor spaces are characterized 
by the presence of an effective $T^2$-action.

The twistor spaces of Joyce metrics
were investigated in \cite{Hon_Cre1} from a
different aspect.
Namely,  for a suitably chosen subgroup
$S^1\subset T^2$ (or equivalently,
a subgroup $\CC^*\subset
\CC^*\times\CC^*$),
a quotient surface of the twistor space
under the $\CC^*$-action was obtained
as the image of the meromorphic map 
associated to (a sub-system of) a pluri-system
$|mF|$,
where $m$ is the number which can be readily determined
from the $T^2$-action and the $S^1$-subgroup.
The quotient surface is called the minitwistor
space, and it is equivalent to an Einstein-Weyl
structure on a 3-manifold.
Defining equations of the minitwistor spaces were explicitly given, and they are birational
to conic bundles over $\CP^1$.
By the quotient maps from the twistor spaces
to the minitwistor spaces,
the twistor spaces of Joyce metrics are bimeromorphic
to the total spaces of conic bundles over the minitwistor spaces.
This description can be regarded as
a generalization of that for
LeBrun twistor spaces described above,
in the sense that the quadric in $\CP^3$,
which is the image of the meromorphic map associated to $|F|$,
 is the simplest example
of a minitwistor space.
Further, it was shown that there are many twistor spaces
which are bimeromorphic to a conic bundle over the above minitwistor spaces,
but which do not admit an effective $T^2$-action.
These `new' twistor spaces are characterized by a presence
of a real irreducible fundamental divisor
which admits an effective $\CC^*$-action that has a pair of
fixed curves,
and conjecturally all the twistor spaces admit a $\CC^*$-action.
From the description of the spaces,
we call these twistor spaces as {\em generalized LeBrun spaces}.
These include the twistor spaces of
LeBrun metrics and Joyce metrics as special cases,
and we have `new' twistor spaces with $\CC^*$-action already in the case $n=3$.
Except the original LeBrun spaces,
generalized LeBrun spaces satisfy $\dim |F|=1$.
We mention that 
each of the minitwistor spaces itself admits
a non-trivial deformation of complex structure
in general,
and in addition we again have a freedom for moving discriminant
divisor for the conic bundle structure.
Moreover, the number of minitwistor spaces
is infinite even up to deformation equivalence.

Other families of Moishezon twistor spaces
on $n\CP^2$, $n\ge 4$, which can be described in a different manner
from the generalized LeBrun spaces,
were found in \cite{Hon_Inv, Hon_Cre2}.
These twistor spaces enjoy the following properties:

\begin{itemize}
\item
$\dim|F|=1$, $\dim|(n-2)F|=n$, and 
the meromorphic map $\Phi:Z\to \CP^n$
associated to the system $|(n-2)F|$ is 
two-to-one over the image.
Moreover, we have $\Phi(Z)=p\inv(\Lmd)$ for the image,
where $p$ is a linear projection 
from $\CP^n$ to $\CP^{n-2}$, and $\Lmd\subset\CP^{n-2}$ is 
a rational normal curve of degree $(n-2)$.
\item
The branch divisor of the degree-two map $\Phi:Z\to\Phi(Z)$
is a cut of the image $p\inv(\Lmd)$ by a quartic hypersurface in $\CP^n$.
The quartic hypersurface 
is defined  by the equation of the form
\begin{align}\label{quartic3}
h_1h_2h_3h_4=Q^2,
\end{align}
where $h_i$ are linear and $Q$ is quadratic.
\end{itemize}

As in the case $n=4$,
the indeterminacy locus of the linear projection $p:\CP^n\to\CP^{n-2}$ 
is a line, and fibers of $p$ are exactly  planes
which contain the line.
The curve $\Lmd\subset\CP^n$ determines
a 1-dimensional non-linear family of such planes,
and the variety $p\inv(\Lmd)$ is the union of the planes in this family.
We again call it
as {\em the scroll of planes} over a rational normal
curve, and the line the {\em ridge} of the scroll.
While the scroll has  cyclic quotient singularities 
of order $(n-2)$ along the ridge,
it is a very simple rational variety, and is uniquely determined from the number $n$. 
Evidently the above twistor spaces are generalization
of the ones over $3\CP^2$ and $4\CP^2$
which are described by a quartic hypersurface.
These twistor spaces are characterized by 
a structure of a real fundamental divisor on them,
but it was not evident to what extent these spaces
occupy in the space of Moishezon twistor spaces.

With these backgrounds, 
the main result of this paper is 
a classification and description of Moishezon twistor spaces on 
$n\CP^2$, $n\ge 4$, which satisfy $\dim |F|=1$:

\begin{theorem}\label{thm:main1}
Let $n\ge 4$ and $Z$ be a Moishezon twistor space on $n\CP^2$.
Suppose $\dim |F|=1$ and $Z$ is not a generalized LeBrun space.
Then there exists $m\ge 2$ such that 
the pluri-system $|mF|$ includes an $(m+2)$-dimensional
sub-system
whose meromorphic map
$\Phi:Z\to\CP^{m+2}$ 
satisfies the following properties.
\begin{itemize}
\item
The image $\Phi(Z)$ is a scroll of planes
over a rational normal curve in $\CP^{m}$.
\item
The meromorphic map $\Phi:Z\to \Phi(Z)$ is
two-to-one, and the branch divisor 
is a cut of the scroll by a quartic hypersurface in $\CP^{m+2}$.
\item
The quartic hypersurface is defined by an equation of the form
\begin{align}\label{quartic4}
h_1h_2h_3h_4=Q^2,
\end{align}
where $h_i$ are linear and $Q$ is quadratic.
\end{itemize}
Further, the integer $m$ necessarily satisfies $m\ge n-2$.
\end{theorem}

We again call these twistor spaces
as {\em \!double solid type}.
The theorem implies that 
if $n\ge 4$, a Moishezon twistor space on $n\CP^2$ which satisfies $\dim |F|=1$
is either a generalized LeBrun  space
or 
of double solid type.
Thus the conclusion is unexpectedly simple, 
and the classification result over $3\CP^2$
stated in Proposition \ref{prop:3} carries over to the case $n\ge 4$
under the assumption $\dim |F|=1$.
Of course, Theorem \ref{thm:main1} is compatible with 
Proposition \ref{prop:4} concerning 
Moishezon twistor spaces on 
$4\CP^2$. 
Actually, twistor spaces belonging to (i) in the proposition 
are either of double solid type or
generalized LeBrun spaces.
For the former spaces we really need a proper sub-system
of $|mF|$ to get the double covering map to the scroll,
because otherwise the map becomes bimeromorphic to the image.
This means that in general we really need to take a sub-system
of $|mF|$ instead of the system $|mF|$ itself to obtain 
the double covering structure as in Theorem \ref{thm:main1}.
Also,  twistor spaces belonging to (iii) in the proposition 
are always generalized LeBrun spaces if 
$\dim |F|=1$ as supposed in Theorem \ref{thm:main1}.
We remark that in the case $n=4$, 
the integer $m$ in Theorem \ref{thm:main1} can be taken as two by Proposition \ref{prop:4}.

Needless to say, the twistor spaces studied in \cite{Hon_Inv, Hon_Cre2} mentioned above
are examples of those in Theorem \ref{thm:main1},
and for them we have $m=n-2$ and the system $|mF|$ itself induces
the double covering map to the scroll of planes in $\CP^{m+2}$.
From these examples, it is tempting to expect that
we can always take the lower bound $(n-2)$ as the number $m$ 
in Theorem \ref{thm:main1}.
We show that this is not correct, and there actually exist
examples of twistor spaces for which the strict inequality $m>n-2$ holds.
Note that, from upper semi-continuity
of dimension of sections of line bundles,
the twistor spaces are more general if
the number $m$ is larger.
This means that the twistor spaces
studied in \cite{Hon_Inv, Hon_Cre2}
are special among all Moishezon twistor spaces of double solid type, despite the title of \cite{Hon_Cre2}!
Also from the aspect of the method of 
analysis, the technique
we employed in \cite{Hon_Inv, Hon_Cre2}
rely too heavily on specific properties of these examples
and it looks difficult to apply it
to general twistor spaces of double solid type.

As a consequence of Theorem \ref{thm:main1}, 
the present understanding of Moishezon twistor spaces
on $n\CP^2$, $n\ge 4$, may be summarized 
as in Table \ref{table0}.

\begin{table}[h]
\begin{tabular}{|c||c||c||c||c||c}
$h^0(F)>4$ & $h^0(F)=4$ & 
{$h^0(F) = 3$} &
$h^0(F)=2$ & $h^0(F)<2$\\
\hline
\multirow{2}{*}{$\not\exists$}
 & 
 \multirow{2}{*}{LeBrun} &
 \multirow{2}{*}{Campana-Kreussler} &
 generalized LeBrun, &
 \multirow{2}{*}{$\exists?$}\\
 &&&double solid type &
\end{tabular} 
\medskip
\caption{Classification of Moishezon twistor spaces on $n\CP^2$ when $n\ge 4$}
\label{table0}
\end{table}

Thus
what is still missing is to understand Moishezon twistor spaces
which satisfy $h^0(F)<2$.
Up to now  such a twistor space seems not known.
If such a twistor space does not exist,
the present result completes a classification 
of Moishezon twistor spaces without any assumption.
But there seems no enough reason to expect such non-existence.

\subsection{Outline of the proof
of Theorem \ref{thm:main1}}
In the rest of this section
we briefly describe basic ideas for analyzing the present twistor spaces.
By a fundamental result of Pedersen-Poon
\cite{PP94},
if a twistor space $Z$ on $n\CP^2$ satisfies 
$\dim |F|\ge 1$, general members of the pencil $|F|$ are smooth rational surfaces
which satisfy $K^2 = 8-2n$.
Such members, especially real ones,
are always denoted by $S$.
By a result of Kreussler \cite{Kr98},
when $\dim |F|=1$,
the base locus of the pencil $|F|$ is
 an anti-canonical curve on the surface
 $S$,
 and it is 
 either a smooth elliptic curve or a cycle of smooth rational curves.
If the twistor space $Z$ is moreover Moishezon,
only the latter can happen.
The cycle, which is the base curve of
the pencil $|F|$, will be always denoted by $C$
throughout this paper.
This is real with respect to the real structure on $Z$
and consists of even number of components.
If $2k$ denotes the number of components,
the pencil $|F|$ has exactly $k$ reducible members
by Kreussler \cite{Kr98}.
Let $\hat Z\to Z$ be the blow-up at the cycle $C$.
This eliminates the base locus of the pencil $|F|$
and so we obtain 
a morphism $\hat Z\to\CP^1$.
This has exactly $k$ reducible fibers,
and on each of them 
the space $\hat Z$ has two ordinary double points.
We take small resolutions
for them which preserve the real structure.
The resulting smooth space will  always be denoted
by  $\tilde Z$, and the composition morphism $\tilde Z\to\hat Z\to\CP^1$
is denoted by $\tilde f$.
We shall make full use of this fibering structure $$\tif:\tiZ\to\CP^1$$ to 
analyze the structure of the twistor space $Z$.
This fibration also has exactly $k$ reducible fibers.
Throughout this paper we denote them by 
$$
S\upone, S\uptwo,\dots,S\upk,
$$
and the corresponding points on $\CP^1$ are
denoted by $\lmd\upone,\lmd\uptwo,\dots,\lmd\upk$ respectively.
These are real points since each reducible fiber is real.

From the assumption that the twistor space
$Z$ is Moishezon,
the anti-canonical bundle of any smooth fiber
$S$ of $\tif$ is big.
Since $K_S^2=8-2n$ as above and this is non-positive
if $n\ge 4$,
from the bigness, 
the pluri-anti-canonical system $\big|mK_S\inv\big|$ has a fixed component for any $m>0$.
As the first step in proving Theorem \ref{thm:main1},
we show from the bigness that 
the surface $S$ admits a non-trivial $\CC^*$-action,
or for some integer $m>0$
a movable part of the system $\big|mK_S\inv\big|$ induces a generically two-to-one
morphism $\phi:S\to\CP^2$
whose branch divisor is a quartic (Proposition \ref{prop:clsf1}).
If the former is the case, 
a result in \cite{Hon_Cre1} means that
the twistor space $Z$ has to be
a generalized LeBrun space.
Thus we need to show that 
if the latter is the case,
$Z$ has to be of double solid type.
We note that 
the above two possibilities for the surface $S$ 
are {\em not} alternative but
$S$ can satisfy both properties.
This implies that there is overlap
for generalized LeBrun spaces
and twistor spaces of  double solid type.
This fact will play some role in proving the existence
of the present twistor spaces (considered in Theorem \ref{thm:main1}).

The movable part of the pluri-system
$\big|mK_S\inv\big|$ which induces the degree-two morphism
$\phi:S\to\CP^2$ has a distinguished member
which is characterized by the property
that the support is exactly the cycle $C$.
This member is always denoted by $D$.
Since $\phi\inv(\bm l)=D$ for some line $\bm l\subset\CP^2$,
the image $\phi(C)$  is the line $\bm l$.
This line will be common for all smooth members
of the pencil $|F|$,
and will be the ridge of the scroll.
Among components of the cycle $C$,
there exist precisely two components
which are mapped isomorphically to this line by
the morphism $\phi:S\to\CP^2$.
These components will be called 
the {\em line components} of the cycle $C$.
In the cycle $C$
these are placed at the opposite sides.
The remaining components in the cycle form two chains which are
mutually conjugate.
These chains are contracted to points by the morphism $\phi:S\to\CP^2$.
We write these points by $\bm q$ and $\ol {\bm q}$.
These are distinct points on the line $\bm l$.
We show that the branch quartic of
the degree-two morphism $\phi:S\to\CP^2$ satisfies one of the 
following four properties.
\begin{itemize}
\item[$\Azero$.]
It is tangent to the line $\bm l$
at the points $\bm q$ and $\ol{\bm q}$.
\item[$\Aone$.]
It has $\Aone$-singularities at $\bm q$ and $\bm {\ol q}$.
\item[$\Atwo$.]
It has $\Atwo$-singularities at $\bm q$ and $\bm {\ol q}$.
\item[$\Athree$.]
It has $\Athree$-singularities at $\bm q$ and $\bm {\ol q}$.
\end{itemize}

\begin{figure}
\includegraphics{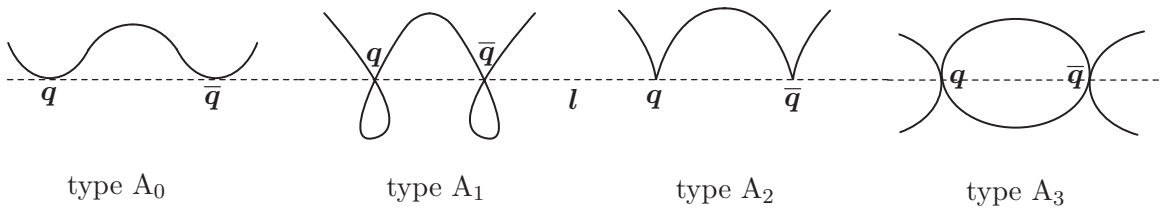}
\caption{
branch quartic and the ridge $\bm l$}
\label{fig:bq}
\end{figure}

\noindent
See Figure \ref{fig:bq} for these.
Thus  possible structure of the surface
$S$ can be classified into four types 
by the singularities of the branch quartic.
Since the type is independent of a choice
of a smooth fundamental divisor in the pencil $|F|$,
the type makes sense also for 
the twistor space (Definition \ref{def:type}),
and it reflects singularities of the 
branch divisor on the scroll at the two points
$\bm q$ and $\bm {\ol q}$
on the ridge $\bm l$.

Thus structure of smooth fibers of the fibering $\tif:\tiZ\to\CP^1$ has become evident,
and the next problem is to relate or extend
the morphism $\phi:S\to\CP^2$ from smooth fibers
to a map from the space $\tiZ$ to a scroll of planes.
For this, we write
the anti-canonical cycle $C$ in the form
$$
C =\big(C_1+ \dots + C_k\big) + 
\big(\ol C_1 + \dots + \ol C_k\big)
$$
respecting the real structure.
Then since the support of the above effective divisor $D$ is equal to $C$ as above,
we can write 
\begin{align}\label{D000}
D =
\big(d_1C_1+ \dots + d_kC_k\big) + 
\big(d_1\ol C_1 + \dots + d_k\ol C_k\big)
\end{align}
for some positive integers $d_1,\dots,d_k$.
We will always let $C_1$ and $\ol C_1$
to be the line components.
Since all other components are contracted to 
points by $\phi$ as above, we have $D.\,C_i=D.\,\ol C_i=0$
for any $i>1$, while $D.\,C_1= D.\,\ol C_1 = 1$.
Hence the divisor $D$ is not ample,
but {\em nef and big}.
These numerical properties will be used a number of times
throughout this paper.
For the purpose of obtaining an effective divisor on $\tilde Z$
whose restriction to a smooth fiber of $\tif$ is equal to $D$,
letting $E_i$ and $\ol E_i$ be
the exceptional divisors of the components
$C_i$ and $\ol C_i$ respectively,
we define a divisor $\bm D$ on $\tiZ$ by
$$
\bm D =
\big(d_1E_1+ \dots + d_kE_k) + 
\big(d_1\ol E_1 + \dots + d_k\ol E_k\big).
$$
This is an extension of the divisor $D$ on a fiber $S$ to the whole of $\tilde Z$.
But the linear system $|\bm D|$
consists of a single member $\bm D$ basically because 
all components of the divisor $\bm D$ are exceptional divisors 
of blowing up. Therefore the restriction map
$H^0\big(\tilde Z,\ms O_{\tilde Z}(\bm D)\big)
\to H^0\big(S,\ms O_{S}(D)\big)$ cannot be surjective.

So instead, 
writing $\ms D:=\ms O_{\tilde Z}(\bm D)$
for the associated invertible sheaf,
to any integer $l$, we  consider the invertible sheaf
\begin{align}\label{}
\ms D(l):=\ms D \otimes \tilde f^*\ms O_{\CP^1}(l)
\end{align}
and compute the dimension $h^0\big(\tilde Z,\ms D(l)\big)$.
To this end, we take the direct image sheaf
$\tilde f_*\ms D(l)$ under the morphism
$\tif:\tiZ\to\CP^1$.
By the projection formula,
this is isomorphic to 
$(\tilde f_*\ms D)\otimes\ms O_{\CP^1}(l)$.
Moreover, since the morphism $\tilde f$ is proper and flat,
the direct image sheaf $\tilde f_*\ms D$ is 
torsion free \cite{Tel15}.
Since any torsion free sheaf over a smooth curve
is locally free, the sheaf $\tilde f_*\ms D$ is 
locally free, and hence, isomorphic to 
a direct sum of invertible sheaves (i.e.\,line bundles). 
The rank of the sheaf $\tilde f_*\ms D$ is  three 
since $h^0(S,D)=3$ as above for any smooth fiber $S$
of $\tilde f$.

Since $\tif$ is proper, we have $h^0\big(\tiZ,\ms D(l)\big)
= h^0\big(\CP^1,(\tif_*\ms D)\otimes\ms O(l)\big)$,
and to compute $h^0\big(\tiZ,\ms D(l)\big)$
it suffices to determine the locally free sheaf $\tif_*\ms D$.
For this purpose,
by restricting the line bundle
$\ms D$ on $\tiZ$ to the exceptional divisor
$$E:=(E_1+\dots + E_k)+ (\ol E_1+\dots+ \ol E_k)$$
of the cycle $C$,
we get the standard exact sequence
\begin{align}\label{ses:o3ye}
0 \lras \ms D(-E)
\lras \ms D\lras \ms D|_E\lras 0.
\end{align}
By taking the direct image under $\tif$, we obtain the long exact sequence
\begin{align}\label{ses:o3yf}
0 \lras \tilde f_*\ms D(-E)
\lras \tilde f_*\ms D\stackrel{\rho}\lras 
\tilde f_*\big(\ms D|_E\big)\lras
R^1\tilde f_*\ms D(-E)\lras \dots.
\end{align}
It is not difficult to see that $\tilde f_*\ms D(-E)\simeq\ms O_{\CP^1}$
(Proposition \ref{prop:d_im4}).
Further by applying Kodaira-Ramanujan vanishing theorem 
to the divisor $D$ (which is nef and big as above) on
any smooth fiber $S$ of $\tif$, we obtain
$H^1\big(S,\ms D(-E)|_S\big)=H^1(S,D+K)=0$.
This means that the sheaf
$R^1\tilde f_*\ms D(-E)$ is a torsion sheaf.
Furthermore, it is possible to show that 
the support of the torsion sheaf $R^1\tilde f_*\ms D(-E)$ is
included in the locus $\{\lmd\upone,\lmd\uptwo,\dots,\lmd\upk\}$
over which the reducible fibers of $\tif$ lie
(Proposition \ref{prop:d_im4}).
This implies that,
if 
\begin{align}\label{restred}
\rho\upi:\big(\tif_*\ms D\big)_{\lmd\upi}\lras
\big(\tif_*\ms D|_E\big)_{\lmd\upi}
\end{align}
denotes the restriction  of the homomorphism $\rho$
in the exact sequence \eqref{ses:o3yf} to the stalks over the point $\lmd\upi$,
we have an exact sequence
\begin{align}\label{der003}
0 \lras \ms O \lras 
\tilde f_*\ms D \stackrel{\rho}{\lras} 
\tif_*\big(\ms D|_E\big)
{\lras} 
\bigoplus _{i=1}^k \Coker \,\rho\upi
{\lras} 
0
\end{align}
over $\CP^1$.
This sequence will be called the derived sequence,
and  of fundamental significance in analyzing 
the structure of the present twistor spaces as we next discuss.

Over the locus $\CP^1\minus\{\lmd\upone,\dots,\lmd\upk\}$,
 the restriction $\tif|_E:E\to\CP^1$
is a fiber bundle whose fibers are isomorphic to the 
anti-canonical cycle $C$.
The fibers of $\tif|_E$ over the points $\lmd\upone,\dots,\lmd\upk$ are 
not isomorphic to $C$ because the exceptional curves
of the small resolutions $\tiZ\to\hat Z$ are inserted
at the two nodes of $\hat Z$
on each reducible fiber of the fibration $\hat Z\to\CP^1$.
We will denote these `singular fibers' of $\tif|_E$ by 
$C\upone,\dots,C\upk$ respectively.
In Section \ref{ss:sbc}, we investigate fixed components 
of the line bundle of the form $\ms D(l)|_E$,
by calculating the 
intersection number of $\ms D$ with components of the 
`singular' cycles $C\upi$.
For each $i$,
the fixed components of $\ms D(l)|_E$ consist of two chains in $C\upi$
which are mutually conjugate if they are not empty.
We call them the {\em stable base curves} of the line bundle $\ms D$.
By using them, we describe the space
$H^0\big(\ms D(l)|_E\big)$ in a concrete form for any $l$
(Proposition \ref{prop:rest3}).
From this we show 
\begin{align}\label{d_im001}
\tilde f_*(\ms D|_E)\simeq \ms O(-e)^{\oplus 2}
\end{align}
holds for an integer $e>0$ (Proposition \ref{prop:e2}).
The number $e$ can be concretely expressed 
in terms of the coefficients $d_1,\dots, d_k$
for $D$ in the expression \eqref{D000}.
From the expression we obtain an inequality $e\ge n-2$.

Next in Section \ref{ss:dim1}, by using the description of the space
$H^0\big(\ms D(l)|_E\big)$,
we show that there exists an integer $m\ge e$ such that
\begin{align}\label{d_im000}
\tilde f_*\ms D\simeq \ms O(-m)^{\oplus 2}
\oplus \ms O
\end{align}
holds (Proposition \ref{prop:d_im06}).
Then by using this we show that the image of
 the meromorphic map
associated to the linear system $\big|\ms D(m)\big|$
is a scroll of planes of degree $m$,
and also that the map is of degree two over the scroll
(Theorem \ref{thm:dc1}).
We also show that the map contracts the two components
$E_1$ and $\ol E_1$ to the ridge $\bm l$, and 
the remaining components $E_2\cup\dots\cup E_k$
and $\ol E_2\cup\dots\cup \ol E_k$ to distinct points
$\bm q$ and $\ol{\bm q}$ on $\bm l$
(Proposition \ref{prop:pim}).
Next by investigating the linear system $\big|\ms D(m+1)\big|$,
we show that the branch divisor of the degree two map from $\tiZ$
to the scroll of planes is a cut of the scroll by a quartic hypersurface
(Theorem \ref{thm:dc2}).
Since there is a natural inclusion 
$H^0\big(\ms D(m)\big)\subset H^0(mF)$, 
at this stage we obtain the conclusion of 
Theorem \ref{thm:main1} except the defining equation \eqref{quartic4} of the quartic hypersurface.

From \eqref{d_im001} and \eqref{d_im000}, the derived sequence \eqref{der003} becomes
\begin{align}\label{der004}
0 \lras \ms O \lras 
\ms O(-m)^{\oplus 2}\oplus\ms O\lras
\ms O(-e)^{\oplus 2}
{\lras} 
\bigoplus _{i=1}^k \Coker \,\rho\upi
{\lras} 
0.
\end{align}
From the exactness of this sequence, we obtain the relation
\begin{align}\label{me0}
m = e + \frac12\sum_{i=1}^k\dim \Coker \,\rho\upi.
\end{align}
The number $e$ is concretely known as above, and the next problem is
to determine the dimension of the vector space
$\Coker \,\rho\upi$ for each index $i=1,\dots,k$.
In this respect, by using the description of the space
$H^0\big(\ms D(l)|_E\big)$, we find that 
the linear system on the exceptional divisor $E$ associated to
the image of the restriction homomorphism
$H^0\big(\ms D(m)\big)\to H^0\big(\ms D(m)|_E\big)$,
has fixed components of the form
\begin{align}\label{cbci}
\sum_{i=1}^k\mu\upi C\upi,\quad\mu\upi\ge 0,
\end{align}
and that 
we have 
\begin{align}\label{dimc}
\dim \Coker \,\rho\upi = 2\mu\upi
\end{align}
 for each index $i=1,\dots,k$
(Propositions \ref{prop:cbc1} and \ref{prop:mvo}). 
Because the support of \eqref{cbci} consists
of whole of the cycles $C\upi$, 
we call the fixed components \eqref{cbci} 
the {\em cyclic base curves} (Definition \ref{def:cbc}).
We note that in contrast with the stable base curves,
they are not determined from the line bundle
$\ms D|_E$ itself.
Rather the multiplicity $\mu\upi$
of the cyclic base curve is determined from 
an infinitesimal neighborhood of 
the reducible fiber $S\upi$ in $\tiZ$
(Proposition \ref{prop:in}).
This seems to be an interesting invariant
associated to the present degeneration 
of smooth fibers into the reducible fiber $S\upi$.
From \eqref{me0} and \eqref{dimc}, it is possible to say that
the failure of the equality $m=e$,
or equivalently, the failure of the surjectivity
for the homomorphism $\rho$ in the derived sequence,  comes exactly from presence of
the cyclic base curves.

In Section \ref{ss:ex1}, we discuss various examples of the present twistor spaces.
In particular we see that the twistor spaces investigated
in \cite{Hon_Inv} and \cite{Hon_Cre2}
satisfies $m=e=n-2$.
This implies that they satisfy $\mu\upi = 0$ for any index $i$.
We also give an example of a twistor space on $n\CP^2$
for which $\mu\upi$ is strictly positive for some index $i$. Hence the strict inequality $m>e$ holds in general.
Obtaining such an example requires a deep analysis on the base locus of the line bundle
$\ms D(m)$.
We also given an example of a twistor space
for which the number $e$ satisfies $e=F(n)+1$, where
$F(n)$ is the $n$-th Fibonacci number.
This example shows that the number $m$ 
in Theorem \ref{thm:main1} can be much larger than the number $(n-2)$ in general.

For completing the proof of Theorem \ref{thm:main1},
it remains to show that a defining equation of 
the quartic hypersurface in $\CP^{m+2}$ which
cuts out the branch divisor  from the scroll has to be of the form
\begin{align}\label{qtc}
h_1h_2h_3h_4 = Q^2,
\end{align}
where $h_i$ are linear and $Q$ is quadratic.
In the cases of $3\CP^2$ and $4\CP^2$,
as we have already explained, 
the particular form of the equation
is a consequence of the existence of real {\em reducible} members
of the systems $|F|$ and $|2F|$ respectively.
Same is the case for the twistor spaces on $n\CP^2$
studied in \cite{Hon_Cre2}.
In all these cases, the precise cohomology classes
to which irreducible components
of the real reducible members should belong were given,
and the existence of the real reducible members was proved
by showing that these cohomology classes 
are represented by  an effective divisor.

We take a similar strategy for proving the existence
of the real reducible members of the system $|mF|$ 
on the present twistor spaces, but 
we work on the blown up space $\tiZ$ rather than the twistor space $Z$ itself.
In Section \ref{ss:rm}, 
we first observe that such a member
is directly related to real bitangents of the branch
quartics on real planes of the scroll.
More concretely, we observe that if $H\subset\CP^{m+2}$ is a real
hyperplane which corresponds to such a reducible member on $\tiZ$,
then the cut of real planes of the scroll $Y_m$ by $H$ 
is a real bitangent of the quartic on the plane.
This means that the existence of a real bitangent
of the branch quartic on real planes of the scroll is 
a necessarily condition for the system 
of $|mF|$ or $\big|\ms D(m)\big|$ to admit a real reducible member.
We will show by using the Pl\"ucker formula that,
the branch quartic on real planes
 has exactly three, two, or one real bitangent(s) if the twistor space $Z$ is of type $\Azero,\Aone$ or $\Atwo$ respectively
 (Propositions \ref{prop:btan1} and \ref{prop:btan2}).
The inverse image of a real bitangent
consists of two $(-1)$-curves on a fundamental divisor $S$ and they are mutually conjugate.
The union of these $(-1)$-curves
formed by moving members of the pencil
is expected to constitute a component of 
a real reducible member of $|mF|$.
In order to show that this is really the case,
we express the cohomology classes of these
$(-1)$-curves on $S$ concretely, and next,
by using them, we define non-real line bundles
over the space $\tiZ$, 
denoted by $\msDha$, which satisfy the property
\begin{align}\label{dha}
\msDha \otimes \omsDha \simeq \ms D.
\end{align}
Here, the superscript `h' stands for 
the `half' of the line bundle $\ms D$,
 and the index $\aaa$ is an integer parameter
whose range depends on the type of the twistor space.
These line bundles are defined in such a way
that their restriction to a smooth fiber $S$
of $\tif:\tiZ\to\CP^1$ gives one of the classes
of two $(-1)$-curves over a real bitangent.
A key for finding such a line bundle $\msDha$ is a formula
that expresses the Chern classes of the $(-1)$-curves
over a real bitangent in a concrete form 
(Proposition \ref{prop:cc01}).

It is not difficult to show that, for some integer
$\tau>0$, the line bundle $\msDha(\tau):=\msDha\otimes \tilde f^*\ms O(\tau)$ admits a non-zero section $\theta_{\aaa}$ which is unique up to constant,
and that such $\tau$ satisfies $2\tau\ge m$.
We denote the divisor defined by the section $\theta_{\aaa}$
by $\Theta_{\aaa}$.
From the relation \eqref{dha} we have
$\Theta_{\aaa} + \ol\Theta_{\aaa}\simeq 
\msDha(\tau) + \omsDha(\tau)\simeq
\ms D(2\tau)$.
Therefore if $2\tau=m$ holds, the divisor 
$\Theta_{\aaa} + \ol\Theta_{\aaa}$ 
is a real reducible member of the system $\big|\ms D(m)\big|$,
and its image to the twistor space $Z$ is the real reducible member of the system $|mF|$ we are searching for.
However this is not always correct.
For example the integer $m$ is not necessarily even
and in that case the equality $2\tau=m$ cannot hold.
The divisor $\Theta_{\aaa} + \ol\Theta_{\aaa}$ gives
a real reducible member of $\big|\ms D(2\tau)\big|$ but
we have $2\tau> m$ in general.

We show that, if $2\tau>m$, the divisor $\Theta_{\aaa}$
necessarily contains some components in the reducible fibers of $\tif$ as irreducible components.
In order to show this,
we make use of a non-real section $s_{\bm q}$ of the line bundle 
$\ms D(m)$ which corresponds to any non-real hyperplane in $\CP^{m+2}$
that passes the special point $\bm q$ on the ridge $\bm l$.
Since the divisors $E_2,E_2,\dots,E_k$ are contracted to 
the point $\bm q$ (Proposition \ref{prop:pim}), 
the section $s_{\bm q}$
vanishes identically on the exceptional divisors
$E_2,E_2,\dots,E_k$.
Then for any $l\ge 0$, over the polynomial ring
on $\CP^1$, 
the space $H^0\big(\ms D(m+l)\big)$ is generated by
the two sections $s_{\bm q}, \ol s_{\bm q}$ and
a defining section $s_{\sst{\bm D}}$ of the effective divisor $\bm D$
(Proposition \ref{prop:fvr}).
Hence, the product $\ta\ota\in H^0\big(\ms D(2\tau)\big)$
can be expressed as
\begin{align}\label{red60}
\ta\ota = g s_{\sst{\bm D}} + hs_{\bm q} + \ol {h s_{\bm q}},
\end{align}
for some $g\in H^0\big(\ms O_{\CP^1}(2\tau)\big)$ and 
$h \in H^0\big(\ms O_{\CP^1}(2\tau-m)\big)$,
where $g$ is real.
Then by making use of the above vanishing property for the section $s_{\bm q}$,
we can deduce that 
the zeroes of $h$ are contained in the set $\{\lmd\upone,\dots,\lmd\upk\}$,
and $g$ is divisible by $h$ (Propositions \ref{prop:h} and \ref{prop:g}).
Theses imply the desired conclusion about the divisor $\Theta_{\aaa}$.
If we write $\tilde T_{\aaa}$ for the divisor on $\tiZ$
obtained by removing all components of the reducible fibers from the divisor 
$\Theta_{\aaa}$, then the sum
 $\tilde T_{\aaa} + \sigma(\tilde T_{\aaa})$ 
gives the required real reducible member of $\big|\ms D(m)\big|$.
By projecting this divisor to $Z$, we finally obtain the 
real reducible member of the system $|mF|$. 

Once the existence of the real reducible members
of the system $\big|\ms D(m)\big|$ is obtained,
it is not difficult to derive the quartic equation
\eqref{quartic3}.
This will be done in the final section.
We also obtain that the type of the twistor space
is reflected to some constraint for the four linear polynomials
$h_i$ in the quartic polynomial.
Conversely, we can read off the type of the singularities
$\Aone, \Atwo$ and $\Athree$ on the plane quartics from
the last constraint for the polynomials $h_i$. 
But the author does not know whether it is possible to derive the quartic equation \eqref{quartic4} from the information on the singularities of the plane quartics.

\section{Moishezon twistor spaces
and the fundamental system}
\label{s:fs}
\subsection{Structure of divisors with small degree}
\label{ss:1and2}
Let $Z$ be a compact twistor space and
$D$ an effective divisor on $Z$.
By a \!{\em degree} of $D$, we mean 
the intersection number $D.\,l$
where
$l$ is a twistor line.
In this section we first recall briefly general results on the structure of divisors on compact twistor spaces
whose degrees are one or two.
All results are due to Poon \cite{P86,P92}
and Pedersen-Poon \cite{PP94}.
These result will be leverages when we
investigate the
structure of twistor spaces.

First, the degree of any effective divisor on $Z$ is
always positive since there always exists 
a twistor line intersecting the divisor
and we can  move it  fixing any one of the 
intersection points in a way that the 
curve is not contained in the divisor.
Such moving is possible because $N_{l/Z}\simeq
\ms O(1)^{\oplus 2}$ for the normal bundle
of a twistor line $l$.
This argument moreover means that 
any effective divisor of degree one is smooth.
Also it means that  such a divisor is irreducible and non-real.
%
%
When $Z$ is a twistor space on $n\CP^2$,
the structure of degree-one divisors is 
understood in much more concrete form:

\begin{proposition}
\cite{P92} 
\label{prop:poon1}
Any degree-one divisor $D$ on a twistor space
on $n\CP^2$ contains exactly one twistor line, and
it is a $(+1)$-curve on $D$.
Let $l$ be this twistor line,
and $\phi:D\to\CP^2$ the birational 
morphism induced by the linear system $|l|$.
Then $\phi$ is identified with $n$ points blowing-up of $\CP^2$,
where the points can be infinitely near.
If $e_1,\dots, e_n$ denote the exceptional
curves of  $\phi$,
we have a linear equivalence
\begin{align}\label{nD}
D|_D\simeq l - \sum_{i=1}^n e_i.
\end{align}
\end{proposition}

\noindent
{\em Outline of a proof.}
Let $\pi:Z\to n\CP^2$ be the twistor projection.
If the divisor $D$ does not contain any twistor line,
from the property $D.\,l=1$, 
the restriction $\pi|_D$ gives an orientation reversing
diffeomorphism $D\simeq n\CP^2$.
This means that $n\ol{\CP}^2$ (where $\ol\CP^2$
is the complex projective plane whose complex orientation
is reversed) admits a
complex structure, which cannot happen.
Indeed, as $b_1(n\ol\CP^2)=0$, the complex surface
admits a K\"ahler metric.
This contradicts that 
the intersection form of $n\ol{\CP}^2$ is 
negative definite.
So there is at least one twistor line 
$l$ which is contained in $D$.
It readily follows from $\deg D=1$ that 
$l$ is a $(+1)$-curve on $D$.
If $D$ contained another twistor line,
$D$ would have two $(+1)$-curves which are mutually disjoint.
This contradicts Hodge index theorem.

Now by the projection $\pi$ we have an orientation reversing diffeomorphism
$D\minus l\simeq n\CP^2\minus\{\pi(l)\}$. 
Hence checking the second Betti number, we obtain that 
the birational morphism $\phi$ is identified with
$n$ points blowing up of $\CP^2$.
Hence we have $K_D\inv \simeq 3l
-\sum_{1\le i\le n} e_i$.
On the other hand, noting that 
$D+\ol D\in |F|$ and $\ol D|_D=l$,
from adjunction formula, we readily obtain 
$D|_D\simeq K_D\inv - 2l$.
These mean the linear equivalence \eqref{nD}.
\proofend

\medskip
Note that the relation \eqref{nD} means that 
if $n\ge 2$, we always have
$\dim|D|\le 1$ and the equality can hold
only when the $n$ points  on $\CP^2$ to be blown up are in 
collinear configuration. 

An effective divisor whose degree is greater than one can have a singularity even if it is irreducible, as is already the case on $\CP^3$ (which is the twistor space of $S^4$).
But for a {\em real} irreducible  divisor of degree two,
we have the following smoothness result,
which we will make use of later.
Let (a.p.) be the anti-podal map $(z:w)\mapsto
(-\ol w:\ol z)$ on $\CP^1$, and (c.c.)  the 
complex conjugation $(z:w)\mapsto(\ol z:\ol w)$.

\begin{proposition}\label{prop:PP}
(Pedersen-Poon \cite{PP94})
Any real, irreducible divisor $S$ of
degree two on a compact twistor space
$Z$ is non-singular.
If $Z$ is a twistor space on $n\CP^2$,
for any such a divisor $S$, there is a sequence of blowing downs
\begin{align}\label{seqbr}
S = S_n \lras S_{n-1} \lras \dots \lras S_{1}\lras S_0 = \qdr,
\end{align}
where each map contracts a real pair of $(-1)$-curves, with
the following properties.
\begin{itemize}
\item  the induced real structure on $S_0$ is 
identified with the product 
(c.c.)$\times$(a.p.),
\item
the sequence \eqref{seqbr} does not blow up points
on the locus 
$\RP^1\times\CP^1\subset S_0$,
\item 
the strict transforms of real $(1,0)$-curves
on $S_0$ into $S$ are twistor lines.
\end{itemize}
%
\end{proposition}

\noindent
{\em Outline of a proof.}
If the divisor $S$
had a singularity, then the reality of $S$ and
absence of real points on $Z$ mean
that $S$ contains the twistor line which passes the singularity.
Then by blowing-up this twistor line
and looking at the strict transform of $S$,
it can be seen that the surface $S$ has  double points 
along the twistor line.
Hence the intersection of 
the strict transform of $S$ with the exceptional divisor
of the blowing up
is a curve which is two-to-one over the twistor line.
But it can be seen that this curve consists
of two disjoint $(+1)$-curves on the strict transform of $S$.
This contradicts Hodge index theorem,
so $S$ has to be non-singular.

If $S$ does not contain any twistor line,
the projection to $n\CP^2$ defines
unramified double cover $S\to n\CP^2$, which cannot happen as $n\CP^2$ is simply connected.
So $S$ contains at least one twistor line.
By adjunction formula, the self-intersection number of this twistor
line in $S$ equals zero.
From Serre duality, the vanishing of $H^1(\ms O_Z)$
and connectedness of $S$,
we obtain $H^1(\ms O_S)=0$.
Hence $S$ is rational and the linear system $|l|$ on $S$
defines a morphism $S\to\CP^1$.
All real fibers of this morphism has to be twistor lines in $Z$.
Thus $S$ has an $S^1$-family of twistor lines,
where $S^1=\RP^1\subset\CP^1$.
By repeating blow-downs of a real pair of $(-1)$-curves
contained in fibers of the above morphism $S\to\CP^1$,
one reaches at a relatively minimal model of $S$.
From the absence of real points on $S$,
this model has to be $\CP^1\times\CP^1$.
Moreover, since we have $K^2_S = (F|_S)^2 
= F^3 = 8-2n$ from \cite{Hi81},
the blow-down of a conjugate pair
of $(-1)$-curves can be done
exactly $n$ times.
The birational morphism $S\to S_0=\qdr$ obtained this way
clearly satisfies all the three property in the proposition.
\proofend

\medskip
Again some of the $2n$ points on $S_0=\qdr$ to be blown up
in the sequence \eqref{seqbr} can be infinitely near,
and such configurations of $2n$ points are a source of diversity of Moishezon twistor spaces.
On the other hand, unlike the birational morphism $\phi:D\to\CP^2$
for a degree-one divisor $D$ in Proposition \ref{prop:poon1},
the sequence \eqref{seqbr} is not canonical.
We note that if a fundamental divisor $S$ in a twistor space $Z$ admits a smooth rational curve
$G$ which satisfies $G^2\le -2$, then 
the curve $G$ is an obstruction for the twistor space to admit a K\"ahler metric
since we have $F.\,G = K_S\inv.\,G = G^2+2\le 0$.

\begin{remark}{\em
\label{rmk:cfg}
(i) The banned locus $\RP^1\times\CP^1\subset\CP^1$ for the blowing up
in the sequence \eqref{seqbr}
disconnects $S_0=\CP^1\times\CP^1$ into
two connected components, and they are interchanged by the real structure.
Hence the $2n$ points on $S_0=\qdr$ to be blown up are naturally 
divided into two groups.
(ii) 
From the three conditions in Proposition \ref{prop:PP},
at most half of the $2n$ points on $S_0=\qdr$ can lie on the same $(1,0)$-curve  or on the same $(0,1)$-curve.
}\end{remark}

\begin{remark}\label{rmk:Hopf}{\em
A {\em non-real} irreducible fundamental divisor
on a twistor space of $n\CP^2$ is
not necessarily a rational surface.
Indeed, if $n\ge 4$, there is an example of a twistor space
on $n\CP^2$ which has a non-real irreducible fundamental divisor  
that is birational to a Hopf surface
\cite{Hon_CM}.
In these examples, such a divisor has self-intersection along a curve.
Of course, these twistor spaces are not Moishezon,
but they are obtained as a small deformation
of Moishezon twistor spaces.
}\end{remark}

\subsection{The anti-canonical system of a fundamental divisor and structure of a twistor space}\label{ss:acsfd}
Let $S$ be a smooth fundamental divisor on a twistor space $Z$ of $n\CP^2$.
Then since $2F\simeq -K_Z$ for the anti-canonical class of $Z$, 
we obtain from adjunction formula 
$
K_S\simeq K_Z + S\,|_S\simeq -2F +F\,|_S\simeq -F|_S.
$
Hence we obtain a basic relation
$F|_S\simeq K_S\inv$,
as well as an exact sequence
\begin{align}\label{ses:000}
0 \lras \ms O_Z \lras F \lras K_S\inv \lras 0.
\end{align}
Then since $H^1(\ms O_Z)=0$ from
simply connectedness of $n\CP^2$ as mentioned in 
Section \ref{ss:I},  we obtain from this a simple formula 
\begin{align}\label{fml001}
h^0(F) = h^0\big(K_S\inv\big)+1.
\end{align}
The surjectivity of the restriction homomorphism
$H^0(F)\to H^0\big(K_S\inv\big)$ also implies  the
coincidence of the base locus $\Bs\,|F|=\Bs\,\big|K_S\inv\big|$,
and the restriction of the meromorphic map from $Z$ induced by
the fundamental system $|F|$ to the divisor $S$ is
identified with the anti-canonical map from $S$.
Thus the fundamental system of a twistor space
is closely related to the anti-canonical system
of a fundamental divisor.
In this subsection, based on 
Proposition \ref{prop:PP}, 
we consider some configurations of the $2n$ points
on $S_0=\qdr$
for which the anti-canonical system
of the blowing up have enough members,
and discuss how they bring information on the structure of
twistor spaces.
The existence of a real irreducible 
fundamental divisor on a twistor space will be discussed
in the next subsection.

If $S$ is a real irreducible fundamental divisor
on a twistor space $Z$ over $n\CP^2$, 
it is smooth and isomorphic to 
$2n$ points blowing up of $S_0=\qdr$
by Proposition \ref{prop:PP}.
When $n=2$, by using Remark \ref{rmk:cfg} (ii), it is not difficult to see that 
the anti-canonical system $\big|K_S\inv\big|$ {\em always} satisfies
\begin{align}\label{p2}
\Bs\,\big|K_S\inv\big| = \emptyset,\quad \dim \big|K_S\inv\big| = 4\qandq
K_S^2 = 4.
\end{align}
One also obtains that the image of the anti-canonical morphism $S\to\CP^4$ is
an intersection of two quadrics and the morphism is
 birational over this complete intersection.
For the twistor space $Z$ on $2\CP^2$, from \eqref{p2}, we obtain
$$
\Bs\,|F| = \emptyset,\quad \dim |F| = 5\qandq
F^3 = 4.
$$
Poon \cite{P86} further deduced that
the image of the morphism
$Z\to\CP^5$ induced by the fundamental system is 
an intersection of two quadrics,
and that the morphism is birational
onto the complete intersection in $\CP^5$,
which contracts
precisely four smooth rational curves
to ordinary double points. (See Proposition \ref{prop:Poon}.)

The case $n=3$ is much more complicated, but we can 
determine the structure of the anti-canonical map 
from the blowing up, for all possible configurations of $6$ points
on $S_0=\qdr$ in Proposition \ref{prop:PP} in the following manner.
Although the argument will be a bit long and
does not yield a new result on twistor spaces
on $3\CP^2$, 
it will be significant in the present study
of twistor spaces.

Recall from Remark \ref{rmk:cfg} (ii) that 
at most half of the 6 points on $S_0=\qdr$ can lie on 
a common $(1,0)$-curve
and on a common $(0,1)$-curve.
Suppose that a half of the 6 points are on a common $(1,0)$-curve
or on a common $(0,1)$-curve.
Let $C_0\subset S_0$ be this curve.
The remaining 3 points are on the conjugate curve $\ol C_0$.
The strict transforms of $C_0$ and $\ol C_0$ into the blowing up $S$ are $(-3)$-curves, and hence
they are base curves of the anti-canonical system of $S$.
Further we readily obtain that after removing
these two base curves the system is base point free
and is composed with a  pencil 
generated by pull-backs of $(0,1)$-curves (resp.\,$(1,0)$-curves).
This means $\dim \big|K_S\inv\big|=2$ and
that the anti-canonical map
of $S$ is a morphism $S\to \CP^2$ whose image is
a conic.
This surface $S$ admits an obvious non-trivial $\CC^*$-action
which fixes all points on the two base curves.


In the following we assume that 
no three points among the 6 points on $S_0=\qdr$
belong to a $(1,0)$-curve nor a $(0,1)$-curve.
To discuss it is convenient to introduce the following

\begin{definition}\label{def:acc}
{\em
Let $S$ be a real irreducible fundamental divisor
on a twistor space over $n\CP^2$.
By a {\em real anti-canonical cycle} on $S$,
we mean a real, reducible anti-canonical curve on the surface $S$, which is of the form
\begin{align}\label{cycle20}
\big(C_1+C_2+\dots + C_k\big) + 
\big(\ol C_1+ \ol C_2+\dots + \ol C_k\big),\quad k\ge 1,
\end{align}
where 
\begin{itemize}
\item 
any components $C_i$ and $\ol C_i$ are
smooth rational curves,
\item
two adjacent components in the presentation \eqref{cycle20}, as well as
the two components $\ol C_k$ and $C_1$,
intersect transversally at one point,
\item 
any other two components do not intersect.
\proofend
\end{itemize}
}
\end{definition}

\begin{figure}
\includegraphics{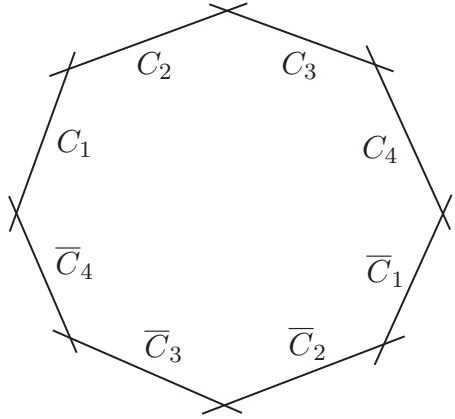}
\caption{
A real anti-canonical cycle in the case $k=4$.}
\label{fig:octagon0}
\end{figure}

The case $k=4$ is illustrated in Figure \ref{fig:octagon0}.
We note that from the definition,
the real structure
on the cycle interchanges the opposite components
in the cycle.
In particular,
a real component is not allowed to be included
in the cycle.
(This excludes $S$ which is contained 
in the twistor spaces of Campana-Kreussler.)

As in Remark \ref{rmk:cfg} (i),  on each of the two connected components of 
the complement of  $\RP^1\times\CP^1$
in $S_0=\qdr$,
precisely half of the 6 points are belonging,
and now assume that these 3 points are not on the same
$(1,0)$-curve nor the same $(0,1)$-curve.
This condition means that {\em the 3 points
on the same connected component
uniquely determine a $(1,1)$-curve on $S_0=\qdr$
by passing requirement}.
Hence we get two $(1,1)$-curves on $S_0$.
If these two curves were equal,
it has to be real.
But since the real structure 
on $S_0 =\qdr$ was given by
 (c.c.)$\times$(a.p.) as in Proposition \ref{prop:PP}, 
there is no real
$(1,1)$-curve on $S_0$.
Therefore the two $(1,1)$-curves are distinct.
Let $C_0$ and $\ol C_0$ be these $(1,1)$-curves on $S_0$.
Each of these can be reducible.
But it is not difficult to see that even in that case
we can change the choice of the blowing down $S\to S_0=\qdr$ in \eqref{seqbr} in such a way
that the curves $C_0$ and $\ol C_0$ are irreducible.

Hence in the following we suppose that the $(1,1)$-curves
$C_0$ and $\ol C_0$ are 
irreducible, and investigate
all possible structure for the anti-canonical map
of $S$ using these curves. 
Since there is no real point on $S_0=\qdr$,
the two irreducible $(1,1)$-curves $C_0$ and $\ol C_0$
intersect transversally at two points.
So the curve $C_0+\ol C_0$ is a real 
anti-canonical cycle on $S_0$ in the sense
of Definition \ref{def:acc}.
Then the following four situations $\Azero,\,\Aone,\,\Atwo$
and $\Athree$ will exhaust all
possible configurations of the 6 points on $S_0$ to be blown up.
The reason for these notations will be evident in
the following argument.
\begin{figure}
\includegraphics{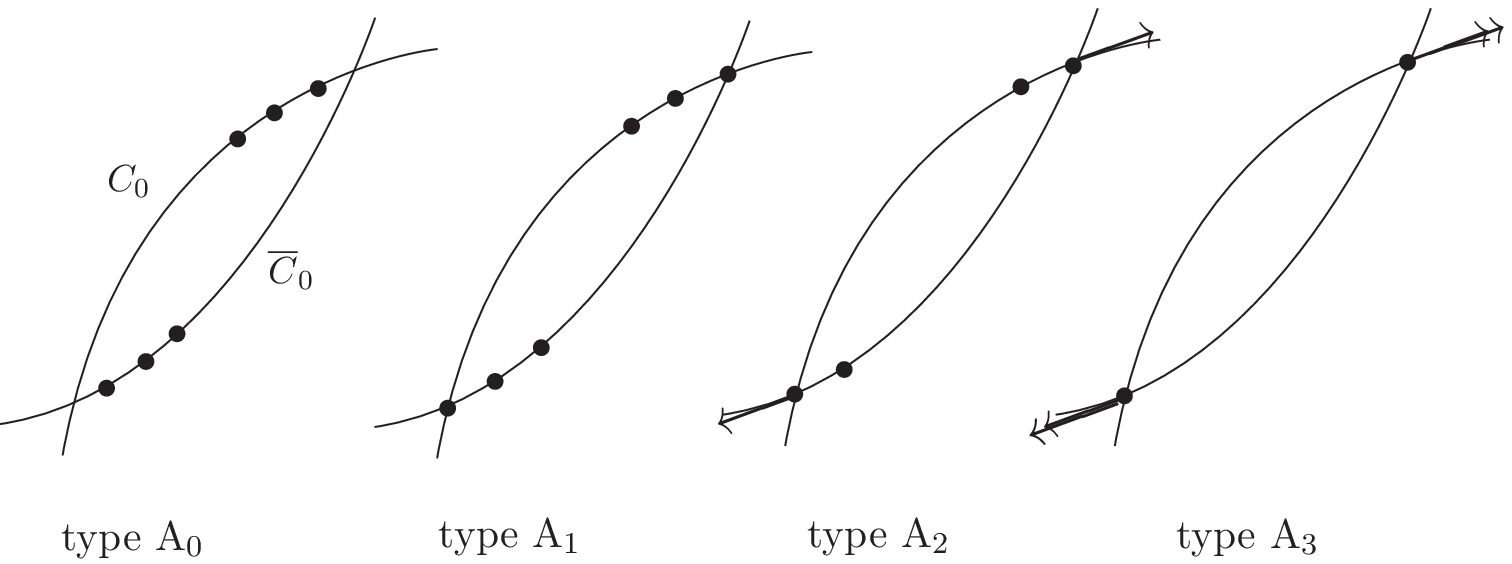}
\caption{
Configurations of the 6 points on $S_0=\qdr$ which yield 
the surface $S$ satisfying $\Bs\,\big|K_S\inv\big|=\emptyset$ in the case $n=3$.
}
\label{fig:bu_cfg1}
\end{figure}

\medskip\noindent
{{\bf A}$_0$.}
Assume that none of the 6 points belongs to the 
intersection $C_0\cap \ol C_0$ (see Figure \ref{fig:bu_cfg1}, type $\Azero$).
If we denote the strict transforms of $C_0$ and $\ol C_0$ into $S$ by $C_1$ and $\ol C_1$ respectively,
these are $(-1)$-curves on $S$, and 
the curve $C:=C_1+\ol C_1$ on $S$ 
is a real anti-canonical cycle on $S$,
still consisting of two components.
We readily have $K_S\inv.\, C_1 = K_S\inv.\,\ol C_1 = 1$,
and  $h^0\big(K_S\inv|_C\big) = 2$.
From the exact sequence
\begin{align}\label{sesS}
0 \lras \ms O_S \lras K_S\inv \lras K_S\inv|_C \lras 0,
\end{align}
this means $h^0\big(K_S\inv\big)=1+2=3$ and $\Bs\,\big|K_S\inv\big|=\emptyset$.
Since $K_S^2=8-6=2$, this implies that the anti-canonical morphism
$\phi:S\to\CP^2$ is
generically two-to-one.
From the adjunction formula a generic anti-canonical curve on $S$ is a
smooth elliptic curve,
and this means that the branch divisor of the morphism $\phi$
is a quartic.
The morphism $\phi$ maps the two curves $C_1$ and $\ol C_1$
to the same line isomorphically.
Let $\bm l$ be this line. This is real.
Then the branch quartic
intersects $\bm l$ tangentially at the images of the two points $C_1\cap \ol C_1$. 
We denote by $\bm q$ and $\bm {\ol q}$ for these two tangent points.(See Figure \ref{fig:bq}.)
This surface $S$ is a del Pezzo surface of degree two.

\medskip\noindent
{{\bf A}$_1$.}
Second suppose that 
among the 6 points precisely two belong to the intersection 
$C_0\cap\ol C_0$ (see Figure \ref{fig:bu_cfg1}, type $\Aone$).
Then on the curve $C_0$ exactly 4 points are belonging among the 6 points.
Let $C_1$ and $\ol C_1$ be the exceptional curves over 
the two points $C_0\cap \ol C_0$,
and $C_2$ and $\ol C_2$ the strict transforms
of $C_0$ and $\ol C_0$ respectively.
Put
$C:=C_1+C_2 + \ol C_1 + \ol C_2$ on $S$.
This is a real anti-canonical cycle
on $S$ and we have $C_1^2 = \ol C_1^2 = -1$ and $C_2^2 = \ol C_2^2=-2$.
Again we have the exact sequence \eqref{sesS},
and from it we deduce that  the anti-canonical map $\phi$ of 
$S$ is a morphism onto $\CP^2$ which is generically two-to-one
and whose branch is a quartic.
The components $C_2$ and $\ol C_2$ are 
mapped to mutually distinct points by $\phi$.
Let $\bm q$ and $\bm {\ol q}$ be these two points.
If $\bm l$ denotes the real line through these points,
the exceptional curves $C_1$ and $\ol C_1$
are mapped isomorphically to $\bm l$.
The branch quartic passes $\bm q$ and $\bm {\ol q}$,
and they are ordinary double points (i.e.\,$\Aone$-singularities) 
of the quartic.
(See Figure \ref{fig:bq}.)
The genuine double cover of $\CP^2$ branched along this quartic is
a Del Pezzo surface of degree two with two $\Aone$-singularities, 
and the surface $S$ is obtained from this double cover by 
taking the minimal resolution.

\medskip\noindent
{{\bf A}$_2$.}
Next suppose that
among the 6 points exactly 4 belong to the intersection $C_0\cap \ol C_0$
(see Figure \ref{fig:bu_cfg1}, type $\Atwo$).
This implies that, at each of these intersection points,
two points are placed as infinitely near points.
Hence,  one of these two points is a tangent vector
at the point, and from our choice of the $(1,1)$-curves $C_0$ and $\ol C_0$ on $S_0$,
the vector has to be tangent to $C_0$ or $\ol C_0$.
So there are two possibilities according to 
which curve the vector is tangent, 
but without loss of generality we may suppose
that the vectors are as in Figure \ref{fig:bu_cfg1}.
From each of the two points $C_0\cap \ol C_0$
there arise two exceptional curves, 
and consequently on the surface $S$ 
we obtain a real anti-canonical cycle $C$ consisting of 6 components.
Exactly two of them are $(-1)$-curves,
and the remaining 4 components are $(-2)$-curves.
These form two chains in $C$.
For the cycle $C$ we still have the exact sequence
\eqref{sesS}, and again the anti-canonical map is 
a generically two-to-one morphism $S\to\CP^2$ whose branch 
is a quartic.
Among the 6 components of the cycle $C$, 
the two chains of $(-2)$-curves are mapped to 
distinct points.
Let $\bm q$ and $\ol{\bm q}$ be these points,
and $\bm l$ the real line through these points.
Then the two $(-1)$-curves in $C$
are mapped isomorphically to $\bm l$,
and the branch quartic has ordinary cusps (i.e.\,$\Atwo$-singularities) at the points
 $\bm q$ and $\ol{\bm q}$.
The genuine double cover of $\CP^2$ branched along this quartic is
a Del Pezzo surface of degree two with two $\Atwo$-singuarities,
and the surface $S$ is obtained from this double cover by 
taking the minimal resolution.

\medskip\noindent
{{\bf A}$_3$.}
Finally suppose that 
all the 6 points belong to the intersection $C_0\cap \ol C_0$.
This means that exactly 3 points are placed
as infinitely near points at each of 
the two points $C_0\cap \ol C_0$.
Hence 3 exceptional curves appear
from each of the two points, and 
we get a real anti-canonical cycle on the surface $S$
consisting of 8 components.
We again write $C$ for this cycle.
At each of the two points of $C_0\cap \ol C_0\subset S_0$,
two directions are specified by 
two among the three infinitely near points.
From the passing requirement, these two directions have to be  
the tangent direction of the $(1,1)$-curve $C_0$ or $\ol C_0$ on $S_0$.
Hence the two directions are either the same or
 different at each point.

If the two directions are the same as in Figure
\ref{fig:bu_cfg1}, type $\Athree$,
exactly two components among the 8 components of the cycle $C$
are $(-1)$-curves, and all other components are
$(-2)$-curves.
Thus we again obtain two chains of $(-2)$-curves in the cycle $C$,
and this time each chain consists of three components.
From these, by using the exact sequence \eqref{sesS},
we again obtain that $\dim \big|K_S\inv\big|=2$,
$\Bs\,\big|K_S\inv\big|=\emptyset$ and the anti-canonical
morphism $\phi:S\to \CP^2$ is  generically two-to-one 
with branch being a quartic curve.
The two chains of $(-2)$-curves are mapped to distinct points
on $\CP^2$ by $\phi$.
Let $\bm q$ and $\ol{\bm q}$ be these points,
and $\bm l$ the real line through these points.
The two $(-1)$-curves in $C$ are mapped isomorphically
to the  line $\bm l$.
The branch quartic
has $\Athree$-singularities at $\bm q$ and $\bm {\ol q}$.
An important remark is that in this situation,
the surface $S$ admits a non-trivial $\CC^*$-action
which is the lift of the $\CC^*$-action 
on $S_0=\qdr$ preserving the two $(1,1)$-curves
$C_0$ and $\ol C_0$.
The cycle $C$ on $S$ is of course $\CC^*$-invariant,
and points on the middle components of the two chains
of $(-2)$-curves are fixed by the $\CC^*$-action.
The genuine double cover of $\CP^2$ branched along the above quartic is
a Del Pezzo surface of degree two with two $\Athree$-singuarities,
and the surface $S$ is obtained from this double cover by 
taking the minimal resolution.

Alternatively, if the two directions at each of the
intersection points $C_0\cap \ol C_0$ are  different,
two among the 6 exceptional curves of the 
blowing up $S\to S_0=\qdr$ are $(-3)$-curves,
and hence they are base curves of the anti-canonical system of $S$.
It is not difficult to see that the surface $S$ obtained in this way is 
isomorphic to the one
determined by 6 points on $S_0$ whose half belong to the same $(1,0)$-curve or the same $(0,1)$-curve.
This is the surface we discussed right before the case A$_0$,
and in the following we exclude this configuration
from the case of type $\Athree$.

\medskip
Thus all configurations $\Azero,\Aone,\Atwo$ and $
\Athree$ of the 6 points on $S_0=\qdr$
generate a surface $S$ whose anti-canonical system
is base point free,
inducing
a degree-two morphism $\phi:S\to\CP^2$
with the branch being a quartic.
Combined with the configuration discussed
right before the case $\Azero$
(which is excluded from the $\Athree$-configuration),
we have obtained the following conclusion about 
structure of the anti-canonical system of 
a real irreducible fundamental divisor on 
a twistor space of $3\CP^2$.
\begin{proposition}\label{prop:S2}
Let $S$ be a real irreducible fundamental divisor
on a twistor space over $3\CP^2$.
Then we always have $\dim \big|K_S\inv\big|=2$,
and  only the following two situations occur:
\begin{itemize}
\item[\em (i)]
The system $\big|K_S\inv\big|$ has
exactly two base curves, and 
after removing them, the system 
is composed with a pencil of rational curves which is 
base point free.
\item[\em (ii)]
The system $\big|K_S\inv\big|$ is base point free,
and the induced morphism $S\to\CP^2$ is 
generically two-to-one whose branch is 
a quartic.
\end{itemize}

\end{proposition}

From the above argument, if $S$ satisfies the item (i)
of the proposition,
one can always take the sequence \eqref{seqbr}
in a way that precisely a half (i.e.\,three) of the blown up points on $S_0=\qdr$ are on the same $(1,0)$-curve or on
the same $(0,1)$ curve.
The strict transform of 
this curve and its conjugation
are  base curves of $\big|K_S\inv\big|$.
Moreover, the surface $S$ admits a non-trivial $\CC^*$-action,
and  all points are fixed by the $\CC^*$-action.

Also from the above argument,
the divisor $S$ satisfying the property (ii) in the proposition
can be classified into four types $\Azero,\Aone,\Atwo$ and $\Athree$
once we fix the sequence \eqref{seqbr}
whose $(1,1)$-curves on $S_0$ determined by the passing
requirement are irreducible.
These two $(1,1)$-curves determine
the real anti-canonical cycle $C$ on $S$,
and the number of components of $C$ depends on the type.
The type can also be realized  as the type of singularities of the branch quartic.
Note that a single surface $S$ can have two real anti-canonical cycles whose numbers of components
are distinct,
and the type is well-defined not for the surface $S$ itself but
rather than for the pair $(S,C)$.
This ambiguity will disappear when we consider a fundamental divisor
on a twistor space on $n\CP^2$, $n>3$.

If the divisor $S$ in Proposition \ref{prop:S2}
satisfies the property (i) in the proposition,
the twistor space is a LeBrun space
which is the situation (i) 
in Proposition \ref{prop:3},
and if $S$ satisfies the property (ii)
in Proposition \ref{prop:S2},
the twistor space is of double
solid type as in (ii) of Proposition \ref{prop:3}.
%
This way one can derive the structure
of a twistor space from
the structure of the anti-canonical system
of a real irreducible fundamental divisor.

%
As for mutual relationship between these four kinds of surfaces, from the configurations of the 6 points on $S_0=\qdr$
as presented in Figure \ref{fig:bu_cfg1},  
they
are under the adjacent relations
\begin{align}\label{adr}
{\rm A}_0 \lras {\rm A}_1 \lras {\rm A}_2 \lras {\rm A}_3,
\end{align}
where the arrows indicate specialization.
Type $\Azero$ is most general, and type $\Athree$ is most special.
The surfaces of type A$_3$ admit a non-trivial
$\CC^*$-action.

As an easy consequence of Proposition \ref{prop:S2},
we have the following.

\begin{proposition}\label{prop:LB}
If $n\ge 4$ and $Z$ is a twistor space on $n\CP^2$,
any real irreducible fundamental divisor $S$ satisfies an inequality
\begin{align}\label{inLB}
\dim \big|K_S\inv\big|\le 2.
\end{align}
Further, the equality holds
if and only if we can take the sequence \eqref{seqbr}
in such a way that a half of the blowing up points on $S_0=\qdr$ are on the same $(1,0)$-curve or on the same $(0,1)$-curve.
\end{proposition}

\proof
By Proposition \ref{prop:S2}, the surface $S_3$ in 
the sequence \eqref{seqbr} always satisfies $h^0\big(K\inv\big)= 3$,
as well as the property (i) or (ii) in the proposition.
Since $h^0\big(K\inv\big)$ does not increase under blowing up,
we obtain the inequality \eqref{inLB}.
Suppose that the surface $S_3$ satisfies the property (i).
If $C_0$ and $\ol C_0$ are the base curves of $\big|K_{S_3}\inv\big|$, we may suppose that  
the image of $C_0$ to $\qdr$ by
the sequence \eqref{seqbr} is a $(1,0)$-curve
or a $(0,1)$-curve, and then
a half (i.e.\,three) of the blown up points
on $S_0$ are on this curve.
Then $h^0\big(K_S\inv\big)=3$ holds
if and only if all the remaining points to be
blown up lie on the curves $C_0$ and $\ol C_0$.
This condition implies that a half of blowing up points
are on
the same $(1,0)$-curve or on the same $(0,1)$-curve.

If the surface $S_3$ satisfies (ii) in 
Proposition \ref{prop:S2},
since $\big|K_S\inv\big|$ is base point free,
the surface $S$ always satisfies an inequality $h^0\big(K_S\inv\big)<3$
regardless of where the remaining points
to be blown up are located.
\proofend

\medskip
From the proof, the surface $S$ which attains the equality in \eqref{inLB}
has a non-trivial $\CC^*$-action which is a lift of the $\CC^*$-action
in the case $n=3$.
The $\CC^*$-action has exactly two curves whose points are fixed,
and these curves are the base locus of the anti-canonical system
of $S$.
The subtraction of these curves from the anti-canonical system
is composed with a pencil of rational curves.
If a twistor space $Z$ on $n\CP^2$, $n>3$, admits
a real irreducible fundamental divisor $S$ of this kind,
$Z$ has to be a LeBrun space \cite{P92,KK92}.

Later we will generalize Proposition \ref{prop:S2}
to the case $n\ge 3$.
For this purpose we rephrase the proposition
in the following form.

\begin{proposition}
\label{prop:S2'}
Let $S$ be as in Proposition \ref{prop:S2}.
Then it satisfies at least one of the following  two properties.
\begin{itemize}
\item[\em (i)]
The surface $S$ admits a 
non-trivial $\CC^*$-action which preserves 
a real anti-canonical cycle $C$ on $S$,
and there is a real pair of components of $C$
whose points are fixed by the $\CC^*$-action.
\item[\em (ii)]
The anti-canonical system $\big|K_S\inv\big|$ is 
base point free and 
the associated morphism is a generically two-to-one morphism $S\to\CP^2$ whose branch
divisor is a quartic curve.
\end{itemize}
\end{proposition}

In this form the surface $S$ can satisfy
both (i) and (ii),
and such surfaces are exactly the ones
in the case of type A$_3$.
If $S$ satisfies (i) in Proposition \ref{prop:S2},
the real anti-canonical cycle $C$ in  (i) of the present proposition
is formed by the two base curves of the system $\big|K_S\inv\big|$ and
a real member of the movable part of the
anti-canonical system.

We end this subsection by presenting a 
(well-known) configuration of  $2n$ points
on $S_0=\qdr$
which yield surfaces satisfying 
$\dim \big|K_S\inv\big|=1$ for arbitrary $n>3$.
Together with the above characterization of LeBrun spaces,
this configuration provides a typical example
where the structure of a twistor space
may be described using the fundamental system.

Recall that the real structure on 
$S_0=\qdr$ is given by (c.c.)$\times$(a.p),
and 
let $C_0\subset\qdr$ be a real irreducible $(1,2)$-curve.
This is automatically smooth.
Next take $2n$ points on $C_0$ 
which are real as a whole and which are
not on the locus $\RP^1\times\CP^1$.
Let $S$ be the surface obtained from $S_0$ by blowing up these $2n$ points,
and use the same letter to 
mean the strict transform of $C_0$ into $S$.
Then if $n>3$, we have $C_0^2 = 4 -2n<-2$
and hence $C_0\subset S$ is a base curve of 
the anti-canonical system $\big|K_S\inv\big|$.
The system $\big|K_S\inv-C_0\big|$ is a pencil and 
base point free.
In particular we have $\dim\big|K_S\inv\big|=1$.
Conversely if $n>4$ and $\dim \big|K_S\inv\big|=1$ hold,
$S$ has to be obtained from $S_0=\qdr$ in this way.
(When $n=4$ this does not hold and 
we need to add an assumption $\Bs\,\big|K_S\inv\big|\neq\emptyset$
or $\kappa\inv(S)=2$ to obtain the same conclusion.)
If a twistor space $Z$ on $n\CP^2$, $n>3$, admits
this surface $S$ as a real  fundamental divisor, then we have 
$$\dim |F|=2\qandq \Bs\,|F|=C_0,$$
and $Z$ has to be a twistor space
studied by Campana-Kreussler \cite{CK98}.
But unlike LeBrun spaces, the corresponding
self-dual structure on $n\CP^2$ is not
explicitly known yet.

\subsection{Existence of a real irreducible fundamental divisor}
In the last subsection we discussed
how structure of the anti-canonical system
of a real irreducible fundamental divisor 
brings information on the structure of a twistor space.
In this subsection we mainly discuss the existence
of a real irreducible fundamental divisor.

\begin{proposition}\label{prop:irr}
Let $Z$ be any twistor space on $n\CP^2$.
If $h^0(F)\ge 2$, generic members of the fundamental system $|F|$ are irreducible.
\end{proposition}

By Proposition \ref{prop:PP}, this means that
generic fundamental divisors are smooth if 
$h^0(F)\ge 2$.

\medskip\noindent
{\em Proof of Proposition \ref{prop:irr}.}
It is enough to show the case $n\ge1$.
Since the degree of any effective divisor is positive
as in the beginning of Section \ref{ss:1and2},
if all members of $|F|$ are reducible,
there exists a degree-one divisor $D$
which satisfies $\dim |D|\ge 1$.
Choose  sections $s,t\in H^0\big(\ms O_Z(D)\big)$
which are linearly independent
and consider the 4 sections  $s\ol s, s\ol t, t\ol s$
and $t\ol t$ of $F$.
These satisfy an obvious quadratic relation.
Hence if these sections were linearly dependent,
the sub-system of $|F|$ spanned by these 4 sections would be composed with a pencil.
But since $n\ge 1$, the two pencils
spanned by $s,t$ 
and $\ol s,\ol t$ are not included in the same linear system because $D$ and $\ol D$ are not 
linearly equivalent. 
This means that the sub-system of $|F|$ generated by the
above 4 sections cannot be composed with a pencil.
Hence the 4 sections  are linearly independent,
and define a meromorphic map
from $Z$ onto a smooth quadric in $\CP^3$.
Then members of $|F|$ which correspond to
generic hyperplane sections of the quadric
are irreducible, since otherwise we would 
obtain a 3-dimensional family of degree-one divisors, which contradicts the relation \eqref{nD}.
\proofend

\medskip
By the Riemann-Roch formula \eqref{RR1}
and Hitchin's vanishing theorem \eqref{Hv}, 
if $Z$ is a twistor space on $n\CP^2$ whose self-dual structure
is of positive type,
we have an inequality
\begin{align}\label{RR7}
h^0(F)\ge 10-2n.
\end{align}
Hence we have $h^0(F)\ge 2$ if $n\le 4$.
Therefore when $n\le 4$,
there always exists a real irreducible fundamental divisor.
Hence the results on the structure of twistor spaces
on $2\CP^2$ and $3\CP^2$ given in the last subsection 
which were derived under the assumption on
the presence of a real irreducible fundamental divisor
are not restrictive at all.
This is a possible explanation as to why the classification results
for the twistor spaces on $2\CP^2$ and $3\CP^2$
stated in Propositions \ref{prop:Poon} and \ref{prop:3} hold.
Also, from the discussion
in the last subsection,
Proposition \ref{prop:irr} leads us to the conclusions that 
when $n\ge 4$, a twistor space
on $n\CP^2$ satisfying $\dim|F| = 3$ has to be a LeBrun space,
and when $n> 4$, a twistor space
on $n\CP^2$ satisfying $\dim|F| = 2$ has to be a Campana-Kreussler space.

We will also need the following property which holds for any fundamental divisor.

\begin{proposition}\label{prop:S03}
Any fundamental divisor $S$ on a twistor space $Z$ on $n\CP^2$ satisfies
$$H^1(\ms O_S)=H^2(\ms O_S)=0.$$
\end{proposition}

\proof
Since we always have
$
H^q(\ms O_Z) = 0
$ for any $q>0$,
from the exact sequence
\begin{align}\label{}
0 \lras \ms O_Z(-S)
\lras \ms O_Z 
\lras \ms O_S
\lras 0,
\end{align}
we have $H^q(\ms O_S)\simeq H^{q+1}
\big(\ms O_Z(-S)\big)$  for any $q>0$.
As $\ms O_Z(-S)\simeq -F$,
when $q=1$, we have
$$H^1(\ms O_S)\simeq H^2\big(\ms O_Z(-S)\big)\simeq H^2(-F)
\simeq H^1(K_Z+F)^* \simeq H^1(-F)^*,$$
and the last one vanishes by 
Hitchin's vanishing theorem \eqref{Hv}.
Similarly, we have 
$H^2(\ms O_S)\simeq H^0(-F)^*$, which 
again vanishes from the degree.
So we get $H^1(\ms O_S)=
H^2(\ms O_S)=0$.
\proofend

\medskip
As in Remark \ref{rmk:Hopf},
even if the divisor $S$ is irreducible, it is not necessarily rational.

\subsection{Pencil of fundamental divisors
and a real anti-canonical cycle}
\label{ss:pfd}
In the last two subsections we have seen that 
basic structure of some twistor spaces can be obtained
from that of the anti-canonical system of a real irreducible
fundamental divisor.
Evidently this method works effectively
only when the anti-canonical system of the fundamental divisor
is of positive dimension since otherwise
the fundamental system of a twistor space is at most a pencil by \eqref{inLB}.
We have also seen that  such a situation
is quite restrictive and only LeBrun spaces
and Campana-Kreussler spaces satisfy this property
when $n>4$.

However, 
there are numerous twistor spaces on $n\CP^2$
whose fundamental system is a pencil.
These twistor spaces generically satisfy $a(Z)=1$,
and the pencil induces
an algebraic reduction of the space.
It seems difficult to obtain further results on the structure of such twistor spaces.
When $n>4$, the algebraic dimension cannot be two
if the fundamental system is a pencil \cite{HK17}.
In contrast, {\em there are a lot of Moishezon twistor 
spaces on $n\CP^2$, $n\ge 4$, whose fundamental system is a pencil.}
For example, as mentioned in Section \ref{ss:I}, 
the twistor spaces associated to Joyce metrics
\cite{J95} and also some small deformations of them
satisfy these properties.
As a matter of fact, when $n\ge 4$, 
all known Moishezon twistor spaces on $n\CP^2$ satisfy the property $\dim |F|=1$ 
except the LeBrun spaces
(which satisfy $\dim |F|=3$) and Campana-Kreussler spaces (which satisfy $\dim |F|=2$).

As mentioned in Section \ref{ss:I}, 
we always have $a(Z)=\kappa(F)$
for the algebraic dimension of $Z$, and therefore 
$F$ is big if the twistor space is Moishezon.
Thus a natural idea for studying Moishezon twistor spaces satisfying $\dim |F|=1$
is to investigate pluri-system $|lF|$
for sufficiently large $l$
instead of $|F|$ itself.
By Proposition \ref{prop:irr},
these twistor spaces always admit
a real irreducible fundamental divisor $S$.
So we have an exact sequence
\begin{align}\label{ses:002}
0 \lras (l-1)F \lras lF \lras lK_S\inv
\lras 0,
\end{align}
and we expect that the dimension $h^0(lF)$ is computable from the cohomology exact sequence
of this sequence.
But as explained in Section \ref{ss:I},
when $n>4$, 
we have $H^1(lF)\neq 0$ for any $l>0$ from 
Riemann-Roch formula \eqref{RR1} and Hitchin's vanishing theorem \eqref{Hv},
so it seems difficult to calculate 
$h^0(lF)$ from the cohomology exact sequence
of \eqref{ses:002}.
%
%

But when $n=4$,
from Riemann-Roch \eqref{RR1}, we have $\chi(F) = 2$.
So from the vanishing theorem \eqref{Hv},
if $\dim |F|=1$, we have $H^1(F)=0$.
Hence the cohomology exact sequence associated to \eqref{ses:002} gives an exact sequence
\begin{align*}
0 \lras H^0(F) \lras H^0(2F) \lras
H^0 (2K_S\inv) \lras 0.
\end{align*}
This allows us \cite{Hon_JAG2} to investigate
the anti-canonical system $|2F|$ on $Z$
from the structure of the system $\big|2K_S\inv\big|$ on $S$,
and it turns out that $|2F|$ is always
{\em not} composed with the pencil $|F|$
if $Z$ is Moishezon (and $\dim |F|=1$).
Moreover, 
it is not difficult to investigate 
the system $\big|2K_S\inv\big|$.
Thus detailed structure of $Z$  on $4\CP^2$ can be deduced
using $|2F|$.
This is how we obtained the structure of all
Moishezon twistor spaces on $4\CP^2$
as in Proposition \ref{prop:4}.

Although the exact sequence \eqref{ses:002} seems not useful 
for computing $h^0(lF)$ when $n>4$ as above,
it has the following important implication.

\begin{proposition}\cite{C91-2}\label{prop:Ca}
If $S$ is a smooth fundamental divisor,
the inequality $a(Z)\le 1+\kappa\inv(S)$ holds,
where $\kappa\inv(S)=\kappa\big(S,K\inv\big)$,
the anti-Kodaira dimension of $S$.
In particular, if $Z$ is Moishezon, 
$\kappa\inv(S)=2$ holds for any smooth fundamental 
divisor $S$.
\end{proposition}
												
\proof
From \eqref{ses:002} we obtain
$h^0(lF)-h^0\big((l-1)F\big)\le h^0\big(lK_S\inv\big)$ for any $l$.
This means an inequality $\kappa(F)\le \kappa\big(K_S\inv\big)+1$.
Since $a(Z)=\kappa(F)$, the assertion follows.
\proofend

\medskip

Before starting investigation on the structure of 
Moishezon twistor spaces satisfying $\dim |F|=1$,
we discuss some properties of a real irreducible fundamental
divisor on such twistor spaces,
which are mostly due to Kreussler.
The first one is the following.


\begin{proposition}\label{prop:cycle}
\cite[Prop.\,3.6]{Kr98}
Let $Z\to n\CP^2$, $n\ge 4$, be a twistor space
and suppose $\dim |F|=1$.
Then the base locus of the pencil $|F|$
is either a smooth elliptic curve
or a real anti-canonical cycle (see Definition \ref{def:acc}) on 
any real smooth fundamental divisor $S$.
If $Z$ is Moishezon, it
has to be the latter,
and the self-intersection numbers of
components of the cycle in $S$ are
independent of a choice of a smooth fundamental 
divisor $S$.
\end{proposition}

\noindent
{\em Outline of a proof.}
Let $C$ be the unique anti-canonical curve on $S$, and $\psi:S\to S_0=\CP^1\times\CP^1$ the 
birational morphism \eqref{seqbr} in Proposition \ref{prop:PP}.
Then the direct image $\psi_*(C)$ (namely, the image which takes the multiplicities of components
into account) is a real anti-canonical curve
on $S_0$.
By using $h^0\big(K_S\inv\big) = 1$ and 
the real structure on $S_0=\CP^1\times\CP^1$
which is (complex conjugation)$\times$(anti-podal)
as in Proposition \ref{prop:PP},
it can be readily seen that $\psi_*(C)$
does not have a multiple component,
and that it is either a smooth elliptic curve,
or a cycle of rational curves
consisting of two or four components.
Hence the original curve $C$  is also
a smooth elliptic curve or a cycle of rational curves.
Moreover, from the real structure on $S_0$,
if the image $\psi_*(C)$ is reducible,
the real structure on the cycle $\psi_*(C)$
interchanges the opposite components.
This readily implies that 
the same is true for the original curve $C$.

If $C$ is an elliptic curve,
then we obtain $\kappa\inv(S)\le 1$
since
$C^2=8-2n \le 0$ on $S$ as $n\ge 4$.
Hence from Proposition \ref{prop:Ca} we obtain $a(Z)\le 2$ for 
the algebraic dimension.
Therefore if $a(Z)=3$ the curve $C$ has to be
a cycle of rational curves.
For the final assertion, 
by adjunction formula, we have, for any component
$C_i$ of the cycle $C$,
\begin{align*}
\deg K_{C_i}&= (K_{S},{C_i})_{S} +(C_i,
C_i)_{S}\\
&=(-F,C_i)_Z +(C_i,C_i)_{S}.
\end{align*}
This means that
the self-intersection number $C_i^2$ on $S$ is independent
of a choice of a smooth member $S$ of the pencil $|F|$.
\proofend

\medskip
We remark that if $n>4$ and the unique anti-canonical curve $C$ is an elliptic curve,  the meromorphic map $Z\to\CP^1$ induced 
by the pencil $|F|$ gives an algebraic reduction of $Z$.
Generic fibers of the reduction are
rational surfaces by Proposition \ref{prop:PP}.\
These are the  most typical examples of twistor spaces on $n\CP^2$ which satisfy $a(Z)=1$.
We also note that even if $\dim |F|=1$ and
the base locus of the pencil $|F|$ is a real anti-canonical cycle, the twistor space $Z$ is not necessarily Moishezon.
This is because a rational surface $S$ with 
an anti-canonical cycle often
satisfies $\kappa\inv(S)\le 1$,
and so $a(Z)<3$ by Proposition \ref{prop:Ca}.
Moreover, such a twistor space always satisfies $a(Z)=1$
by the result of  \cite{HK17}.

\begin{figure}
\includegraphics{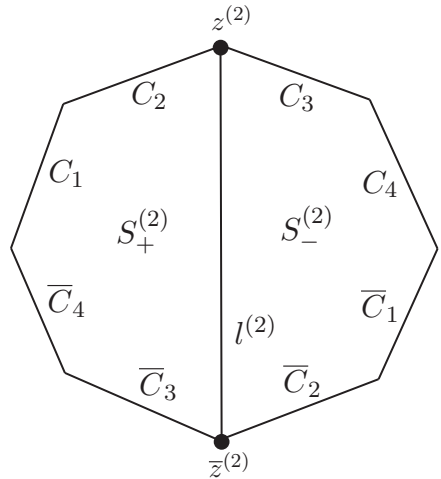}
\caption{
The twistor line $l\upi$ and the reducible 
member $S\upi_++S\upi_-$}
\label{fig:octagon1}
\end{figure}

From Proposition \ref{prop:cycle},
if $n\ge 4$ and $Z$ is a Moishezon twistor space
on $n\CP^2$ that satisfies $\dim |F|=1$,
the base locus of $|F|$ is a real anti-canonical cycle
on a smooth member of $|F|$.
In the remainder of this paper
we always denote this cycle by $C$.
This will be a key geometric object for 
analyzing the structure of these twistor spaces.
For later use, under the presentation 
\begin{align}\label{cycle4}
C=\big(C_1+C_2+\dots + C_k\big) + 
\big(\ol C_1+ \ol C_2+\dots + \ol C_k\big)
\end{align}
as before,
we denote the nodes of $C$ by
\begin{align}\label{z_i}
z\upi = C_i \cap C_{i+1}
\quad{\text{and}}\quad
\ol z\upi = \ol C_i \cap \ol C_{i+1},
\quad
i\in\{1,\dots,k\},
\end{align}
where $C_{k+1}$ means $\ol C_1$
and $\ol C_{k+1}$ means $C_1$.
Furthermore, we denote by $l\upi$ for the twistor
line passing  the points $z\upi$ and $\ol z\upi$.
Thus each twistor line $l\upi$ divides the cycle $C$ into halves.
(See Figure \ref{fig:octagon1} for the case $k=4$ and $i=2$.)
So recalling that  a degree-one divisor
$D$ always contains a unique twistor line
as in Proposition \ref{prop:poon1}
and it satisfies $D+\ol D\in |F|$,
it would be natural to expect that 
for each index $i\in\{1,\dots,k\}$,
there exists a degree-one divisor $D_i$
which satisfies $D_i\cap \ol D_i=l\upi$,
and which divides the cycle $C$ into halves.
The next proposition means that this is correct.
%

\begin{proposition}\label{prop:d1}
  \cite[Prop.\,3.7]{Kr98}
 Under the above situation,
 the pencil $|F|$ has exactly $k$ reducible members.
Each reducible member is real and of the form
  $S\upi_++S\upi_-$ ($1\le i\le k$),
  where $S\upi_+$ and $S\upi_-$ are degree-one divisors
  satisfying $\ol S\upi_+ = S\upi_-$.
  Furthermore, the divisor $S\upi_++S\upi_-$
  splits the cycle $C$ into  halves in the following manner:
  \begin{itemize}
  \item  $S\upi_+\cap S\upi_-=l\upi$, 
  \item the intersections $S\upi_+\cap C$ and $S\upi_-\cap C$ are 
    connected.
  \end{itemize}
\end{proposition}

See Figure \ref{fig:octagon1} for the situation.

\medskip
\noindent
{\em Outline of a proof.}
Fix any $i\in\{1,\dots,k\}$ and take any point
$w\in l\upi\minus\{z\upi,\ol z\upi\}$.
Then $w$ is not on the base curve $C$ of the pencil $|F|$, and hence there exists a unique member
$S_w\in |F|$ which contains the point $w$.
Then $l\upi\subset S_w$ since
otherwise the intersection number $\big(S_w,l\upi\big)_Z$ would be at least 3.
Further the divisor $S_w$ is real since otherwise
$\Bs\,|F|=S_w\cap \ol S_w\supset l\upi$ 
which contradicts
$l\upi\not\subset C=\Bs\,|F|$.
Furthermore if $S_w$ were irreducible,
it would be smooth by Proposition \ref{prop:PP},
and this means $\big(F,\,l\upi\big)_Z = \big(K_{S_w}\inv,l\upi\big)_{S_w}
= \big(C,l\upi\big)_{S_w}$.
But since $C$ is singular at $z\upi$ and $\ol z\upi$,
the last intersection number is at least 4.
This again contradicts $\big(F,l\upi\big)_Z = 2$.
Hence the divisor $S_w$ is reducible.
So we can  write 
$S_w = S\upi_+ + S\upi_-$.
It is easy to see that  $S\upi_+$ and $S\upi_-$ satisfy the required properties.
\proofend

\medskip
In the rest of this paper we use the notations
$S\upi_+$ and $S\upi_-$ to mean the degree-one divisors as in the proposition, 
and write $S\upi:=S\upi_+ + S\upi_-$
for the reducible members of the pencil $|F|$.
(We do not specify which one is $S\upi_+$ or $S\upi_-$ until Section \ref{ss:rb}.)
Hence the pencil $|F|$ has  $k$ reducible members
$$
S\upone, S\uptwo,\dots,S\upk,
$$
and all these are real.
The letter $D$ will be used to mean some other 
divisor on a real irreducible fundamental divisor.

\begin{remark}
{\em 
In Proposition \ref{prop:d1}, it is assumed that $\dim |F|=1$.
But from the above proof, 
the same conclusion holds
for a real pencil of fundamental divisors
whose base locus is a real anti-canonical cycle.
(Here, the real pencil means a pencil
which is invariant under the real structure.) 
This remark is meaningful for the case $n=3$
because in this case we have $\dim |F|=3$
but it is possible to choose a pencil
of the above kind  by fixing a real anti-canonical cycle
on $S$ as in Section \ref{ss:acsfd}.
\proofend
}
\end{remark}

We will also need the following easy but useful property on singularities of non-real irreducible members
of the pencil $|F|$.

\begin{proposition}\label{prop:sm1}
Under the situation as in Proposition \ref{prop:d1},
any non-real \,{\em irreducible}\! fundamental divisor is non-singular
at any point of the real anti-canonical cycle $C$.
\end{proposition}

\proof
Let $S$ be any real irreducible member of the pencil $|F|$ 
and $S'$  any non-real irreducible member of $|F|$.
The former surface $S$ is smooth by Proposition \ref{prop:PP}.
Since $S'|_S=C$, the divisor $S'$ is smooth at any smooth point
of the cycle $C$.
Hence it suffices to show that  $S'$ is smooth 
as well
at any point of $\Sing C
=\{z\upi,\ol z\upi\set 1 \le i\le k\}$.
The member $S'$ does not contain $l\upi$ for any $i\in\{1,\dots,k\}$ since $\Bs\,|F|=S'\cap \ol S'=C$ 
and $l\upi\not\subset C$.
So the intersection number between $S'$
and $l\upi$ at any intersection point makes sense
and it is positive.
Further, since $\{z\upi,\ol z\upi\}\subset S'\cap l\upi$, we have an inequality 
$$
\big(S',l\upi\big)_Z \ge \big(S',l\upi\big)_{z\upi} + \big(S',l\upi\big)_{\ol z\upi}.
$$
On the other hand $\big(S',l\upi\big)_Z = F.\,l\upi = 2$.
This means  
$\big(S',l\upi\big)_{z\upi} = \big(S',l\upi\big)_{\ol z\upi} = 1$.
This implies that the divisor $S'$ is smooth at
the points $z\upi$ and $\ol z\upi$
for any index $i\in\{1,\dots,k\}$.
\proofend

\subsection{Structure of a real irreducible fundamental divisor
which has a real anti-canonical cycle}
\label{ss:strS}
In Section \ref{ss:acsfd} we have seen that 
in the case of $3\CP^2$,
there are two possible structures for  real irreducible 
fundamental divisors  
and they are distinguished by the structure of
the anti-canonical system
(Proposition \ref{prop:S2}).
This was the first possible step for proving the structural result
for  twistor spaces on $3\CP^2$ stated in Proposition \ref{prop:3}.
In this subsection we first see that,
if we use a pluri-anti-canonical system instead of the 
anti-canonical system,
a similar result for real fundamental divisors
holds in the case of $n\CP^2$ for arbitrary $n\ge 3$
when the twistor space is Moishezon and satisfies $\dim |F|=1$.
This result will be the first step for proving our main result
in this paper (Theorem \ref{thm:main1}). 
Next we intensively study the surfaces which have
a structure of double covering over $\CP^2$ by a pluri-anti-canonical map.

%
%
%
From Proposition \ref{prop:PP},
let $S$ be a rational surface which is obtained from $S_0=\qdr$ 
by blowing up $2n$ points as in the proposition.
Suppose that the surface $S$ has a real anti-canonical cycle
$C$ (see Definition \ref{def:acc}) and satisfies $\kappa\inv(S)=2$.
Then as a generalization of Proposition \ref{prop:S2'}, we have the following proposition.

\begin{proposition}
\label{prop:clsf1}
Let $S$ and $C$ be as above.
If $n\ge 3$, 
the surface $S$ satisfies  at least one of the following 
two properties.
\begin{itemize}
\item[\em (i)]
The surface $S$ admits a 
non-trivial $\CC^*$-action which preserves the cycle $C$,
and there is a real pair of components of $C$
whose points are fixed by the $\CC^*$-action.
\item[\em (ii)]
There exists a real effective divisor $D$ on $S$
satisfying $\Supp D=C$, 
$\Bs\,|D|=\emptyset$ and\, $\dim |D|=2$,
such that the associated morphism $\phi:S\to\CP^2$
is generically two-to-one with branch
divisor being a quartic.
\end{itemize}
\end{proposition}

In the proposition, the surface $S$ can satisfy both
of (i) and (ii). 
This issue will be discussed later.

\medskip
\noindent
{\em Proof of Proposition \ref{prop:clsf1}.}
%
%
%
The case $n=3$ follows immediately from 
Proposition \ref{prop:S2'} by taking the anti-canonical cycle $C$ itself
for the effective divisor $D$ in (ii).
We first assert that if the pair $(S,C)$ satisfies the 
property (i) (resp.\,(ii)),
and if $(S',C')$ is a pair which is obtained
from $(S,C)$ by blowing up a real pair
of double points of the cycle $C$,
then the pair $(S',C')$ also satisfies (i) (resp.\,(ii)).
This is obvious for the property (i)
because any intersection points of components
of the cycle $C$ are necessarily fixed points of 
the $\CC^*$-action from the $\CC^*$-invariance of $C$.
If $(S,C)$ satisfies the property (ii),
let $D'$ be the pullback of $D$
under the blowing up $S'\to S$.
Then evidently we have $\Supp D'=C'$, $\Bs\,|D'|=\emptyset$
and $\dim|D'|=2$.
Further,  the  morphism induced by $|D'|$
is exactly the composition of the blowing down
$S'\to S$ and the morphism $S\to\CP^2$ induced by $|D|$.
These mean that $(S',C')$ satisfies the property (ii)
by taking $D'$ for $D$.
Thus we obtain the assertion.

Suppose that 
the surface $S$ admits a sequence 
\begin{align}\label{bsp}
S = S_n \lras S_{n-1} \lras \dots \lras S_3,
\end{align}
where each map is the  blowing down
of a real pair of $(-1)$-components of the
real anti-canonical cycle.
Let $C^3\subset S_3$ be the image of the cycle $C$
on $S$.
This is also a real anti-canonical cycle.
Therefore by Proposition \ref{prop:S2'}, the  pair $(S_3, C^{3})$ satisfies at least one of (i) and (ii).
Hence by what we have just proved, the 
original pair $(S,C)$ also satisfies the property (i) or (ii).
So assume that the contraction 
process of $(-1)$-components
of the cycle always stops 
before reaching the surface $S_3$,
and let $(S_{\nu}, C^{\nu})$ be the resulting pair.
We have ${\nu}>3$,  and 
the cycle $C^{\nu}$ has no $(-1)$-component.

First suppose that all components of the cycle $C^{\nu}$ have negative self-intersection numbers.
Since $C^{\nu}$ has no $(-1)$-component,
this implies that
 the self-intersection number of any component of $C^{\nu}$
is at most $(-2)$.
If all components have self-intersection number $(-2)$,
the anti-canonical cycle $C^{\nu}$ is nef,
and   the self-intersection number of the cycle $C^{\nu}$ is zero.
Hence, by \cite[Theorem 3.4]{Sa84},
we obtain $\kappa\inv( S_{\nu})\le 1$.
Since the anti-Kodaira dimension does not increase
by  blowing up, we  obtain $\kappa\inv(S)\le 1$. This contradicts the assumption $\kappa\inv(S)=2$.
Hence some component of $C^{\nu}$
has self-intersection number which is less than $(-2)$.
Then this time $C^{\nu}$ itself becomes
the negative part of the Zariski decomposition
of $C^{\nu}$ since the intersection matrix
formed by components of the cycle $C^{\nu}$
becomes negative definite.
Again by \cite[Theorem 3.4]{Sa84}, this implies $\kappa\inv(S)=0$.
This is again a contradiction, and
we obtain that some component of the cycle $C^{\nu}$ has
a non-negative self-intersection number.

From the real structure, there are at least two such components.
If some of the self-intersection numbers is positive,
by moving the two components, we obtain that the surface $S_{\nu}$ satisfies $\dim |K\inv| > 2$.
This cannot happen by Proposition \ref{prop:LB}
since $\nu>3$.
Hence the above non-negative self-intersection number is zero,
and it implies $\dim |K\inv|=2$ for the surface $S_{\nu}$.
Then again by Proposition \ref{prop:LB}, the pair
$(S_{\nu},C^{\nu})$ admits a non-trivial $\CC^*$-action
which has a pair of components whose points are fixed.
In particular, the pair
$(S_{\nu},C^{\nu})$ satisfies the property (i) in the present
proposition.
This implies that the original pair $(S,C)$ also satisfies the same property.
\proofend

\medskip

The structure of twistor spaces on $n\CP^2$ which have the surface $S$
satisfying the property (i) in
Proposition \ref{prop:clsf1} as a real  fundamental divisor
was studied in detail in the paper \cite{Hon_Cre1}.
At the end of this subsection, for reader's convenience, 
we will briefly describe how the structure of such twistor spaces
was analyzed.

As in the proof of the proposition, if the surface $S$ can reach 
a surface $S_3$ (which satisfies $K^2=2$)
as in \eqref{bsp} by repeating 
pairwise blowing down of $(-1)$-components of the real anti-canonical cycle,
and if the resulting pair $(S_3,C^3)$ satisfies the property (ii)
in the proposition,
then the divisor $D$ on $S$ is obtained
as the pull-back of the real anti-canonical cycle $C^3$ on $S_3$.
But from the proof,
there is still a possibility that
the blowing down process \eqref{bsp} always
stops at a surface $S_{\nu}$ with $\nu>3$,
but nevertheless $(S,C)$ satisfies the property (ii).
If such a situation would happen, 
 the surface $S$ could  be a new type of surface
having the double covering structure as in (ii).
Namely such a surface would not be obtained as
a blowing up of $(S_3,C^3)$  by the blowing up 
of the above kind.
To deny this possibility, as in the proof, let $C^{\nu}$ be
the image of $C$ into $S_{\nu}$, $\nu>3$.
This is a real anti-canonical cycle of $S_{\nu}$
and has no $(-1)$-component.
Of course, we are allowing the situation $(S,C) = (S_{\nu},C^{\nu})$.

\begin{proposition}\label{prop:clsf2}
Under the above situation, 
the surface $S$ cannot satisfy the property
(ii) in Proposition \ref{prop:clsf1}.
\end{proposition}

\proof
As in the proof of Proposition \ref{prop:clsf1},
the surface $S_{\nu}$ satisfies $\dim|K\inv|=2$.
By Proposition  \ref{prop:LB}, this implies that
$S_{\nu}$ is obtained from $S_0=\qdr$ by 
blowing up $2\nu\,(\ge 8)$ points whose half belong to a $(1,0)$-curve
or a $(0,1)$-curve.
We denote $C_0$ for the strict transform of this curve
into $S_{\nu}$.
Points on the curves $C_0$ and $\ol C_0$ are fixed by 
the non-trivial $\CC^*$-action on $S_{\nu}$.
From the choice of the sequence \eqref{bsp},
the original surface $S$
is obtained from the surface $S_{\nu}$ by repetition of 
blowing up at double points of the
real anti-canonical cycle $C^{\nu}$ on $S_{\nu}$.
Hence $S$ also admits a non-trivial $\CC^*$-action.
We denote the strict transforms of $C_0$ and $\ol C_0$
into $S$ by the same letters.
As above we have $\dim |K\inv|=2$ on $S_{\nu}$,
and the curves $C_0$ and $\ol C_0$ are
fixed components of the anti-canonical system of $S_{\nu}$.
Hence the cycle $C^{\nu}$ includes these curves
as components.
Therefore so is the original real anti-canonical cycle $C$ on $S$.
Moreover, since the cycle $C^{\nu}$ on $S_{\nu}$
is $\CC^*$-invariant, so is the cycle $C$ on $S$.

Suppose that $S$ satisfies the property (ii)
in  Proposition \ref{prop:clsf1}.
Then, since the divisor $D$ on $S$ is supposed to consist
of components of $C$ as in (ii) of Proposition \ref{prop:clsf1}, 
the divisor $D$ is also $\CC^*$-invariant.
Hence the morphism $\phi:S\to\CP^2$ induced by the system $|D|$ is  $\CC^*$-equivariant.
Let $l_1$ be the line which satisfies $\phi\inv(l_1)=D$.
This is real and $\CC^*$-invariant, and
we have $\phi(C)=l_1$ since $\Supp D = C$
as supposed.

Assume that $\CC^*$ acts trivially on the line $l_1$.
Then by the equivariancy of $\phi$,
we have $\phi(C_0)=\phi(\ol C_0)=l_1$.
Besides the line $l_1$ there exists precisely one fixed point
of the $\CC^*$-action on $\CP^2$,
and 
the closure of a generic orbit in $\CP^2$ is a line
through this point,
and it intersects the line $l_1$.
But the closure of a generic orbit on $S$ passes the two components $C_0$ and $\ol C_0$
of $C$.
This contradicts $\CC^*$-equivariancy of 
the morphism $\phi$.
Therefore $\CC^*$ has to act non-trivially on  the line $l_1$.
Since $\phi$ preserves the real structure, 
this means that, in homogeneous coordinates, 
%
the induced $\CC^*$-action on $\CP^2$
may be supposed to be of the form
\begin{align}\label{Cstar2}
(x:y:z)\stackrel{t}\longmapsto
(tx:t\inv y:z),\quad t\in\CC^*,
\end{align}
where the locus $\{z=0\}$ is the line $l_1=\phi(C)$ and 
the induced real structure 
on $\CP^2$ is given by $(x:y:z)\mapsto
(\ol y:\ol x:\ol z)$.
In particular, the closure of a generic orbit of the $\CC^*$-action on $\CP^2$
is a conic passing the two points $(1:0:0)$ and $(0:1:0)$.
Also, since the action \eqref{Cstar2}
has only isolated fixed points,
 the components $C_0$ and $\ol C_0$
 of the cycle $C$
are mapped to points by $\phi$, and they are distinct.
By changing $C_0$ and $\ol C_0$ if necessary,
we may suppose $\phi(C_0) = (1:0:0)$
and $\phi(\ol C_0) = (0:1:0)$.

Let $S\to T\to \CP^2$ be the Stein factorization of 
the morphism $\phi:S\to\CP^2$.
The morphism $T\to \CP^2$ is finite and of degree two.
Since $\phi$ is $\CC^*$-equivariant, so is the morphism $T\to\CP^2$.
Further, the branch divisor of $\phi:S\to\CP^2$ equals
that of $T\to\CP^2$ and hence it is a quartic by the assumption in (ii).
Therefore, from the $\CC^*$-action on $\CP^2$
as in \eqref{Cstar2},
the branch of $T\to\CP^2$ consists of two $\CC^*$-invariant conics,
and these conics touch each other at the two points
$(1:0:0)$ and $(0:1:0)$.
Let $\tilde T\to T$ be the minimal resolution of $T$.
From the
universal property of the minimal resolution,
the morphism $S\to T$ factors through as $S\to\tilde T\to T$,
and the morphism $S\to\tilde T$
is $\CC^*$-equivariant and birational.
Then from the concrete form of the branch divisor,
the surface $\tilde T$ satisfies $K^2 = 2$.
Further, since the morphism $S\to\tilde T$ is 
birational and $\CC^*$-equivariant, the image of the cycle $C$ to $\tilde T$ is a $\CC^*$-invariant anti-canonical cycle on $\tilde T$.

We show that  
the morphism $S\to\tilde T$ is a repetition of
contraction of real pair of $(-1)$-components
in the real anti-canonical cycle.
From $\CC^*$-equivariance of the morphism $S\to\tilde T$, each factor
of it can blow up only a fixed point of 
the $\CC^*$-action.
If some factor blows up a real pair of smooth points of the real anti-canonical cycle, 
the exceptional curves of such points
have to be included in the divisor $D$,
because $\phi\inv(l_1)=D$ and the morphism
$\phi$ and the composition $S\to \tilde T\to T\to \CP^2$ are equal.
%
%
This implies that the support of the divisor $D$ is not a cycle,
which contradicts the assumption that
$\Supp\,D = C$ in (ii) of Proposition \ref{prop:clsf1}.
Therefore the morphism $S\to\tilde T$ is a repetition of
contractions of a real pair of $(-1)$-components
in the real anti-canonical cycle $C$ on $S$.
By taking $\tilde T$ as the surface $S_3$,
this implies that the morphism $S\to \tilde T$
provides the sequence \eqref{bsp},
which contradicts our assumption that the contraction process always stops 
before reaching the surface $S_3$.
Hence the surface $S$ cannot fulfill
the property (ii) in Proposition \ref{prop:clsf1}.
\proofend\medskip

Thus a pair $(S,C)$ satisfying the property (ii)
in Proposition \ref{prop:clsf1} can always be obtained from
the one in the case $n=3$ (namely the pair
$(S_3,C^{3})$ in the previous notation)
by repetition of blowing up at
a real pair of double points of the real anti-canonical cycle.
The effective divisor $D$ on $S$ which induces
the degree-two morphism $\phi:S\to\CP^2$
is obtained from the real anti-canonical cycle $C^{3}$ on $S_3$
as the total transform.
Therefore we have the factorization
of the morphism $\phi:S\to\CP^2$ into
the composition of a sequence
\eqref{bsp} and the anti-canonical map
from the surface $S_3$:
\begin{align}\label{diagram:fac1}
 \xymatrix{ 
S \ar^{}[r] \ar[dr]_{\phi} &S_3 \ar^{}[d]\\
&\CP^2
 }
\end{align}
As above the morphism $S\to S_3$ is a repetition
of contractions of $(-1)$-components in the 
real anti-canonical cycle.
From this, for a pair $(S,C)$ satisfying the property (ii)
in Proposition \ref{prop:clsf1},
the type can be defined as the type of 
the pair $(S_3,C^3)$,
which was one of $\Azero,\Aone,\Atwo$ 
and $\Athree$ as seen in Section \ref{ss:acsfd}.
This definition of the type for $(S,C)$
does not depend on the choice of the 
sequence \eqref{bsp} from dependence on the number of 
the components of $C$ o the type.
Also, whether a pair $(S,C)$ 
 satisfies the property
(ii) in the proposition is determined from 
the self-intersection numbers of components of the cycle $C$ in $S$.
(The same is true for the property (i).)
We also note that from Proposition \ref{prop:clsf2},
among these four types of surfaces,
only type A$_3$ admits a non-trivial $\CC^*$-action.
In other words, {\em the surface $S$ (and $C$) in Proposition
\ref{prop:clsf1} can satisfy both of the two properties (i) and (ii)
if and only if it is of type $\Athree$.}

%
%

Since the anti-canonical system
of $S$ consists of a single member $C$ when $n>3$,
the type is well-defined not only for the pair 
$(S,C)$ but also for the surface $S$ itself.
Since the type of the surface is determined from
the self-intersection numbers of the components as above,
and since these numbers are independent of a choice of a (real) irreducible member of $|F|$ as in Proposition \ref{prop:cycle},
 we may define the type of twistor spaces
as follows.

\begin{definition}
\label{def:type}
Suppose $n>3$ and let $Z$ be a twistor space on $n\CP^2$.
We call $Z$ is of type 
$\Azero,\Aone,\Atwo$ or $\Athree$
if $Z$ contains a real irreducible fundamental divisor $S$ which is of type 
$\Azero,\Aone,\Atwo$ or $\Athree$ respectively.
\end{definition}

These twistor spaces are of our principal interest
in this paper.

We next discuss existence of these twistor spaces.
First, from the realization of a surface of type $\Athree$ in the case $n=3$
as a 6 points blowing up of $S_0=\qdr$ as described in Figure \ref{fig:bu_cfg1},
it is easy to see that a surface $S_3$ of type $\Athree$ is obtained 
from a toric surface with real structure by small ($\CC^*$-equivariant) deformation preserving the real structure.
From this and the diagram \eqref{diagram:fac1},
for any $n\ge 3$, a surface $S$ of type $\Athree$ satisfying $K^2 = 8-2n$
can be obtained as a small ($\CC^*$-equivariant) deformation
of a toric surface with real structure, regardless
of the realization $S\to S_3$ 
in 
the diagram \eqref{diagram:fac1}.
Since the toric surface is realized
as a real fundamental divisor on the twistor space of a Joyce metric
on $n\CP^2$ \cite{F00},
in the same way to Section 5.3 in \cite{Hon_Cre2},
a deformation theoretic argument shows that
there exists a twistor space on $n\CP^2$ which has
a surface $S$ of type $\Athree$ as a real fundamental divisor,
regardless of the realizing sequence $S\to S_3$
for $S$ in \eqref{diagram:fac1}.
Next, for a surface $S$ of type $\Atwo$,
from the configuration of the 6 points on $S_0=\qdr$ as in 
Figure \ref{fig:bu_cfg1},
the surface $S$ is always obtained from a surface of type $\Athree$
by small deformation preserving the real structure, for any choice of the blowing up sequence
$S\to S_3$ in \eqref{diagram:fac1}.
Hence again by deformation theoretic argument
we obtain the existence of twistor spaces on $n\CP^2$ which 
have a surface of type $\Atwo$ as a real fundamental divisor, for any choice of the blowing up sequence
$S\to S_3$ in \eqref{diagram:fac1}.
Similarly, surfaces of type $\Aone$ and $\Azero$ are
realized as real fundamental divisors on twistor spaces on $n\CP^2$ respectively,
for any realization from $(S_3,C^3)$.

In brief, similarly to the adjacent relations
\eqref{adr}
for the surfaces,
we have the following diagram
for mutual relationship for the twistor spaces
\begin{align}\label{adr2}
{\rm A}_0 \lras {\rm A}_1 \lras {\rm A}_2 \lras {\rm A}_3 \lras {\text{Joyce}},
\end{align}
where the arrows still indicate specialization.

Next we investigate
these surfaces (satisfying the property (ii)
in Proposition \ref{prop:clsf1}) more closely.
The number of components of the 
real anti-canonical cycle $C$ is readily seen
to be given as in Table \ref{table5},
depending on the type.%
\begin{table}[h]
\begin{center}
\begin{tabular}{c||c|c|c|c}
type & $\Azero$ & $\Aone$ & $\Atwo$ & $\Athree$\\
\hline
$2k$ & $2(n-2)$   & $2(n-1)$ & $2n$ & $2(n+1)$
\end{tabular}
\end{center}
\caption{The number of components of the cycle $C$.}
\label{table5}
\end{table}
As before, we are using the letter $k$ to mean the half of the
number of components of $C$.
We note that in the case of the twistor spaces
of Joyce metrics on $n\CP^2$,
the real anti-canonical cycle $C$ consists
of $2(n+2)$ components.

Since the divisor $D$
 inducing the degree-two morphism
$\phi:S\to\CP^2$ is real and 
satisfies $\Supp\,D = C$,
we may write 
\begin{align}\label{cycle32}
D=\sum_{i=1}^k d_iC_i + 
\sum_{i=1}^k d_i\ol C_i,
\end{align}
where any $d_i$ is positive.
For later purpose we define an integer $d$ by
\begin{align}\label{d}
d=\max\{d_1,d_2,\dots, d_k\}.
\end{align}
Then the divisor $D$ is a sub-divisor of $dC
\in \big|dK_S\inv\big|$.

Next we give the list of 
all the surfaces in the case $n=3,4,5$,
by presenting the sequences of the self-intersection numbers
of components of the anti-canonical cycle $C$
arranging in order according to 
the presentation $C = (C_1 + \dots + C_k)
+ (\ol C_1 + \dots + \ol C_k)$. 

When $n=3$, from the result in Section \ref{ss:acsfd},
the type uniquely 
determines not only the number of components
of the cycle $C$ but also
the self-intersection numbers of the components.
So each type consists of one possibility.
Moreover 
the divisor $D$ is the anti-canonical cycle $C$ itself.
They are summarized as in Table \ref{table1}.
\begin{table}[h]
\begin{center}
\begin{tabular}{c||c|c}
type & the sequence for $D\,(=C)$ & $d$\\
\hline\hline
$\Azero$  & $\bold{-1;-1}$ & 1\\
\hline
$\Aone$ & $\bold{-1},-2;\bold{-1},-2$ & 1  \\
\hline
$\Atwo$  & $\bold{-1},-2,-2;\bold{-1},-2,-2$  & 1    \\ 
\hline
$\Athree$  & $\bold{-1},-2,-2,-2;\bold{-1},-2,-2,-2$   & 1   
\end{tabular}
\end{center}
\caption{Concrete forms of the divisor $D\,(=C)$ in the case $n=3$}
\label{table1}
\end{table}
For example, in the case of type $\Atwo$,
the cycle $C$ consists of 6 components,
and the sequence $\bold{-1},-2,-2;\bold{-1},-2,-2$
is the self-intersection numbers of these components put in order.
The bold figures represent special components
(the line components) which will be soon defined.


By the diagram \eqref{diagram:fac1},
all possibilities for the surface $S$ in the case $n=4$
can be obtained by considering all blowing ups
at double points of the cycle $C$ in the case $n=3$,
and we obtain the list as in Table \ref{table:2}.
\begin{table}[h]
\begin{center}
\begin{tabular}{c||c|c}
Type & the sequence for $D$ & $d$\\
\hline\hline
\rule[0mm]{0mm}{4mm}
 $\Azero$   & $\bold{-3}^1,-1^2;\bold{-3}^1,-1^2$ & 2\\
\hline\hline
\rule[0mm]{0mm}{4mm}
$\Aone$ & $\bold{-2}^1,-1^2,-3^1;\bold{-2}^1,-1^2,-3^1$   & 2\\
\hline\hline
\rule[0mm]{0mm}{4mm}
\multirow{2}{*}{$\Atwo$}  & $\bold{-2}^1,-1^2,-3^1 ,-2^1;
\bold{-2}^1,-1^2,-3^1 ,-2^1$     & 2 \\ \cline{2-3}
\rule[0mm]{0mm}{4mm}
 & $\bold{-1}^1,-3^1,-1^2,-3^1; \bold{-1}^1,-3^1,-1^2,-3^1 $   & 2  \\ 
\hline\hline
\rule[0mm]{0mm}{4mm}
\multirow{2}{*}{$\Athree$}  & $\bold{-2}^1,-1^2,-3^1,-2^1,-2^1;\bold{-2}^1,-1^2,-3^1,-2^1,-2^1  $ & 2     \\ \cline{2-3}
\rule[0mm]{0mm}{4mm}
  & $\bold{-1}^1,-3^1,-1^2, -3^1,-2^1;\bold{-1}^1,-3^1,-1^2, -3^1,-2^1  $  & 2   \\ 
\end{tabular}
\end{center}
\caption{Concrete forms of the divisor $D$ in the case $n=4$}
\label{table:2}
\end{table}
In the table the digits
on upper right 
represent the multiplicities for
the components in the divisor $D$.
So for instance, in the case of type $\Aone$,
 the divisor $D$ can be written as
$$
D = C_1 + 2 C_2 + C_3 + \ol C_1 + 2 \ol C_2 + \ol C_3,
$$
where the self-intersection numbers of
the components are $-2,-1,-3,-2,-1,-3$ 
respectively in order.
Note that by reality the multiplicities of 
$C_i$ and $\ol C_i$ are always equal.
At this stage the types $\Atwo$ and  $\Athree$ include	
multiple situations.
On the other hand the number $d$ is always two. 

The list in the case $n=5$ is again easily
obtained from that of the case $n=4$
and it is as in Table \ref{table:3}.
As one sees, the number of all possibilities is 
quite large.
Also it turns out the number $d$ is not determined from $n$ nor type.

\begin{table}[h]
\begin{center}
\begin{tabular}{c||c|c}
Type & the sequence for $D$ & $d$\\
\hline\hline
\rule[0mm]{0mm}{4mm}
$\Azero$   & $\bold{-4}^1,-1^3,-2^2;\bold{-4}^1,-1^3,-2^2,$ & 3\\
\hline\hline
\rule[0mm]{0mm}{4mm}
\multirow{3}{*}{$\Aone$}   & $\bold{-3}^1,-1^3,-2^2,-3^1;
\bold{-3}^1,-1^3,-2^2,-3^1$   & 3\\ \cline{2-3}
\rule[0mm]{0mm}{4mm}
& $\bold{-2}^1,-2^2,-1^3,-4^1;
\bold{-2}^1,-2^2,-1^3,-4^1$  & 3 \\ \cline{2-3}
\rule[0mm]{0mm}{4mm}
 & $\bold{-3}^1,-1^2,-4^1,-1^2;
\bold{-3}^1,-1^2,-4^1,-1^2$  & 2 \\
\hline\hline
\rule[0mm]{0mm}{4mm}
\multirow{8}{*}{$\Atwo$} & $\bold{-3}^1,-1^3,-2^2 ,-3^1, -2^1;
\bold{-3}^1,-1^3,-2^2 ,-3^1, -2^1$ & 3\\ \cline{2-3}
\rule[0mm]{0mm}{4mm}
 & $\bold{-2}^1,-2^2,-1^3,-4^1,-2^1;
\bold{-2}^1,-2^2,-1^3,-4^1,-2^1$ & 3\\ \cline{2-3}
\rule[0mm]{0mm}{4mm}
 & $\bold{-2}^1,-1^2,-4^1,-1^2,-3^1;
\bold{-2}^1,-1^2,-4^1,-1^2,-3^1$ & 2\\ \cline{2-3}
\rule[0mm]{0mm}{4mm}
  & $\bold{-3}^1,-1^2,-3^1,-3^1,-1^2;
\bold{-3}^1,-1^2,-3^1,-3^1,-1^2$ & 2\\ \cline{2-3}
\rule[0mm]{0mm}{4mm}
 & $\bold{-2}^1,-1^2,-4^1,-1^2,-3^1;
\bold{-2}^1,-1^2,-4^1,-1^2,-3^1$ & 2\\ \cline{2-3}
\rule[0mm]{0mm}{4mm}
  & $\bold{-1}^1,-4^1,-1^3,-2^2,-3^1;
\bold{-1}^1,-4^1,-1^3,-2^2,-3^1$ & 3\\ \cline{2-3}
\rule[0mm]{0mm}{4mm}
  & $\bold{-1}^1,-3^1,-2^2,-1^3,-4^1;
\bold{-1}^1,-3^1,-2^2,-1^3,-4^1$ & 3\\ \cline{2-3}
\rule[0mm]{0mm}{4mm}
  & $\bold{-2}^1,-3^1,-1^2,-4^1,-1^2;
\bold{-2}^1,-3^1,-1^2,-4^1,-1^2$ & 2\\ \cline{2-3}
\hline\hline
\rule[0mm]{0mm}{4mm}
\multirow{10}{*}{$\Athree$}  & $\bold{-3}^1,-1^3,-2^2,-3^1,-2^1,-2^1;
\bold{-3}^1,-1^3,-2^2,-3^1,-2^1,-2^1 $  & 3   \\ 
\cline{2-3}
\rule[0mm]{0mm}{4mm}
  & $\bold{-2}^1,-2^2,-1^3, -4^1,-2^1,-2^1;
\bold{-2}^1,-2^2,-1^3, -4^1,-2^1,-2^1  $  & 3  \\ 
\cline{2-3}
\rule[0mm]{0mm}{4mm}
  & $\bold{-2}^1,-1^2,-4^1, -1^2,-3^1,-2^1;
\bold{-2}^1,-1^2,-4^1, -1^2,-3^1,-2^1 $   & 2  \\ 
\cline{2-3}
\rule[0mm]{0mm}{4mm}
  & $\bold{-2}^1,-1^2,-3^1, -3^1,-1^2,-3^1;
\bold{-2}^1,-1^2,-3^1, -3^1,-1^2,-3^1 $  & 2  \\ 
\cline{2-3}
\rule[0mm]{0mm}{4mm}
  & $\bold{-3}^1,-1^2,-3^1, -2^1,-3^1,-1^2;
\bold{-3}^1,-1^2,-3^1, -2^1,-3^1,-1^2 $ & 2  \\ 
\cline{2-3}
\rule[0mm]{0mm}{4mm}
  & $\bold{-2}^1,-1^2,-4^1, -1^2,-3^1,-2^1;
\bold{-2}^1,-1^2,-4^1, -1^2,-3^1,-2^1$   & 2  \\ 
\cline{2-3}
\rule[0mm]{0mm}{4mm}
  & $\bold{-1}^1,-4^1,-1^3, -2^2,-3^1,-2^1;
\bold{-1}^1,-4^1,-1^3, -2^2,-3^1,-2^1$   & 3  \\ 
\cline{2-3}
\rule[0mm]{0mm}{4mm}
  & $\bold{-1}^1,-3^1,-2^2, -1^3,-4^1,-2^1;
\bold{-1}^1,-3^1,-2^2, -1^3,-4^1,-2^1$   & 3  \\ 
\cline{2-3}
\rule[0mm]{0mm}{4mm}
  & $\bold{-1}^1,-3^1,-1^2, -4^1,-1^2,-3^1;
\bold{-1}^1,-3^1,-1^2, -4^1,-1^2,-3^1$ & 2  \\ 
\cline{2-3}
\rule[0mm]{0mm}{4mm}
  & $\bold{-2}^1,-3^1,-1^2, -3^1,-3^1,-1^2;
\bold{-2}^1,-3^1,-1^2, -3^1,-3^1,-1^2$  & 2   
\end{tabular}
\end{center}
\caption{Concrete forms of the divisor $D$ in the case $n=5$}
\label{table:3}
\end{table}

Next we define the special components
in the anti-canonical cycle $C$ 
which are denoted by bold numbers in 
the above tables.
In the case $n=3$, 
the anti-canonical cycle $C$
always has exactly one pair of $(-1)$-components,
and all other components of $C$ 
are, if any, $(-2)$-curves.
Further, as we have seen in Section \ref{ss:acsfd},
the degree-two morphism $\phi:S\to\CP^2$
(induced by $|D| = |C| = \big|K_S\inv\big|$) maps
the $(-1)$-curves  isomorphically 
to the same line $\bm l$, and  $(-2)$-curves
 to points (denoted by $\bm q$ and $\bm {\ol q}$)
 on the line $\bm l$.
Because the morphism $\phi:S\to\CP^2$
for general $n$ factors through the anti-canonical
morphism in the case $n=3$ as in the diagram \eqref{diagram:fac1}, 
among components of the anti-canonical cycle $C$,
there exists a unique pair of components
which are mapped isomorphically to the line $\bm l$
in $\CP^2$.
These are mutually conjugate.

\begin{definition}\label{def:lc}{\em
We call these two components of the 
real anti-canonical cycle $C$ as the
{\em line components}.
By a cyclic permutation, we always assume that 
$C_1$ and $\ol C_1$ are
these components.
Also we write $\bm l$ for 
the line which is the image of the line components
under the morphism $\phi:S\to\CP^2$.
This is real.\proofend
}
\end{definition}

The line components will play a special
role in the reminder of this paper.
The line $\bm l$ will be the ridge of the scroll
of planes.
Note that
{\em the divisor $D$ 
has the line components (i.e.\,$C_1$ and $\ol C_1$) as 
multiplicity-one components\,} because
this is the case when $n=3$
and $D$ is the pull-back of the cycle  $C$
in the case $n=3$.
The bold figures in Tables
\ref{table1}--\ref{table:3} are exactly the self-intersection numbers of the line components.
If we write the cycle $C$ as
\begin{align}\label{cycle22}
\big(C_1+C_2+\dots + C_k\big) + 
\big(\ol C_1+ \ol C_2+\dots + \ol C_k\big) 
\end{align}
as before, the two chains 
$C_2 + \dots + C_k$ and $\ol C_2 + \dots + \ol C_k$
are contracted to points by the morphism $\phi:S\to\CP^2$.
We write these points by $\bm q$ and $\bm {\ol q}$.
Further, from Section \ref{ss:acsfd}, relative to the real line $\bm l$, the branch quartic satisfies the following properties.
\begin{itemize}
\item
In the case of type $\Azero$, the branch quartic is tangent to the line $\bm l$ at the two points $\bm q$ and $\bm {\ol q}$.
\item
In the case of type $\Aone$
(resp.\,$\Atwo$ and $\Athree$), the branch quartic 
has $\Aone$-singularities
(resp.\,$\Atwo$- and $\Athree$-singularities) at 
the points $\bm q$ and $\bm {\ol q}$.
\end{itemize}
See Figure \ref{fig:bq} for these.

The intersection numbers between
the effective divisor $D$ and  components
of the anti-canonical cycle $C$ can be readily 
calculated as follows.

\begin{proposition}\label{prop:bn1}
Let $S,C, C_1$ and $D$ be as above.
Then the intersection numbers between
the divisor $D$ and components of $C$ are given by
\begin{align}\label{int0011}
D.\,C_i = D.\,\ol C_i= 
\begin{cases}
1 & {\text{if $i=1$}},\\
0 & {\text{otherwise.}}
\end{cases}
\end{align}
Moreover, $D^2=2$. In particular the divisor
$D$ is nef and big.
\end{proposition}

\proof
The degree-two morphism $\phi:S\to\CP^2$
is induced by the linear system $|D|$ which
is base point free.
If $i\neq 1$, the components $C_i$ and $\ol C_i$ are mapped to points by the morphism.
Hence we have $D.\,C_i = D.\,\ol C_i= 0$
for these $i$.
Therefore, recalling that the multiplicities of 
$C_1$ and $\ol C_1$ in $D$ are one,
we have 
$
D^2 = D.\,(C_1+\ol C_1).
$
On the other hand, since $D$ is the pull-back 
of the anti-canonical cycle in the case $n=3$
by the sequence of blow-ups,
and since the self-intersection 
number of the last cycle
is $2\,(=8-6)$, we have $D^2=2$.
Therefore $D.\,C_1=D.\,\ol C_1=1$
as $D.\,C_1=D.\,\ol C_1$ from the real structure.
\proofend

\medskip
We will use the intersection numbers \eqref{int0011}
countless times and call them 
the {\em basic intersection numbers}.

Next we prove some vanishing
property for cohomology groups on irreducible 
fundamental divisors.
This will be used in the next section for analyzing the structure
of the present twistor spaces.
\begin{proposition}\label{prop:van8}
Let $Z\to n\CP^2$ be any one of the twistor spaces in Definition \ref{def:type}, and
$S$ any irreducible member of the pencil $|F|$
which is not necessarily real.
Then we have
$$H^q(S,D-C)=
\begin{cases}
0 &  q>0,\\
\CC & q=0.
\end{cases}$$
\end{proposition}

\proof
If the irreducible divisor $S$ is real, this follows immediately from Proposition \ref{prop:bn1} by applying
Kodaira-Ramanujan vanishing theorem to the divisor $D$ since $C=K_S\inv$ and $D$ is nef and big.
This argument works without any change as long as
$S$ is smooth.
When $S$ is not necessarily smooth, 
we need some circuitous argument as follows.
We recall that $S$ is smooth at least at any point
of $C$ by Proposition \ref{prop:sm1} and 
the self-intersection numbers of components
of the cycle $C$ in $S$ are independent of a choice
of the irreducible member $S$ of $|F|$ 
from the adjunction formula.
Therefore as in the case for a real member,
we can successively blow down real
pairs of components of $C$
$(n-3)$ times.
As before let
\begin{align}\label{sbu}
S=S_n\stackrel{\psi_n}\lras
S_{n-1}\stackrel{\psi_{n-1}}\lras
\dots \stackrel{\psi_{5}}\lras S_4
\stackrel{\psi_{4}}\lras S_3
\end{align}
be the sequence of the pairwise blowing down obtained this way.
Let $D^j$ be the divisor on $S_j$
which is inductively obtained from 
the divisor $D=D^n$ on $S$ by letting
$D^{j-1}  =(\psi_j)_*D^{j}$.
In other words, $D^j = \psi_j^*D^{j-1}$
for any index $j$.
Let $C^j\subset S_j$ be the 
anti-canonical cycle which is the image
of the original cycle $C=C^n$ on $S$,
and we write $E^j$ and $\ol E^j$ 
for the pair of $(-1)$-curves on $S_j$
which are contracted by $\psi_j$.
Note that $D^{3} = C^{3}$.
Then we  have the relation
\begin{align}
\psi^*_j(D^{j-1} - C^{j-1}) &=
D^j - \psi_j^* C^{j-1}\notag\\
& = D^j - C^j - E^j - \ol E^j
\label{pbck1}
\end{align}
and 
\begin{align}
\big(D^{j} - C^{j}, E^j\big)_{S_j}
&=
\big(D^{j},E^j\big)_{S_j} - \big(C^{j}, E^j\big)_{S_j}\notag\\
&= 0 - (2-1) = -1
\label{pbck2}
\end{align}
for any index $j$.
From these we obtain, for any $q\ge 0$,
\begin{align*}
H^q\big(S_{j}, D^{j}-C^{j}\big)&\simeq 
H^q\big(S_j, D^j-C^j-E^j - \ol E^j\big)
\quad(\because \eqref{pbck2})\\
&\simeq 
H^q\big(S_{j-1},\psi_j^*(D^{j-1}-C^{j-1})\big)
\quad(\because \eqref{pbck1})\\
&\simeq
H^q\big(S_{j-1},D^{j-1}-C^{j-1}\big).
\end{align*}
Connecting these isomorphisms we get an isomorphism
$$
H^q\big(S,D-C\big)\simeq H^q\big(S_3, D^{3} - C^{3}\big)
$$
for any $q\ge 0$.
Hence, as $D^{3} = C^{3}$, we obtain 
$H^q(S,D-C)\simeq H^q(S_3,\ms O_{S_3})$.
For the right-hand side, pulling back by 
the above sequence \eqref{sbu},
we get $H^q(S_3,\ms O_{S_3})
\simeq H^q(S,\ms O_S)$ for any $q\ge 0$.
If $q>0$, this is zero by Proposition \ref{prop:S03}.
\proofend

\medskip
Finally in this subsection, as we promised,
we briefly explain how one can obtain the structure
of twistor spaces on $n\CP^2$ 
which admit a real fundamental divisor $S$
that satisfies the property (i) in 
Proposition \ref{prop:clsf1},
following \cite{Hon_Cre1}.

So let $Z\to n\CP^2$, $n\ge 3$, be a twistor space,
$S$ such a divisor on $Z$,
and $C$ the real anti-canonical cycle on $S$.
(In the following argument
 we do not need to assume that the $\CC^*$-action on $S$
 extends to that on $Z$;
namely we do not need to assume that $Z$ itself has a non-trivial
$\CC^*$-action.
Also we do not need to assume that $Z$ is Moishezon.)
We write the cycle $C$ in the form \eqref{cycle22}.
By renumbering the indices, we may suppose that $C_1$
and $\ol C_1$ are the components whose points
are fixed by 
the $\CC^*$-action.
Then we have a quotient morphism $S\to\CP^1$
by the $\CC^*$-action,
and the curves $C_1$ and $\ol C_1$ are 
sections of this morphism.
On the other hand, the two chains formed
by the remaining components
$C_2,\dots, C_k$ and $\ol C_2,\dots,\ol C_k$
of the cycle $C$
are the supports of two fibers of the quotient morphism $S\to\CP^1$.
We can readily represent these two fibers
as positive linear combinations of the components 
$C_2,\dots, C_k$ and $\ol C_2,\dots,\ol C_k$
respectively.
Let $f$ and $\ol f$ be these divisors, 
and $m$ the highest multiplicity
for the components in this representation.
Then the divisor $$mC + f - \ol f$$ is clearly effective
and belongs to the system $\big|mK_S\inv\big|$ since
the divisors $f$ and $\ol f$ are linearly equivalent.
Moreover, by making use of Proposition \ref{prop:d1}, we can concretely find a
collection of degree-one divisors on $Z$
which satisfies the following property:
if we write $Y$ for
the sum of the divisor belonging to the collection,
then the restriction $Y|_S$ is exactly
the divisor $mC + f - \ol f$.
Here, the collection can have the same degree-one divisor
 as elements in general.
In other words, the divisor $Y$ is non-reduced in general.
Since the restriction map
$H^2(Z,\ZZ)\to H^2(S,\ZZ)$ is always
injective, this
implies that the reducible divisor $Y$ belongs to
the pluri-system $|mF|$.
Moreover, the divisor $Y$ does not belong to the
sub-system generated by the pencil $|F|$.
Then the $(m+2)$-dimensional sub-system of $|mF|$ generated by $Y, \ol Y$
and this sub-system induces a meromorphic map
to a projective space, whose image is a (possibly singular)
rational surface.
This rational surface is the {\em minitwistor space}
described in Section \ref{ss:I}.
The images of generic twistor lines to the minitwistor spaces
are of course rational, but {\em they have
precisely $(m-1)$ ordinary nodes.}
Hence, if $m>1$, these are not minitwistor spaces in the sense of Hitchin \cite{H82-1},
but can be regarded as a natural generalization of them.
See \cite{HN11} in this respect.


As in Section \ref{ss:I}, we call these twistor spaces
{\em generalized LeBrun spaces}.
They are bimeromorphic to the total spaces of conic
bundles over the minitwistor spaces, 
and
are characterized by the property that it has a real irreducible
fundamental divisor $S$ satisfying the property (i)
in Proposition \ref{prop:clsf1}.
(As above we do not need to assume it to be Moishezon.)

%
%
%
%
%
%
%
%

\section{Study on the pluri-system $|lF|$ via
relativization}
\label{s:dim}
\subsection{The derived sequence}
\label{ss:dim}
We were interested in the structure
of Moishezon twistor spaces on $n\CP^2$
which satisfy $\dim |F|=1$.
By Propositions \ref{prop:Poon} and \ref{prop:3},
the condition $\dim|F|=1$ implies $n>3$.
Real irreducible members of the 
pencil $|F|$ on such a twistor space satisfy at least one of the two 
properties (i) and (ii) in Proposition \ref{prop:clsf1},
and we are left to understand 
the structure of twistor spaces which have a divisor
satisfying the property (ii)
in the proposition.
As in Definition \ref{def:type}, these twistor spaces are classified into 
four types $\Azero,\Aone,\Atwo$ and $\Athree$
according to the structure of the branch quartic
for real irreducible members of the pencil $|F|$

In Sections \ref{s:dim} and \ref{s:sd}, we show that 
for such a twistor space there always exists
an integer $m$ such that the pluri-system $|mF|$ 
on the twistor space contains an $(m+2)$-dimensional sub-system
which induces a degree-two map
to a scroll of planes over a rational normal curve
of degree $m$.
The restriction of this map to a smooth 
fundamental divisor $S$ will be the degree-two morphism $\phi:S\to\CP^2$
induced by the linear system $|D|$.
Therefore, from the relation $F|_S\simeq K_S\inv$, 
one may expect that if $d$ denotes the highest multiplicity
for the components of the divisor $D$ as in the previous section,
then we can take $d$ for the value of $m$.
Indeed, this was actually the case for the twistor spaces
studied in \cite{Hon_Inv} and \cite{Hon_Cre2}.
But later we will see that this is not always correct
and  in general we really need $m$ which is strictly greater than $d$.

So assume $n>3$ and let $Z$ be a twistor space on $n\CP^2$
which satisfies $\dim |F|=1$
and suppose that real irreducible fundamental divisors
satisfy the property (ii) of Proposition \ref{prop:clsf1}.
(We do not need to suppose that 
$Z$ is Moishezon.)
Then $S$ is one of the four types $\Azero$--$\Athree$ as in 
Definition \ref{def:type}.
As before, let $2k$ be the number of components
of the real anti-canonical cycle $C$ on $S$, and write it as
\begin{align}\label{cycle23}
C=\big(C_1+C_2+\dots + C_k\big) + 
\big(\ol C_1+ \ol C_2+\dots + \ol C_k\big). 
\end{align}
The number $k$ is determined by $n$ and the type (Table \ref{table5}).
We keep the convention that $C_1$ and $\ol C_1$
are the line components (see Definition \ref{def:lc}).
Recall from \eqref{z_i} that we denoted the double points of 
$C$ by $z\upi=C_i\cap C_{i+1}$ and 
$\ol z\upi=\ol C_i\cap \ol C_{i+1}$,
where $z^{(k+1)} = \ol z\upone$ and $\ol z^{(k+1)} = z\upone$.
We have $F|_S\simeq K_S\inv$ and $\Bs\,|F|=C$.
Also as before let $D$ be the effective divisor on $S$
whose support is exactly the cycle $C$, 
and whose linear system induces the generically two-to-one morphism $\phi:S\to\CP^2$.
The branch divisor of $\phi$ is a quartic.
As before we write
\begin{align}\label{cycle24}
D=\sum_{i=1}^k d_iC_i + 
\sum_{i=1}^k d_i\ol C_i,
\end{align}
and write $d=\max\{d_1,d_2,\dots,d_k\}$.
By the above convention we have $d_1=1$.

We want to show that the morphism $\phi:S\to\CP^2$ 
extends to
a map from the twistor space $Z$ to a scroll of planes.
But on $Z$, since the support of the divisor $D$
is exactly the base locus $C$ of the pencil $|F|$,
 there  exists no divisor on $Z$
whose restriction to a smooth member of $|F|$ is $D$.
A natural way to resolve this is to 
blow up the twistor space at the cycle $C$
to separate members of the pencil.
Let
\begin{align}\label{bbb}
\bbb:\hat Z\lras Z
\end{align}
be the blowing up at $C\subset Z$, and
$E_i$ and $\ol E_i$ the exceptional divisors 
over the components $C_i$ and $\ol C_i$
of $C$ respectively.
The divisors $E_i$ and $\ol E_i$ are
$\CP^1$-bundles over $C_i$ and $\ol C_i$ 
respectively.
Further, since the blowing up $\bbb$ separates members
of the pencil $|F|$, 
the intersections
of $E_i$ with the strict transforms
of distinct members of the pencil are
mutually disjoint sections of these bundles.
Therefore these are trivial bundles. Namely 
we have isomorphisms $E_i\simeq C_i\times\CP^1$ and 
$\ol E_i\simeq \ol C_i\times\CP^1$ for any index $i$.

\begin{figure}
\includegraphics{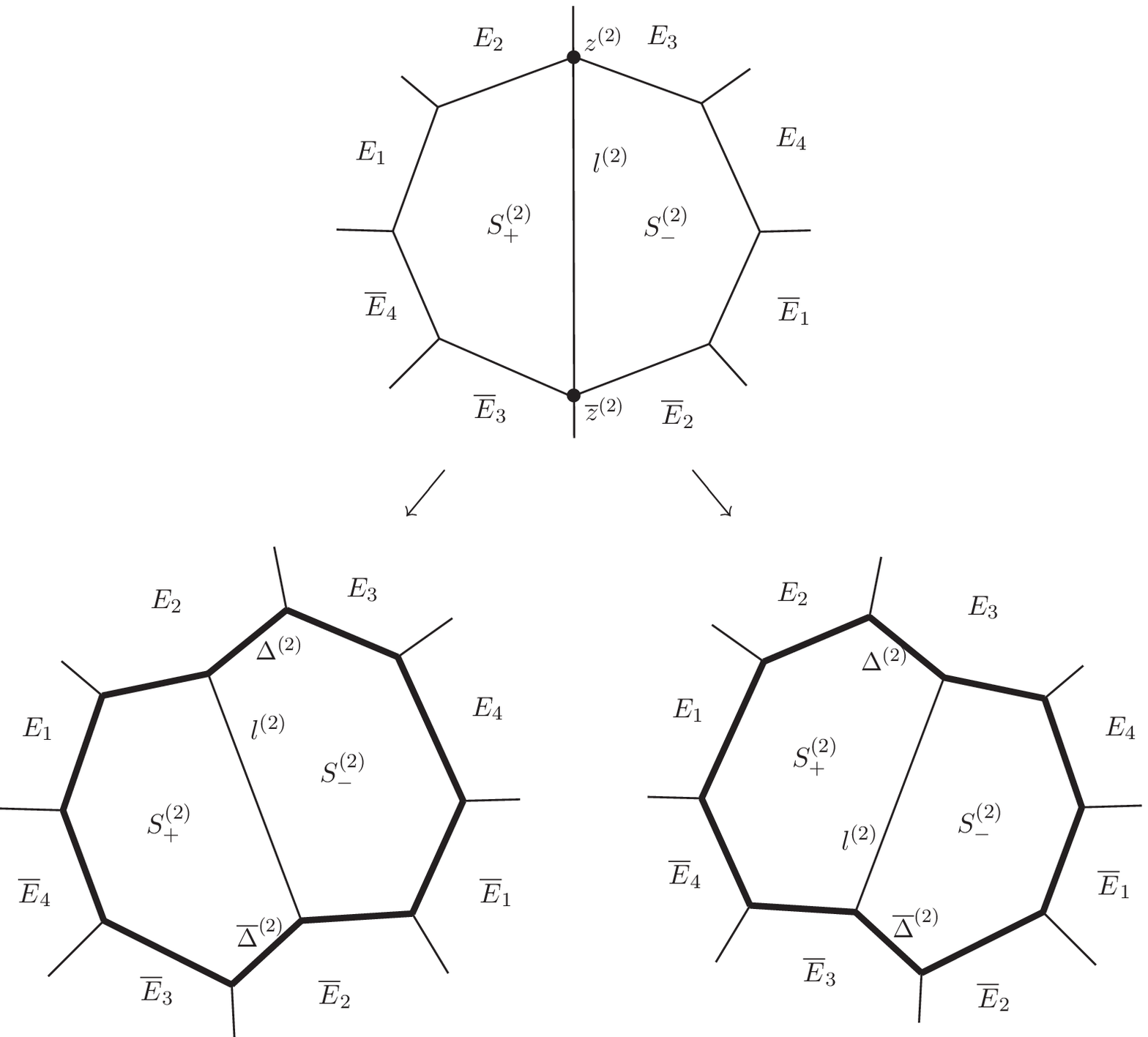}
\caption{
A reducible fiber $S\uptwo$ of $\hat Z\to\CP^1$ (up) and 
two small resolutions (down).
The two points $z\uptwo$ and $\ol z\uptwo$ are ordinary double points
of $\hat Z$, and $\DDD\uptwo$ and $\ol\DDD\uptwo$ are
the exceptional curves.
}
\label{fig:octagon2}
\end{figure}

The blowing up $\bbb$ in \eqref{bbb} eliminates the base locus
of the pencil $|F|$, so we obtain 
a morphism $\hat Z\to\CP^1$.
Any fiber of this morphism is isomorphic to the corresponding
member of the pencil $|F|$, and since the pencil
has exactly $k$ reducible members
$S\upi_+ + S\upi_-$
by Proposition \ref{prop:d1}, the morphism
$\hat Z\to\CP^1$ has exactly
$k$ reducible fibers.
We use the same symbols $S\upi_+$
and $S\upi_-$ for the strict transforms
of the degree-one divisors into $\hat Z$,
and again write  $S\upi=S\upi_++S\upi_-$.
We denote
\begin{align}\label{lmdi}
\lmd\upone,\lmd^{(2)},\dots ,\lmd\upk
\end{align}
for the points of $\CP^1$
over which the  reducible fibers 
$S\upone,S^{(2)},\dots,S\upk$ lie respectively.
These are real points as
the fibers $S\upi$ are real.
On each of these fibers, there are two points which correspond
to the points $z\upi$ and $\ol z\upi$ by the morphism $\bbb$.
We denote these points by the same symbols respectively.
As in the upper picture of Figure \ref{fig:octagon2},
these points are shared by four divisors,
and are ordinary double points of $\hat Z$.
Thus the space
$\hat Z$ has $2k$ ordinary double points
$$\big\{z\upi,\ol z\upi\set 1\le i\le k\big\}.$$


Each of these points
admits two small resolutions,
but a choice of a small resolution at the node $z\upi$
uniquely determines that of 
the conjugate node $\ol z\upi$ by the real structure.
Further, the resolution at the node $z\upi$
is distinguished by which one of $E_i$ and  $E_{i+1}$
is blown up
(see Figure \ref{fig:octagon2}).
We use the notation $\zeta\upi$ to mean 
small resolution at 
$z\upi$ and $\ol z\upi$ which preserves the real structure.
So there are two possible choices for $\zeta\upi$
for each $i\in\{1,\dots, k\}$.
Further we write 
\begin{align}\label{zeta}
\zeta:=\zeta\upone\circ\zeta\uptwo
\circ\dots\circ\zeta\upk:\tilde Z\lras \hat Z.
\end{align}
This resolves all the $2k$ nodes
of $\hat Z$ by small resolution,
and preserves the real structure.
The space $\tiZ$ is  smooth and equipped with 
the natural real structure.
There are $2^k$ possibilities for
the total small resolution $\zeta$.
For each $i$,
we write $\Delta\upi$ and $\ol \DDD\upi$ for the exceptional curves
over the nodes $z\upi$ and $\ol z\upi$ respectively.
These are exceptional curves of the factor $\zeta\upi$,
and isomorphic to $\CP^1$.
We denote the composition 
$
\tilde Z\to \hat Z\to \CP^1
$
by 
$\tilde f$.
We have the commutative diagram
\begin{align}\label{diagram000}
 \xymatrix{ 
\tilde Z \ar[r]^{\bbb\circ\zeta} \ar[d]_{\tilde f} &  Z
\ar@{.>}[dl]\\
\CP^1 &
   }
\end{align}
The morphism $\tif:\tiZ\to\CP^1$ is of fundamental significance
throughout this paper.

We use the same letters $E_i$ and $\ol E_i$,
$S\upi_+$ and $S\upi_-$
for the strict transforms
into $\tilde Z$ of the exceptional divisors of $\bbb$
and components of the reducible fibers
in the singular space $\hat Z$.
For a reason that will be apparent later,
we call the two components $E_1$ and $\ol E_1$
the {\em components over the ridge} or simply the {\em ridge components}.
Like the line components $C_1$ and $\ol C_1$ of the real anti-canonical cycle $C$,
these components play special role
in the reminder of this paper.
We put
\begin{align}\label{E}
E:=\sum_{i=1}^kE_i + \sum_{i=1}^k\ol E_i.
\end{align}
This is a real divisor on $\tiZ$.
We have $E|_S\simeq C$ for any smooth fiber $S$ of $\tif$.
This divisor will be critical for analyzing 
the structure of the present twistor spaces.

From the diagram \eqref{diagram000}, we have 
the basic relation 
\begin{align}\label{rel00}
(\bbb\circ\zeta)^*F \simeq \tilde f^*\ms O_{\CP^1}(1) + E.
\end{align}
The restriction $\tilde f|_E:E\to\CP^1$
is a trivial bundle
whose fibers are isomorphic to the cycle $C$,
except over the $k$ points $\lmd\upone,\lmd\uptwo,
\dots,\lmd\upk$ corresponding to the reducible fibers of $\tif$.
Over each of these $k$ points, 
the exceptional curves $\DDD\upi$ and $\ol\DDD\upi$ are inserted
at the nodes $z\upi$ and $\ol z\upi$ respectively
in the fiber, and the number of components of the cycle
increases by two over these points.
For each $i\in\{1,\dots,k\}$, we denote
\begin{align}\label{cyclei}
C\upi:=S\upi\cap E
\end{align}
for the cycle (in $\tiZ$).
So each $C\upi$ consists of $(2k+2)$ components.
In lower pictures in Figure \ref{fig:octagon2},
this cycle is indicated by bold cycles.


%
%

With these preliminaries, we now define an effective divisor $\bm D$ on
the smooth space $\tilde Z$ by
\begin{align}\label{bmD}
\bm D = \sum_{i=1}^k d_i E_i + 
\sum_{i=1}^k d_i \ol E_i.
\end{align}
We write $\ms D$ for the line bundle or the invertible sheaf 
associated to the divisor $\bm D$.
Clearly for a smooth fiber $S$ of $\tif$,
we have an isomorphism $\bm D|_S\simeq D$, 
and the divisor $\bm D$  is an enlargement of 
the divisor $D$ to the space $\tilde Z$.
But unlike the linear system
$|D|$ on $S$, the system $|\bm D|$ on $\tilde Z$ consists
of the divisor $\bm D$ itself.
Indeed, from negativity of the degree of the line bundle
$\ms D$ on a generic fiber of the projections 
$E_i\to C_i$ and $\ol E_i\to\ol C_i$,
we can subtract components of the divisor $\bm D$
one by one as fixed components
until it reaches the zero divisor.
Hence $|\bm D| = \{\bm D\}$ holds.
So in order to obtain a linear system on $\tilde Z$
whose restriction to a smooth fiber of $\tilde Z\to\CP^1$  equals
the linear system $|D|$,
for any integer $l\ge0$, we define 
a line bundle or an invertible sheaf $\ms D(l)$ on $\tilde Z$ by 
\begin{align}\label{}
\ms D(l):= \ms D\otimes \tilde f^*\ms O(l). 
\end{align}
We also write $\bm D(l)$ instead of $\ms D(l)$.
For the relation between the linear system
$\big|\ms D(l)\big|$ on $\tilde Z$ 
and the pluri-system $|lF|$ on $Z$,
since $\bm D\le lE$ (i.e.\,$\bm D$ is a sub-divisor of $lE$) if $l\ge d
=\max\{d_1,\dots, d_k\}$, 
from the basic relation \eqref{rel00}, 
\begin{align*}
(\bbb\circ\zeta)^*(lF)&\simeq
lE + \tif^*\ms O(l)\\
&\ge \bm D + \tif^*\ms O(l)\\
& = \bm D(l).
\end{align*}
Hence if $l\ge d$ we have a natural inclusion 
$H^0\big(\ms D(l)\big)\subset
H^0\big((\bbb\circ\zeta)^*(lF)\big)$
which is obtained by multiplying a defining 
section of the effective divisor $lE-\bm D$.
Therefore connecting with the natural identification 
$H^0\big((\bbb\circ\zeta)^*(lF)\big)\simeq H^0(lF)$, 
there is a natural inclusion
\begin{align}\label{ni}
H^0\big(\tilde Z,\ms D(l)\big)\subset H^0(Z,lF),\quad l\ge d.
\end{align}
Hence when $l\ge d$, there is a natural identification
between the linear system $\big|\ms D(l)\big|$ on 
the blown-up space $\tiZ$ and
a sub-system of the pluri-system $|lF|$ on the 
original space $Z$.
Thus our problem is transferred to the study on 
the linear system $\big|\ms D(l)\big|$ on 
the blown up space $\tiZ$ when $l\ge d$.
But we first see that,
when $l<d$, the linear system $\big|\ms D(l)\big|$ is composed with 
a pencil:
%


\begin{proposition}
\label{prop:l<m}
If $0<l<d$, 
the effective divisor $\bm D$ is 
a fixed component of
the linear system $\big|\bm D(l)\big|$.
Namely if $l<d$ we have 
$$\big|\bm D(l)\big| = \bm D+ \big|\tilde f^*\ms O(l)\big|.$$
\end{proposition}

\proof
The degree of the line bundle $\ms D(l)$
over a general fiber of the projections
$E_i\to C_i$ and $\ol E_i\to\ol C_i$ 
is $(l-d_i)$.
Hence, if $l<d=\max\{d_1,\dots, d_k\}$, this is negative
for some $i$, and therefore 
the components $E_i$ and $\ol E_i$
for such an index $i$
are fixed components of $\big|\ms D(l)\big|$.
Further, on any smooth fiber $S$ of $\tilde f$,
we have $\bm D|_S = D$ and by making use of the morphism $\phi:S\to\CP^2$,
the system $\big|D-C_i-\ol C_i\big|$ consists of 
the divisor $D-C_i-\ol C_i$ itself.
This implies that the divisor $\bm D- E_i - \ol E_i$
is fixed component of $\big|\ms D(l) - E_i - \ol E_i\big|$.
Hence the divisor $\bm D$ is a fixed component
of $\big|\ms D(l)\big|$ if $l<d$.
\proofend

\medskip

From the proposition, 
if the restriction homomorphism
$$
H^0\big(\tiZ, \ms D(l)\big)\lras H^0(S,D)
$$
is surjective for a smooth fiber $S$ of $\tilde f$ (as we wish),
we have $l\ge d$.
The problem is whether the converse also holds.
For this, the main idea we employ in this paper is to take
the direct image by the morphism $\tif$.
First, by the projection formula, we have
\begin{align}\label{}
\tilde f_*\big(\ms D(l)\big)\simeq \big(\tilde f_*\ms D\big)
\otimes\ms O_{\CP^1}(l).
\end{align}
Moreover, since the morphism $\tilde f$ is proper and flat, 
the direct image sheaf $\tilde f_*\ms D$ is
torsion-free \cite[Appendix A]{Tel15}.
Hence, since a torsion-free sheaf over a smooth curve
is locally free, the sheaf $\tilde f_*\ms D$ is locally free.
Therefore, as this is a sheaf over $\CP^1$, 
the sheaf $\tilde f_*\ms D$ is isomorphic to a direct sum
of invertible sheaves.
Furthermore, for any smooth fiber $S$ of $\tif$,
we have $h^0\big(\ms D(l)|_S\big)=h^0(S,D)=3$.
Hence the  sheaf
$\tilde f_*\ms D$ is of rank $3$.
From these, 
we obtain
\begin{align}\label{m123}
\tilde f_*\ms D\simeq
\ms O(m_1)\oplus\ms O(m_2)\oplus\ms O(m_3)
\end{align}
for some integers $m_1,m_2$ and $m_3$.
In order to calculate the dimension $h^0\big(\ms D(l)\big)$, it is enough
to determine these three numbers.

Since the space $\tiZ$ is determined from 
the small resolutions $\zeta\upone,\dots,\zeta\upk$,
the direct image
sheaf $\tif_*\ms D$ may depend on 
these.
But we have

\begin{proposition}\label{prop:indep0}
The sheaf $\tilde f_*\ms D$
is independent of a choice of the small resolutions $\zeta\upone,\dots,\zeta\upk$.
\proofend
\end{proposition}

\proof
Choose and fix the small resolutions $\zeta\upone,\dots,\zeta\upk$ in arbitrary way,
and let $\zeta:\tiZ\to \hat Z$ be
obtained from these resolutions.
It is enough to show that for any index $i=1,\dots,k$,
a replacement of $\zeta\upi$ 
by another small resolution does not affect
the form of the direct image.
Let $\tiZ'$ be the smooth space obtained
as a result of the replacement of $\zeta\upi$.
Then the spaces $\tiZ$ and $\tiZ'$ 
 are related by a flop 
at the exceptional curves $\DDD\upi = \zeta\inv(z\upi)$ and $\ol\DDD\upi=\zeta\inv(\ol z\upi)$:
\begin{align}\label{diagram3}
 \xymatrix{ 
&W \ar[dl]_{\nu} \ar[dr]^{\nu'}&\\
 \tilde Z \ar[dr]_{\zeta} & & \tilde Z'\ar[dl]
 ^{\zeta'}\\
&\hat Z&
 }
\end{align}
Here $\nu$ is the blow up at the 
exceptional curves
$\DDD\upi\cup\ol\DDD\upi$,
$\nu'$ is the blowing down of $Q_i:=\nu\inv(\DDD\upi)\simeq\DDD\upi\times\CP^1$
and $\ol Q_i=\nu\inv(\ol\DDD\upi)\simeq\ol \DDD\upi\times\CP^1$ to the other direction,
and $\zeta'$ is the contraction of the 
images $\DDD\upi\,\!':=\nu'(Q_i)$ and 
$\ol\DDD\upi\,\!':=\nu'(\ol Q_i)$
to the ordinary double points
$ z\upi$ and $\ol z\upi$ of $\hat Z$.
Let $\ms D'$ be the transformation of $\ms D$ to
$\tiZ'$, and $\tif':\tiZ'\to\CP^1$ be the composition $\tiZ'\to \hat Z\to\CP^1$.
It is enough to show 
$\tilde f_*\ms D
\simeq
\tilde f'_*\ms D'$.

Because either the small resolution $\zeta\upi$ or its replacement blows up the component $E_i$,
without loss of generality we may suppose that
 $\zeta\upi$ blows up
the component $E_i$ (among $\{E_i,E_{i+1}\}$).
This implies $\DDD\upi\subset E_i$ in $\tiZ$
(see the picture in lower left in Figure \ref{fig:octagon2}).
Hence we have
\begin{align}\label{pb1}
\nu^*\ms D
\simeq
\bm D + d_i(Q_i+\ol Q_i),
\end{align}
where in the right-hand side,
$\bm D$ means the divisor \eqref{bmD} 
with each component being understood as the
strict transform into $W$.
In a similar way, noting that
$\DDD\upi\,\!'\subset E_{i+1}$ in another resolution $\tilde Z'$
(see the picture in lower right in Figure \ref{fig:octagon2}),
we have
\begin{align}\label{pb2}
(\nu')^*\ms D'
\simeq
\bm D + d_{i+1} (Q_i+\ol Q_i).
\end{align}
Suppose $d_i=d_{i+1}$. Then
 \eqref{pb1} and \eqref{pb2} mean
$\nu^*\ms D\simeq
(\nu')^*\ms D'$.
So taking the direct image under
$\tilde f\circ\nu =\tilde f'\circ\nu'$,
we obtain
$$
(\tilde f\circ\nu)_*
\big(\nu^*\ms D\big)
\simeq
(\tilde f'\circ\nu')_*
(\nu^{'*}\ms D'\big).
$$
This implies 
$
\tilde f_*
\ms D
\simeq
\tilde f'_*
\ms D',
$
obtaining the assertion in the case
$d_i=d_{i+1}$.

To prove the case $d_i\neq d_{i+1}$,
we do some computation.
Since components of $\bm D$ which intersect
$\DDD\upi$ are $E_{i}$ and $E_{i+1}$ only, 
we have
$\bm D.\,\DDD\upi= d_{i}E_{i}.\,\DDD\upi 
+ d_{i+1} E_{i+1}.\,\DDD\upi$ in $\tiZ$.
Here, $d_{k+1}$ means $d_1$ and
$E_{k+1}$ means $\ol E_1$.
As $E_{i+1}$ and $\DDD\upi$ intersect 
transversally at one point, we have
$E_{i+1}.\,\DDD\upi =1$.
On the other hand, since the exceptional curve
$\DDD\upi$ has normal bundle $\ms O(-1)^{\oplus 2}$ in $\tiZ$,
we have  $E_i.\,\DDD\upi = -1$.
From these we obtain
\begin{align}\label{int040}
\bm D.\,\DDD\upi = d_{i+1} - d_i.
\end{align}
If use the symbol $\ms O(a,b)$ to represent
 a second cohomology class on the divisor $Q_i$ following
the canonical isomorphism with the product
$Q_i\simeq\DDD\upi\times \DDD\upi\,\!'$,
from \eqref{int040},
we obtain 
\begin{align*}
(\nu^*\ms D)|_{Q_i}
&\simeq\ms O(d_{i+1}-d_i,0).
\end{align*}
%
From this, by restricting 
$\nu^*\ms D$ to the divisors $Q_i$ and $\ol Q_i$,
 we obtain an exact sequence
\begin{multline*}\label{}
0 \lra \ms O_W \big(\bm D + (d_i-1)(Q_i+\ol Q_i)
\big)
\lras\ms O_W \big(\bm D + d_i(Q_i+\ol Q_i)
\big)\\
\lras \ms O_{Q_i}(d_{i+1}- d_i,0)
\oplus
\ms O_{\ol Q_i}(d_{i+1}- d_i,0)
\lras 0.
\end{multline*}

Now suppose $d_1> d_{i+1}$.
Taking the direct image under $\tilde f\circ\nu:
W\to\CP^1$, as $d_{i+1}-d_i<0$ from the assumption,
we obtain an isomorphism
\begin{align}\label{}
(\tilde f\circ\nu)_*\ms O_W \big(\bm D + (d_i-1)(Q_i+\ol Q_i)
\big)
\simeq
(\tilde f\circ\nu)_*
\ms O_W \big(\bm D + d_i(Q_i+\ol Q_i)
\big).
\end{align}
Continuing this process $(d_i-d_{i+1})$ times
to reduce the coefficient of $(Q_i+\ol Q_{i})$
one by one,
we eventually obtain 
\begin{align}\label{}
(\tilde f\circ\nu)_*
\ms O_W \big(\bm D + d_i(Q_i+\ol Q_i)
\big)
\simeq
(\tilde f\circ\nu)_*\ms O_W \big(\bm D + d_{i+1}(Q_i+\ol Q_i)
\big).
\end{align}
Replacing $\tilde f\circ \nu$ on the right-hand side
by $\tilde f'\circ \nu'$, from \eqref{pb1} and \eqref{pb2}, this means the desired isomorphism
$
\tilde f_*
\ms D
\simeq
\tilde f'_*
\ms D'.
$

The case $d_i<d_{i+1}$ can be shown
by replacing the roles of $\tilde Z$ and 
$\tilde Z'$ in the above argument.
\proofend
%


\medskip
Hence the three numbers $m_1,m_2$ and $m_3$ in the isomorphism \eqref{m123} 
are independent of a choice of the small resolution.
In order to derive a simple constraint for these numbers, 
we make use of the exceptional divisor $E$
in \eqref{E}
and consider  the exact sequence
\begin{align}\label{ses0002}
0 \lras \ms D(-E)
\lras \ms D
\lras \ms D|_E\lras 0.
\end{align}
If we take the direct image under the morphism $\tilde f:
\tiZ\to \CP^1$
of this exact sequence, we obtain the long exact sequence
\begin{multline}\label{der01}
0 \lras \tilde f_*\ms D(-E) \lras 
\tilde f_*\ms D
 \stackrel{\rho}{\lras} 
\tilde f_*\ms D|_{E}
{\lras} 
R^1\tilde f_*\ms D(-E) 
{\lras} 
R^1\tilde f_*\ms D
 \lras \dots.
\end{multline}
We call this as the {\em derived sequence}.
Let $\rho:\tif_*\ms D\to \tif_*\ms D|_E$
be the induced homomorphism of $\ms O_{\CP^1}$-modules as indicated in \eqref{der01}.
For the purpose of investigating the direct image
$\tilde f_*\ms D$,
we determine the two sheaves
$\tilde f_*\ms D(-E)$ and $\tif_*\ms D|_E$,
and also investigate the image of the homomorphism
$\rho$.
For the first sheaf and the support of the cokernel 
sheaf of the homomorphism $\rho$, 
it is not difficult to derive the following conclusion.
Recall that
$\lmd\upone,\lmd^{(2)},\dots ,\lmd\upk$
are points on $\CP^1$ over which 
the real reducible fibers
$S\upone,S^{(2)},\dots,S\upk$
lie respectively.

\begin{proposition}\label{prop:d_im4}
We have
\begin{align}\label{d_im4}
\tilde f_*\ms D(-E)
\simeq \ms O_{\CP^1}.
\end{align}
Moreover, for any point $\lmd\in\CP^1\minus\{
\lmd\upone,\lmd^{(2)},\dots ,\lmd\upk\}$,
the restriction of the homomorphism $\rho$
to the stalk over the point $\lmd$ is surjective.
\end{proposition}

\proof
By the same reason to the direct image
$\tilde f_*\ms D$,
the sheaf $\tilde f_*\ms D(-E)$
is locally free.
From Proposition \ref{prop:van8},
we have $h^0(S, D-C)=1$ for any smooth fiber $S$
of $\tif$.
Therefore the sheaf $\tilde f_*\ms D(-E)$ is 
an invertible sheaf.
%
%
So we can write 
$\tilde f_*\ms D(-E)\simeq
\ms O_{\CP^1}(\nu)$ for some $\nu\in\ZZ$.
On the other hand, we have an obvious isomorphism
\begin{align}\label{isom1}
H^0\big(\CP^1,\tilde f_*\ms D(-E)\big)
\simeq
H^0\big(\tilde Z,\ms D(-E)\big).
\end{align}
Now the right-hand side is 1-dimensional,
since the linear system $|\ms D|$ consists of
a single element $\bm D$ as remarked before, and 
the divisor $\bm D$ includes the divisor $E$ as a sub-divisor.
(In other words the system $\big|\ms D(-E)\big|$ consists of a single member
$\bm D-E$.)
From the isomorphism \eqref{isom1}, this means $\nu=0$, which implies 
$\tilde f_*\ms D(-E)
\simeq \ms O_{\CP^1}$.

For the second assertion, 
let $S$ be {\em any} irreducible fiber of $\tif$.
From Proposition \ref{prop:van8},
we have $H^1(S,D-C)=0$.
By Grauert's theorem about direct image sheaves, 
this implies that the support of the higher direct image
$R^1 \tif_*\ms D(-E)$ is contained in 
the set $\big\{\lmd\upone,\lmd^{(2)},\dots ,\lmd\upk\big\}$.
Since the derived sequence \eqref{der01} is exact, this means that $\rho$ is surjective over any 
point $\lmd\neq \lmd\upone,\dots,\lmd\upk$.
\proofend

\medskip
From the second assertion in the previous proposition,
 the third sheaf
 $\tif_*\ms D|_E$ in \eqref{der01}
is of rank two.
Moreover by the same reason for the sheaf $\tilde f_*\ms D$,
the sheaf  $\tif_*\ms D|_E$ is also locally free.
So we can write  $\tif_*\ms D|_E\simeq\ms O(-e)\oplus\ms O(-e')$ for some
integers $e$ and $e'$.
In the following two subsections,
by investigating common base curves of
the line bundle of the form $\ms D(l)$, $l\in\ZZ$, 
 we will 
show that the two numbers $e$ and $e'$ are equal and positive,
and express it in terms of the numbers
$d_1,\dots,d_k$ in \eqref{cycle32} or \eqref{bmD}.

\subsection{Stable base curves of the line bundle $\ms D$}
\label{ss:sbc}
In order to express the above base curves
concretely, 
we first return to the singular space $\hat Z$
and  prepare some notations.
Recall that the singular space $\hat Z$ is obtained
from the twistor space $Z$ by blowing up 
the cycle $C$, and as a result, we obtain
a morphism $\hat Z\to \CP^1$.
This fibering has exactly $k$ reducible fibers. 
As before, we write 
the reducible fiber over the point $\lmd\upi$
as
$$S\upi=S\upi_+ + S\upi_-,\quad 1\le i\le k.$$
All these are real.
As before let $E=\sum_{i=1}^k E_i + \sum_{i=1}^k\ol E_i$ be the exceptional divisor
of the blowing up.
Then for any index $i\in\{1,\dots,k\}$ the intersection $S\upi\cap E$ in 
the singular space $\hat Z$ is a cycle of rational curves
which is isomorphic to the original cycle $C$ in $Z$.
By using this isomorphism,
we write this cycle as
\begin{align}\label{Cupi}
\sum_{j=1}^k C\upi_j + \sum_{j=1}^k \ol C\upi_j,
\end{align}
where $C\upi_j$ and $\ol C\upi_j$ are the components
which correspond to $C_j$ and $\ol C_j$ respectively.
We have $C\upi_j = S\upi\cap E_j$ and $\ol C\upi_j = S\upi\cap \ol E_j$ for any indices $i,j\in\{1,\dots,k\}$
in the singular space $\hat Z$.


For studying the base curves, we need to take and fix a total  
small resolution.
For this, let $\zeta:\tiZ\to \hat Z$ be the total small resolution such that any factor 
$\zeta\upi$ blows up the components $E_i$ and $\ol E_i$ (as in the picture in lower left in Figure \ref{fig:octagon2}).
As before, let $\DDD\upi\simeq\CP^1$ be the exceptional
curve over the node $z\upi$.
Also we do not change notations for
the strict transforms of the curves $C\upi_j$ and the reducible fibers
$S\upi=S\upi_+ + S\upi_-$ into $\tiZ$.
So we have $S\upi=\tif\inv(\lmd\upi)$.
As in \eqref{cyclei}, for each $i\in\{1,\dots,k\}$, we write
$$
C\upi= S\upi\cap E.
$$
This is a cycle of rational curves consisting of $2k+2$ components and includes the two curves
$\DDD\upi$ and $\ol\DDD\upi$.
Some chains in these $k$ cycles
$C\upone,\dots,C\upk$ are of our main interest in this subsection.
In addition to these, in Section \ref{ss:cbc}, we will find that some of these cycles themselves
will be  base curves of the system 
of the form $\big|\ms D(l)\big|$.

\begin{figure}
\includegraphics{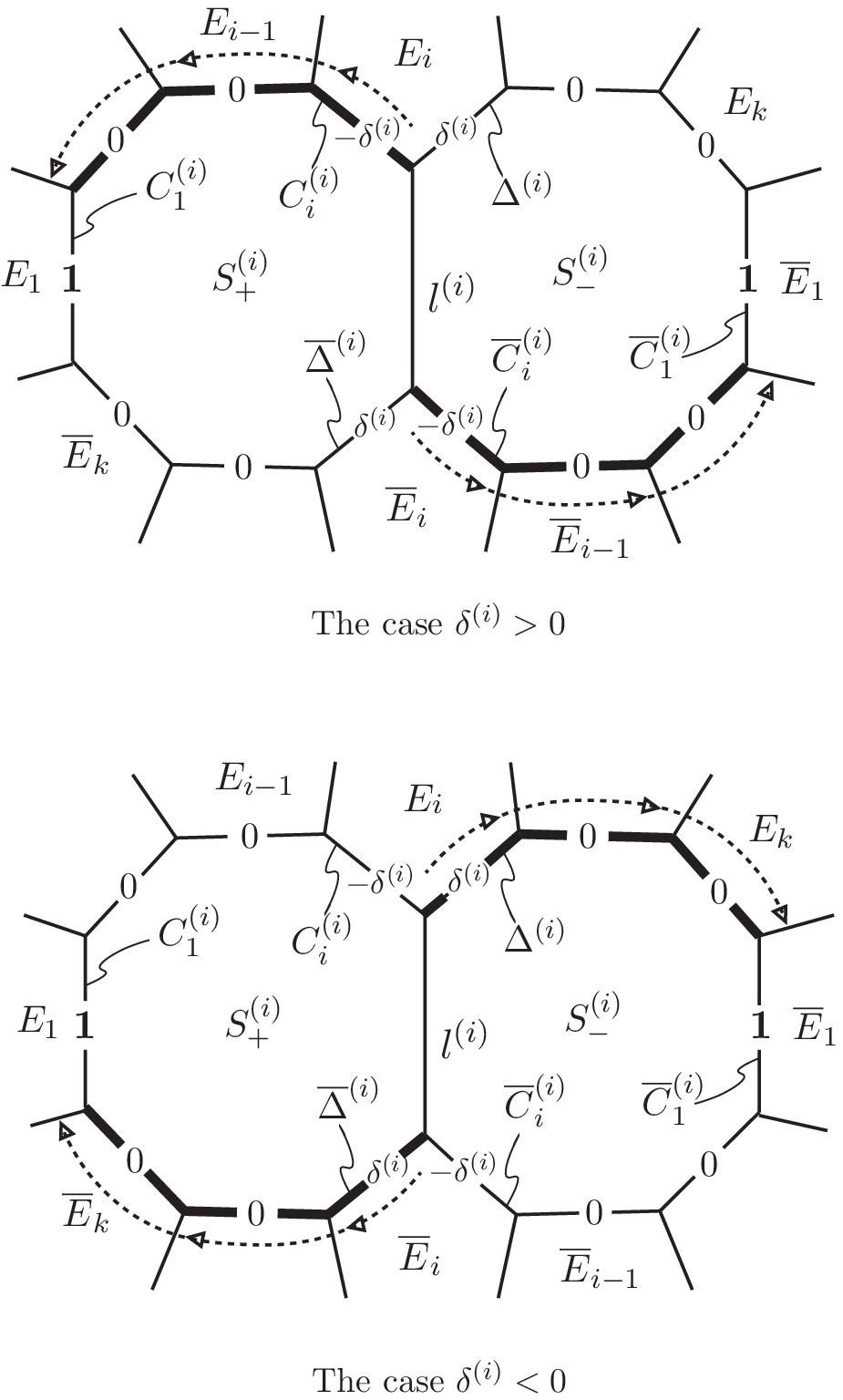}
\caption{The intersection numbers with the divisor $\bm D$
(the numbers on the segments),
and the stable base curves (bold chains) 
when $i\neq 1$}
\label{fig:rf1}
\end{figure}

From the above choice of each factor $\zeta\upi$,
the exceptional curve $\DDD\upi$ is inserted in the component $E_{i}$ for any $i$, 
and we have
\begin{align}\label{sf}
S\upi|_{E_i} = \DDD_{i} + C\upi_i,\quad i\in\{1,\dots,k\}.
\end{align}
(See Figure \ref{fig:rf1}.)
As in \eqref{bmD} let $\bm D=\sum_{i=1}^k d_iE_i + \sum_{i=1}^k d_i\ol E_i$ be the effective divisor on 
$\tiZ$ which is the enlargement of the divisor 
$D$ on a smooth fiber $S$ to the whole of $\tiZ$,
and $\ms D$ the invertible sheaf associated to $\bm D$.

First,
as we obtained in \eqref{int040}, we have
\begin{align}\label{int04}
\bm D.\,\DDD\upi = d_{i+1} - d_i
\end{align}
for any index $i$.
In the following for simplicity we denote
$$\ddd\upi:= d_{i+1} - d_i,$$
so that $\bm D.\,\DDD\upi = \ddd\upi$.

Next, let $S$ be a general fiber of $\tilde f$
and write the curve $S\cap E_i$ simply by $C_i$.
Then in the divisor $E_i$, the reducible curve $\DDD_{i}+ C\upi_i$
is homologous to the curve $C_i$ from \eqref{sf}.
Therefore
$$\bm D.\,\big(\DDD_{i}+ C\upi_i\big)
= \bm D.\,C_i=(D,C_i)_S
=\begin{cases}
0& {\text{ if $i\neq 1$}},\\
1& {\text{ if $i=1$,}}
\end{cases}$$
where the final equality is due to 
the basic intersection numbers
\eqref{int0011}.
Hence by \eqref{int04} and the real structure, we obtain 
\begin{align}\label{int05}
\bm D.\,C\upi_i = \bm D.\,\ol C\upi_i =
\begin{cases}
-\ddd\upi & {\text{ if $i\neq 1$}},\\
-\ddd\upi +1 = 2-d_2 & {\text{ if $i= 1$}}.
\end{cases}
\end{align}
Here, for the final equality we used $d_1=1$.

\begin{figure}
\includegraphics{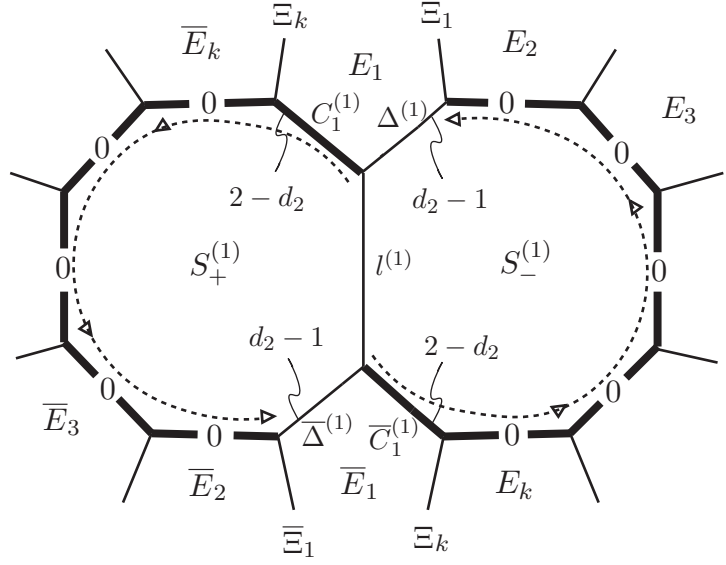}
\caption{The case $i= 1$}
\label{fig:rf2}
\end{figure}
The intersection numbers of the divisor $\bm D$ 
with other components of the cycle $C\upi$ can be
obtained in  a similar way from the basic intersection numbers,
and 
we have, when $i\neq 1$,
\begin{align}\label{int06}
\bm D.\,C\upi_j = \bm D.\,\ol C\upi_j =
\begin{cases}
0, & {\text{ if $j\neq 1,i$}},\\
1, & {\text{ if $j= 1$}},\\
\end{cases}
\end{align}
and 
\begin{align}\label{int07}
\bm D.\,C\upone_j = \bm D.\,\ol C\upone_j = 
0\,\,  {\text{ for any $j\neq 1$}}.
\end{align}
In Figures \ref{fig:rf1} and \ref{fig:rf2}
these intersection numbers are put on the segments
which represent the components respectively. 

From these intersection numbers
\eqref{int05}--\eqref{int07}, 
we immediately find 
base curves of the line bundle $\ms D(l)$ on the present $\tiZ$,
that are included in the cycles $C\upi$
as follows.
We promise $d_{k+1} := d_1$ as before.
We need to distinguish the cases $i=1$ and $i\neq 1$
because of the discrepancy in \eqref{int05}.
We first discuss the case $i\neq 1$.

\begin{proposition}\label{prop:bc1}
(i) If $i\neq 1$ and $\ddd\upi>0$ ,
the chain of rational curves
\begin{align}\label{bc01}
C\upi_i+ C\upi_{i-1}+ C\upi_{i-2}+
\dots+ C\upi_2\subset C\upi,
\end{align}
is base curve of the line bundle 
$\ms D(l)$ for any $l\in\ZZ$.
Moreover, this chain is fixed component 
of the restriction $\ms D(l)|_{E}$ with 
multiplicity $\ddd\upi$ for any $l\in\ZZ$.\\
(ii) If $i\neq 1$ and $\ddd\upi<0$,
the chain of rational curves 
\begin{align}\label{bc02}
\DDD\upi+ C\upi_{i+1}+ C\upi_{i+2}+ C\upi_{i+3}+ 
\dots+  C\upi_{k}\subset C\upi,
\end{align}
is  base curve of the bundle $\ms D(l)$ for any $l\in\ZZ$.
(If $i=k$, this reads $\DDD\upi$.)
Moreover, this chain is  fixed component 
of the restriction $\ms D(l)|_{E}$ with 
multiplicity $\big(-\ddd\upi\big)$ for any $l\in\ZZ$.
\end{proposition}

From the real structure, the conjugate chains to \eqref{bc01} and \eqref{bc02}
are also base curves of $\ms D(l)$ for any $l$.
They do not intersect the original ones.
See Figure \ref{fig:rf1} for these base curves
when $\ddd\upi>0$ and $\ddd\upi<0$ respectively.

\medskip
\noindent
{\em Proof of Proposition \ref{prop:bc1}.}
We only prove (i) because (ii) can be shown in a similar way.
From \eqref{int05} and \eqref{int06},
noting that 
the curves $C\upi_j$ are contained in a fiber of $\tilde f$,
we have, for any $l\in\ZZ$,
\begin{align*}
\ms D(l).\,C\upi_j
= \bm D.\,C\upi_j
=
\begin{cases}
-\ddd\upi & {\text{if $j=i$}},\\
0 & {\text{if $j\neq 1,i$}}.
\end{cases}
\end{align*}
Since we have $-\ddd\upi<0$ from the assumption,
these imply that
the chain \eqref{bc01} is a base curve of $\ms D$.
Note that the chain stops when it hits the component $\ol E_1$ because $\bm D.\,C\upi_1 = 1>0$.

Next we consider the restriction $\ms D(l)|_E$ to the exceptional divisor $E$ 
and see that any section of 
the line bundle $\ms D(l)|_E$ vanishes
along the chain \eqref{bc01} 
with multiplicity $\ddd\upi$ as in the proposition.
First, for the  end component $C\upi_{i}$
of the chain \eqref{bc01}, we have, on the component $E_i$
of $E$,
\begin{align}\label{int08}
\big(\bm D|_{E_i}- C\upi_i, C\upi_{i}\big)_{E_i}
= \big(\bm D,\,C\upi_{i}\big)_{E_i} - \big(C\upi_{i}\big)^2_{E_i}
= \big(\bm D,C\upi_{i}\big)_{\tilde Z}+1
\end{align}
since $C\upi_i$ is a $(-1)$-curve on $E_i$.
By \eqref{int05} this is still negative if 
$\ddd\upi>1$, and hence $C\upi_{i}$ is yet a
base curve of the system $\big|\bm D|_{E_i}- C\upi_{i}\big|$ if $\ddd\upi>1$.
Moreover, as $\big(C\upi_{j}\big)^2=0$ on $E_{j}$
for any $j$ with $i>j>1$,
we have
\begin{align}\label{int09}
\big(\bm D|_{E_{j}}- C\upi_{j}\big).\,C\upi_{j}
= \bm D.\,C\upi_{j},
\end{align}
and all these are zero as above.
From \eqref{int08} and \eqref{int09}, 
we obtain that any section of the line bundle
$\ms D|_E$ vanishes along the chain \eqref{bc01} 
with multiplicity at least two if $\ddd\upi>1$.
%
Repeating this process, the chain  \eqref{bc01} can be
removed $\ddd\upi$ times.
Thus we have seen that the chain \eqref{bc01} is fixed component of $\ms D(l)|_E$ with multiplicity
$\ddd\upi$.
\proofend

\medskip
If $d_i = d_{i+1}$, a base curve does not appear on $S\upi$.
As seen above, each of the chains \eqref{bc01}
and \eqref{bc02} has a unique component (which is 
$C\upi_i$ and $\DDD\upi$ respectively)
whose intersection number with $\bm D$ is negative,
and this is always an end component of the chain.
This will be called the fixed end of the chain.


In a similar way, in the remaining case $i=1$,
we find base curves on 
the reducible fiber $S\upone =S^{(1)}_+ + S\upone_-$ as follows. 

\begin{proposition}\label{prop:bc2}
If $i=1$ and $d_2>2$, 
the chain of rational curves 
\begin{align}\label{bc03}
C\upone _{1}+ \ol C\upone _{k}+ \ol C\upone _{k-1}+ 
\dots +  \ol C\upone _{2}\subset C\upone 
\end{align}
is base curve of the line bundle
$\ms D(l)$ for any $l\in\ZZ$.
Moreover this chain is  fixed component of the restriction $\ms D(l)|_E$ with multiplicity $(d_2-2)$ for any $l\in\ZZ$.
\end{proposition}

See Figure \ref{fig:rf2} for these base curves.

\medskip
\proof
From \eqref{int04} we have $\ms D(l).\,\DDD\upone = \bm D.\,\DDD\upone = d_2 - 1$ for any $l$ and this is non-negative since any $d_i$ is always positive.
By \eqref{int05} we have $\ms D(l).\,C\upone_1 = 2- d_2$
and hence $C\upone_1$ is a base curve if $d_2>2$.
Then by \eqref{int07} the chain \eqref{bc03} is
a base curve of $\ms D(l)$ for any $l$.
The assertion about the multiplicity when restricted
to $E$ can be obtained in a similar way to 
the previous proof.
\proofend


\medskip
If $d_2 = 1,2$, no base curve appears on $S\upone$.

Thus we have obtained chains on the exceptional divisor $E$, which are
base curves of the linear system $\big|\ms D(l)\big|$ for arbitrary $l$,
for the total small resolution $\zeta:\tilde Z\to \hat Z$ for which each factor $\zeta\upi$
blows up the components $E_i$ and $\ol E_i$.

\begin{definition}\label{def:sbc}{\em 
We call these chains
\eqref{bc01}, \eqref{bc02} and \eqref{bc03} and 
their conjugations 
as the \!{\em stable base curves} of the line bundle $\ms D$.
By the {\em fixed end} of the stable base curves, 
we mean the end component of the chain
whose intersection number with $\ms D$ is negative.
\proofend}
\end{definition}

We also prepare notations for the stable base curves.
For each index $i=1,\dots,k$, we denote $\CCC\upi_E$ to mean the sum of the stable base curves
which are included in the cycle $C\upi$, with the multiplicities on the divisor $E$
being taken into account.
For example, if $i\neq 1$ and $\ddd\upi=d_{i+1}-d_i>0$,
we have
$$
\CCC\upi_E = \ddd\upi
\Big\{
\big(C\upi_i+ C\upi_{i-1}+ \dots+ C\upi_2\big)
+
\big(\ol C\upi_i+ \ol C\upi_{i-1}+ \dots+ \ol C\upi_2\big)
\Big\}.
$$
We further denote
\begin{align}\label{tg}
\CCC_E = \sum_{i=1}^k \CCC\upi_E.
\end{align}
The subscript $E$ is put in order to emphasize that
they are considered as divisors on $E$.
But we note that this is {\em not} a Cartier divisor on $E$,
because 
each of the chains \eqref{bc01}, \eqref{bc02} and \eqref{bc03}
stops at the intersection when it
hits the ridge component $E_1$ or $\ol E_1$.

So far we have fixed a simultaneous small resolution 
$\zeta:\tiZ\to\hat Z$ as the particular one and obtained the base curves of $\ms D(l)$
on $\tiZ$ for this resolution.
The base curves of the line bundle  $\ms D(l)$ can  be obtained in a similar way
for any total  small resolution of $\hat Z$
preserving the real structure.
We also call them the stable base curves of $\ms D$, and use the same 
notations for them as above.

\subsection{Explicit form of the sheaf
$\tilde f_*\ms D|_E$
}
\label{ss:e}
In this subsection, 
by utilizing the stable base curves of $\ms D$,
we express the direct image sheaf
$\tilde f_*\ms D|_E$ in terms of the multiplicities
$d_1,\dots, d_k$ appearing in the divisor $D$ on $S$ or $\bm D$ on $\tiZ$.

As above, the stable base curve $\CCC_E$ of $\ms D$ is not a Cartier divisor on $E$.
But the restriction $\CCC_{E-E_1 - \ol E_1}:=\CCC_E|_{E-E_1-\ol E_1}$ (i.e.\,the collection
of the components which lie on $E-E_1-\ol E_1$) is clearly a Cartier divisor on 
$$E-E_1-\ol E_1 = (E_2\cup E_3\cup\dots\cup E_k)\sqcup (\ol E_2\cup \ol E_3\cup\dots\cup \ol E_k).
$$
Therefore, the subtraction
\begin{align}\label{lsb}
\ms D|_{E-E_1-\ol E_1} - \CCC_{E-E_1-\ol E_1}
\end{align}
makes sense as a line bundle or an invertible sheaf over $E-E_1-\ol E_1$.
From the intersection numbers \eqref{int04}--\eqref{int07}
and the multiplicities of the chains in $\CCC_E$,
the intersection numbers of the line bundle \eqref{lsb}
with any irreducible component of fibers of 
$\tif:E_2\cup E_3\cup \dots \cup E_k\to\CP^1$
and 
$\tif:\ol E_2\cup\ol E_3\cup \dots \cup\ol E_k\to\CP^1$
are all zero, regardless of a choice of 
the small resolution $\zeta:\tiZ\to\hat Z$
preserving the real structure.
(Here we are omitting the subscripts for $\tif$
that indicate the restriction.)
This implies that we always have 
\begin{align}\label{lsb2}
\ms D|_{E-E_1-\ol E_1} - \CCC_{E-E_1-\ol E_1}
\simeq\tif^*\ms O_{\CP^1}(-e)
\end{align}
for some integer $e$.
(Minus is put here because $(-e)$ will be positive.)
Namely we can write 
\begin{align}\label{lsb3}
\ms D(e)|_{E-E_1-\ol E_1} \simeq 
\CCC_{E-E_1-\ol E_1}
\end{align}
A priori the number $e$ depends on the choice of 
the small resolution $\zeta$, and it will turn
out that this is actually the case in general.
For simplicity of notation, 
we write $\CCC_{E_2\cup\dots \cup E_k}$
for the restriction of the stable base curves
$\CCC_E$ to the union $E_2\cup\dots \cup E_k$.
Because the divisor $\CCC_E$ is 
fixed component of the line bundle $\ms D(l)|_E$ for any $l$, 
\eqref{lsb3} means that for any non-zero section
$t\in H^0\big(\ms D(e)|_E\big)$
that does not vanish identically on $E_2\cup
\dots\cup E_k$, 
we have the coincidence
\begin{align}\label{1dim0}
(t)|_{E_2\cup\dots\cup E_k} =
\CCC_{E_2\cup\dots\cup E_k}.
\end{align}
Further, for any $l\ge e$ and 
any section $t\in H^0\big(\ms D(l)|_E\big)$
that does not vanish identically on $E_2\cup
\dots\cup E_k$, 
we have an inequality
\begin{align}\label{1dim1}
(t)|_{E_2\cup\dots\cup E_k} \ge
\CCC_{E_2\cup\dots\cup E_k},
\end{align}
as effective divisors,
and the difference
\begin{align}\label{dif0}
(t)|_{E_2\cup\dots\cup E_k} -
\CCC_{E_2\cup\dots\cup E_k}
\end{align} is always
a sum of $(l-e)$ fibers of $\tif|_{E_2\cup
\dots\cup E_k}:E_2\cup
\dots\cup E_k\to\CP^1$.
Of course, these fibers are not necessarily distinct.
By the real structure, analogous properties hold
on another union $\ol E_2\cup\dots \cup \ol E_k$.


If we consider the  subspace $\Omega(l)$ of 
$H^0\big(\ms D(l)|_{ E_2\cup\dots \cup  E_k}\big)$ defined by
\begin{align}\label{om1}
\Omega(l):=\Big\{t\in H^0\big(\ms D(l)|_{ E_2\cup\dots \cup  E_k}\big)
\,\big|\,
(t)\ge \CCC_{E_2\cup\dots \cup E_k}
\Big\}
\end{align}
when $l\ge e$, 
then from \eqref{1dim1} the image of the restriction homomorphism
$$
 H^0\big(\ms D(l)|_{ E}\big)
\lras
H^0\big(\ms D(l)|_{ E_2\cup\dots \cup  E_k}\big)
$$
is contained in the subspace $\Omega(l)$.
Similarly  the subspace
$\ol\Omega(l)$ of $ H^0\big(\ms D(l)|_{\ol E_2\cup\dots \cup \ol E_k}\big)$
is defined and it satisfies an analogous property by the real structure.

If the total small resolution $\zeta:\tiZ\to \hat Z$
does not blow
up the ridge components $E_1$ and $\ol E_1$,
no component of $\CCC_E$ is on $E_1$ nor $\ol E_1$,
and this means that 
the condition $(t)\ge \CCC_{E_2\cup\dots \cup E_k}$ is
satisfied for any $t\in H^0\big(\ms D(e)|_{E_2\cup\dots\cup E_k}\big)$.
So for such a small resolution,
the subspace $\Omega(l)$ is nothing but 
$H^0\big(\ms D(l)|_{ E_2\cup\dots \cup  E_k}\big)$.
If the resolution $\zeta$ blows up the components $E_1$ and $\ol E_1$
(namely if at least one of the two factors $\zeta\upone$ and $\zeta\upk$ blows up the ridge components),
the stable base curves $\CCC_E$ can have a component 
on $E_1$ and $\ol E_1$, and in that case 
$\Omega(l)$ is strictly smaller than $H^0\big(\ms D(l)|_{ E_2\cup\dots \cup  E_k}\big)$
because
the condition $(t)\ge \CCC_{E_2\cup\dots \cup E_k}$ gives 
a non-trivial constraint for 
$t\in H^0\big(\ms D(l)|_{ E_2\cup\dots \cup  E_k}\big)$
 coming from such a component.
 
By \eqref{1dim0} and \eqref{1dim1},
if $l=e$, a section of $\ms D(e)$ over $E_2\cup\dots \cup E_k$ 
whose zero is exactly $\CCC_{E_2\cup\dots \cup E_k}$
gives a generator of the space $\Omega(e)$,
and we have $\dim\Omega(e)=1$.
Similarly, by the above property on the difference
\eqref{dif0},
 if $l>e$, we have
\begin{align}\label{om1d}
\dim \Omega(l) = l-e+1.
\end{align}
%
%
%

From \eqref{lsb2}, the number $e$ can be calculated 
if the intersection number of the left-hand side
of \eqref{lsb2} with a section of $\tif|_E:E\to\CP^1$
is known.
The intersection of any two adjacent components of the exceptional divisor $E$ 
provides such a section, and we prepare
a notation for these sections
by putting
$$
\Xi_i:=E_i\cap E_{i+1},\quad
1\le i\le k.
$$
Here, $E_{k+1}$ means $\ol E_1$ as before.
We also define, for any integer $r$,  
\begin{align}\label{rplus}
r_+ := 
\max\{r,0\}.
\end{align}

\begin{proposition}\label{prop:e}
The number $e$ in \eqref{lsb2} or \eqref{lsb3}
is concretely given as follows.
\begin{itemize}
\item If both of the small resolutions $\zeta\upone$
and $\zeta\upk$ 
do not blow up
the ridge components $E_1$ and $\ol E_1$,
\begin{align}\label{case1}
e =  d_2 + \sum_{1<j< k} (d_{j+1} - d_j)_+.
\end{align}
\item
If we replace the factor $\zeta\upone$
by another small resolution from the first item,
$(d_2-2)_+$ is added to \eqref{case1}.
\item 
If we replace the factor $\zeta\upk$
by another small resolution from the first item,
$(d_k-2)_+$ is added to \eqref{case1}.
\item 
If we replace the factors $\zeta\upone$
and $\zeta\upk$ by another resolutions respectively
 from the first item,
$(d_2-2)_++(d_k-2)_+$  is added to \eqref{case1}.
\end{itemize}
\end{proposition}

\proof
Taking the intersection number with a section $\Xi_1\subset E_1$,
from \eqref{lsb2}, we have
\begin{align}\label{mie}
-e = \bm D.\,\Xi_1 - \big(\CCC_E,\Xi_1\big)_{E_1}.
\end{align}
Both intersection numbers $\bm D.\,\Xi_1$ and $\big(\CCC_E,\Xi_1\big)_{E_1}$
depend on the choice of the small resolution $\zeta:\tiZ\to\hat Z$
preserving the real structure,
and we need case by case calculation as follows.

If the small resolution $\zeta\upone$ does not blow up the components $E_1$ and $\ol E_1$, we have 
$\Xi_1^2 = 0$ on $E_1$ and $\Xi_1^2 = -1$ on $E_2$
for the self-intersection numbers.
Hence we have
\begin{align}\label{int101}
\bm D.\,\Xi_1 &= (d_1 E_1 + d_2 E_2).\,\Xi_1\notag\\
&= d_1E_1|_{E_{2}}.\,\Xi_1
+d_{2}E_{2}|_{E_{1}}.\,\Xi_1\notag\\
&= d_1\cdot (-1) +d_{2}\cdot 0\notag\\
&= -1 \quad(\because d_1=1).
\end{align}
Alternatively, if $\zeta\upone$ 
blows up the component $E_1$
we have 
$\Xi_1^2 = -1$ on $E_1$ and $\Xi_1^2 = 0$ on $E_2$
for the self-intersection numbers.
Hence by similar calculation, we obtain
\begin{align}\label{int100}
\bm D.\,\Xi_1 &= d_1\cdot 0 + d_2\cdot(-1) =-d_2.
\end{align}

\begin{figure}
\includegraphics{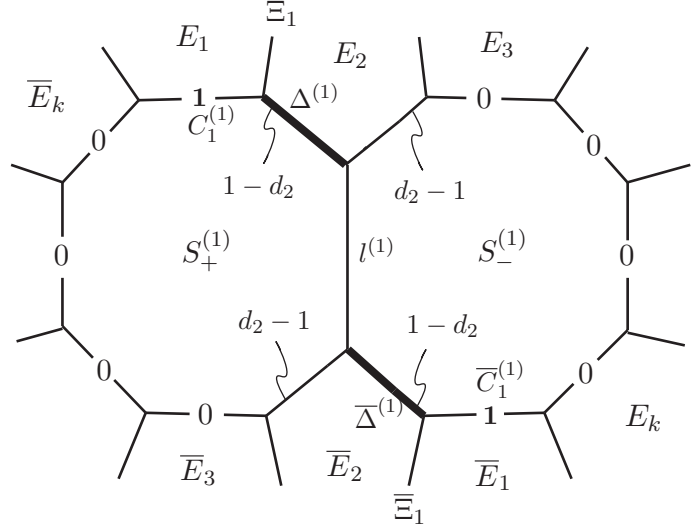}
\caption{Stable base curves (bold segments) when 
resolutions at $z\upone$ and $\ol z\upone$ are replaced}
\label{fig:rf3}
\end{figure}

Next we calculate another intersection number $\big(\CCC_E,\Xi_1\big)_{E_1}$, 
by gathering all components of  $\CCC_E$ which intersect $\Xi_1$.
First, from the concrete description of 
the stable base curves,
if $i\neq 1,k$, the support of $\CCC\upi_E$ intersects
$\Xi_1$ iff $\ddd\upi=d_{i+1} - d_i >0$,
and in this situation the multiplicity of $\CCC\upi_E$ is $\ddd\upi$.
Second, for the multiple chain $\CCC\upone_E$,
if the small resolution $\zeta\upone$ blows up $E_1$ and $\ol E_1$,
the support of $\CCC\upone_E$ always intersects $\Xi_1$
and its multiplicity is $(d_2-2)_+$.
(See Figure \ref{fig:rf2}.)
Alternatively if $\zeta\upone$ does not up the components $E_1$ nor $\ol E_1$,
we have $\CCC\upone_E = (d_2-1)
\big(\DDD\upone
+ \ol\DDD\upone\big)$ (see Figure \ref{fig:rf3}),
and this always intersects $\Xi_1$.
Third,
if the small resolution $\zeta\upk$ blows up $E_1$ and $\ol E_1$,
the support of $\CCC\upk_E$ always intersects 
$\Xi_1$, and the multiplicity is $(d_k-2)_+$.
Alternatively, if the small resolution $\zeta\upk$ does not blow up $E_1$ and $\ol E_1$, we have $\CCC\upk_E = (d_k-1)
\big(\DDD\upk + \ol\DDD\upk\big)$
but this does not intersect $\Xi_1$.
From these, we obtain the following.
\begin{itemize}
\item
If the small resolution $\zeta$ does not blow up $E_1$ and $\ol E_1$, we have 
\begin{align}\label{cg1}
\big(\CCC_E,\Xi_1\big)_{E_1} = (d_2-1) + \sum_{i=2}^{k-1} \ddd\upi_+,
\end{align}
where we are using the notation \eqref{rplus}.
\item
If we replace the factor $\zeta\upone$
from the first item,
$(d_2-1)$ in \eqref{cg1} is replaced by $(d_2-2)_+$.
\item
If we replace the factor $\zeta\upk$
from the first item,
$(d_k-2)_+$ is added to \eqref{cg1}.
\item
If we replace the factor
$\zeta\upone$ and $\zeta\upk$ from the first item,
$(d_2-1)$ in \eqref{cg1} is replaced by $(d_2-2)_+$,
and $(d_k-2)_+$ is added.
\end{itemize}
The number $e$ is obtained from these
and the intersection numbers  \eqref{int101}
and \eqref{int100}
by using \eqref{mie}.
\proofend

%

\medskip
Since $d_2$ is positive, 
the number 
\eqref{case1} is positive.
Therefore, the number $e$ is  positive
for any total small resolution $\zeta:\tiZ\to\hat Z$.

Next we identify the space
$H^0\big(E, \ms D(l)|_{E}\big)$ when $l\ge e$.

\begin{proposition}\label{prop:rest3}
For any integer $l\ge e$, the restriction homomorphism from the space
$
H^0\big(E, \ms D(l)|_{E}\big)
$ 
to the space of sections over the union
\begin{align}\label{union1}
\big(E_2\cup\dots\cup E_k\big)\amalg
\big(\ol E_2\cup\dots\cup \ol E_k\big)
\end{align}
gives an isomorphism
\begin{align}\label{isdl}
H^0\big(E, \ms D(l)|_{E}\big)\simeq \Omega(l)\oplus\ol\Omega(l).
\end{align}
In particular, 
$h^0\big(E, \ms D(l)|_{E} \big)
=2(l-e+1)$ by \eqref{om1d}.
\end{proposition}

\proof
That the image of the restriction homomorphism is included
in the direct sum $\Omega(l)\oplus\ol\Omega(l)$ is obvious
from the definition of the subspace $\Omega(l)$.

For the isomorphicity (or rather surjectivity),
again we need to investigate case by case depending on 
whether the small resolutions $\zeta\upone$
and $\zeta\upk$ blow up the ridge component
$E_1$ or not.
We begin with the easiest case where 
both $\zeta\upone$ and $\zeta\upk$ do not blow up
the components $E_1$ and $\ol E_1$
 so that we have  
isomorphisms $E_1\simeq C_1\times\CP^1$ and
$\ol E_1\simeq \ol C_1\times\CP^1$ as they were in the singular space $\hat Z$.
Following these isomorphisms, we denote $(1,0)$ 
(resp.\,$(0,1)$) for the
fiber class of the projection to the first (resp.\,second) factor.
Then letting $S$ be a smooth fiber of $\tif$,
again from the basic intersection number \eqref{int0011},
we have
$$
\ms D|_{E_1}.\,(0,1) = \bm D|_S.\,C_1 = 1,
$$ 
where $C_1$ is the line component $S\cap E_1$.
Hence, with \eqref{int101}, we have, for any $l\in\ZZ$
\begin{align}\label{}
\ms D(l)|_{E_1}\simeq (1,l-1).
\end{align}
Therefore the kernel sheaf of the restriction homomorphism
from the invertible sheaf $\ms D(l)|_{E_1}$
to the  curves $\Xi_1\cup\ol \Xi_k$ is isomorphic to
the class $(-1,l-1)$.
As any cohomology group of this class vanishes, 
any section of $\ms D(l)|_{E_1}$ defined on the curves
$\Xi_1\cup\ol\Xi_k$ uniquely extends to 
the component $E_1$.
By the real structure, analogous property holds
on another component $\ol E_1$.

Let $(t,t')$ be any element of $\Omega(l)\oplus\ol\Omega(l)$.
By restriction this determines 
a section of 
$\ms D(l)$ defined on $\Xi_1\cup\ol\Xi_k$,
as well as a section of $\ms D(l)$ defined on $\ol\Xi_1\cup\Xi_k$.
From the above extendability, the pair $(t,t')$ admits
an extension to $E_1$ and $\ol E_1$
and it is unique.
Thus we obtain the required isomorphicity for these small resolutions.

Next we replace a resolution $\zeta\upone$,
so that the exceptional curve $\DDD\upone$ is 
inserted in the component $E_1$.
(See Figure \ref{fig:rf2}.)
Again let $(t,t')$ be any element of $\Omega(l)\oplus\ol\Omega(l)$.
Then $t'$ vanishes along the chain
$\ol C\upone_k\cup \ol C\upone_{k-1}\cup
\dots \cup \ol C\upone_2$ by multiplicity $(d_2-2)_+$.
In particular the restriction $t'|_{\ol\Xi_k}$
vanishes at the intersection point $\ol\xi\upone_k:=\ol\Xi_k\cap S\upone$
by multiplicity $(d_2-2)_+$.
Let $u'$ be a section which is obtained from 
$t'|_{\ol\Xi_k}$ by dividing by a defining section
of the divisor $(d_2-2)_+\ol\xi\upone_k$ on $\ol\Xi_k$.
Then $u'$ can be regarded as a section of 
the line bundle $\ms D(l)|_{E_1}-(d_2-2)_+C\upone_1$
defined over the curve $\ol\Xi_k$.
On the other hand, the restriction 
$t|_{\Xi_1}$ can be regarded as a section of 
the same line bundle $\ms D(l)|_{E_1}-(d_2-2)_+C\upone_1$
defined over the curve $\Xi_1$
because the curves $\Xi_1$ and $C\upone_1$ are disjoint.




We show that any section of the line bundle $\ms D(l)|_{E_1} - (d_2-2)_+ C\upone_1$ defined on the curves
$\Xi_1\cup\ol \Xi_k$ extends to $E_1$ 
in a unique way.
The three classes $(1,0),(0,1)$ and $\DDD\upone$ generate $\Pic\,E_1$.
For the intersection number of $\bm D$ with these classes are
respectively given by 
\begin{gather*}
\bm D.\,(1,0) = \bm D.\,\ol\Xi_k = -1 {\text{ by \eqref{int101}}},\\
\bm D.\,(0,1) = 1\qandq \bm D.\,\DDD\upone= d_2-1.
\end{gather*}
From these we obtain 
\begin{align}\label{DE1}
\ms D(l)|_{E_1}\simeq (1,l-1)  - (d_2 - 1) \DDD\upone.
\end{align}
Since $C\upone_1\simeq (0,1)-\DDD\upone$ , we have
\begin{align}\label{DE2}
\ms D(l)|_{E_1} -  (d_2-2)_+ C\upone_1
\simeq
\big (1,l-1-  (d_2-2)_+\big)  - \eee_1\DDD\upone
\end{align}
where $\eee_1$ is either $1$ or $0$ depending on whether
$d_2>1$ or $d_2=1$ respectively.
Subtracting from this the class
$$
\Xi_1 + \ol \Xi_k \simeq (2,0) - \DDD\upone,
$$
we obtain that the kernel sheaf of the restriction homomorphism
from $\ms D(l)|_{E_1} - (d_2-2)_+ C\upone_1$ to
$\Xi_1\cup\ol \Xi_k$ is 
$$
\big (-1,l-1-  (d_2-2)_+\big)  +(1-\eee_1)\DDD\upone.
$$
By using that $1-\eee_1$ is either 0 or 1,
we obtain that any cohomology group of this class vanishes.
Therefore again the unique extendability
to $E_1$ holds.

To finish the proof of the isomorphicity for the present small resolution, 
we consider the pair $\big(t|_{\Xi_1}, u'\big)$.
As above, this can be regarded as an element
of the line bundle 
$\ms D(l)|_{E_1}-(d_2-2)_+C\upone_1$
defined over the curves $\Xi_1\cup\ol\Xi_k$.
By the unique extendability just proved,
this pair extends to $E_1$ in a unique way.
Multiplying a defining section of 
the divisor $(d_2-2)_+C\upone_1$ to the extension,
we obtain an element of $H^0\big(\ms D(l)|_{E_1}\big)$.
On the curves $\ol\Xi_k$ and $\Xi_1$, this section agrees with the original ones $t'|_{\ol\Xi_k}$
and $t|_{\Xi_1}$ respectively.
Thus we have obtained an extension of 
the section $(t,t')$
to the ridge component $E_1$, and it is unique.
By reality, we also obtain a unique extension to 
another component $\ol E_1$.
Hence we have shown the isomorphicity of 
the restriction homomorphism in the proposition
when the factor $\zeta\upone$ blows up 
the component $E_1$.

If the small resolution $\zeta\upk$ blows up $E_1$ and $\ol E_1$
but $\zeta\upone$ does not blow up $E_1$ and $\ol E_1$,
we obtain the isomorphicity in the same way to the case of 
the last resolution, exchanging the role of 
$\CCC\upone_E$ and $\CCC\upk_E$.

Finally, if the small resolutions $\zeta\upone$
and $\zeta\upk$ blow up $E_1$ and $\ol E_1$,
the two exceptional curves $\DDD\upone$ and $\ol\DDD\upk$ are inserted
in $E_1$, and among sub-divisors of $\CCC_E$,
the parts $(d_2-2)_+C\upone_1$
and $(d_k-2)_+  C\upk_1$ are included in $E_1$.
The idea for proving the unique extendability
is the same as above, and for the present resolution
it is enough to show that 
any section of the line bundle $\ms D(l)|_{E_1} - 
(d_2-2)_+ C\upone_1 - (d_k-2)_+  C\upk_1$ 
defined on the curves
$\Xi_1\cup\ol \Xi_k$ extends to $E_1$.
This time, the four classes $(1,0),(0,1),\DDD\upone$ and $\ol\DDD\upk$ 
generate $\Pic\,E_1$, and we have
\begin{gather*}
C\upone_1\simeq (0,1)-\DDD\upone,\quad
\ol C\upone_k\simeq (0,1)-\ol\DDD\upk,\\
\Xi_1 \simeq (1,0)-\DDD\upone,\quad
\ol\Xi_k \simeq (1,0)- \ol\DDD\upk.
\end{gather*}
Also we have
\begin{gather*}
\bm D.\,(1,0)=-1,\quad
\bm D.\,(0,1)=1,\\
\bm D.\,\DDD\upone = d_2 -1,\quad
\bm D.\,\ol\DDD\upk = d_k -1.
\end{gather*}
From these, we obtain 
\begin{multline}\label{DE3}
\ms D(l)|_{E_1} - 
(d_2-2)_+ C\upone_1 - (d_k-2)_+  C\upk_1\\
\simeq
\big (1,l-1-  (d_2-2)_+
 -  (d_k-2)_+\big)  - \eee_1\DDD\upone
 - \eee_k\ol\DDD\upk,
\end{multline}
where 
\begin{align}\label{e1e2}
\eee_1 = \begin{cases}
1 & d_2>1\\
0 & d_2=1
\end{cases}
\qandq
\eee_k = \begin{cases}
1 & d_k>1\\
0 & d_k=1.
\end{cases}
\end{align}
Subtracting $\Xi_1+\ol\Xi_k\simeq (2,0) - \DDD\upone - \ol\DDD\upk$
from \eqref{DE3},
we obtain that the kernel sheaf of the restriction homomorphism
from $\ms D(l)|_{E_1} - (d_2-2)_+ - (d_k-2)_+ C\upone_1$ to
$\Xi_1\cup\ol \Xi_k$ is 
$$
\big (-1,l-1-  (d_2-2)_+\big)  +(1-\eee_1)\DDD\upone
+  (1-\eee_k)\ol\DDD\upk.
$$
Again by using that  $1-\eee_1$ and $1-\eee_k$ are either 0 or 1,
any cohomology groups of this class are easily seen to vanish.
This implies the required unique extendability.
\proofend

\medskip



Now we are able to determine the sheaf $\tif_*\ms D|_E$.
Recall that the number $e$ is defined by 
the property $\ms D|_{E-E_1-\ol E_1} - \CCC_{E-E_1-\ol E_1}
\simeq\tif^*\ms O_{\CP^1}(-e)$ as in \eqref{lsb2}.

\begin{proposition}\label{prop:e2}
We have an isomorphism
$$\tif_*\ms D|_E\simeq\ms O_{\CP^1}(-e)^{\oplus 2}.$$
\end{proposition}

%
\proof
As we showed at the end of Section \ref{ss:dim}, 
the direct image 
$\tif_*\ms D|_E$ is of rank two and locally free.
Hence so is the sheaf $\tif_*\ms D(e)|_E\simeq \big(\tif_*\ms D|_E\big)
\otimes\ms O(e)$.
We write
$\tif_*\ms D(e)|_E\simeq\ms L_1\oplus\ms L_2$
with invertible sheaves $\ms L_1$ and $\ms L_1$
over $\CP^1$.
Since $h^0\big(\ms D(e)|_E\big)=2$ from
Proposition \ref{prop:rest3}, by exchanging $\ms L_1$ and $\ms L_2$
if necessary,
either $\ms L_1\simeq\ms L_2\simeq\ms O_{\CP^1}$,
or $\ms L_1\simeq\ms O_{\CP^1}(1)$ and $\deg\ms L_2<0$
holds.
Since $\ms D$, $\tif$ and $E$ are real, the sheaf
$\tif_*\ms D(e)|_E$ is equipped with a real structure,
and taking a composition with the complex conjugation acting on each fiber,
we obtain a holomorphic involution on $\tif_*\ms D(e)|_E\simeq\ms L_1\oplus\ms L_2$.
We show that $\ms L_1\simeq\ms L_2\simeq\ms O_{\CP^1}$ holds and the isomorphism
$\tif_*\ms D(e)|_E\simeq\ms O\oplus\ms O$
 can be taken in a way that the involution is exactly switching involution
of the two factors.

Let $\omega_0\in H^0\big(\ms D(e)|_E\big)$ be
a non-zero element which vanishes identically on 
the union $E_2\cup\dots\cup E_k$.
By \eqref{om1d} and Proposition \ref{prop:rest3},
such an element exists and is unique up to constant.
By abuse of notation we write $\ol\omega_0$ for the element
$\ol{\sigma^*\omega_0}\in H^0\big(\ms D(e)|_E\big)$.
This vanishes identically on $\ol E_2\cup\dots\cup \ol E_k$.
The two sections $\omega_0$ and $\ol\omega_0$ are clearly linearly independent
and form a basis
of $H^0\big(\ms D(e)|_E\big)$ from the isomorphism \eqref{isdl}.
Via the  isomorphisms
$H^0\big(\ms D(e)|_E\big)\simeq 
H^0\big(\tif_*\ms D(e)|_E\big)
\simeq H^0(\ms L_1\oplus\ms L_2)$,
the sections $\omega_0$ and $\ol\omega_0$
can be regarded as sections of $\ms L_1\oplus\ms L_2$.

For a point $\lmd\in\CP^1$,
let $\big(\tif_*\ms D(e)|_E\big)_{\lmd}$ be
the stalk of the sheaf $\tif_*\ms D(e)|_E$
at the point $\lmd$, and write $S_{\lmd}=\tif\inv(\lmd)$.
From definition of the direct image sheaf,
elements of this stalk are represented
by sections of the line bundle $\ms D(e)|_E$,
defined in neighborhoods in $E$ of $S_{\lmd}\cap E$.
Let $\lmd\in\CP^1$ be a point for which the fiber $S_{\lmd}$ is smooth.
Then the intersection $S_{\lmd}\cap E$ is
the anti-canonical cycle $C$ on the surface $S_{\lmd}$,
and from the basic intersection numbers \eqref{int0011},
we have $D.\,C_1=D.\,\ol C_1 = 1$ 
and $D.\,C_i=D.\,\ol C_i=0$ for any other index $i$.
These imply that, 
the line bundle $\ms D$ is trivial over
the two chains $C-C_1-\ol C_1$, and moreover, sections defined
over these two chains uniquely extend to whole of the cycle $C$.
Namely we have a natural identification
$H^0\big(\ms D|_{S_{\lmd}\cap E}\big)\simeq\CC^2$.
From this, for such a point $\lmd$, we obtain a natural isomorphism
\begin{align}\label{gst}
\Big(\tif_*\ms D(e)|_E\Big)_{\lmd}
\simeq\ms O_{\lmd}\oplus\ms O_{\lmd},
\end{align}
where $\ms O_{\lmd}$ denotes the local ring (i.e.\,the space of convergent
power series)
at the point $\lmd\in\CP^1$.

We regard $\omega_0$ and $\ol\omega_0$ as sections of 
the direct image sheaf $\tif_*\ms D(e)|_E$,
and let $(\omega_0)_{\lmd}$ and $\big(\ol\omega_0\big)_{\lmd}$
be the germs represented by these sections at the point $\lmd$
for which $S_{\lmd}$ is smooth.
These are elements of the stalk
$\big(\tif_*\ms D(e)|_E\big)_{\lmd}$.
From the choice of $\omega_0$,
since the real structure switches the two 
connected components of $C-C_1-\ol C_1$,
the two germs
$(\omega_0)_{\lmd}$ and $\big(\ol\omega_0\big)_{\lmd}$ 
can be identified with two elements
$(\varphi,0)$ and $(0,\varphi)$ respectively under the isomorphism
\eqref{gst},
where $\varphi$ is the germ of a holomorphic function which does not vanish 
at the point $\lmd$.
In particular, the two germs
$(\omega_0)_{\lmd}$ and $\big(\ol\omega_0\big)_{\lmd}$ 
generate
the stalk $\big(\tif_*\ms D(e)|_E\big)_{\lmd}$.

Suppose that 
$\ms L_1\simeq\ms O_{\CP^1}(1)$ and $\deg\ms L_2<0$.
This implies $H^0(\ms L_2)=0$.
Hence the above germs $(\omega_0)_{\lmd}$ and $\big(\ol\omega_0\big)_{\lmd}$ have to belong to the sub-module $(\ms L_1)_{\lmd}$.
Therefore they cannot generate the stalk
$\big(\tif_*\ms D(e)|_E\big)_{\lmd}$.
This is a contradiction.
Hence $\ms L_1\simeq\ms L_2\simeq\ms O$.
Namely $\tif_*\ms D(e)|_E\simeq\ms O^{\oplus 2}$.
This is equivalent to $\tif_*\ms D|_E\simeq\ms O(-e)^{\oplus 2}$.
Further, the holomorphic involution on $\tif_*\ms D(e)\simeq\ms O\oplus\ms O$ switches of
the two factors under the isomorphism
\eqref{gst}.
\proofend

\medskip
From Proposition \ref{prop:e},
this proposition in particular means that 
the sheaf $\tif_*\ms D|_E$  depends on the small resolution $\zeta:\tiZ\to\hat Z$.
Also, the dependency occurs only from the 
resolutions $\zeta\upone$ and $\zeta\upk$,
 and those for all other resolutions 
$\zeta\upi$, $1<i<k$,
are irrelevant.
Notice also that the dependency
occurs only if
$d_2>2$ or $d_k>2$ holds,
because otherwise $(d_2-2)_+=(d_k-2)_+=0$.

As in the proof of Proposition \ref{prop:e2}, 
let $\omega_0$ be a section of
$\ms D(e)|_E$ which vanishes on $E_2\cup\dots\cup E_k$.
This is unique up to a constant and the zero divisor on $\ol E_2\cup\dots\cup \ol E_k$ coincides with
$\CCC_{\ol E_2\cup\dots\cup \ol E_k}$ (see \eqref{lsb3}).
The sections $\omega_0$ and $\ol\omega_0$
form a basis of $H^0\big(\ms D(e)|_E\big)$.
Then we have the following identification 
of the space $H^0\big(
\ms D(l)|_E
\big)$ when $l\ge e$.
\begin{proposition}\label{prop:DEl}
If $l\ge e$, we have
\begin{align}\label{DEl}
H^0\big(
\ms D(l)|_E
\big)
= \big\langle \omega_0,\ol\omega_0\big\rangle
\otimes 
H^0\big(
\tif^*\ms O_{\CP^1}(l-e)
\big).
\end{align}
\end{proposition}

\proof
Since $\tif_*\ms D|_E\simeq\ms O(-e)^{\oplus 2}$
by Proposition \ref{prop:e2},
we have
\begin{align*}
H^0\big( \ms D(l)|_E \big)
\simeq
H^0\big(\tif_* \ms D|_E \otimes\ms O_{\CP^1}(l)\big)
\simeq
H^0\big( \ms O_{\CP^1}(l-e)^{\oplus 2} \big)
\simeq\CC^{2(l-e+1)}.	
\end{align*}
The right-hand side of \eqref{DEl} is also 
$2(l-e+1)$-dimensional. 
Hence the conclusion follows.
\proofend

\medskip
By using the concrete forms for the value of $e$,
we can give a lower bound for it.

\begin{proposition}\label{prop:lbd}
For any total small resolution $\zeta:\tiZ\to\hat Z$ preserving the 
real structure, we have  equalities
\begin{align}\label{ineqb0}
e\ge n-2\qandq e\ge d.
\end{align}
\end{proposition}

\proof
For the first estimate in \eqref{ineqb0},
the case $n=3$ is obvious since $e$ is always positive
by Proposition \ref{prop:e}.
Among the numbers for $e$ in the proposition,
the number \eqref{case1} is the smallest.
Hence, for any small resolution, we have an inequality
\begin{align}\label{ineqd0}
e\ge d_2+ 
\sum_{1<j< k} (d_{j+1} - d_j)_+.
\end{align}
If $n=4$, from Table \ref{table:2}, 
we easily obtain that the right-hand side
of this inequality 
is always two. Hence the former estimate
in \eqref{ineqb0} holds.
Obviously the right-hand side of \eqref{ineqd0} increases
by at least one if we raise $n$ by one.
Hence the inequality $e\ge n-2$ follows by induction on $n$.

For the second estimate in \eqref{ineqb0}, 
let $\nu\in\{1,2,\dots, k-1\}$ be any number
for which $d_{\nu+1} = d = \max \{d_1,\dots, d_k\}$ holds.
Then we have the estimate
\begin{align*}
{\text{\big(RHS of \eqref{ineqd0}\big)}} &\ge d_2  +
\sum_{1<j< \nu} (d_{j+1} - d_j)_+\\
&\ge d_2  +
\sum_{1<j< \nu} (d_{j+1} - d_j)\\
& = d_{\nu+1} = d.
\end{align*}
So we obtain $e\ge d$.
\proofend

\medskip
Finally in this subsection,
we make some remark on the formula for the number $e$ given 
in Proposition 
\ref{prop:e}.
Let $\zeta:\tiZ\to\hat Z$ be a total small resolution,
which satisfies the property
in the second item in Proposition \ref{prop:e}.
If we change  the order of the numbering for 
the components of the cycle $C$ by 
the replacement
\begin{align}\label{oc}
(d_1,d_2,d_3,\dots, d_{k-1},d_k)\longmapsto (d_1,d_{k}, d_{k-1},\dots,d_3, d_2)
\end{align}
(namely fix $d_1$ and reverse the order of
the numbering for 
the remaining components)
then the two resolutions $\zeta\upone$ and $\zeta\upk$ are exchanged and consequently the total 
resolution $\zeta$ becomes to satisfy
the property in the third item in Proposition 
\ref{prop:e}.
Therefore, the two numbers for $e$ in the second and the third
items have to be exchanged
under the change \eqref{oc}.
Similarly, if the total resolution $\zeta$ satisfies the property
in the first (resp.\,forth) item in Proposition \ref{prop:e} under the change of the numbers
\eqref{oc},
it still satisfies the same property.
Therefore the two numbers for $e$ in these two items have
to be invariant under the change \eqref{oc}.
All these follows from the relation
\begin{align}\label{}
d_2 + \sum_{1<j< k} (d_{j+1} - d_j)_+
=
d_k + \sum_{1<j< k} (d_{j} - d_{j+1})_+,
\end{align}
which holds for arbitrary sequence $d_1,d_2,\dots, d_k$
of numbers.

\section{Study on the direct image sheaf}
\label{s:sd}
In the previous two subsections, we mainly worked on 
the exceptional divisor $E$ and determined the direct image
sheaf $\tif_*\ms D|_E$.
In this section, we work on the space $\tiZ$,
and investigate the direct image $\tif_*\ms D$.
As a result we will prove Theorem \ref{thm:main1}
except the concrete form of the quartic hypersurface
that cuts out the branch divisor of the double covering map
from the scroll of planes.
We also discuss various concrete examples.

\subsection{The sheaf $\tif_*\ms D$}
\label{ss:dim1}
For each index $i\in\{1,\dots,k\}$,
we write
\begin{align}\label{rhoi2}
\rho\upi:\big(\tif_*\ms D\big)_{\lmd\upi}\lra
\big(\tif_*\ms D|_E\big)_{\lmd\upi}
\end{align}
for the restriction of
the homomorphism $\rho:\tif_*\ms D\to\tif_*\ms D|_E$
to the stalks over the point $\lmd\upi$.
This is a homomorphism from a
free $\ms O_{\lmd\upi}$-module
of rank three to that of rank two,
where as before we are writing $\ms O_{\lmd}$
for the local ring at a point $\lmd\in\CP^1$.
Since the analogous homomorphism is surjective 
over any $\lmd\neq\lmd\upi$ by Proposition \ref{prop:d_im4},
for the cokernel sheaf
of $\rho$, we have
\begin{align}\label{}
\ms Coker(\rho)\simeq
\bigoplus _{i=1}^k \Coker\, \rho\upi.
\end{align}
Therefore,
the derived exact sequence \eqref{der01} becomes the following form.
\begin{align}\label{der024}
0 \lras \ms O_{\CP^1} \lras 
\tilde f_*\ms D \stackrel{\rho}{\lras} 
\ms O_{\CP^1}(-e)^{\oplus 2} 
{\lras} 
\bigoplus _{i=1}^k \Coker\, \rho\upi
{\lras} 
0.
\end{align}

Next
we derive the following constraint for the 
direct image sheaf $\tif_*\ms D$.

\begin{proposition}\label{prop:d_im06}
We have 
\begin{align}\label{d_im06}
\tilde f_*\ms D
\simeq
\ms O_{\CP^1}\oplus \ms O_{\CP^1}(-m)^{\oplus 2}
\end{align}
for a positive integer $m$ which is not smaller than the number $e$
in Proposition \ref{prop:e2}.
\end{proposition}

Proposition \ref{prop:indep0} means that
the integer $m$ does not depend on the choice
of the total small resolution $\zeta:\tiZ\to\hat Z$
preserving the real structure.

\medskip\noindent
{\em Proof of Proposition \ref{prop:d_im06}.}
Let $\ms F$ be the image sheaf of the
homomorphism  $\rho$.
This is a coherent subsheaf of 
$\ms O_{\CP^1}(-e)^{\oplus 2}$, and is real because $\tif, E$ and $\ms D$
are real.
Obviously $\rank\,\ms F=2$.
So we can write 
$
\ms F\simeq\ms O(-m_1)\oplus\ms O(-m_2)
$
for some integers $m_1$ and $m_2$.
From the inclusion $\ms O(-m_1)\oplus \ms O(-m_2)
\hookrightarrow\ms O(-e)^{\oplus 2}$,
we have $m_1\ge e$ and $m_2\ge e$. 
Since $e>0$, it follows $m_1,m_2>0$. Therefore
the short exact sequence
$0 \to \ms O
\to \tilde f_*\ms D
\to  \ms O(-m_1)\oplus\ms O(-m_2)
\to 0$ splits, and 
we obtain $\tif_*\ms D\simeq\ms O(-m_1)\oplus\ms O(-m_2)\oplus \ms O$.
It seems difficult to conclude $m_1=m_2$ at this stage even with the real structure, and we 
derive a contradiction supposing $m_1>m_2$.

We consider the exact sequence
\begin{align*}
0 \lras \ms D(m_1)\otimes\ms O_{\tiZ}(-E)
\lras  \ms D(m_1)
\lras  \ms D(m_1)|_E
\lras 0
\end{align*}
and the associated exact sequence
\begin{align}\label{lesE}
0 \lras  H^0\big(\ms D(m_1)\otimes\ms O_{\tiZ}(-E)
\big)
\lras  H^0\big(\ms D(m_1)\big)
\stackrel{r_E}\lras  H^0\big(\ms D(m_1)|_E\big).
\end{align}
By the same reason to $\big|\ms D(-E)
\big|=\{\bm D-E\}$ shown in the proof of Proposition \ref{prop:d_im4}, 
we have an equality
$\big|\ms D(m_1)\otimes\ms O_{\tiZ}(-E)
\big|= (\bm D - E) + \big|\tif^*\ms O(m_1)\big|$.
Hence $h^0\big(\ms D(m_1)\otimes\ms O_{\tiZ}(-E)
\big)= m_1 + 1$.
We also have isomorphisms
\begin{align}\label{}
H^0\big(\ms D(m_1)\big)
\simeq
H^0\big(\tif_*\ms D(m_1)\big)
\simeq
H^0\big(\ms O_{\CP^1}\oplus\ms O_{\CP^1}(m_1-m_2)\oplus\ms O_{\CP^1}(m_1)\big),
\end{align}
and this means
$h^0\big(\ms D(m_1)\big) = 2m_1-m_2+3$.
Let $\tilde\Phi_{m_1}$ be the meromorphic map from 
$\tiZ$ induced by the system $\big|\ms D(m_1)\big|$.
Then from the sub-system $\bm D+\big|\tif^*\ms O(m_1)\big|$ of $\big|\ms D(m_1)\big|$, 
we have the commutative diagram of meromorphic maps
\begin{align}\label{diagram44}
 \xymatrix{ 
\tilde Z \ar^{\hspace{-20pt}\tilde\Phi_{m_1}}[r] \ar[d] &\CP^{2m_1-m_2+2} \ar^{p}[dl]\\
\CP^{m_1}
 }
\end{align}
Here, $p$ is the linear projection induced 
by the last sub-system, and  $\tiZ\to\CP^{m_1}$
is the meromorphic map associated to
the sub-system. 
The image of the last map is a rational normal curve
of degree $m_1$.
From the exact sequence \eqref{lesE},
the center of the projection $p$ is 
identified with the dual projective
space of the vector space Image\,$(r_E)$,
and we have
\begin{align}\label{dimrE}
\dim {\rm Image}\,(r_E) &= 
h^0\big(\ms D(m_1)\big)- 
h^0\big(\ms D(m_1)\otimes\ms O_{\tiZ}(-E)
\big)\notag\\
&= (2m_1-m_2+3) - (m_1 + 1) \notag\\
&= m_1 - m_2 + 2.
\end{align}
The image $\tilde\Phi_{m_1}(E)$
is included in the dual projective space
of ${\rm Image}\,(r_E)$;
namely the center of the projection $p$.

So far we have not used the assumption $m_1>m_2$.
Now if $m_1>m_2$, we have
$\dim {\rm Image}\,(r_E)>2$ from \eqref{dimrE}.
By using this we prove that the map $\tilde\Phi_{m_1}$ separates
fibers of $\tif|_{E_2\cup\dots\cup E_k}: E_2\cup
\dots\cup E_k\to\CP^1$
(and so does for fibers of $\tif|_{\ol E_2\cup\dots\cup \ol E_k}:\ol E_2\cup\dots\cup\ol E_k\to\CP^1$).
This property contradicts the commutativity 
of the diagram \eqref{diagram44}
since $\tilde\Phi_{m_1}$ maps the divisor $E$ 
to the center of the projection $p$ as above. 
Therefore the equality $m_1=m_2$ has to
hold, and the proof of the proposition will be over by putting $m:=m_1\,(=m_2)$.

We write $r_{\Xi_1}$ and 
$r_{\ol\Xi_1}$ for the restriction homomorphism
from $H^0\big(\ms D(m_1)\big)$ to the curves
$\Xi_1=E_1\cap E_2$ and $\ol\Xi_1=\ol E_1\cap \ol E_2$ respectively.
If $\dim{\rm Image}\,(r_{\Xi_1})=0$,
we have $\dim{\rm Image}\,(r_{\ol\Xi_1})=0$
from the real structure.
From the basic intersection numbers \eqref{int0011},
these mean that whole of the divisor $E$ is  fixed component
of the system $\big|\ms D(m_1)\big|$ and 
this contradicts $h^0\big(\ms D(m_1)\big)=
2m_1-m_2+3\ge 2m_1-m_1+3=m_1+3$
which implies that $E$ is not  fixed component of 
$\big|\ms D(m_1)\big|$.
Hence ${\rm Image}\,(r_{\Xi_1})\neq 0$.
Suppose $\dim{\rm Image}\,(r_{\Xi_1})=1$,
so that $\dim{\rm Image}\,(r_{\ol\Xi_1})=1$.
Since the restriction homomorphism
$$
H^0\big(\ms D(l)|_{E_2\cup\dots\cup E_k}\big)\lras
H^0\big(\ms D(l)|_{\Xi_1}\big)
$$
is injective for any $l$ as $D.\, C_i=0$ for any $i\neq 1$,
this implies 
that the restriction homomorphism
$H^0\big(\ms D(m_1)\big)\to
H^0\big(\ms D(m_1)|_{E_2\cup\dots\cup E_k}\big)$
is also 1-dimensional.
Moreover, the image of this homomorphism
is in the subspace $\Omega(m_1)$
of $H^0\big(\ms D(m_1)|_{E_2\cup\dots\cup E_k}\big)$.
(See \eqref{om1} for the definition of 
this subspace.)
By Proposition \ref{prop:rest3},
this means that 
$\dim{\rm Image}\,(r_{E})=1+1=2$.
If $m_1>m_2$, this contradicts \eqref{dimrE}.
Therefore $\dim{\rm Image}\,(r_{\Xi_1})>1$.
This implies that the images of the restriction maps
from $H^0\big(\ms D(m_1)\big)$ to the
divisors $E_2\cup \dots \cup E_k$
and $\ol E_2\cup \dots \cup \ol E_k$ are
more than $1$-dimensional.
This means that the map $\tilde\Phi_{m_1}$
separates
fibers of $\tif|_{E_2\cup\dots\cup E_k}: E_2\cup
\dots\cup E_k\to\CP^1$
and $\tif|_{\ol E_2\cup\dots\cup \ol E_k}:
\ol E_2\cup
\dots\cup \ol E_k\to\CP^1$ respectively,
and this is what we need to show.
\proofend

%

\medskip
Thus the direct image $\tilde f_*\ms D$ takes a simple form as in \eqref{d_im06}.
From now on, the letter $m$ always means the one in 
\eqref{d_im06}.
This will be the number $m$ in Theorem \ref{thm:main1}.
From Proposition \ref{prop:d_im06},
the derived exact sequence becomes
\begin{align}\label{der07}
0 \lras \ms O_{\CP^1} \lras 
\ms O_{\CP^1}\oplus \ms O_{\CP^1}(-m)^{\oplus 2}
\stackrel{\rho}{\lras} 
\ms O_{\CP^1}(-e)^{\oplus 2} 
{\lras} 
\bigoplus _{i=1}^k \Coker\, \rho\upi
{\lras} 
0
\end{align}
with $0<e\le m$.
From this we obtain an equality
\begin{align}\label{mec2}
2(m-e) = \sum_{i=1}^k\dim\Coker\,\rho\upi.
\end{align}

By the inequality $m\ge e$,
the  estimate $e\ge n-2$ in \eqref{ineqb0}
means the inequality
\begin{align}\label{ineqc}
m\ge n-2.
\end{align}
This will prove the inequality in Theorem \ref{thm:main1}.

From Proposition \ref{prop:d_im06},
the dimension formula for
the space of global sections of the line bundle $\ms D(l)$
can be easily obtained:

\begin{proposition}\label{prop:df7}
Let $m$ be the positive integer as in Proposition \ref{prop:d_im06}.
Then for any integer $l\ge 0$, we have
\begin{align}\label{dimm0}
h^0\big(\tilde Z, \ms D(l)\big)
= 
\begin{cases}
l+1 & {\text{if $l<m$}},\\
3l - 2m + 3 & {\text{if $l\ge m$}}.
\end{cases}
\end{align}
In particular, we have
\begin{align}\label{dimm}
h^0\big(\tilde Z, \ms D(m)\big)
= 
m + 3.
\end{align}
\end{proposition}

\proof
For any $l$,
taking the direct image under $\tilde f$, we obtain
\begin{align*}
H^0\big(\tilde Z, \ms D(l)\big)
&\simeq
H^0\big(\CP^1, \tilde f_*\ms D \otimes \ms O(l)\big)\\
&\simeq H^0\big(\CP^1, \ms O(l)\oplus
\ms O(l-m)^{\oplus 2}\big),
\end{align*}
where the second isomorphism is due to 
Proposition \ref{prop:d_im06}.
From this, the conclusion \eqref{dimm0} easily follows.
\proofend

\medskip
We note that by \eqref{dimm}, the system
$\big|\ms D(m)\big|$ is strictly larger
than the system $\bm D + |\tif^*\ms O(m)|$,
and from Proposition \ref{prop:l<m}, this means the inequality 
\begin{align}\label{dm}
d\le m.
\end{align}

By using Proposition \ref{prop:df7},
we next describe the space of
sections of the line bundle $\ms D(m+l)$
for arbitrary $l\ge 0$.
It will be used in Section \ref{ss:redpf}
for proving the existence of  real 
reducible members
of the linear system $\big|\ms D(m)\big|$.
To state it, let 
$s_{\sst{\bm D}}$ be a section
of the line bundle $\ms D$ which defines
the effective divisor $\bm D$.
Then 
\begin{align}\label{pg}
H^0\big(\tiZ,\tilde f^*\ms O(m)\big)\otimes s_{\sst{\bm D}}
=\Big\{s\otimes s_{\sst{\bm D}}\,\Big|\,
s\in H^0\big(\tiZ,\tilde f^*\ms O(m)\big)\Big\}
\end{align}
is an $(m+1)$-dimensional linear subspace
of $H^0\big(\tiZ,\ms D(m)\big)\simeq\CC^{m+3}$.
The corresponding sub-system is 
$\bm D + \big|\tif^*\ms O(m)\big|$.

\begin{proposition}\label{prop:fvr} 
Let $w_0$ and $w_1$ be
elements of $H^0\big(\ms D(m)\big)$ which generate,
together with the $(m+1)$-dimensional subspace
\eqref{pg}, the space $H^0\big(\ms D(m)\big)$.
Then for any non-negative integer $l$, 
we have
\begin{align}\label{fvr}
H^0\big(
\ms D(m+l)\big)
=
H^0\big(\tilde f^*\ms O(m+l)\big)\otimes s_{\sst{\bm D}}
\quad\bigoplus\quad
H^0\big(\tilde f^*\ms O(l)\big)\otimes
\langle
w_0,w_1
\rangle.
\end{align}
\end{proposition}

\proof
By Proposition \ref{prop:df7} we have
$$
h^0\big(
\ms D(m+l)\big) = 
m+3l+3.
$$
On the other hand we have
$h^0\big(\tilde f^*\ms O(m+l)\big) = 
h^0\big(\ms O(m+l)) = m+l+1,$ and
$$
\dim H^0\big(\tif^*\ms O(l)\big)\otimes
\langle
w_0,w_1
\rangle
= 2h^0\big(\tilde f^*\ms O(l)\big) = 2(l+1).
$$
From these it follows that, if the intersection of the two direct summands in the right-hand
side of \eqref{fvr} is zero, then the equality 
\eqref{fvr} holds.
If the intersection were non-zero, 
there would exist $g_0,g_1\in H^0\big(\tif^*\ms O(l)\big)\simeq H^0\big(\ms O_{\CP^1}(l)\big)$
and $g\in H^0\big(\tif^*\ms O(m+l)\big)
\simeq H^0\big(\ms O_{\CP^1}(m+l)\big)$
with $(g_0,g_1)\neq (0,0)$ and $g\neq 0$,
such that a linear relation
\begin{align}\label{lrel2}
g_0 w_0 + g_1 w_1 = g s_{\sst{\bm D}}
\end{align}
holds in $H^0\big(\bm D(m+l)\big)$.

Since $\deg g >\deg g_1=\deg g_2$ from the 
inequalities $m>0$ and $l\ge 0$,
by dividing both sides of \eqref{lrel2}
by common factors if necessary,
we may suppose that there is a point $\lmd\in\CP^1$
such that $g(\lmd)=0$ but $
\big(g_1(\lmd),g_2(\lmd)\big)\neq (0,0)$.
By restricting the relation \eqref{lrel2}
to the fiber $S$ over such a point $\lmd$,
we obtain the relation
$$
aw_0|_S + bw_1|_S = 0
$$
in $H^0(S,\bm D|_S)$, where $a$ and $b$ are constants that satisfy $(a,b)\neq (0,0)$.
This means $(aw_0+bw_1)|_S = 0$.
Hence 
\begin{align}\label{lrel1}
aw_0+bw_1\in
H^0\big(\bm D(m-1)\big)\otimes H^0\big(\tif^*\ms O(1)\big).
\end{align}
By Proposition \ref{prop:df7}, we have
$H^0\big(\bm D(m-1)\big) = H^0\big(\tif^*\ms O(m-1)\big)\otimes s_{\sst{\bm D}}$,
and hence the right-hand side of \eqref{lrel1} equals 
$H^0\big(\tif^*\ms O(m)\big)\otimes s_{\sst{\bm D}}$.
So
$aw_0+bw_1\in H^0\big(\tif^*\ms O(m)\big)\otimes s_{\sst{\bm D}}$.
As $(a,b)\neq (0,0)$, this contradicts the choice of $w_0$ and $w_1$.
Therefore if the relation \eqref{lrel2}
holds, we have $g_0=g_1=g=0$.
Namely, the intersection of the two direct summands in 
\eqref{fvr} is zero.
\proofend

\subsection{Double covering structure of the twistor spaces}
\label{ss:dc}
Let $\tilde\Phi_m$ be the meromorphic map associated to 
the linear system $\big|\ms D(m)\big|$.
By letting $m_1=m_2=m$ in the diagram \eqref{diagram44}, we obtain the 
commutative diagram
\begin{align}\label{diagram4}
 \xymatrix{ 
\tilde Z \ar^{\hspace{-10pt}\tilde\Phi_m}[r] \ar[d] &\CP^{m+2} \ar^{p}[dl]\\
\CP^{m}
 }
\end{align}
Here, as before, $p$ is the linear projection 
induced by the sub-system $\bm D+\big|\tif^*\ms O(m)\big|$ (or the subspace
$H^0\big(\tif^*\ms O(m)\big)\otimes s_{\sst{\bm D}}$
in $H^0\big(\ms D(m)\big)$).
The center of the projection $p$ is 1-dimensional.
We denote this line by the bold letter $\bm l$.
The image of the map $\tiZ\to\CP^m$ is a rational normal curve of degree $m$. We denote it by $\Lmd_m$.
This is canonically identified with 
the parameter space of the pencil $|F|$.
By using Proposition \ref{prop:d_im06} and 
the diagram \eqref{diagram4} we are able to show that the present twistor spaces have  structure of double covering 
over the scroll of planes,
for which we give a precise definition.

\begin{definition}\label{def:scroll}{\em 
Let  $m$ be a positive integer.
By a {\em scroll of planes over a rational curve
in $\CP^m$}, we mean the 3-dimensional variety
$$
Y_m:=p\inv(\Lmd_m),
$$
where $p:\CP^{m+2}\to\CP^m$ is a linear projection
and $\Lmd_m\subset\CP^m$ is a rational normal curve
of degree $m$.
The scroll $Y_m$ contains the line $\bm l$
which is the center of the projection $p$.
This line is called the {\em ridge} of the scroll.
\proofend
}
\end{definition}

The space $\CP^m$,  the target space of the projection $p$, can be regarded as the space of 
planes in $\CP^{m+2}$ which contain the ridge $\bm l$.
The scroll $Y_m$ is a union of planes parametrized by 
the curve $\Lmd_m$.
Then we have the following result,
which, together with Theorem \ref{thm:dc2} that will be presented later,
is the main result in the first half of this paper.

\begin{theorem}\label{thm:dc1}
Let $m$ be the positive integer  in Proposition \ref{prop:d_im06},
and 
\begin{align}\label{}
\tilde\Phi_m:
\tilde Z\to\CP^{m+2}
\end{align}
the meromorphic map associated to 
the linear system $\big|\ms D(m)\big|$ on $\tilde Z$ as above.
Then the image $\tilde\Phi_m(\tiZ)$ is the scroll\, $Y_m$, and 
the map $\tilde\Phi_m$ is of degree-two over $Y_m$.
\end{theorem}

Before giving a proof of this proposition,
we give an immediate consequence followed by 
the proposition.
From \eqref{ni} and the inequality $m\ge d$
obtained in \eqref{dm},
we have the natural inclusion 
\begin{align}\label{ni2}
H^0\big(\tilde Z,\ms D(m)\big)\subset H^0(Z,mF).
\end{align}
Hence the system $\big|\ms D(m)\big|$ on $\tiZ$ is identified
with a sub-system of $|mF|$ on the twistor space $Z$, 
and if we denote $\Phi_m:Z\to\CP^{m+2}$ for
the meromorphic map associated to 
this sub-system of $|mF|$
and $\eta:\tiZ\to Z$ for the composition 
of the total small resolution $\zeta:\tiZ\to\hat Z$
preserving the real structure
and the blowing up $\bbb:\hat Z\to Z$ at the cycle $C$, 
we have the commutative diagram
\begin{align}\label{diag}
 \xymatrix{ 
\tilde Z \ar_{\eta}[d] \ar^{\tilde\Phi_m}[dr] \\
Z  \ar_{\Phi_m}[r] & Y_m.
}
\end{align}
Therefore Theorem \ref{thm:dc1} means that 
the system $|mF|$ on $Z$ has an $(m+2)$-dimensional sub-system whose meromorphic map 
is of degree-two over the scroll $Y_m$
as in Theorem \ref{thm:main1}.

\medskip\noindent
{\em Proof of Theorem \ref{thm:dc1}.}
Let $S$ be any smooth fiber of the morphism $\tif:\tiZ\to\CP^1$.
We have isomorphisms
$$
H^0\big(\ms D (m)|_{S}\big)\simeq
H^0(\ms D|_{S})\simeq
H^0(S,D),
$$
and the last space is 3-dimensional.
By Proposition \ref{prop:df7},
we obtain
$$
h^0\big(
\ms D(m)
\big) -
h^0\big(
\ms D (m-1)
\big) = (m+3) - m =3.
$$
These mean that the restriction homomorphism
\begin{align}\label{rest002}
H^0\big(\tilde Z, \ms D (m)\big)
\lras 
H^0\big(S, \ms D(m)|_S\big)
\end{align}
is surjective.



From the commutativity of the diagram \eqref{diagram4},
the image $\tilde{\Phi}_m(\tiZ)$ is contained in the scroll 
$Y_m=p\inv(\Lmd_m)$.
Moreover, from the surjectivity
of the restriction map \eqref{rest002},
the restriction of the map $\tilde\Phi_m$ to
smooth fibers of $\tilde Z\to\Lmd_m$
can be identified with the morphism associated
to the system $|D|$ on the fibers.
The latter morphism is generically two-to-one.
Hence the image of the meromorphic map $\tilde \Phi_m$ is precisely
the scroll $p\inv(\Lmd_m)$,
and $\tilde\Phi_m$ is of degree-two
over the scroll.
\proofend

\medskip
We also obtain the following properties on 
the map $\tilde\Phi_m$.

\begin{proposition}
\label{prop:pim}
The map $\tilde\Phi_m:\tiZ\to Y_m$ maps 
the ridge components $E_1$ and $\ol E_1$ to 
the ridge $\bm l$ of $Y_m$,
and the unions $E_2\cup\dots\cup E_k$
and $\ol E_2\cup\dots\cup\ol E_k$
to distinct points on $\bm l$.
\end{proposition}

\proof
As in the proof of Proposition \ref{prop:d_im06},
from the exact sequence \eqref{lesE},
the image $\tilde\Phi_m(E)$ is included in the 
center $\bm l$ of the projection $p:\CP^{m+2}\to\CP^m$.
By letting $m_1=m_2$ in \eqref{dimrE}, 
the image of the restriction homomorphism
$$
r_E:H^0\big(\ms D(m)\big)
\lras H^0\big(\ms D(m)|_E\big)$$
is 2-dimensional. Hence the image  $\tilde\Phi_m(E)
$ cannot be a point and we obtain $\tilde\Phi_m(E)=\bm l$.
As in the  proof of Theorem \ref{thm:dc1},
for a smooth fiber $S$ of $\tif$, 
the restriction $\tilde\Phi_m|_S$ is identified
with the morphism $\phi:S\to\CP^2$
induced by the system $|D|$ on $S$.
This implies that the line components
$C_1 = E_1\cap S$ and $\ol C_1 = \ol E_1\cap S$
are mapped to the line $\bm l$.
Therefore $\tilde\Phi_m(E_1) = \tilde\Phi_m(\ol E_1) = \bm l$.
Finally, the morphism $\phi$ contracts 
the two chains $C_2\cup\dots\cup C_k$
and $\ol C_2\cup\dots\cup \ol C_k$
to distinct points on the line.
Because $E_1$ and $\ol E_1$ are contracted to $\bm l$,
this implies that the unions $E_2\cup\dots\cup E_k$
and $\ol E_2\cup\dots\cup\ol E_k$
are mapped to points on $\bm l$.
These two points have to be distinct since
$\tilde\Phi_m$ maps $C_1$ and $\ol C_1$ isomorphically to the line.
\proofend

\medskip
In accordance with the notation in Section \ref{s:fs},
we denotes the images of the unions $E_2\cup\dots\cup E_k$
and $\ol E_2\cup\dots\cup\ol E_k$ by $\bm q$ and $\ol{\bm q}$
respectively.
As above, these are distinct points on the ridge $\bm l$.


Next we discuss a relation between two numbers $d$ and $m$.
Recall from \eqref{d} that the number $d$ is the smallest one such that 
the system $\big|dK_S\inv\big|$ on a real irreducible fundamental divisor
$S$ has a sub-system which induces a degree-two morphism 
$\phi:S\to\CP^2$.
For the number $m$,
we have $h^0\big(\ms D(l)\big)=l+1$ if $l<m$ by Proposition \ref{prop:df7},
and by \eqref{ni} this means that 
the linear system $|lF|$ on the original twistor space $Z$ is composed
with the pencil $|F|$.
Hence if $l<m$ the linear system $|lF|$ cannot
induce a degree-two map to the scroll of planes.
This implies that $m$ is the smallest number such that 
the system $|mF|$ on the twistor space has a sub-system
which induces a degree-two map to the scroll of planes.
We note that when we fix the number $n$,
{\em the twistor space is more general if the number $m$ is larger,}
by upper semi-continuity of dimension of the space of sections
of line bundles.

From these, since $F|_S\simeq K_S\inv$, it would be tempting to expect
that 
the equality $m=d$ always holds.
Indeed, the twistor spaces investigated in \cite{Hon_Inv} and 
\cite{Hon_Cre2} satisfy $m=d=n-2$.
In this respect, from the second estimate in \eqref{ineqb0},
we obtain a relation
\begin{align}\label{dem}
d\le e\le m.
\end{align}
Hence discrepancy between $m$ and $d$ can 
occur from a gap between $m$ and $e$,
and also from a gap between $e$ and $d$.
The number $(n-2)$ is not greater than $e$, but any 
one of the three situations
$n-2>d$, $n-2=d$ and $n-2<d$ can happen.
We can readily find concrete examples for which the 
strict inequality $e>d$ holds.
For example, when $n=5$, the number $d$ is two for some $S$
as in Table \ref{table:3}, 
but from the first estimate in \eqref{ineqb0}
we have $e\ge n-2=3$.
So  $e>d$ holds,
and the strict inequality $m>d$ can really happen.
Namely, although a sub-system of $\big|dK_S\inv\big|$ on $S$ induces the double covering map
$\phi:S\to\CP^2$, the system $|dF|$
on the twistor space may be composed with the pencil $|F|$
so that it does not  induce
a degree-two map to the scroll.
On the other hand, it is harder to find examples such that
the strict inequality $m>e$ holds.
We give such examples in Section \ref{ss:ex1},
after investigating the base curves of the system $\big|\ms D(m)\big|$
more closely.

\subsection{The branch divisor on the scroll}
From the proof of Theorem \ref{thm:dc1},
the branch divisor of the map 
$\tilde\Phi_m:\tiZ\to Y_m$ is of degree four
on generic planes of the scroll $Y_m$.
In this subsection we show that the branch divisor is 
 a cut of the scroll $Y_m$
by a quartic hypersurface in $\CP^{m+2}$.
For this purpose, 
we need some preliminary considerations 
about the scroll $Y_m$.

When $m>1$, the scroll $Y_m$ has quotient singularities of order $m$ along the ridge $\bm l$,
and it is resolved by blowing up $\bm l$.
We write $\tilde Y_m$ for the resolution.
We have
$$
\tilde Y_m \simeq \mathbb P\big(\ms O_{\CP^1}(m)^{\oplus 2}\oplus \ms O_{\CP^1}\big).
$$
We write $\Sigma\subset \tilde Y_m$ for the divisor 
$\mathbb P\big(\ms O(m)^{\oplus 2}\big)$ under the above identification.
This is exactly the exceptional divisor
of the resolution $\tilde Y_m \to Y_m$, and
there is a canonical isomorphism
\begin{align}\label{is}
\Sigma\simeq\Lmd_m\times \bm l.
\end{align}
Let $\mathfrak f$ be the fiber class of the $\CP^2$-bundle map
$\tilde Y_m\to\Lmd_m$.
Similarly to the notation $\ms D(l)$ on the space $\tiZ$, we write 
a divisor class  on $\tilde Y_m$ by
$$
\Sigma(l):= \Sigma + l\,\mathfrak f,
\quad l\in\ZZ.
$$
The rational maps induced by the linear systems
$\big|\Sigma(m)\big|$ and $\big|\Sigma(m+1)\big|$ on $\tilde Y_m$ have the following properties.
 
\begin{proposition}
\label{prop:scroll}
(i) 
The linear system $\big|\Sigma(m)\big|$ is base point free and
$(m+2)$-dimensional as
a linear system. The induced morphism $\tilde Y_m\to\CP^{m+2}$ 
is identified with the composition of the contraction map 
$\tilde Y_m\to Y_m$ of $\Sigma$ to $\bm l$ and the embedding $Y_m\hookrightarrow\CP^{m+2}$.
(ii) The linear system $\big|\Sigma(m+1)\big|$ is base point free
and $(m+5)$-dimensional as a linear system.
The induced morphism $\tilde Y_m\to\CP^{m+5}$ is an embedding.
\end{proposition}

\proof
By Leray spectral sequence,
for any $q\ge 0$ and $l\in\ZZ$, we have an isomorphism
\begin{align}\label{Ls}
H^q\big(\tilde Y_m,\ms O(l\,\mathfrak f)\big)
\simeq
H^q\big(\CP^1,\ms O(l)\big).
\end{align}
If we write $\ms O_{\Sigma}(1,0)$ for the class $\mathfrak f|_{\Sigma}$,
then we have $\Sigma|_{\Sigma}\simeq \ms O_{\Sigma}(-m, 1)$.
This implies 
\begin{align}\label{Qm}
\Sigma(m)|_{\Sigma}
\simeq \ms O_{\Sigma}(0,1).
\end{align}
Hence by restricting to $\Sigma$,
we have the exact sequence
\begin{align}\label{exQ1}
0\lras \ms O_{\tilde Y_m}( m\mathfrak f) \lras 
\ms O_{\tilde Y_m}\big(\Sigma(m)\big) \lras 
\ms O_{\Sigma}(0,1) \lras 0.
\end{align}
From \eqref{Ls} and the cohomology exact sequence 
associated to this sequence, we obtain the base point freeness of the system
$\big|\Sigma(m)\big|$ and the desired dimension.
We write $Y'_m\subset \CP^{m+2}$ for the image of the 
morphism induced by $\big|\Sigma(m)\big|$.
Then from the injection 
$H^0\big(\ms O_{\tilde Y_m}( m\mathfrak f)\big) \subset
H^0 \big(\ms O_{\tilde Y_m}(\Sigma + m\mathfrak f)\big)$,
we obtain a projection $Y'_m\to \Lmd_m$.
Hence we
have an inclusion $Y'_m\subset p\inv(\Lmd_m) = Y_m$.
Moreover the restriction of the morphism $\tilde Y_m\to Y'_m$ induced by $\big|\Sigma(m)\big|$
to any fiber of $\tilde Y_m\to\Lmd_m$ is easily seen to be 
an isomorphism.
These imply that the morphism $\tilde Y_m\to Y'_m$
 is 
generically one-to-one.
Furthermore, from \eqref{Qm}, this morphism contracts
the divisor $\Sigma$ to a line.
Hence we obtain $Y'_m=Y_m$, and the morphism $\tilde Y_m\to\CP^{m+2}$
is identified with the composition map as in the proposition.
Thus we obtain the assertion (i).

For the second assertion (ii), the base point freeness of $\big|\Sigma(m+1)\big|$
is obvious because $\big|\Sigma(m)\big|$ 
is base point free and the difference $
\Sigma(m+1) - \Sigma(m) = \mathfrak f$ can move
without base point.
Let $\Psi$ be the morphism induced by 
$\big|\Sigma(m+1)\big|$.
Since this includes $\big|\Sigma(m)\big|$ as
a linear sub-system (by fixing a fiber $\mathfrak f$) and this sub-system induces embedding
at least outside the divisor $\Sigma$,
$\Psi$ is also embedding at least outside $\Sigma$.
It remains to show that $\Psi$ is embedding also
on the divisor  $\Sigma$ and $\dim \big|\Sigma(m+1)\big| = m+5$.

For these,
by multiplying $\mathfrak f$ to \eqref{exQ1}, we have an exact sequence
\begin{align}\label{exQ2}
0\lras \ms O_{\tilde Y_m}\big((m+1)\mathfrak f\big) \lras 
\ms O_{\tilde Y_m}\big(\Sigma  (m+1)\big) \lras 
\ms O_{\Sigma}(1,1) \lras 0.
\end{align}
Again from \eqref{Ls} the cohomology group $H^1$ of the first term vanishes.
Hence we get an exact sequence
\begin{align}\label{exQ3}
0\lras H^0\Big(\ms O_{\tilde Y_m}\big((m+1)\mathfrak f\big)\Big) \lras 
H^0\Big(\ms O_{\tilde Y_m}\big(\Sigma  (m+1)\big)\Big) \lras 
H^0\big(\ms O_{\Sigma}(1,1)\big) \lras 0.
\end{align}
So we obtain $h^0\big(\ms O_{\tilde Y_m}(\Sigma  (m+1))\big)
= (m+2) + 4 = m+6$.
This exact sequence also means that 
the morphism $\Psi$ embeds 
the divisor $\Sigma$ to a smooth quadric in a 3-dimensional
linear subspace of $\CP^{m+5}$.
Hence $\Psi$ is one-to-one over the image.
Therefore to show that $\Psi$ is an embedding,
it is enough to prove that the image $\Psi(\tilde Y_m)$ is non-singular.

We show this by identifying the image
$\Psi(\tilde Y_m)$ concretely.
First, the zero section of 
the vector sub-bundle $\ms O(m)^{\oplus 2}
\subset \tilde Y_m$ is embedded by $\Psi$
to $\CP^{m+1}$ as a rational normal
curve, say $\Lmd_{m+1}$, in a linear subspace $\CP^{m+1}$ in $\CP^{m+5}$.
From the exact sequence \eqref{exQ3}, the last $\CP^{m+1}$ is disjoint from the linear subspace
$\CP^3$ in which the divisor $\Sigma$ is embedded.
Also, there is a natural identification $\Lmd_{m+1}\simeq\Lmd_m$
because both are the base space of the $\CP^2$-bundle map $\tilde Y_m\to\CP^1$.
Each fiber of $\tilde Y_m\to\CP^1$ is 
embedded (by $\Psi$) as a plane in $\CP^{m+5}$.
Then the rational normal curve $\Lmd_{m+1}$ is canonically identified with
one of the two factors 
of the quadric $\Psi(\Sigma)\simeq\qdr$
through the canonical isomorphism $\Sigma\simeq \Lmd_m\times \bm l$ in \eqref{is}
and the above identification $\Lmd_{m+1}\simeq\Lmd_m$.
From these, we obtain that the image $\Psi(\tilde Y_m)$ is 
obtained from the chordal variety
spanned by points on the curve $\Lmd_{m+1}$ and points of the quadric $\Psi(\Sigma)$,
by taking only lines of the following kind:
for any point $\lmd$ of $\Lmd_{m+1}$,
we consider only lines through $\lmd$ and 
 points on the fiber of the projection $\Psi(\Sigma)\simeq\Lmd_m\times \bm l\to\Lmd_m$ over the point $\lmd$ under the isomorphism
$\Lmd_{m+1}\simeq\Lmd_m$.
From this identification of the image,
we obtain that $\Psi(\tilde Y_m)$ is non-singular.
Hence $\Psi$ is  embedding.
\proofend

\medskip
We will also need the following result which gives
a basis of $H^0\Big(\tilde Y_m,\ms O\big(\Sigma(m+1)\big)\Big)$ using a basis of $H^0\Big(\tilde Y_m,\ms O\big(\Sigma(m)\big)\Big)$.
For this, let $u_0,u_1$ be a basis
of $H^0\big(\ms O_{\CP^1}(1)\big)$.
We use the same notations for the pullbacks
of $u_0,u_1$ to $\tilde Y_m$ by the projection
$\tilde Y_m\to\Lmd_m$.

\begin{proposition}
Let $\varsigma\in H^0\big(
\ms O_{\tilde Y_m}(\Sigma)\big)$ 
be a defining section of the divisor $\Sigma$,
and suppose that $v_0$ and $v_1$ are elements of 
$H^0\Big(\tilde Y_m,\ms O\big(\Sigma(m)\big)\Big)$
such that $$\varsigma u_0^m,\varsigma u_0^{m-1}u_1,\dots,\varsigma u_1^m\qandq
v_0,v_1$$ generate $H^0\Big(\tilde Y_m,\ms O\big(\Sigma(m)\big)\Big)$.
Then the collection
\begin{align}\label{basis2}
\varsigma u_0^{m+1},\varsigma u_0^{m}u_1,\dots,\varsigma u_1^{m+1}
\qandq u_0v_0,u_0v_1,u_1v_0,u_1v_1
\end{align}
form a basis of 
$H^0\Big(\tilde Y_m,\ms O\big(\Sigma(m+1)\big)\Big)$.
\end{proposition}

\proof
We use the notation in the previous proof
to represent cohomology classes on $\Sigma\simeq\qdr$.
From the cohomology exact sequence  of the sequence \eqref{exQ1},
the pair of the restrictions $v_0|_{\Sigma}$ and $v_1|_{\Sigma}$
form a basis of $H^0\big(\ms O_{\Sigma}(0,1)\big)$.
On the other hand, another pair of 
the restrictions $u_0|_{\Sigma}$ and $u_1|_{\Sigma}$
form a basis of $H^0\big(\ms O_{\Sigma}(1,0)\big)$. 
These imply that the restrictions of the four sections
$u_0v_0,u_0v_1,u_1v_0,u_1v_1$
to the divisor $\Sigma$ form a basis of 
$H^0\big(\ms O_{\Sigma}(1,1)\big)$.
From the  exact sequences \eqref{exQ3},
this implies that the collection
\eqref{basis2} is a basis of 
$H^0\Big(\tilde Y_m,\ms O\big(\Sigma(m+1)\big)\Big)$.
\proofend

\medskip
In order to identify the branch divisor
of the map $\tilde\Phi_m:\tiZ \to Y_m\subset\CP^{m+2}$ induced
by the system $\big|\ms D(m)\big|$, we take a look at the map 
induced by $\big|\ms D(m+1)\big|$.
We have $\dim \big|\ms D(m+1)\big| = m+5$ from 
Proposition \ref{prop:df7}.
 
\begin{proposition}
Let $\tilde\Phi_{m+1}:\tiZ\to\CP^{m+5}$ be the meromorphic map
associated to the system $\big|\ms D(m+1)\big|$.
Then the image $\tilde\Phi_{m+1}(\tiZ)$ is
isomorphic to the resolution $\tilde Y_m$,
and $\tilde\Phi_{m+1}$ is of degree-two
over $\tilde Y_m$.
Further, the branch divisor of this map
belongs to the linear system 
$\big|4(\Sigma + m\mathfrak f)\big|$ on $\tilde Y_m$.
\end{proposition}

\proof
Let $w_0,w_1\in H^0\big(\ms D(m)\big)$
be as in Proposition \ref{prop:fvr}.
By letting $l=1$ in the proposition, the collection
\begin{align}\label{basis3}
s_{\sst{\bm D}} u_0^{m+1},s_{\sst{\bm D}} u_0^{m}u_1,\dots,s_{\sst{\bm D}} u_1^{m+1}
\qandq u_0w_0,u_0w_1,u_1w_0,u_1w_1
\end{align}
is a basis of $H^0\big(\ms D(m+1)\big)$.
This is the same form as the basis 
\eqref{basis2} of $H^0\Big(\tilde Y_m,\ms O\big(\Sigma(m+1)\big)\Big)$.
Hence by Proposition \ref{prop:scroll} (ii),
the image of the meromorphic map $\tilde \Phi_{m+1}$
is isomorphic to the resolution $\tilde Y_m$.

Since the restriction of $\tilde \Phi_{m+1}$
to each smooth fiber  $S$ of $\tif$ is identified with 
the map $\phi:S\to\CP^2$ induced by the system $|D|$,
the map $\tilde\Phi_{m+1}:\tiZ\to\tilde Y_m$ is again of degree two.
Let $\tilde B\subset\tilde Y_m$ be the branch divisor.
The cohomology group $H^2(\tilde Y_m,\ZZ)$
is free of rank two, and is generated by 
the classes of $\Sigma$ and $\mathfrak f$.
Moreover under  the restriction homomorphism
 $H^2(\tilde Y_m,\ZZ)\to H^2(\Sigma,\ZZ)$
these generators are mapped to 
\begin{align}\label{genm}
\Sigma \longmapsto (-m,1)
\qandq
\mathfrak f \longmapsto (1,0).
\end{align}
From this it follows that 
the last homomorphism is an isomorphism.
Recall that the branch quartic of the 
double covering map $\phi:S\to\CP^2$
induced by $|D|$ is tangent to a (real) line
at two points
or has double points at two points (see Figure \ref{fig:bq}).
Further by Proposition \ref{prop:pim} the line is naturally identified with 
the ridge $\bm l$ of $Y_m$,
and the two points on the line are common for any fiber of $\tif$;
the points $\bm q$ and $\ol{\bm q}$.
These imply that the restriction $\tilde B|_{\Sigma}$
belongs to the class $(0,4)$.
(Recall that $(0,1)$ is the fiber class of 
the projection $\Sigma\simeq\Lmd_m\times \bm l\to \bm l$.)
From \eqref{genm}, the class $\Sigma + m\mathfrak f$
is mapped to the class $(0,1)$.
Hence, since the restriction homomorphism 
$H^2(\tilde Y_m,\ZZ)\to H^2(\Sigma,\ZZ)$ is isomorphic, the divisor 
$\tilde B$ belongs to the class 
$4(\Sigma + m\mathfrak f)$, as desired.
\proofend

\medskip
Thus the relation between  two meromorphic maps $\tilde \Phi_m$
and $\tilde \Phi_{m+1}$ (associated to $\big|\ms D(m)\big|$
and $\big|\ms D(m+1)\big|$ respectively) may be summarized 
in the commutative diagram
\begin{align}\label{diagram5}
 \xymatrix{ 
\tilde Z \ar_{\tilde\Phi_{m+1}}[d] \ar^{\tilde\Phi_m}[dr] \\
\tilde Y_m  \ar[r] & Y_m,
}
\end{align}
where $\tilde Y_m\to Y_m$ is the resolution.
Now we are able to show 

\begin{theorem}\label{thm:dc2}
The branch divisor of the degree-two map
$\tilde\Phi_m:\tiZ\to Y_m\subset\CP^{m+2}$
induced by $\big|\ms D(m)\big|$
 is a cut of the scroll $Y_m$ by
a quartic  hypersurface.
\end{theorem}

\proof
We use the same symbol $\mathfrak f$ to denote
the class on $Y_m$ represented by a plane of the scroll
$Y_m$.
Then the cohomology group $H^2(Y_m,\ZZ)$ is generated
by the class $m\mathfrak f$.
If $h\in H^2(\CP^{m+2},\ZZ)$ is the hyperplane class,
we have $h|_{Y_m} = m\mathfrak f$.
Let $\psi:\tilde{\CP}^{m+2}\to \CP^{m+2}$
be the blowing up of $\CP^{m+2}$ at the ridge $\bm l$,
and $\pi:\tilde{\CP}^{m+2}\to\CP^m$ the induced
$\CP^2$-bundle projection.
Write temporary $E$ for the exceptional divisor of $\psi$,
and $h'$ be the hyperplane class of $\CP^m$.
By taking a hyperplane $H\subset\CP^{m+2}$
which contains $\bm l$,
we obtain $\psi^*h = E + \pi^* h'$.
Restricting to $\tilde Y_m = \pi\inv(\Lmd_m)$,
we obtain
$(\psi^*h)\big|_{\tilde Y_m} = \Sigma + m\mathfrak f$.

By passing $H^2(Y_m,\ZZ)$, this implies that
the generator $m\mathfrak f$ of $H^2(Y_m,\ZZ)$ is 
mapped to the class $\Sigma + m\mathfrak f$
under the resolution $\tilde Y_m\to Y_m$.
On the other hand, from the previous proposition,
we have $\tilde B\in \big|
4(\Sigma + m\mathfrak f)\big|$ for the branch divisor
$\tilde B$ of $\tilde\Phi_{m+1}:\tiZ\to \tilde Y_m$.
Therefore from the diagram \eqref{diagram5} the branch divisor of $\tilde\Phi_m:\tiZ\to Y_m$
belongs to the class $4h|_{Y_m}$.
%
%
Hence the conclusion follows.
\proofend

\medskip
As in the commutative diagram \eqref{diag},
the map $\tilde\Phi_m:\tiZ\to Y_m$ factors through
the map $Z\to Y_m$ via the birational morphism $\eta:\tiZ\to Z$.
Further, as we obtained in \eqref{ineqc}, the inequality $m\ge n-2$ holds.
Hence combination of Theorems \ref{thm:dc1} and \ref{thm:dc2}
gives Theorem \ref{thm:main1} except on the assertion on the equation of the quartic.

As in the proof of Theorem \ref{thm:dc2},
let $\tilde B\subset\tilde Y_m$ denote the 
branch divisor of the map $\tilde\Phi_{m+1}:\tiZ\to\tilde Y_m$.
According to the type of the twistor space,
the divisor $\tilde B$ has $\Aone,\Atwo$ or $\Athree$-singularities
respectively
along the two fibers of the projection $\Sigma\to \bm l$
over the points $\bm q$ and $\ol{\bm q}$.
If the twistor space is of type $\Azero$, 
the divisor $\tilde B$ is smooth along the last two fibers,
but it is tangent to $\Sigma$ along the two fibers.
Let $Z'\to\tilde Y_m$ be the (genuine) double cover of $\tilde Y_m\simeq
\mathbb P(\ms O(m)^{\oplus 2}\oplus\ms O)$
with branch $\tilde B$.
This is possible because we have $\tilde B\simeq
4(\Sigma + m\mathfrak f)$ as above and this is divisible by two.
Then the space $Z'$ is projective and bimeromorphic to the twistor space $Z$. 
We have the commutative diagram of morphisms:
\begin{align}\label{prj}
\xymatrix{ 
Z' \ar[r] \ar[dr] & \tilde Y_m \ar[d]\\
& \CP^1
}
\end{align}
Generic fibers of the morphism
$Z'\to \CP^1$ are Del Pezzo surfaces
of degree two.
Thus the present twistor space is 
bimeromorphic to the total space of 
a Del Pezzo fibration over $\CP^1$.
The Del Pezzo surfaces have $\Aone,\Atwo$ or $\Athree$-singularities
respectively
at two points according to the type of the twistor space.

\subsection{Cyclic base curves}
\label{ss:cbc}
In this subsection, we take a closer look at
the linear system $\big|\ms D(m)\big|$ on $\tiZ$
and find that, besides the stable base curves,
the system has a particular type of base curves
on the reducible fibers of $\tif:\tiZ\to\CP^1$
in general.

As we have seen in Section \ref{ss:dc},
for the image of the restriction homomorphism
$$r_E:H^0\big(\ms D(m)\big)
\lras  H^0\big(\ms D(m)|_E\big),$$
we have $\dim {\rm Image}\,(r_E) = 2$.
Next we show the following important property
about generators of the image.
Let $\lmd$ be a non-homogeneous coordinate
on $\CP^1$ which is the parameter space of 
the pencil $|F|$.
We may suppose that the induced real structure on $\CP^1$
is given by the complex conjugation in this coordinate.
We recall that in the proof of Proposition \ref{prop:e2},
we have chosen a special basis 
$\omega_0$ and $\ol\omega_0$ of $H^0\big(\ms D(e)|_E\big)$.
These satisfy the properties
$$
\omega_0|_{E_2\cup\dots\cup E_k} = 0
\qandq
\ol\omega_0|_{\ol E_2\cup\dots\cup \ol E_k} = 0.
$$

\begin{proposition}\label{prop:rE}
There exists a polynomial $P=P(\lmd)$ of degree
$(m-e)$ with real coefficients, which satisfies
\begin{align}\label{imrE}
{\rm Image}\,(r_E) = \big\langle
\omega_0,\ol\omega_0\big\rangle
\otimes P(\lmd).
\end{align}
Moreover the zeroes of $P(\lmd)$ are contained 
in the set $\{\lmd\upone,\dots,\lmd\upk\}$.
\end{proposition}

\proof
From Proposition \ref{prop:pim},
we have $\tilde\Phi_m (E_2\cup\dots\cup E_k) = \bm q$
and $\tilde\Phi_m (\ol E_2\cup\dots\cup \ol E_k) = \ol{\bm q}$.
Take any non-real hyperplane $H_{\bm q}$ which passes 
the point $\bm q$, and let $s_{\bm q}
\in H^0\big(\ms D(m)\big)$ be  a non-real element 
whose divisor $(s_{\bm q})$ corresponds to the hyperplane $H_{\bm q}$.
By abuse of notation, we write 
$s_{\ol{\bm q}}$ for $\ol{\sigma^*s_{\bm q}}\in H^0\big(\ms D(m)\big)$.
Then the section $s_{\bm q}$ (resp.\,$s_{\ol{\bm q}}$)
 vanishes identically on 
$E_2\cup \dots\cup E_k$ (resp.\,$\ol E_2\cup \dots\cup \ol E_k$),
and the restrictions 
$s_{\bm q}|_E$ and $s_{\ol{\bm q}}|_E$
form a basis of the 2-dimensional space ${\rm Image}\,(r_E)$.
Note that the non-reality of the hyperplane $H_{\bm q}$ 
is equivalent to $H_{\bm q}\not\supset \bm l$.
By Proposition \ref{prop:pim},
this implies that the divisor $(s_{\bm q})$ does not
includes $E_1$ nor $\ol E_1$.

From the explanation right after \eqref{dif0}, the difference
$(s_{\bm q})|_{\ol E_2\cup\dots\cup \ol E_k} -
\CCC_{\ol E_2\cup\dots\cup\ol E_k}$
is a sum of $(m-e)$ fibers of  $\tif|_{\ol E_2\cup
\dots\cup \ol E_k}:\ol E_2\cup
\dots\cup \ol E_k\to\CP^1$ which are not necessarily distinct.
Suppose that some fiber over a point $\lmd\in\CP^1$,
$\lmd\not\in\{\lmd\upone,\dots,\lmd\upk\}$, is included in this sum.
Write $S=\tif\inv(\lmd)$, and $C=S_{\lmd}\cap E$.
Since the subspace $H^0\big(\tiZ,\tilde f^*\ms O(m)\big)\otimes s_{\sst{\bm D}}$
and two sections $s_{\bm q}$, $s_{\ol{\bm q}}$ generate
the space $H^0\big(\ms D(m)\big)$,
by using $D.\,C_1=D.\,\ol C_1 = 1$, it follows that 
 the cycle $C$ on $S$ is a base curve of 
the line bundle $\ms D(m)$.
On the other hand,  by Proposition \ref{prop:van8},
we have $H^0(S,D-C) \simeq \CC$.
Hence the image of the restriction homomorphism
$H^0\big(\ms D(m)\big)\to H^0\big(\ms D(m)|_S\big)$
is at most 1-dimensional.
This contradicts that 
$h^0\big(
\ms D(m)
\big) -
h^0\big(
\ms D (m-1)
\big) = (m+3) - m =3$ which follows from 
Proposition \ref{prop:df7}.
Therefore the difference 
$\big(s_{\bm q}|_{\ol E_2\cup\dots\cup \ol E_k}\big)
-
\CCC_{\ol E_2\cup\dots\cup \ol E_k} $ consists of
 fibers over the points
$\lmd\upone,\dots,\lmd\upk$.
Namely we can write 
\begin{align}\label{sqE}
\big(s_{\bm q}|_{\ol E_2\cup\dots\cup \ol E_k}\big)
=
\CCC_{\ol E_2\cup\dots\cup \ol E_k} + 
\sum_{i=1}^k \mu\upi C\upi\,|_{\ol E_2\cup\dots\cup \ol E_k} 
\end{align}
for some non-negative integers $\mu\upone,\dots,\mu\upk$
that satisfy
\begin{align}\label{musum0}
\sum_{i=1}^k\mu\upi = m-e.
\end{align}

We define a polynomial $P(\lmd)$ by
\begin{align}\label{Pl}
P(\lmd ) = \prod_{i=1}^k 
\big(
\lmd-\lmd\upi
\big)^{\mu\upi}.
\end{align}
Since any $\lmd\upi$ is a real number from the choice of the coordinate $\lmd$, 
the polynomial $P(\lmd)$ is of real coefficients,
and by \eqref{musum0} the degree is $(m-e)$.
We compare the two elements
$s_{\bm q}|_E$ and $\omega_0\otimes P(\lmd)$ of 
$H^0\big(\ms D(m)|_E\big)$.
Both vanish identically on the union $E_2\cup E_3\dots\cup  E_k$.
Further, on another union $\ol E_2\cup\ol E_3\cup\dots\cup \ol E_k$,
both of the zero divisors are the restriction of the divisor
$\CCC_E+\sum_{i=1}^k\mu\upi C\upi$ to the union.
Therefore by multiplying a constant,
we may suppose that the two elements
$s_{\bm q}|_E$ and $\omega_0\otimes P(z)$
 agree on $E-E_1-\ol E_1$.
By Proposition \ref{prop:rest3},
this implies coincidence on the whole of $E$.
Namely we may suppose
\begin{align}\label{szE}
s_{\bm q}|_E = \omega_0\otimes P(\lmd).
\end{align}
Hence by using reality of the polynomial $P(\lmd)$,
we obtain $s_{\ol{\bm q}}|_E = \ol\omega_0\otimes P(\lmd)$.
The two sections $\omega_0\otimes P(\lmd)$
and $\ol\omega_0\otimes P(\lmd)$
are clearly linearly independent.
Therefore we obtain \eqref{imrE}.
\proofend

\medskip
The section $s_{\bm q}$  
obtained in this proof plays some role later.
So here we summarize its basic properties.

\begin{proposition}\label{prop:szero}
The element $s_{\bm q}\in H^0\big(\ms D(m)\big)$ 
satisfies the following properties.
\begin{itemize}
\item The restrictions $s_{\bm q}|_E$ and $s_{\ol{\bm q}}|_E$ generate
the image of the  restriction homomorphism
$r_E:H^0\big(\ms D(m)\big)\to H^0\big(\ms D(m)|_E\big)$.
\item It satisfies the following.
\begin{align}\label{szr0}
s_{\bm q}|_{ E_2\cup \dots\cup E_k} =0,
\quad
\big(s_{\bm q}|_{\ol E_2\cup \dots\cup \ol E_k}\big) =
\Big(\CCC_E + \sum_{i=1}^k\mu\upi C\upi
\Big)\,\Big|_{\ol E_2\cup \dots\cup \ol E_k}.
\end{align}
\end{itemize}
\end{proposition}

As above, the divisor $(s_{\bm q})$
is a pull-back of
a non-real hyperplane in $\CP^{m+2}$ which passes the point $\bm q$
on the ridge $\bm l$.
(So there are many choices for $s_{\bm q}$.)

\medskip
Proposition \ref{prop:rE} immediately implies the following
property on the base locus of the system $\big|\ms D(m)\big|$.

\begin{proposition}\label{prop:cbc1}
As in \eqref{Pl}, let $\mu\upi$
be the multiplicity of the zero
at the point $\lmd\upi$ of the polynomial $P(\lmd)$ obtained in Proposition \ref{prop:rE}.
Then the divisor
\begin{align}\label{cbc1}
\sum_{i=1}^k \mu\upi C\upi
\end{align}
on $E$ is fixed component of the linear system (pencil)
associated to the image of the restriction homomorphism $r_E:
H^0\big(\ms D(m)\big)\to  H^0\big(\ms D(m)|_E\big)$.
\end{proposition}

\begin{definition}
\label{def:cbc}{\em
We call the divisor \eqref{cbc1} on $E$
as the {\em cyclic base curves}.
\proofend
}
\end{definition}

Obviously this is base curve of the 
line bundle $\ms D(m)$.
Thus the linear system $\big|\ms D(m)\big|$ 
has two kinds of base curves,
the stable base curves and the cyclic base curves.
The former base curves consist of chains in 
the cycle $C\upi$, and are determined from the line bundle $\ms D|_E$.
On the other hand, the latter base curves
consist of whole of the cycle $C\upi$,
and are not determined from the line bundle $\ms D|_E$
itself, in the sense that they are relevant to 
the space of sections over the 3-dimensional space
$\tiZ$.
Later we show that $2\mu\upi=\dim\Coker\,\rho\upi$ holds for any index $i=1,\dots,k$.
This implies that the number $\mu\upi$ is a kind of 
local invariant determined from an infinitesimal neighborhood of
the reducible fiber $S\upi$.

From the proof of Proposition \ref{prop:rE},
these two kind of base curves are mutually independent;
namely, while these have common components, the sum of
the stable base curves and 
the cyclic base curves is really fixed component
of the pencil on $E$ in Proposition \ref{prop:cbc1}.

Although the cyclic base curves are never determined from
the line bundle $\ms D|_E$ as above,
this line bundle itself has the following interesting relevant property.

\begin{proposition}\label{prop:vsn1}
Suppose $i\neq 1,k$, and let $s$ be any section of $\ms D|_E$,
which is defined on a neighborhood in $E$ of the cycle $C\upi$.
If $s$ vanishes identically on the line components $C\upi_1$ and $\ol C\upi_1$,
then $s$ vanishes identically on the whole of the cycle $C\upi$.
Moreover, if $s$ vanishes on the two line components 
by multiplicity $\mu$, then $s$ vanishes on all
components of $C\upi$ by multiplicity
at least $\mu$.
\end{proposition}

\begin{figure}
\includegraphics{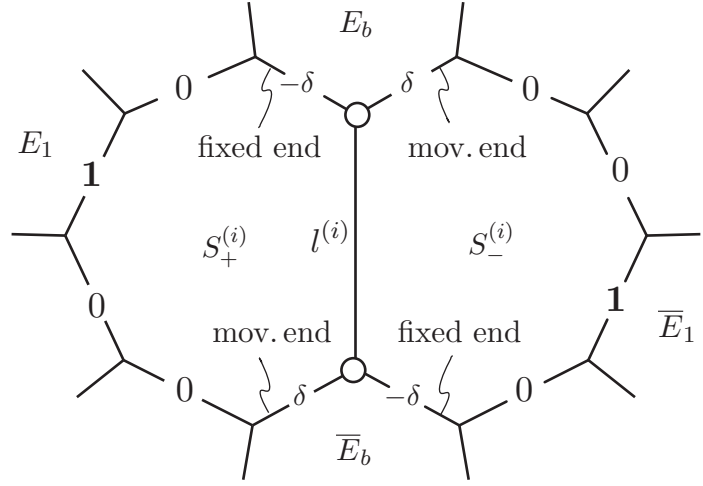}
\caption{Proof of Proposition \ref{prop:vsn1}}
\label{fig:rf4}
\end{figure}

\proof
In this proof, for simplicity of presentation we reset the notation
for the components of $C\upi_+=S\upi_+\cap E$ and we write
$$
C\upi_+ = C_1 + \dots + C_l + \dots + C_{k+1}
$$
in a way that $C_l$ is the line component $C\upi_1$ or $\ol C\upi_1$.
Then since $i\neq 1,k$, we have $1<l<k+1$.
Without loss of generality we may suppose that 
$\ms D.\,C_1 = -\ddd$ for some $\ddd\ge 0$.
Then $\ms D.\, C_{k+1} = \ddd\ge 0$
from the basic intersection numbers \eqref{int0011}.
If $\ddd>0$, the first component $C_1$ is the fixed end
of the stable base curves.
Let $s$ be as in the proposition and suppose $s|_{C_1}=
s|_{\ol C_l} =0$.
Then since $\ms D.\,C_j=0$ for any $j\not\in\{1,l,k+1\}$ from the basic intersection numbers,
we immediately obtain $s|_{C_j}=s|_{\ol C_j}=0$ for $j\neq 1,k+1$.
In particular, $s$ vanishes at the point $C_k\cap C_{k+1}$.
Since $\ms D.\,\ol C_1=-\ddd$ and 
the section $s$ is defined in a neighborhood of $C\upi$ in $E$,
the restriction $s|_{C_{k+1}}$ vanishes at the point
$C_{k+1}\cap \ol C_1$
(indicated by one of the circled points
in Figure \ref{fig:rf4}) by the order 
at least $\ddd$.
Therefore, because $\ms D.\,C_{k+1} = \ddd$,
we have $s|_{C_{k+1}}=0$.
Similarly we have $s|_{\ol C_{k+1}}=0$.
Thus $s$ vanishes on all components of the cycle $C\upi$.

For the latter assertion,
the cycle $C\upi = S\upi\cap E$ is a Cartier divisor on $E$,
and $[C\upi]$ makes sense as the line bundle
over the reducible surface $E$.
Let $c\upi\in H^0\big(E,[C\upi]\big)$
be an element whose zero is $C\upi$.
Then as $s|_{C\upi} = 0$,
the quotient $s/c\upi$ is a holomorphic
section of the line bundle $\ms D|_E\otimes
[-C\upi]$.
The degree of this bundle over components of $C\upi$
is the same as for the original $\ms D$ except possibly for the end components $C_1, \ol C_1, C_{k+1}$ and $\ol C_{k+1}$,
 since the self-intersection numbers
of non-end components, counted in the components of the divisor $E$,
are all zero.
But for the end components, the degree does not change as well.
In fact, let $E_b$ be the component on which 
the two curves $C_1$ and $\ol C_{k+1}$ lie.
(The subscript comes from `bridge'.
See Figure \ref{fig:rf4}.)
$b$ is either $i$ or $i+1$ depending on the choice of 
the small resolution $\zeta\upi$. Then
 we have
\begin{align*}
\big(\ms D|_E\otimes[- C\upi],C_1\big)_{E_b} &= \ms D.\,C_1 - \big(C_1 + \ol C_{k+1},C_1\big)_{E_b}\\
&=\ms D.\,C_1 -\{ (-1) + 1\}\\
&=\ms D.\,C_1,
\end{align*}
and the same for the other end components $\ol C_1, C_{k+1}$ and $\ol C_{k+1}$.
Hence we can apply the first assertion
in the present proposition 
to the holomorphic section $s/c\upi$ of the line bundle $\ms D|_E\otimes[-C\upi]$,
as long as the quotient $s/c\upi$ is zero on the line components
$C_l$ and $\ol C_l$.
This is the case exactly when $\mu>1$, where $\mu$ is
the integer as in 
the proposition.
Hence we can repeat the subtraction process  $\mu$ times.
This means that the original section $s$ vanishes 
along all components of the cycle $C\upi$ by multiplicity at least $\mu$.
\proofend

\subsection{Infinitesimal study on the direct image sheaves}
In this subsection, for the purpose of identifying the image
of the homomorphism $\rho:\tif_*\ms D
\to \tif_*\ms D|_E$, 
we investigate the stalks of these sheaves
at the points $\lmd\upi$, $1\le i\le k$.
In particular, we give a precise relationship between
the two numbers $e$ and $m$ in terms of the numbers $\mu\upi$
for the cyclic base curves,
and also determine the two numbers $\mu\upone$ and $\mu\upk$.

Elements of the domain
$(\tif_*\ms D\big)_{\lmd\upi}$ of $\rho\upi:
\big(\tif_*\ms D\big)_{\lmd\upi}\to
\big(\tif_*\ms D|_E\big)_{\lmd\upi}
$ are
represented by sections of the line bundle
$\ms D$, which are defined on neighborhoods of 
the reducible fiber $S\upi$.
Since the direct image $\tif_*\ms D$ is locally free of rank three,
we have 
\begin{align}\label{lf3}
\big(\tif_*\ms D\big)_{\lmd\upi}\simeq\ms O_{\lmd\upi}^{\oplus 3}.
\end{align}
We also have a natural homomorphism
$$
\big(\tif_*\ms D\big)_{\lmd\upi}
\lras H^0\big(\ms D|_{S\upi}\big),
$$
which assigns the restriction onto the reducible fiber $S\upi$ to each representative section of $\ms D$.
If $V\upi$ means 
the image of this homomorphism, from \eqref{lf3},
$V\upi$ is a 3-dimensional subspace of $H^0\big(\ms D|_{S\upi}\big)$.
This is a subspace consisting of sections of $\ms D|_{S\upi}$
which can be extended to a neighborhood  of $S\upi$.
The fact $\dim V\upi=3$ will be significant
in computing the number $\mu\upi$.

In a similar way to the stalk $\big(\tif_*\ms D\big)_{\lmd\upi}$, since the direct image sheaf
$\tif_*\ms D|_E$ is a locally free of rank two,
the stalk $\big(\tif_*\ms D|_E\big)_{\lmd\upi}$ is a locally free 
$\ms O_{\lmd\upi}$-module of rank two,
and we have a natural homomorphism
$$
\big(\tif_*\ms D|_E\big)_{\lmd\upi}
\lras H^0\big(\ms D|_{C\upi}\big).
$$
Let $W\upi$ be the image of this homomorphism.
This is a 2-dimensional subspace of 
$H^0\big(\ms D|_{C\upi}\big)$,
and consists of sections of $\ms D|_{C\upi}$
which extend to neighborhoods of $C\upi$ in the exceptional divisor $E$.

From the definition of the subspaces $V\upi$ and $W\upi$,
we have a commutative diagram
\begin{equation}\label{CDVW}
\begin{CD}
\big(\tif_*\ms D\big)_{\lmd\upi} @>{\rho\upi}>> 
\big(\tif_*\ms D|_E\big)_{\lmd\upi}\\
@VVV @VVV\\
V\upi @>{r\upi}>> W\upi,
\end{CD}
\end{equation}
where $r\upi$ is the restriction 
of the restriction homomorphism
$H^0\big(\ms D|_{S\upi}\big)\to H^0\big(\ms D|_{C\upi}\big)$
to the subspace $V\upi$.
Obviously if $\rho\upi$ is surjective, so is $r\upi$.
By Nakayama's lemma, the converse holds.
On the other hand, 
even if $r\upi$ is a zero map,
$\rho\upi$ can be a non-zero map,
and in such a case, the image of $\rho\upi$ is
more difficult to determine.

While the former subspace $V\upi$ seems difficult to identify
as a subspace of $H^0\big(\ms D|_{S\upi}\big)$ in a concrete
form in full generality,
it is not difficult to identify the latter subspace $W\upi$.
We first discuss situations where the stable base curve on $S\upi$
is absent. 
Suppose first $i\neq 1,k$.
This means that the line components $C\upi_1$ and $\ol C\upi_1$
are not the end components of the chains $C\upi_+=S\upi_+\cap E$ and $C\upi_-=S\upi_-\cap E$.
Absence of stable base curves on $S\upi$ means an equality $d_i=d_{i+1}$.
From \eqref{int04} and \eqref{int05}, this means
$\ms D.\,\DDD\upi = \ms D.\,C\upi_i=0$,
and the line bundle $\ms D$ is trivial over all components
of the cycle $C\upi$ except over the line components $C\upi_1$ 
and $\ol C\upi_1$.
Since $\ms D.\,C\upi_1 = \ms D.\,\ol C\upi_1=1$ from the basic
intersection number \eqref{int0011}, 
this implies that a section of the bundle
$\ms D|_{C\upi}$ is uniquely determined by the 
values at the two connected components of $C\upi- C\upi_1 - \ol C\upi_1$.
This gives a concrete isomorphism $H^0\big(\ms D|_{C\upi}\big)\simeq\CC^2$.
In particular,  since $\dim W\upi = 2$, this means  $W\upi = H^0\big(\ms D|_{C\upi}\big)$ under this situation.
If $i=1,k$, whether there exist stable base curves on $S\upi$ 
depends on the choice of small resolution in general,
but when they are absent,
exactly one conjugate pair of components of the cycle $C\upi$ have non-zero intersection number with $\ms D$,
and the number is one.
(See Figures \ref{fig:rf2} and \ref{fig:rf3}.)
This again gives a concrete identification
$H^0\big(\ms D|_{C\upi}\big)\simeq\CC^2$ as well as
the coincidence $H^0\big(\ms D|_{C\upi}\big) = W\upi$.

\begin{figure}
\includegraphics{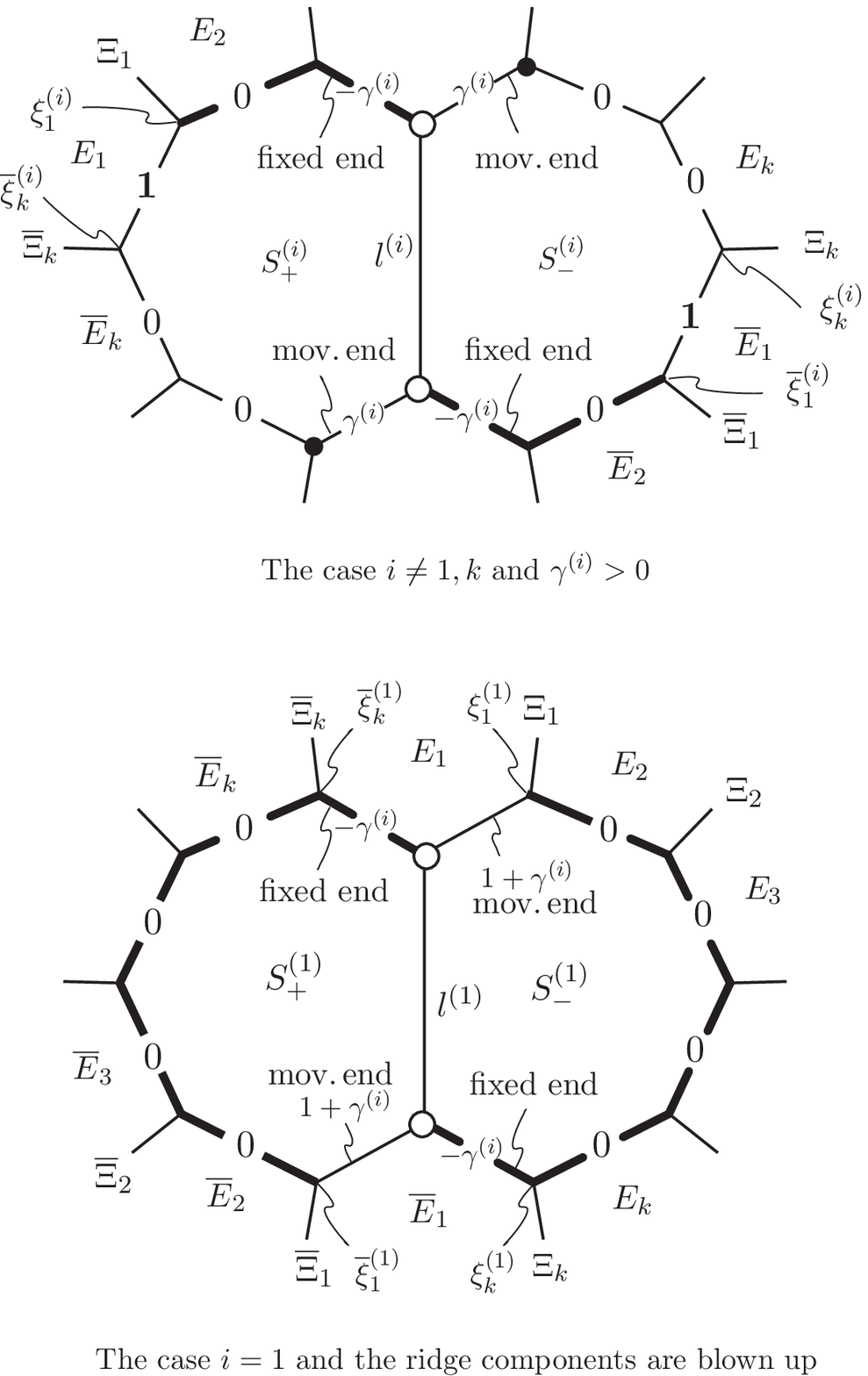}
\caption{Proof of Proposition \ref{prop:Wi}}
\label{fig:rf5}
\end{figure}

Next we discuss the general situation where the stable base curves
of $\ms D$ are present on $C\upi$.
If this is the case, the subspace $W\upi$ is strictly smaller
than the space $H^0\big(\ms D|_{C\upi}\big)$ in general,
and it can be concretely identified.
To state the result, by the {\em movable end} of the chains $C\upi_+$
and $C\upi_-$, we mean an end component of these chains whose
intersection number with $\ms D$ is positive.
Such a component always exists if the stable curves are present on $S\upi$,
and it is the end component which is different from the fixed end
of the stable base curves on $S\upi$ (see Figures \ref{fig:rf5}).
For each index $i\in\{1,\dots,k\}$,
we denote $\ccc\upi$ for the multiplicity of the stable base curve $\CCC\upi_E$
in $C\upi$.
If $i\neq 1,k$, this is equal to $\big|\ddd\upi\big|$
in the notation of the previous section.

\begin{proposition}
\label{prop:Wi}
The above subspace $W\upi$ in 
$H^0\big(\ms D|_{C\upi}\big)$ is identified with 
the subspace consisting of sections of
$\ms D|_{C\upi}$ whose vanishing order along the 
movable end at the intersection point
with the fixed end is at least $\ccc\upi$.
\end{proposition}

See Figures \ref{fig:rf5} where
the last intersection point
is indicated by one of the two small circles.

\medskip
\noindent
{\em Proof of Proposition \ref{prop:Wi}.}
Recall that the subspace $W\upi$ consists of
sections of $\ms D|_{C\upi}$ which can be
extended to a neighborhood of $C\upi$ in $E$.
%
Again we need to distinguish the case $i\neq 1,k$
and the case $i=1,k$.

When $i\neq 1,k$, 
the degree of the line bundle $\ms D$
on each component of the cycle $C\upi$ is
as in Figure \ref{fig:rf5} (or its mirror image).
From this, 
it follows that 
along the movable end, any section of $\ms D|_{C\upi}$
which extends to a neighborhood in $E$
vanishes at the fixed
end by multiplicity $\ccc\upi$ at least.
Hence any element of $W\upi$ satisfies the vanishing
property at the intersection point as in the 
proposition.
Moreover
since $\ms D$ has degree $\ccc\upi$ over the movable end as $i\neq 1,k$, sections of $\ms D|_{C\upi}$ which satisfy
the vanishing property are readily seen to be
determined from its values at the two dotted points
in Figure \ref{fig:rf5}
on the movable ends,
and therefore constitute a 2-dimensional subspace.
Since we know $\dim W\upi=2$,
these mean that $W\upi$ has to coincide with the space of
these sections of $\ms D|_{C\upi}$, as desired.

A proof for the cases $i=1,k$ is completely the same
if the small resolution $\zeta\upi$
does not blow up the 
the ridge components $E_1$ and $\ol E_1$.
Alternatively if the small resolution $\zeta\upi$
blows up the components $E_1$ and $\ol E_1$,
we need to note that the degree of $\ms D$
over the movable end is $\ccc\upi+1$ (instead of $\ccc\upi$), and that $\ms D$ is trivial over
all components of the chains $C\upi_+$ and $C\upi_-$
except the end components.
See lower picture in Figure \ref{fig:rf5}.
The requirement for elements of $W\upi$ at the circled points
is completely the same as the case $i\neq 1,k$,
but this time the stable base curves intersect
the movable ends.
From these, elements of $W\upi$ are uniquely determined by their values along the movable ends,
and over each movable end, all freedom is readily seen to be 1-dimensional
if we take the vanishing order at the circled points  in Figure \ref{fig:rf5} into account.
Thus we again obtain a 2-dimensional subspace of 
$H^0\big(\ms D|_{C\upi}\big)$ as a necessary condition,
and it has to coincide with $W\upi$ since 
$\dim W\upi=2$.
\proofend

\medskip
In order to investigate the image of the homomorphism
$\rho\upi:\big(\tif_*\ms D\big)_{\lmd\upi}\lra
\big(\tif_*\ms D|_E\big)_{\lmd\upi}$, we next identify the stalk $\big(\tif_*\ms D|_E\big)_{\lmd\upi}$.
Recall that for each $j\in\{1,\dots,k\}$, we have written $\Xi_j = E_j\cap E_{j+1}$, with the convention $E_{k+1}=\ol E_1$.
These are sections of the morphism $\tif:\tiZ\to\CP^1$.
For $i,j\in\{1,\dots,k\}$,
we put 
\begin{align}\label{xij}
\xi\upi_j:= S\upi\cap \Xi_j.
\end{align}
These are double points on the intersection cycle $C\upi$.
(See Figures \ref{fig:rf5} for some of these points.)
Then we have a natural homomorphism
\begin{align}\label{homom1}
\chi\upi_j:\big(\tif_*\ms D|_E\big)_{\lmd\upi}\lras
\big(\ms D|_{\Xi_j}\big)_{\xi\upi_j}\oplus \big(\ms D|_{\ol\Xi_j}\big)_{\ol\xi\upi_j}
\end{align}
which may be defined as follows.
Any element of the stalk $\big(\tif_*\ms D|_E\big)_{\lmd\upi}$
is represented by a section of the line bundle
$\ms D|_E$, which is defined on a neighborhood 
in $E$ of 
the cycle $C\upi=S\upi\cap E$.
By first restricting it to the curves $\Xi_j$ and $\ol\Xi_j$
and next taking  germs
at the points $\xi\upi_j$ and $\ol\xi\upi_j$,
we obtain the homomorphism \eqref{homom1}.

\begin{proposition}
\label{prop:homom1}
The homomorphism $\chi\upi_j$
is injective for any indices $i$ and $j$.
\end{proposition}

\proof
Let $s$ be a section of $\ms D|_E$
which is defined on a neighborhood of $C\upi$ in $E$,
and suppose that the restrictions $s|_{\Xi_j}$ 
and $s|_{\ol\Xi_j}$ are identically zero.
Then using the basic intersection numbers \eqref{int0011}
as far as possible,
we  obtain that $s$ vanishes identically on 
any component of $E$ except possibly the 
ridge components $E_1$ and $\ol E_1$
(see Figure \ref{fig:rf5}).
This implies that $s$ vanishes identically on any adjacent
components of $E_1$ and $\ol E_1$.
Therefore, again from the basic intersection numbers,
we further obtain that
the section $s$ vanishes identically also on $E_1$
and $\ol E_1$.
Hence  $s$ is a zero section.
This implies the injectivity of the homomorphism $\chi\upi_j$.
\proofend

\medskip
As before, for a point $\lmd\in\CP^1$, we denote $\ms O_{\lmd}$
for the local ring at the point $\lmd$.
%
Since $\Xi_j$ and $\ol\Xi_j$ are sections of $\tif$,
for any $j\in\{1,\dots,k\}$, we have isomorphisms
\begin{align}\label{psl}
\big(\ms D|_{\Xi_j}\big)_{\xi\upi_j}
\simeq
\big(\ms D|_{\ol\Xi_j}\big)_{\ol\xi\upi_j}
\simeq
\ms O_{\lmd\upi}.
\end{align}
By Proposition \ref{prop:homom1},
the homomorphism $\chi\upi_j$ can be regarded
as an injection to the direct sum 
$\ms O_{\lmd}\oplus\ms O_{\lmd}$ under the isomorphisms \eqref{psl}.

Recall that the symbol 
$\ccc\upi$ means the multiplicity of the stable base curves $\CCC\upi_E$
of $\ms D$ included
in $C\upi$.
In the sequel, for simplicity of presentation,
 for a pair of indices $i,j\in\{1,\dots,k\}$,
we denote 
\begin{align}\label{cij}
\ccc\upi_j:=
\begin{cases}
\ccc\upi, & {\text{if $\xi\upi_j\not\in \Supp\CCC_E$,}}\\
0, & {\text{otherwise.}}
\end{cases}
\end{align}
Let $\ms M_{\lmd\upi}\subset\ms O_{\lmd\upi}$ 
denote the maximal ideal of the local ring 
at the point $\lmd\upi$,
and we promise $\ms M^0_{\lmd\upi}=\ms O_{\lmd\upi}$.
Since any representative of the stalk $\big(\tif_*\ms D|_E\big)_{\lmd\upi}$
vanishes along the stable base curves
by  multiplicity $\ccc\upi$,
under the above identification \eqref{psl},
\begin{align}\label{psl2}
{\rm Image}\,\chi\upi_j\subset
\ms M^{\ccc\upi_j}_{\lmd\upi}
\oplus
\ol{\ms M}^{\ccc\upi_j}_{\lmd\upi}.
\end{align}
(The upper script means the power as an ideal.) 
Our next purpose is to show that this inclusion is always full:

\begin{proposition}
\label{prop:homom2}
Under the isomorphism \eqref{psl},
for any indices $i$ and $j$, 
 we have
\begin{align}\label{psl3}
{\rm Image}\,\chi\upi_j=
\ms M^{\ccc\upi_j}_{\lmd\upi}
\oplus
\ol{\ms M}^{\ccc\upi_j}_{\lmd\upi}.
\end{align}
\end{proposition}

We need some elementary lemmas for the proof.
To state it,
let $U$ be an open disk in $\CC$ which contains 
the origin, and
consider the trivial $\CP^1$-bundle  $U\times\CP^1\to U$
over $U$.
Define $\Xi:=U\times\{0\},\Xi':=U\times\{\infty\}$
and $A:=\{0\}\times\CP^1$. 
Let $(U\times\CP^1)^{\sim}$ be the blowing up
of $U\times\CP^1$ 
at the point $(0,\infty)\in U\times\CP^1$,
and $A'$ the exceptional curve.
We use the same letters $\Xi,\Xi'$ and $A$
to mean the strict transforms of the above curves
on $U\times\CP^1$ into the blow up.
We write the composition $(U\times\CP^1)^{\sim}\to
U\times\CP^1\to U$ by $\pi$.
Assume that  $L$ is a line bundle over $(U\times\CP^1)^{\sim}$
which satisfies
\begin{align}\label{intddd}
L.\,A = \ccc \qandq L.\,A' = -\ccc
\end{align}
for some $\ccc\ge 0$.
Then $L$ is trivial over any fiber of $\pi$ except 
the reducible fiber $A+A'$.
Under this situation we have the following extendability
for sections of $L$.

\begin{lemma}\label{lemma:ext}
Any section of $L$ defined over the curve $\Xi$
can be extended uniquely to the whole of $(U\times\CP^1)^{\sim}$.
Moreover, if the given section has zero at the point $\Xi\cap A$
with multiplicity $m$, then
over the curve $\Xi'$,
the extended section has zero with multiplicity $m+\ccc$ at the point 
$\Xi'\cap A'$.
Conversely, any such a section of $L$ defined over $\Xi'$ extends
to  the whole of $(U\times\CP^1)^{\sim}$,
and over the curve $\Xi$, it vanishes at the point $\Xi\cap A$ with multiplicity $m$.
\end{lemma}

%
\proof
Since the line bundle $L$ is trivial over any irreducible fiber 
of $\pi$ as above,
extendability and uniqueness are obvious over the complement of 
the reducible fiber $A + A'$. 
We show that extension can be taken all over $(U\times\CP^1)^{\sim}$
by using coordinates.

Let $x$ be the standard complex coordinate on $U\subset \CC$,
and $y$ an affine coordinate on $\CP^1$ such that 
$(x,y)=(0,0)$ is the point $\Xi\cap A$
and $(x,y)=(0,\infty)$ is the blown up point $\Xi'\cap A\in U\times\CP^1$.
Putting $u=y\inv$, the pair $(x,u)$ can be used as coordinates around the 
point $(0,\infty)$.
Hence if we put $v =xu\inv$, the pair
 $(u,v)$ can be used as coordinates on the blowing up, around the intersection point $A\cap A'$.
Putting $w = v\inv$, the pair $(x,w)$ can be used as coordinates
around the point $A'\cap \Xi'$.
Thus $(U\times\CP^1)^{\sim}$ can be covered by three
coordinates charts $(x,y)$, $(u,v)$ and $(x,w)$.
The transition laws are given by 
$$
uy=1, \quad xyw=1\qandq v=xy.
$$
Now by fixing a trivialization of the line bundle $L$ around the point
$(x,y)=(0,0)$, 
 any section of $L$ defined over $\Xi$ 
can be represented as a holomorphic function $f=f(x)$,
$x\in U$.
Naturally this section can be regarded as defined over
locus where the coordinates $(x,y)$ are valid.
Then noting that $L$
is isomorphic to the line bundle $[\ccc A']$ from \eqref{intddd},
over the locus where the coordinates $(u,v)$ are valid,
the section can be represented by the function 
$$
\frac{u^{\ccc}}1 f(x) = u^{\ccc}f(uv),
$$
which is holomorphic.
If we rewrite this section over the locus where 
the coordinates $(x,w)$ are valid, we get
$$
\frac{x^{\ccc}}{u^{\ccc}} u^{\ccc} f(uv) = x^{\ccc} f(uv) = x^{\ccc} f(x).
$$
This is  holomorphic in $x$ and $w$, and hence we get the 
desired extendability to the entire space.
The assertion about the multiplicity of zero is obvious from the final expression. The converse is also immediate if we reverse the above argument.
\proofend

\medskip
For the proof of Proposition \ref{prop:homom2},
we will also need the following lemma.

\begin{lemma}\label{lemma:ext2}
Let $(U\times\CP^1)^{\sim}$, 
$A$, $A'$, $\Xi$ and $\Xi'$ be as in Lemma \ref{lemma:ext},
and $L$ be a line bundle over the surface $(U\times\CP^1)^{\sim}$ 
which satisfies 
\begin{align}\label{intddd2}
L.\,A=\ccc+1 \qandq L.\,A'=-\ccc
\end{align}
for some $\ccc\ge 0$
(instead of \eqref{intddd}).
Then any holomorphic section of $L$ which is defined on 
the curves $\Xi\cup\Xi'$ 
extends to $(U\times\CP^1)^{\sim}$ 
in a unique way, provided that the section vanishes at the point $A'\cap \Xi'$
by the order at least $\ccc$.
\end{lemma}

Note that the final condition on the section over $\Xi'$
is a {\em necessary} condition for it to admit an extension to a neighborhood of
$\Xi'\cup A'$ from the latter intersection number in \eqref{intddd2}.
Note also that this time $L$ is of degree-one over irreducible fibers of the projection $\pi$.

\medskip
\noindent
{\em Proof of Lemma \ref{lemma:ext2}.}
We first show that any section of the line bundle $L-\ccc A'$ 
which is defined over the curves $\Xi\cup\Xi'$ extends to the entire surface
in a unique way.
Since $\Xi$ and $\Xi'$ are sections of $\pi$,
the degree of 
the line bundle $(L-\ccc A') - \Xi - \Xi'$ 
over irreducible fibers of $\pi$ is $1-(1+1)=-1$.
Hence this line bundle has no non-zero section.
This means that for proving the unique extension property
for sections of $L-\ccc A'$, it suffices
to show that $H^1\big((L-\ccc A') - \Xi - \Xi'\big)=0$.
To show this we make use of Leray spectral sequence.
From the degree of the line bundle
over irreducible fibers of $\pi$ and torsion freeness
of the direct image under a proper flat morphism, we obtain 
$\pi_*\big((L-\ccc A') - \Xi - \Xi'\big)=0$.
Therefore from the spectral sequence we obtain an isomorphism
$$
E^1=
H^1\big((L-\ccc A') - \Xi - \Xi'\big)
\simeq
H^0\big(R^1\pi_*((L-\ccc A') - \Xi - \Xi')\big)
=
E_2^{0,1}.
$$
Again from the degree on irreducible fibers of $\pi$,
the higher direct image in the right-hand side
is a torsion sheaf whose support can be only at
the origin of $U$.
Now we have
\begin{gather*}
\big((L-\ccc A') - \Xi - \Xi', A\big) = (\ccc+1) - \ccc - 1 = 0,\\
\big((L-\ccc A') - \Xi - \Xi', A'\big) = -\ccc + \ccc - 1 = -1.
\end{gather*}
From these we readily obtain that the first cohomology group of the restriction of
$(L-\ccc A') - \Xi - \Xi'$ to the reducible fiber $A \cup A'$ vanishes.
By Grauert's theorem, this implies
$R^1\pi_*\big((L-\ccc A')-\Xi-\Xi'\big)=0$.
Therefore we have shown the extension property.

To prove the lemma, take any section of $L$ defined over $\Xi\cup\Xi'$,
and write it as the pair $(t,t')$,
where $t$ and $t'$ are defined over the curves $\Xi$ and $\Xi'$ respectively.
Let $\aaa'$ be a section of the line bundle $[A']$
whose zero is $A'$.
Then the pair of  quotients
$$\Big(
\frac{t}{(\aaa'|_{\Xi})^{\ccc}},
\frac{t'}{(\aaa'|_{\Xi'})^{\ccc}}\Big)$$ 
is a holomorphic section
of the line bundle $\big(L-\ccc A'\big)|_{\Xi\cup\Xi'}$ from the condition on the vanishing
order of $t'$ at the point $A'\cap\Xi'$,
and disjointness of $A'$ and $\Xi$.
Hence from the above extendability,
this admits an extension to the whole of the surface.
Then by multiplying $(\aaa')^{\ccc}$ to the extension,
we get the extension of $(t,t')$ to the whole of the surface.
\proofend

\medskip
In the following proof of Proposition  \ref{prop:homom2}
we use these two lemmas by taking as $(U\times\CP^1)^{\sim}$ the components of the divisor $E$
which intersect the (strict transform of) twistor line $l\upi=S\upi_+\cap S\upi_-$
and as $L$ the restriction of the line bundle $\ms D$ to the components.
We recall that for any index $i$ and $j$, the notation $C\upi_j$
means the component of the cycle $C\upi=S\upi\cap E$
which corresponds to the component $C_j$ of the real anti-canonical cycle $C$,
under the birational morphism $\eta:\tiZ\to Z$.
We have $C\upi_j = S\upi\cap E_j$ if $l\upi$ does not intersect $E_j$.
If $E_b$ is a component of $E$ which intersects $l\upi$,
the intersection $S\upi\cap E_b$ consists of 
two components, and 
$C\upi_b$ is the component which is different from the exceptional curve
of the small resolution.
We call the curves $C\upi_1$ and $\ol C\upi_1$  the line components of the cycle $C\upi$.
These are mapped isomorphically to the components $C_1$ and $\ol C_1$ by the birational morphism $\eta:\tiZ\to Z$.

\medskip
\noindent
{\em Proof of Proposition \ref{prop:homom2} when $i\neq 1,k$}.
We first show that $\chi\upi_j$ is surjective
if the index $j$ is chosen in a way that
the point $\xi\upi_j$
(and $\ol\xi\upi_j$) is not on the stable curves of $\ms D$.
Let $(t,t')$ be a representative
of any element of $\big(\ms D|_{\Xi_j}\big)_{\xi\upi_j}\oplus \big(\ms D|_{\ol\Xi_j}\big)_{\ol\xi\upi_j}
$,
So $t$ (resp.\,$t'$) is a section of $\ms D$,
which is defined on a neighborhood in the curve $\Xi_j$ of 
the point $\xi\upi_j$ 
(resp.\,a neighborhood in the curve
$\ol\Xi_j$ of the point $\ol\xi\upi_j$).
Again we write $E_{b}$ and $\ol E_{b}$ for the components of $E$
which intersect the twistor line $l\upi$.
The number $b$ is either $i$ or $i+1$ depending on 
the choice of small resolution at $\zeta\upi$.
The intersection $E_{b}\cap S\upi$ consists
of two $(-1)$-curves on $E_{b}$.

On a neighborhood of $C\upi$ in $E$,
the sections $t$ and $t'$ can be successively extended
to all components of $E$ which can be reached to 
the component $E_{b}$ or $\ol E_{b}$ without 
passing the ridge component $E_1$ or $\ol E_1$,
because over these components the line bundle $\ms D$ is 
trivial  on a neighborhood of $C\upi$ in $E$.
(See the upper picture in Figure \ref{fig:rf5}.)
By using Lemma \ref{lemma:ext}, we can further extend the sections
across the components $E_{b}$ and $\ol E_{b}$.
They can further be extended  to the remaining components of $E$
successively until hitting the component $E_1$ or $\ol E_1$,
in the same way to the first extension we took.
Thus we have taken an extension 
except over neighborhoods of the line components
$C\upi_1$ and $\ol C\upi_1$ in $E_1$ and $\ol E_1$ respectively.
But each of these are sandwiched by other components of $E$.
Since $\ms D.\,C\upi_1 = \ms D.\,\ol C\upi_1=1$
from the basic intersection numbers, this implies that 
the section can be further extended to a neighborhoods
of the curve $C\upi_1$ in $E_1$ and to
a neighborhoods
of the curve $\ol C\upi_1$ in $\ol E_1$ respectively.
Thus we have obtained the desired extension.
Hence the homomorphism \eqref{homom1} is surjective
if $\xi\upi_j\not\in \Supp\,\CCC_E$.

To complete the proof of the proposition,
choose an index $j$ for which the point $\xi\upi_j$ belongs to 
the stable base curves of $\ms D$.
Since we fix $i$ and $j$, we simply write $\ccc$ for 
the multiplicity $\ccc\upi=\ccc\upi_j$.
Let $(t,t')$ be any element of the sub-module
$\ms M^{\ccc}_{\lmd\upi}
\oplus
\ol{\ms M}^{\ccc}_{\lmd\upi}
$
%
Then again the sections $t$ and $t'$ can 
be successively extended
to all components of $E$ which can be reached to 
the component $E_{b}$ or $\ol E_{b}$ without 
passing the ridge component $E_1$ or $\ol E_1$.
(Again see the upper picture in Figure \ref{fig:rf5}.)
Further, this time, these extensions vanish along components of 
the cycle $C\upi$ by the order $\ccc$.
Therefore by the converse part of Lemma \ref{lemma:ext},
the extension can further be extended to the components
$E_{b}$ and $\ol E_{b}$.
Then in the same way to the above case  $\xi\upi_j\not\in
\Supp\,\CCC_E$, we get an
extension to whole of 
a neighborhood of $C\upi$ in $E$.
This means
${\rm Image}\,\chi\upi_j =
\ms M^{\ccc}_{\lmd\upi}
\oplus
\ol{\ms M}^{\ccc}_{\lmd\upi}$, as desired.
\proofend

\medskip
\noindent
{\em Proof of  Proposition \ref{prop:homom2} when $i=1,k$}.
We show  only when $i=1$,
 since the case $i=k$ can be shown in the same way.
First, if the small resolution  $\zeta\upone$ does not blow up
the ridge component $E_1$,
the components
of $E$ which intersect the twistor line $l\upone$
are $E_2$ and $\ol E_2$.
The restrictions of the line bundle 
$\ms D$ to the components $E_{2}$
and $\ol E_{2}$ satisfy the numerical 
property as in Lemma \ref{lemma:ext},
and the previous proof (in the case $i\neq 1,k$)
works without any change.

So in the following we prove the assertion when
the small resolution $\zeta\upone$ blows up
the ridge component $E_1$.
Fix any $j\in\{1,\dots,k\}$
and write $\ccc$ for the multiplicity $\ccc\upone_j$.
Let $(t,t')$ be any element of $
\ms M^{\ccc}_{\lmd\upi}
\oplus
\ol{\ms M}^{\ccc}_{\lmd\upi}$.
From the choice of the small resolution,
except possibly the  end components,
the line bundle $\ms D$ is trivial
over any component of the half chains $C\upone_+$ and $C\upone_-$.
(See lower picture in Figure \ref{fig:rf5}.)
This implies that the line bundle $\ms D$ is trivial over
neighborhoods in $E$ of these components.
Therefore the sections $t$ and $t'$ can be extended
to neighborhoods in $E$ of these components.
In particular, we get a section of $\ms D|_E$ which is defined 
in a neighborhood of the point $\xi\upone_1$ in $\Xi_1$,
and also a section of $\ms D|_E$ which is defined 
in a neighborhood of the point $\ol\xi\upone_k$ in $\ol\Xi_k$.
As we are choosing $(t,t')$ from the sub-module
$$
\ms M^{\ccc}_{\lmd\upi}
\oplus
\ol{\ms M}^{\ccc}_{\lmd\upi},
$$
the section $t'$ vanishes at the point $\ol\xi\upone_k$ by the order $\ccc$.
Therefore applying Lemma \ref{lemma:ext2}
to these sections, they can further be extended
to a neighborhood of the curve $S\upone\cap E_1$
in $E_1$.
This is of course the same over the component  $\ol E_1$.
Thus the pair $(t,t')$ has the desired extension.
\proofend

\medskip
So far in this subsection we have mainly worked on the divisor $E$.
Now we are able to express the image of 
the homomorphism 
$
\rho\upi:\big(\tif_*\ms D\big)_{\lmd\upi}\to
\big(\tif_*\ms D|_E\big)_{\lmd\upi}
$ 
by using the invariant $\mu\upi$ appeared in \eqref{Pl}.

\begin{proposition}\label{prop:mvo}
Under the isomorphism \eqref{psl},
for any indices $i,j\in\{1,\dots,k\}$, we have
\begin{align}\label{imr}
{\rm Image\,}\big(\chi\upi_j\circ\rho\upi\big) = 
\ms M^{\ccc\upi_j+\mu\upi}_{\lmd\upi}
\oplus \ol{\ms M}^{\ccc\upi_j+\mu\upi}_{\lmd\upi}
\end{align}
and
\begin{align}\label{mcok}
2\mu\upi = \dim\Coker\,\rho\upi.
\end{align}
Further, we have the relation
\begin{align}\label{mec}
m=e + \sum_{i=1}^k\mu\upi.
\end{align}
%
\end{proposition}

\proof
First we note that for any $i$ and $j$ there is a commutative diagram
of stalks
\begin{align}\label{diagram51}
 \xymatrix{ 
(\tif_*\ms D)_{\lmd\upi} 
\ar[rr] \ar_{\rho\upi}[dr] && 
\big(\ms D|_{\Xi_j}\big)_{\xi\upi_j}\oplus \big(\ms D|_{\ol\Xi_j}\big)_{\ol\xi\upi_j} \\
&\big(
 \tif_*\ms D|_E
 \big)_{\lmd\upi}
 \ar[ur]_{\chi\upi_j} 
}
\end{align}
where the horizontal arrow
is a map which sends
a germ in the domain to
a pair of germs at the two points
$\xi\upi_j$ and $\ol \xi\upi_j$
after taking restriction to
the curves $\Xi_j$ and $\ol\Xi_j$ respectively.
By Proposition \ref{prop:homom2},
Image\,$\chi\upi_j$ can be identified with the direct sum
$\ms M^{\ccc\upi_j}_{\lmd\upi}
\oplus \ol{\ms M}^{\ccc\upi_j}_{\lmd\upi}$.
Therefore, since $\chi\upi_j$ is always injective
by Proposition \ref{prop:homom1}, if we regard $\chi\upi_j$
as a map to $\ms M^{\ccc\upi_j}_{\lmd\upi}
\oplus \ol{\ms M}^{\ccc\upi_j}_{\lmd\upi}$,
 we have an isomorphism
\begin{align}\label{isosli}
\Coker\,\rho\upi\simeq
\Coker\,\big(\chi\upi_j\circ\rho\upi\big).
\end{align}

We consider the two elements $s_{\bm q}$ and $s_{\ol{\bm q}}$
of $H^0\big(\ms D(m)\big)$ in 
Proposition \ref{prop:szero}.
As in \eqref{szr0}, the former section $s_{\bm q}$ vanishes identically on 
$E_2\cup\dots\cup  E_k$,
and from the property $s_{\bm q}|_E = \omega_0\otimes P(z)$
as in the latter property in \eqref{szr0},
the restriction $s_{\bm q}|_{\ol\Xi_j}$ vanishes at the point
$\ol\xi\upi_j$  by the order exactly $\ccc\upi_j + \mu\upi$,
where $\ccc\upi_j$ is as defined in \eqref{cij}.
This implies that the germ 
$\big(s_{\bm q}|_{\ol\Xi_j}\big)_{\ol\xi\upi_j}$ is a generator of 
the sub-module 
$0\oplus\ms M_{\lmd\upi}^{\ccc\upi_j+\mu\upi}$
of $\ms M^{\ccc\upi_j+\mu\upi}_{\lmd\upi}
\oplus \ol{\ms M}^{\ccc\upi_j+\mu\upi}_{\lmd\upi}$.
By the real structure, the germ 
$\big(s_{\ol{\bm q}}|_{\ol\Xi_j}\big)_{\ol\xi\upi_j}$
is a generator of the sub-module $ \ol{\ms  M}_{\lmd\upi}^{\ccc\upi_j+\mu\upi}\oplus 0$ of $\ms M^{\ccc\upi_j+\mu\upi}_{\lmd\upi}
\oplus \ol{\ms M}^{\ccc\upi_j+\mu\upi}_{\lmd\upi}
$.
%
%
From these, we obtain 
\begin{align}\label{icim2}
\ms M^{\ccc\upi_j+\mu\upi}_{\lmd\upi}
\oplus \ol{\ms M}^{\ccc\upi_j+\mu\upi}_{\lmd\upi}
\subset
{\rm Image}\,\big(
\chi\upi_j\circ\rho\upi
\big).
\end{align} 
Hence there is an exact sequence
\begin{align}\label{sj1}
\Big(\ms M^{\ccc\upi_j}_{\lmd\upi}
\oplus \ol{\ms M}^{\ccc\upi_j}_{\lmd\upi}\Big)
\Big/
\Big(
\ms M^{\ccc\upi_j+\mu\upi}_{\lmd\upi}
\oplus \ol{\ms M}^{\ccc\upi_j+\mu\upi}_{\lmd\upi}
\Big)
\lras
\Coker\big(
\chi\upi_j\circ\rho\upi
\big)
\lras 0.
\end{align}
Obviously the domain of this surjection is isomorphic to 
$\CC^{\mu\upi}\oplus\CC^{\mu\upi} =\CC^{2\mu\upi}$.
Therefore, by \eqref{isosli},
 we obtain an inequality
\begin{align}\label{iem}
2\mu\upi\ge \dim\Coker\,\rho\upi,
\end{align}
and the equality holds exactly when 
the inclusion \eqref{icim2} is an equality.
Summing these up with respect to $i$
and recalling $\sum_{i=1}^k\mu\upi = m-e
$ from \eqref{musum0}, we obtain 
\begin{align}\label{ieem}
2 (m-e) \ge  \sum_{i=1}^k\dim\Coker\,\rho\upi,
\end{align}
and the equality holds exactly when 
the inclusion \eqref{icim2} is equality for any  $i$ and $j$.
On the other hand, as we obtained in \eqref{mec2} from the derived 
exact sequence,
both sides of \eqref{ieem} are equal.
Namely, the inequality \eqref{ieem} is an equality.
Hence the inequality \eqref{iem} is equality for any $i$.
So we obtain \eqref{mcok} and \eqref{mec}.
Therefore 
the inclusion \eqref{icim2} is an isomorphism.
This shows \eqref{imr}.
\proofend

\medskip
The equality $2\mu\upi = \dim\Coker\,\rho\upi$
in \eqref{mcok} implies that
the number $\mu\upi$ is independent of 
a choice of the small resolutions
$\zeta\upj$ as long as $j\neq i$.
Moreover, from the 
identification \eqref{imr}, we obtain 

\begin{proposition}\label{prop:in}
For any index $i\in\{1,\dots,k\}$,
we have an equality
\begin{align}\label{cbc0}
\mu\upi=\max\Big\{
\nu\,\Big|\,
{\text{\em$\big(s|_E\big)\ge \nu C\upi$  for all
$s\in
\big(\tif_*\ms D\big)_{\lmd\upi}$}}
\Big\}.
\end{align}
\end{proposition}

\noindent
In \eqref{cbc0}, for $s\in
\big(\tif_*\ms D\big)_{\lmd\upi}$,
the restriction $s|_E$ is
taken in the following (natural) sense:
any element of the stalk $\big(\tif_*\ms D\big)_{\lmd\upi}$ is represented by 
a section of the line bundle $\ms D$,
which is defined in a neighborhood of the 
fiber $S\upi$ in $\tiZ$.
So the restriction $s|_E$ makes sense
for any $s\in \big(\tif_*\ms D\big)_{\lmd\upi}$.

Thus the number $\mu\upi$ is a 
kind of local invariant, in the sense
that it is determined
from an infinitesimal neighborhood of the reducible fiber $S\upi$ in $\tiZ$.

%


At the end of the next subsection we shall give
an example of twistor spaces
on $n\CP^2$, $n>6$, for which some reducible fiber
$S\upi$ satisfies $\mu\upi>0$.
Since $m = e + \sum\mu\upi$ as in \eqref{mec},
we have the strict inequality $m>e$ for these twistor spaces.
Thus the inequality $m\ge e$ that we obtained in 
Proposition \ref{prop:d_im06} or \eqref{dem} can be really strict.

Next we investigate the image of the homomorphism
$\rho\upi$ when $i=1,k$.
Recall that the number $e$ depends
on the choice of the small resolutions $\zeta\upone$ and $\zeta\upk$
as in Proposition \ref{prop:e}.
On the other hand, the direct image $\tif_*\ms D$
(namely the number $m$ for which $\tif_* \ms D\simeq
\ms O\oplus \ms O(-m)^{\oplus 2}$ holds) is independent of a choice of 
the small resolutions
as in Proposition \ref{prop:indep0}.
This discrepancy fits 
with that of the cokernel of $\rho\upi$, with $i=1,k$,
as follows.

\begin{proposition}\label{prop:cokrho1}
If $i=1,k$, the following hold.
\begin{itemize}
\item
If the small resolution $\zeta\upi$
blows up the ridge components $E_1$ and $\ol E_1$, then 
$\mu\upi=0$ (i.e.\,the homomorphism $\rho\upi$ is surjective
by \eqref{mcok}).
\item
If not, we have $\mu\upi = (d_2-2)_+$ if $i=1$ 
and $\mu\upi = (d_k-2)_+$ if $i=k$.
\end{itemize}
\end{proposition}

\begin{figure}
\includegraphics{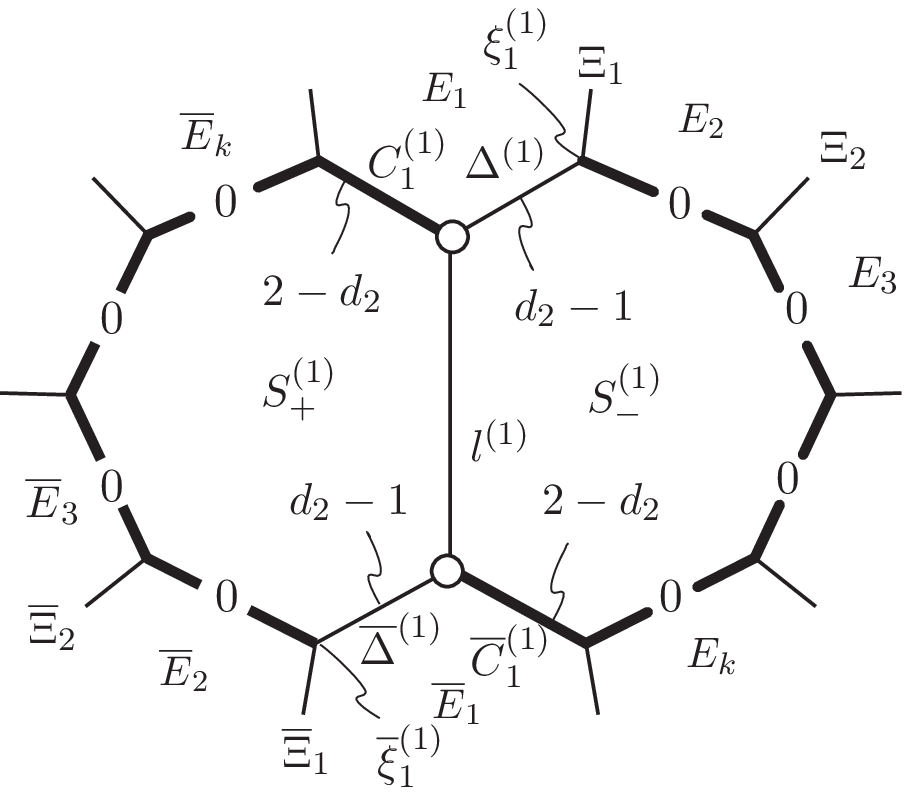}
\caption{
A proof of the first item in Proposition \ref{prop:cokrho1}
}
\label{fig:rf6}
\end{figure}

\noindent
{\em Proof of the first item.}
We prove only in the case $i=1$ since the case $i=k$ can be shown in the same way.
Recall that if the small resolution $\zeta\upone$
blows up the component $E_1$,
 the stable base curves on $S\upone$ are present
precisely when $d_2>2$, and 
if this inequality holds, the multiplicity of
the base curves (namely the number $\ccc\upone$ in 
\eqref{cij}) is $(d_2 - 2)$
(see Figure \ref{fig:rf6}).

First suppose $d_2 = 1$. 
Then $\ccc\upone=\ccc\upone_1 = 0$.
Hence by Propositions \ref{prop:homom1} and \ref{prop:homom2}
we have an isomorphism $\big(\tif_*\ms D|_E\big)_{\lmd\upone}
\simeq \ms O_{\lmd\upone}\oplus\ms O_{\lmd\upone}$ under the 
identification $\big(\ms D|_{\Xi_1}\big)_{\xi\upone_1}
\simeq
\big(\ms D|_{\ol\Xi_1}\big)_{\ol\xi\upone_1}
\simeq\ms O_{\lmd\upone}
$, and 
 by \eqref{imr} 
\begin{align*}
{\rm Image\,}\big(\chi\upone_1\circ\rho\upone\big) = 
\ms M^{\mu\upone}_{\lmd\upone}
\oplus \ol{\ms M}^{\mu\upone}_{\lmd\upone}.
\end{align*}
%
%
This implies that, if $\mu\upone>0$, 
any element of $\big(\tif_*\ms D\big)_{\lmd\upone}$
would vanish at the points $\xi\upone_1$ and $\ol\xi\upone_1$.
From the intersection numbers with $\ms D$ (see Figure \ref{fig:rf6}),
this readily means that any element of $\big(\tif_*\ms D\big)_{\lmd\upone}$ vanishes along the whole of the cycle $C\upone$.
Hence the same would be true for any element
of the 3-dimensional subspace $V\upone$ in 
$H^0\big(\ms D|_{S\upone}\big)$
appearing in the commutative diagram \eqref{CDVW}.
In particular any element of $V\upone$
would vanish identically
on the both ends of the chains $C\upone_+$ and 
$C\upone_-$ respectively.
Therefore 
 we can subtract both ends
from the linear systems $\big|\bm D|_{S\upone_+}\big|$
and $\big|\bm D|_{S\upone_-}\big|$ as base curves.
Then the intersection matrices of the resulting 
systems
(namely the systems
$\big|\bm D|_{S\upone_+} - C\upone_1 - \ol\DDD\upone\big|$
and 
$\big|\bm D|_{S\upone_-} -  \ol C\upone_1 - \DDD\upone\big|$)
are readily seen to be negative definite
 noting that the ingredient effective curves
 do not include the end components
because the line components are of multiplicity one in $\bm D$.
This implies that each of the subtractions consists
of a single member.
Hence the linear system $\big|V\upone\big|$ consists
of a single member.
This contradicts $\dim V\upone = 3$.
Hence $\mu\upi=0$.
The case $d_2 = 2$ can be shown in the same way
by just exchanging the roles of $C\upone_1$
and $\DDD\upone$ in this argument.

The case $d_2>2$ is similarly shown as follows.
This time we have $\ccc\upone_1 = d_2-2>0$.
Hence again by Propositions \ref{prop:homom1} and \ref{prop:homom2}
we have an isomorphism $\big(\tif_*\ms D|_E\big)_{\lmd\upone}
\simeq \ms M_{\lmd\upone}^{d_2-2}
\oplus\ms \ol{\ms M}_{\lmd\upone}^{d_2-2}
$ under the 
identification $\big(\ms D|_{\Xi_1}\big)_{\xi\upone_1}
\simeq
\big(\ms D|_{\ol\Xi_1}\big)_{\ol\xi\upone_1}
\simeq\ms O_{\lmd\upone}
$, and 
 by \eqref{imr} 
\begin{align*}
{\rm Image\,}\big(\chi\upone_1\circ\rho\upone\big) = 
\ms M^{d_2-2+\mu\upone}_{\lmd\upone}
\oplus \ol{\ms M}^{d_2-2+\mu\upone}_{\lmd\upone}.
\end{align*}
Hence if $\mu\upone>0$,
any element of the germ $\big(\tif_*\ms D\big)_{\lmd\upone}$
would vanish at the points $\xi\upone_1$ and $\ol\xi\upone_1$ by the order at least $(d_2-1)$
when restricted to $\Xi_1$ and $\ol\Xi_1$ respectively. 
Then from the degrees of $\ms D$
over components of $C\upone$ as in Figure
\ref{fig:rf6}, any element of $(\tif_*\ms D)_{\lmd\upone}$
would vanish along the stable base curves (namely
the chains indicated by the bold lines)
by the same multiplicity (which is at least
$(d_2-1)$) when restricted to the divisor $E$.
This  implies that any element of 
$\big(\tif_*\ms D\big)_{\lmd\upone}$ vanishes at the
intersection points $\DDD\upone\cap l\upone$
and $\ol \DDD\upone\cap l\upone$,
the circled points in Figure \ref{fig:rf6},
by that order when restricted to $\DDD\upone$ and $\ol\DDD\upone$
respectively.
Moreover, any element of the same stalk
vanishes at the points $\xi\upone_1$ and $\ol\xi\upone_1$.
Hence, since $\ms D.\,\DDD\upone = \ms D.\,\ol\DDD\upone= d_2-1$, any element of the stalk
$\big(\tif_*\ms D\big)_{\lmd\upone}$
would vanish identically on $\DDD\upone$ and $\ol\DDD\upone$
under the assumption $\mu\upone>0$.
Thus we have obtained the same conclusion 
with the case $d_2 = 1$ to the effect that 
any element of $V\upone$ would vanish on
both ends of the chains $C\upone_+$ and $C\upone_-$.
By the same reason to the above proof  in the case $d_2 = 1$,
this means $\dim V\upone=1$, and 
we again obtain a contradiction.
Therefore  $\mu\upone=0$ if $d_2>2$ as well.
\proofend

\medskip
The proof for the case where 
the small resolution $\zeta\upone$ or $\zeta\upk$ does not blow up 
the ridge components is more interesting.

\begin{figure}
\includegraphics{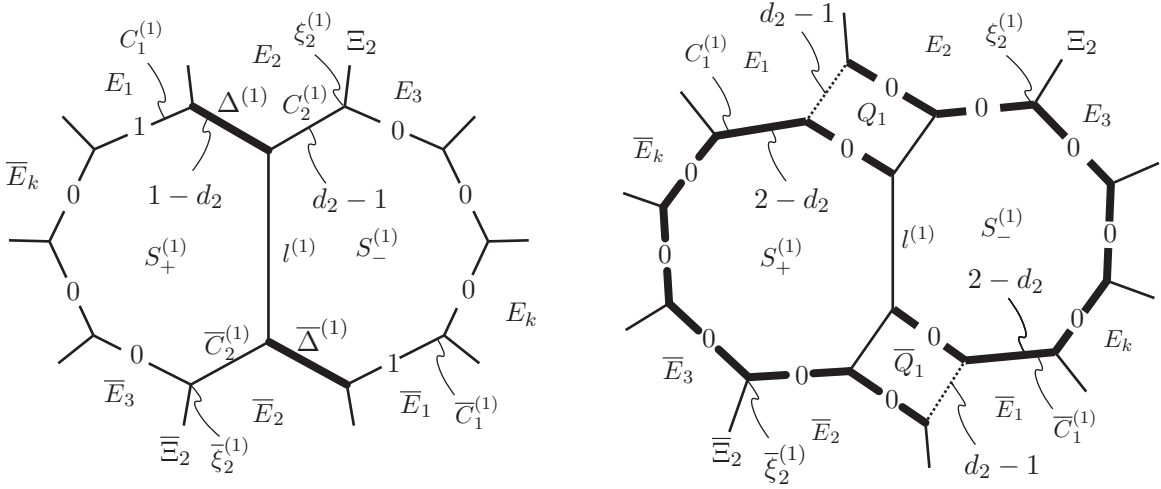}
\caption{A proof of the second item in Proposition \ref{prop:cokrho1}}
\label{fig:rf7}
\end{figure}

\medskip
\noindent
{\em Proof of the second item in Proposition \ref{prop:cokrho1}.}
Again we show the assertion only when $i=1$
since the case $i=k$ can be shown in the same way.
We recall that under the choice of the small resolution,
we have
$$
\ms D.\,\DDD\upone = 1- d_2\qandq
\ms D.\,C\upone_2 = d_2-1,
$$
and the stable base curves on $S\upone$
exist if and only if $d_2 >1$
(see the left picture in Figure \ref{fig:rf7}).
From this, regardless of the value of $d_2$,
the points $\xi\upone_2$ and $\ol\xi\upone_2$
do not belong to the stable base curves on $S\upone$.
Hence by Propositions \ref{prop:homom1}
and \ref{prop:homom2}, 
we obtain an isomorphism
\begin{align}\label{disum}
\big(\tif_*\ms D|_E\big)_{\lmd\upone}
\simeq
\big(\ms D|_{\Xi_2}\big)_{\xi\upone_2}\oplus \big(\ms D|_{\ol\Xi_2}\big)_{\ol\xi\upone_2}.
\end{align}
If $d_2=1$, the stable base curves on $S\upone$ do not appear,
and passing to the 
direct sum in \eqref{disum},
the same argument in the proof of 
the first item of the present proposition
shows $\mu\upone=0$.


Next we suppose $d_2>2$.
(The case $d_2=2$ will be discussed later.)
Again exactly the curves $\DDD\upone$ and $\ol\DDD\upone$
are the stable base curves  on $S\upone$.
This time it seems difficult to get the desired conclusion by just looking at the space $\tiZ$.
Let $\nu:W\to \tiZ$ be the blowing up at these base curves,
and $Q_1$ and $\ol Q_1$ the exceptional divisors.
(See the right picture in Figure \ref{fig:rf7}.)
These are isomorphic to $\qdr$ and are fixed components
of $\nu^*\ms D(l)$ for any $l$
with multiplicity  $(d_2 -1)$.
Namely we have an identification
\begin{align}\label{id0001}
\big|\nu^*\ms D(l)\big| =
(d_2-1)\big(Q_1 + \ol Q_1\big)+
\big|\nu^*\ms D(l) - (d_2 - 1) (Q_1 + \ol Q_1)\big|
\end{align}
for any $l$.
For simplicity in the following we write
$\ms D'(l):= \nu^*\ms D(l) - (d_2 - 1) (Q_1 + \ol Q_1)$.
By using  $N_{Q_1/W}\simeq N_{\ol Q_1/W}\simeq
\ms O(-1,-1)$,
we readily obtain that 
the intersection numbers of $\ms D'(l)$ with
curves in the fiber over $\lmd\upone$ are as illustrated 
in the right picture of Figure \ref{fig:rf7}
on each component.
From this, as $2-d_2<0$, 
the two chains indicated by 
bold segments in the  picture are base curves of 
$\ms D'(l)$ for any $l$.
Moreover, noting that the curves $C\upone_1$
and $\ol C\upone_1$ are $(-1)$-curves
in $E_1$ and $\ol E_1$ respectively,
we can successively subtract these chains
as fixed component of the restriction
$\ms D'(l)|_{E}$,
$(d_2-2)$ times in total.
Namely, these two chains are fixed components
of $\ms D'(l)|_{E}$
with multiplicity $(d_2-2)$ at least.

By transferring to 
the original space $\tiZ$,
this means that any element of 
$\big(\tif_*\ms D\big)_{\lmd\upone}$ vanishes at the points $\xi\upone_2$
and $\ol\xi\upone_2$ by the order $(d_2-2)$
when restricted to the curves $\Xi_2$ and $\ol\Xi_2$ respectively.
Via the isomorphism \eqref{disum},
this implies an equality
\begin{align}\label{ieq}
\mu\upone\ge d_2-2.
\end{align}

Next suppose that the inequality \eqref{ieq} is strict.
This means that any element of
the stalk $\big(\tif_*\ms D\big)_{\lmd\upone}$ vanishes at the points
$\xi\upone_2$ and $\ol\xi\upone_2$ by 
the order at least $(d_2-1)$ when restricted to 
the curves $\Xi_2$ and $\ol \Xi_2$ respectively.
This implies that the above two bold chains 
in the right pictures of Figure \ref{fig:rf7}
are fixed components of the line bundle
$\ms D'(l)|_{E}$ with
multiplicity $(d_2-1)$.
Then since the degree of the line bundle  $\ms D'(l)$ over
fibers of the projection $Q_1\to \DDD_1$ and 
$\ol Q_1\to\ol \DDD_1$ is $(d_2-1)$,
this means that the divisors $Q_1$ and $\ol Q_1$
are fixed components of the above line bundle $\ms D'(l)$ for any $l$.
Therefore 
from \eqref{id0001}, we obtain further identification
\begin{align}\label{id0003}
\big|\nu^*\ms D(l)\big|=
d_2\big(Q_1 + \ol Q_1\big)+
\big|\nu^*\ms D(l) - d_2(Q_1 + \ol Q_1)\big|.
\end{align}

As before let $\bm D$ be the effective divisor on $\tiZ$ which is an enlargement
of the divisor $D$ on $S\in |F|$ to $\tiZ$.
Let $S'$ be any fiber of the composition $W\stackrel{\nu}\to \tiZ 
\stackrel{\tif}\to \CP^1$
which is different from the present fiber $S\upone$.
Then the divisor $\nu^*\bm D - d_2(Q_1 + \ol Q_1)+lS'$ 
is effective and belongs to the
movable part of the right-hand side in \eqref{id0003}.
Moreover it does not have $Q_1$ and $\ol Q_1$ as components any more.
From this it is immediate to see from the 
intersection numbers in the right picture of Figure \ref{fig:rf7}
that the intersection matrix of the restriction of the divisor $\nu^*\bm D - d_2(Q_1 + \ol Q_1)+lS'$ 
to $S\upone_{\pm}\subset W$ is negative definite.
This implies that  intersection of an element
of the right-hand side of \eqref{id0003} and 
the divisors $S\upone_{\pm}$ is independent of a choice of
the element.
Therefore from the identification \eqref{id0003}
the same is true for the intersection of 
an element of $\big|\ms D(l)\big|$ on $\tiZ$ with 
the divisors $S\upone\dpm$.
Hence the image of the restriction homomorphism 
$H^0\big(\ms D(l)\big)\to H^0\big(\ms D(l)|_{S\upone}\big)$ is 1-dimensional.
This means $\dim V\upone = 1$, which  contradicts $\dim V\upone=3$.
Therefore the inequality \eqref{ieq} has to be an equality
(when $d_2>2$).

Finally, suppose $d_2 = 2$.
Then still the curves $\DDD\upone$ and $\ol\DDD\upone$ are
stable base curves of $\ms D$,
and let $\nu:W\to \tiZ$ be the blowing up at these curves as above.
Obviously $Q_1$ and $\ol Q_1$ are
fixed components of the system $|\nu^*\ms D(l)|$
for any $l$, and 
we still have an identification \eqref{id0001}.
But this time we cannot conclude that the two bold chains in the right picture
of Figure \ref{fig:rf7} are base curves since $2-d_2=0$.
However, if we suppose $\mu\upone>0$, 
the two bold chains have to be 
base curves of $\big|\nu^*\ms D(l)-Q_1-\ol Q_1\big|$
for any $l$.
Then the same argument as above implies $\dim V\upone=1$.
Hence we obtain $\mu\upone=0$ also in the case $d_2=2$.
\proofend

\medskip
Thus for the two reducible fibers $S\upi$, $i=1,k$,
we always have $\mu\upi = 0$
if the small resolution $\zeta\upi$ blows up the ridge components.
Alternatively if the resolution $\zeta\upi$
does not blow up
the ridge components,
we have $\mu\upi = (d_2-2)_+$ if $i=1$
and $\mu\upi = (d_k-2)_+$ if $i=k$.
If $i\neq 1,k$, we have
%
%
\begin{proposition}\label{prop:sc}
If $i\neq 1,k$ then the number $\mu\upi 
$ is independent
of a choice the small resolution $\zeta\upi$.
\end{proposition}

\proof
By \eqref{mec}, we have the relation 
\begin{align}\label{me2}
m = e + \sum_{j=1}^k \mu\upj.
\end{align}
As mentioned right after the proof of 
Proposition \ref{prop:mvo},
a replacement of the small resolution $\zeta\upi$ does not affect the value of $\mu\upj$ if $i\neq i$.
Further by Proposition \ref{prop:indep0},
the number $m$ is independent of a choice of
the small resolutions $\zeta\upone,\dots,\zeta\upk$. 
Therefore from \eqref{me2}, $\mu\upi$ cannot change
under the replacement of $\zeta\upi$
if $e$ is invariant under the same replacement.
But this is always the case by  Proposition \ref{prop:e}.\proofend
\medskip


\subsection{Examples}
\label{ss:ex1}

In this subsection we discuss several basic examples of 
the present twistor spaces.
%
They will illustrate how the twistor spaces
on $n\CP^2$ studied in \cite{Hon_Inv, Hon_Cre2} are
special among all the collection of 
the present twistor spaces.


\begin{example}\label{ex:2}
{\em
If $n=4$, 
as in Table \ref{table:2}, any coefficient
of the divisor $D$ on $S$ is either 1 or 2,
so we always have $d\,(=\max\{d_1,\dots,d_k\})=2$.
By Proposition \ref{prop:e}, this means that 
the number $e$ is independent of the total resolution $\zeta$,
and we always have $e=2$.
By results in \cite{Hon_JAG2},
the system $|2F|\,(=\big|K_Z\inv\big|)$ always includes
a 4-dimensional sub-system whose meromorphic map
gives a degree-two map to a scroll of planes over a conic.
Hence we  have $m=2$.
Thus the equality $m=d=e\,(=2)$ always holds
when $n=4$.
 \proofend
}
\end{example}

\begin{example}\label{ex:3}
{\em
If $n=5$, for some $S$, the number $e$ can really depend on
the choice of a total resolution $\zeta:\tilde Z\to\hat Z$.
Let us take $S$ of type $\Azero$ for example.
We have $k=3$ and $(d_1,d_2,d_3) = (1,3,2)$
as in Table \ref{table:3}.
From this, by Proposition \ref{prop:e}, we obtain 
\begin{align}\label{}
e = 
\begin{cases}
4 & {\text{if $\zeta\upone$ blows up $E_1$ and $\ol E_1$,}}\\
3 & {\text{otherwise.}}
\end{cases}
\end{align}
In a similar way, we obtain the same conclusion
for $S$ in the first lines of the cases for type 
$\Aone,\,\Atwo$ and $\Athree$ in Table \ref{table:3}.
Since $m\ge e$ always holds as in \eqref{dem},
from the first small resolution,
this means that the strict inequality $m>n-2\,(=3)$ holds
for these twistor spaces.
%
\proofend
}
\end{example}

\begin{example}
\label{ex:4}
{\em
As we mentioned in Section \ref{ss:I},
in \cite{Hon_Inv} and \cite{Hon_Cre2},
we gave Moishezon twistor spaces on $n\CP^2$
such that the linear system $|(n-2)F|$ induces
a degree-two map onto a
scroll of planes with degree $(n-2)$.
The number $k$, the half of the components
of the real anti-canonical cycle $C$, is $(n+1)$ in \cite{Hon_Inv} and 
$(n-1)$ in \cite{Hon_Cre2}.
By Table \ref{table5}, this means that,
in the terminology of the present paper,
they are of type $\Athree$ and
$\Aone$ respectively.

For the twistor spaces in \cite{Hon_Inv},
by Proposition 2.1 in the paper,
after a cyclic permutation for the components
to make $C_1$ the line component,
the multiplicities for components of the divisor $D$
are given by
\begin{align}\label{Inv}
(d_1,d_2,\dots,d_{n+1}) = 
(1,2,3,\dots, n-2, 1,1,1).
\end{align}
%
Namely $d_i = i$ if $i < n-1$ and $d_{n-1} =
d_n=d_{n+1} = 1$.
(In \cite{Hon_Inv}, the line component was denoted by $C_4$.)
Since $d_2=2$ and $d_k=1$, Proposition \ref{prop:e} means that the number $e$ 
does not depend on a small resolution $\zeta:\tiZ\to \hat Z$.
From the formula \eqref{case1} in Proposition \ref{prop:e}, 
we obtain $e=n-2$. 
From \eqref{Inv}, we have $d=n-2$.
The result in \cite{Hon_Inv} means that
the system $|(n-2)F|$ on the twistor space $Z$
itself induces a degree-two map to the scroll of planes
of degree $n$.
Hence $m=n-2$ holds.
Therefore $m=e=d\,(=n-2)$ holds.
Hence both of the inequalities \eqref{dem} are equality.

For the twistor spaces in \cite{Hon_Cre2},
by equation (2.3) in the paper,
after a cyclic permutation for the components
to make $C_1$ the line component,
the multiplicities for the divisor $D$
are given by
\begin{align}\label{m11}
(d_1,d_2,\dots,d_{n-1}) = 
(1,2,3,\dots, n-2,1).
\end{align}
Namely $d_i=i$ if $i<n-1$ and $d_{n-1} = 1$.
(In \cite{Hon_Cre2},
the line component was denoted by $C_2$.)
%
%
Again Proposition \ref{prop:e} means that the number $e$ 
does not depend on the small resolution,
and again we have $e=n-2$.  
From \eqref{m11}, we have $d=n-2$.
The result in \cite{Hon_Cre2} means that
the system $|(n-2)F|$ on the twistor space $Z$
itself induces a degree-two map to the scroll of planes
of degree $n$. Hence $m=n-2$.
Therefore again we have
equality $m=e=d\,(=n-2)$.

An example which interpolates these two examples
is given by $S$ of type $\Atwo$
whose divisor $D$ satisfies 
\begin{align}\label{m12}
(d_1,d_2,\dots,d_{n}) = 
(1,2,3,\dots, n-2,1,1).
\end{align}
Namely $d_i=i$ if $i<n-1$ and $d_{n-1} = d_n=1$.
These surfaces can be obtained inductively from the one in the case
$n=4$ (the first line in the case of type $\Atwo$ in Table \ref{table:2}) by succession of blowing up at appropriate pair of double points
of the real anti-canonical cycle, which would be easily found.
By the same reason to the above two examples,
we have $e=n-2$, regardless of 
a choice of the small resolution.
A twistor space on $n\CP^2$ which has
this surface $S$ as a real fundamental divisor seems to have not appeared in the literature if $n>4$.
\proofend
}
\end{example}

\begin{example}\label{ex:5}{\em
Next as examples in the opposite direction, we give a twistor space
on $n\CP^2$
whose number $m$ rapidly increases
as $n$ does.
Let $n>4$ and $S$ be the surface satisfying $K^2=8-2n$
which is obtained inductively
from the one in the case $n=4$, type $\Azero$,
by repetition of blowing up 
at the double points $C_i\cap C_{i+1}$ and $\ol C_i\cap\ol C_{i+1}$
which are chosen in such a way that the two numbers $d_i$ and $d_{i+1}$ are
the first two biggest values among the multiplicities $d_1,d_2,\dots, d_k$
of the effective divisor $D$.
\begin{table}[h]
\begin{flushleft}
\begin{tabular}{c||c|c|c|c|c|c}
$n$ & 4 & 5 & 6 & 7 & 8 \\
\hline
$(d_1,\dots, d_k)$ & $(1,2)$ & $(1,3,2)$ & $(1,3,5,2)$ &
$(1,3,8,5,2)$ & $(1,3,8,13,5,2)$ 
\\
\hline
$e$ & 2 & 4 & 6 & 9 & 14 
\end{tabular}
\end{flushleft}
\begin{flushright}
\begin{tabular}{|c|c|c}
 9 & 10 & 11\\
\hline
$(1,3,8,21, 13,5,2)$ & $(1,3,8,21,34, 13,5,2)$ & $\cdots$
\\
\hline
 22 & 35 & 57
\end{tabular}
\end{flushright}
\caption{Example \ref{ex:5}}
\label{table11}
\end{table}
These multiplicities are displayed in Table \ref{table11} for small $n$.
Since $d_2 = 3$ if $n>4$,
the number $e$ depends on the choice of 
the small resolution $\zeta\upone$,
and it becomes larger if $\zeta\upone$ blows up the ridge components $E_1$ and $\ol E_1$.
%
In Table \ref{table11}, the number $e$ is displayed
for such a small resolution.
This again follows from the formula
in Proposition \ref{prop:e}.
In these examples, the value $d=\max\{d_1,\dots, d_{n-2}\}$ is the Fibonacci number
$F(n)$, and we have
\begin{align}\label{Fib}
e= F(n) +1.
\end{align}
Hence $m\ge F(n)+1$ holds, and this is much larger than $(n-2)$.
(Recall that the twistor space is more general if $m$ is larger.)
It is likely that the number $F(n)+1$ is
the maximum for the number $e$ when we consider all the present twistor spaces
on $n\CP^2$.}
\proofend
\end{example}

Examples \ref{ex:4} and \ref{ex:5} are
two extreme cases with respect to the value of the number $e$.
It is likely that most values in between 
$(n-2)$ and $F(n)+1$ are realized by a twistor space
on $n\CP^2$ as the number $e$ for some small resolution.

Finally we give 
an example of twistor space on $n\CP^2$
for which the number $\mu\upi$ (namely the dimension of $\Coker\,\rho\upi$) becomes larger
for some index $i\neq 1,k$.

\begin{example}\label{ex:rho}
{\em
First we give the surface $S$ which will be a member
of the pencil $|F|$ on a twistor space $Z$ over $n\CP^2$, $n>6$,
as a blowing up of 
the one with $n=5$ which is of type $\Azero$.
As in Table \ref{table:3}, there is only one 
kind of such a surface and the self-intersection numbers of components of the real anti-canonical cycle $C$
on the surface
are $-4,-1,-2,-4,-1,-2$. The 
line components are the pair of $(-4)$-curves.
We first blow up the intersection of
the line components and the $(-1)$-curves.
Next we again blow up the surface at the intersection of the line components
and the new $(-1)$-curves.
Repeating this procedure $(n-6)$ times, 
we get a surface whose components of the  anti-canonical cycle
have the self-intersection numbers
$$
-(n-2), -1, {\overbrace{-2,-2,\cdots,-2}^{n-5}},
-(n-2), -1, {\overbrace{-2,-2,\cdots,-2}^{n-5}}.
$$
This surface satisfies $K^2 = 8-2(n-1)$.
Finally, we blow up this surface at the intersection of $(-1)$-curves and the adjacent $(-2)$-curves.
Let $S$ be the resulting surface and $C$ the real anti-canonical cycle on $S$.
$S$ is of type $\Azero$.
This surface satisfies $K^2 = 8-2n$,
and the sequence for the self-intersection numbers of components of $C$ is given by
$$
-(n-2), -2,-1, -3,{\overbrace{-2,-2,\cdots,-2}^{n-6}},
-(n-2), -2,-1, -3,{\overbrace{-2,-2,\cdots,-2}^{n-6}}.
$$
If $D$ denotes the effective divisor on $S$
which induces the double covering map $\phi:S\to\CP^2$
as before, the multiplicities of the components of $C$ in $D$ are given in order by
$$
1, n-3,2n-7, {\overbrace{n-4, n-5,\dots, 3,2}^{n-6}},
1, n-3,2n-7, {\overbrace{n-4, n-5,\dots, 3,2}^{n-6}}.
$$

As we discussed soon after Definition \ref{def:type},
a twistor space $Z$ on $n\CP^2$ 
which contains this surface $S$ as a real member of the pencil $|F|$ exists.
We show that for this twistor space, 
$\mu\uptwo = n-5$ holds.

\begin{figure}
\includegraphics{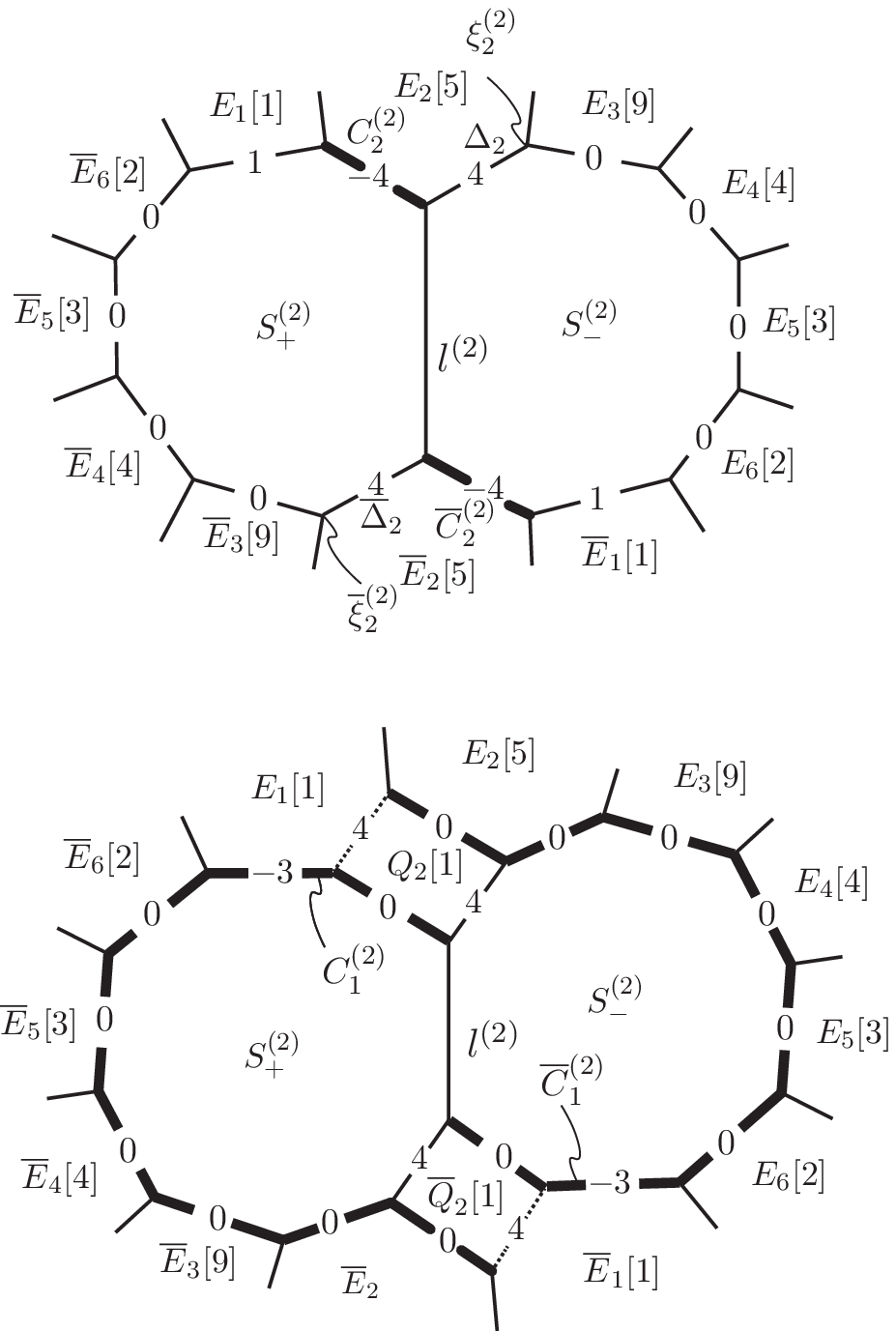}
\caption{The case $n=8$ in Example \ref{ex:rho}.
The numbers in the brackets denote
the multiplicities.}
\label{fig:n8}
\end{figure}
The main idea to show this is similar to the proof for 
the second item in Proposition \ref{prop:cokrho1}.
As a small resolution $\zeta\uptwo$,
we take the one which blows up 
the components $E_2$ and $\ol E_2$.
(By Proposition \ref{prop:sc}, another choice of resolution
does not affect the value of $\mu\uptwo$.)
The case $n=8$ is displayed in the upper picture
in Figure \ref{fig:n8}.
The curves $C\uptwo_2$ and $\ol C\uptwo_2$
are the stable base curves on $S\uptwo$
and the multiplicity is $(n-4)$.
In particular, each of them consists of a single component.
Again let $\nu:W\to \tiZ$ be the blowing up
at these curves, and $Q_2$, $\ol Q_2$
the exceptional divisors.
Since the curves 
$C\uptwo_2$ and $\ol C\uptwo_2$ have normal bundle $\ms O(-1)^{\oplus 2}$ in $\tiZ$,
the divisors $Q_2$ and $\ol Q_2$ are again isomorphic to $\qdr$.
If $\bm D$ denotes the formal extension of $D$ on $S$ to $\tiZ$
as before,
the pull-back $\nu^*\bm D$ includes
the divisors $Q_2$ and $\ol Q_2$ by 
multiplicity $d_2=n-3$,
and they  are fixed components
of $\nu^*\ms D(l)$ with multiplicity $(n-4)$
 for any $l$.
Hence similarly to \eqref{id0001},
we have an identification
\begin{align}\label{id0002}
\big|\ms D(l)\big|=
(n - 4) \big(Q_2 + \ol Q_2\big)+
\big|\nu^*\ms D(l) - (n - 4) (Q_2 + \ol Q_2)\big|
\end{align}
for any $l$.
Again we write the  line bundle $\nu^*\ms D(l) - (n - 4) (Q_2 + \ol Q_2)$
on the right-hand side
by $\ms D'(l)$.
Then the degrees of this line bundle
over curves on the divisor $\nu\inv\big(S\uptwo\big)$ can be easily computed,
 and  in particular we obtain 
the key intersection number
\begin{align}\label{ki}
\ms D'(l).\,C\uptwo_1 = -(n-5).
\end{align}
In Figure \ref{fig:n8} the degrees of the line bundle $\ms D'(l)$ are displayed
 in the lower picture for the case $n=8$.
From \eqref{ki}, the two chains indicated by 
the bold line segments are base curve
of the line bundle $\ms D'(l)|_E$ with multiplicity
 $(n-5)$.
By Proposition \ref{prop:mvo},
this implies an inequality
$
\mu\uptwo \ge n-5.
$
The assertion that this has to be the equality
can be shown in the same way to the second item
in Proposition \ref{prop:cokrho1}
by contradiction using $\dim V\uptwo = 3$.
Hence we obtain $\mu\uptwo = n-5$.}
\proofend
\end{example}

This example shows that, 
the number $\mu\upi$ can be arbitrary large
if we allow $n$ to be large.
Therefore $m$ can be quite larger than $e$.

\section{Existence of real reducible members of $|mF|$}
\subsection{Reducible members of the system $|mF|$ and 
real bitangents}\label{ss:rm}
As before, let $m$ be the number which satisfies
$\tif_*\ms D\simeq\ms O(-m)^{\oplus 2}\oplus \ms O$.
From the results in Sections \ref{s:dim} and \ref{s:sd},
we know that the system $|mF|$ or its sub-system 
on the present twistor space induces a degree-two meromorphic map
$\Phi_m:Z\to Y_m\subset\CP^{m+2}$ where $Y_m$ is  
a scroll of planes over a rational normal curve of degree $m$
(Theorem \ref{thm:dc1}).
Further,  the branch divisor of $\Phi_m$ is a cut
of the scroll by a quartic hypersurface (Theorem \ref{thm:dc2}).


In this section, we investigate the branch divisor
on the scroll of planes
more closely, and show that the quartic 
hypersurface which cuts out the branch divisor
is not generic but must
satisfy a strong constraint.
Needless to say, this is an essential portion
in 
understanding the present twistor spaces.
%
In the case of $3\CP^2$, 
the scroll is nothing but $\CP^3$,
and the branch divisor  is a quartic surface itself.
The constraint for this quartic was intensively 
studied by Kreussler-Kurke \cite{KK92}
and Poon \cite{P92}.
We now briefly describe it
in order to obtain a perspective for
generalizing to the case  $n>3$.

If $B$ is a quartic surface in $\CP^3$,
a conic $\ms C\subset B$ is called a {\em trope} of $B$
if a plane is tangent to $B$ along $\ms C$.
If $H$ is the hyperplane on which $\ms C$ lies,
this means $B|_H=2\ms C$.
In the case of $3\CP^2$, 
the principal role is played by this kind of conics.
Namely, it is shown that there exist exactly 
four tropes $\ms C\upone,\dots, \ms C^{(4)}$
of the branch divisor $B\subset\CP^3$, as well as a quadratic surface
containing all  these tropes.
From these it is shown that the equation of $B$
is of the form 
\begin{align}\label{B001}
Q^2 = h_1h_2h_3h_4,
\end{align}
where $Q$ is an equation of the quadric
and each $h_i$ is that of the hyperplane 
on which the trope $\ms C\upi$ lies.
The equation of the trope $\ms C\upi$ is $Q=h_i=0$.
For any distinct indices $i$ and $j$,
the intersection $\ms C\upi\cap\ms C\upj$
consists of two points, and 
the quartic surface defined by the equation \eqref{B001} has
singularities at these points.
If the twistor space on $3\CP^2$
is generic in the sense that
the four hyperplanes $\{h_i=0\}$ are linearly independent,
all the singularities are ordinary double points.
If they are not linearly independent,
the singularities become worse than 
ordinary double points.
This issue was investigated in \cite{KK92}
in detail.

The existence of the four tropes $\ms C\upone,\dots,
\ms C^{(4)}$ is a consequence
of the fact that the linear system $|F|$ always
has exactly four real reducible members.
Namely if $S\upi_+ + S\upi_-$, $1\le i\le 4$, are such  members
and $H\upi$ are the real hyperplanes corresponding to these four members,
the map $\Phi:Z\to\CP^3$ induced by $|F|$
maps $S\upi_+$ and $S\upi_-$ 
birationally onto $H\upi$,
and the trope $\ms C\upi$ is a double branch divisor
of the restriction $\Phi|_{S\upi_+\cup S\upi_-}\to H\upi$.
The intersection $S\upi_+\cap S\upi_-$ is a twistor line
as in Proposition \ref{prop:poon1},
and $\Phi$ maps it isomorphically onto $\ms C\upi$.
The existence of the four real reducible fundamental
divisors is deduced by 
concretely specifying the cohomology classes
of irreducible components $S\upi_+$ and $S\upi_-$
for each $i$.

Thus existence of reducible members of the system 
$|F|$ was crucial for deriving the equation
of the branch quartic.
In the case of $4\CP^2$, it was shown in
\cite{Hon_JAG2} that the system $|2F|$  contains some 
real reducible members, and again this was done by specifying
the cohomology classes of the components.
Also, in the case of $n\CP^2$, 
analogous result was shown to hold in \cite{Hon_Cre2} for
the system $|(n-2)F|$ in a similar way.
These results were used in order to show that 
the quartic hypersurface which cuts out the branch divisor 
from the scroll is defined by the equation of the form \eqref{B001}.
In view of these results, it would be natural to expect that, 
for any present twistor spaces, the linear system $|mF|$
 has some real reducible members, and as a consequence,
the quartic hypersurface is defined by an equation of the same form
as \eqref{B001}.
However, the cohomology classes obtained in \cite{Hon_Cre2}
were already complicated, 
and it seems difficult (to the author) to generalize it
to the present situation in full generality.

The idea we adapt in this paper is to make use of real 
bitangents of a plane quartic.
Concretely, if there is a real reducible member
of $|mF|$ and $H\subset\CP^{m+2}$ is a 
real hyperplane which corresponds to
this reducible member,
the cut $B\cap H$ of the branch divisor $B$ has to be a double curve.
We call this double curve as a {\em trope
in the generalized sense}.
If $m=1$, this agrees with the trope in the original sense.
If $H$ is generic hyperplane and $\CP^2$ is a generic plane 
of the scroll, the intersection line $H\cap \CP^2$
is a  bitangent of the plane quartic
$B\cap \CP^2$.
Thus  tropes in the generalized sense
can be detected on generic  planes of the scroll
as tangent points on  bitangents.
If the hyperplane $H$ and the plane $\CP^2$ are real,
the line $H\cap \CP^2$ is real, and this is 
a real bitangent of the quartic $B\cap\CP^2$ 
which is also real.
Therefore existence of a real bitangent for
a real plane section of $B$ is a necessary
condition for the system $|mF|$ to
have a real reducible member.
So in the next subsection we investigate
real bitangents to the branch plane quartic,
for the degree-two morphism $\phi:S\to\CP^2$
induced by the linear system $|D|$ 
on a real irreducible member $S$ of the pencil $|F|$.
%
%
%
%

\subsection{Real bitangents of plane quartics}
\label{ss:rb}
\begin{figure}
\includegraphics{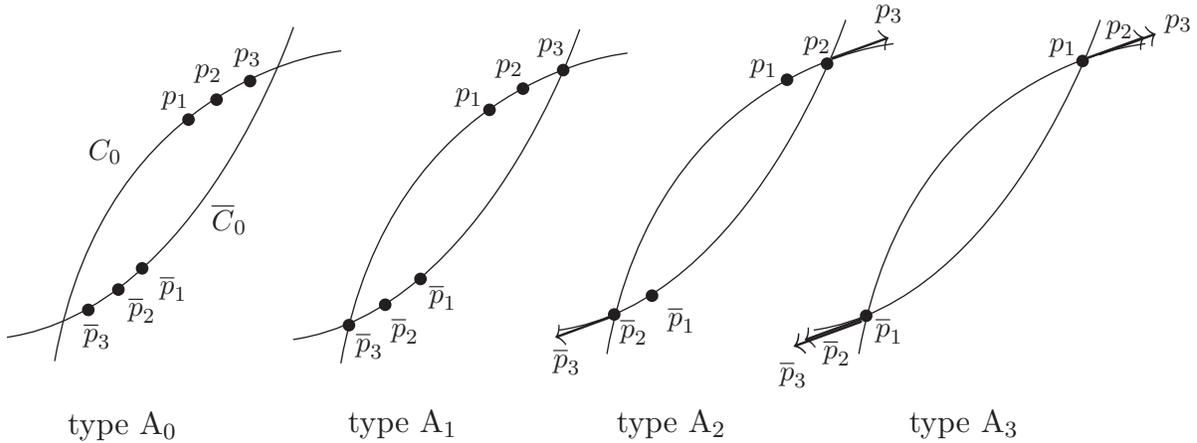}
\caption{
Configuration of the six points on $S_0=\qdr$ which yield the intermediate surface $S_3$.
}
\label{fig:bu_cfg2}
\end{figure}

In the following, for simplicity of presentation,
we always suppose $n>3$.
This implies that the fundamental system $|F|$ of the twistor space is a pencil.
We recall from  \eqref{diagram:fac1} that 
if $S$ is 
a real irreducible member of this pencil,
the morphism $\phi:S\to\CP^2$ induced by the system $|D|$ admits a factorization
\begin{align}\label{diagram1'}
 \xymatrix{ 
S \ar[r] \ar[dr]_{\phi} & S_3  \ar[d]^{\phi_3}\\
&\CP^2 
   }
\end{align}
where 
$S_3$ is a surface which is obtained from 
$S_0=\qdr$ by blowing up 6 points
$$
p_1,p_2,p_3,\ol p_1,\ol p_2, \ol p_3\in\qdr
$$ 
located as in Figure \ref{fig:bu_cfg2} (=Figure \ref{fig:bu_cfg1}) depending on the type,
$S\to S_3$ is iteration of pairwise contractions of $(-1)$-components
in the real anti-canonical cycle,
and $\phi_3:S_3\to \CP^2$ is the anti-canonical morphism from $S_3$.
In particular, the branch quartic of $\phi$
equals that of the anti-canonical map $\phi_3$.
The surface $S_3$ has distinguished
real anti-canonical cycle which is determined 
from the two $(1,1)$-curves $C_0$ and $\ol C_0$
on which the $3+3$ points lie,
and the effective divisor $D$ on $S$ is
the pull-back of this anti-canonical cycle
by the birational morphism $S\to S_3$.
%
%

Recall that the branch quartic of $\phi:S\to\CP^2$
has always two distinguished points
$\bm q$ and $\ol{\bm q}$ which are the images
of the two chains $C-C_1-\ol C_1$ under $\phi$.
The line passing the points $\bm q$ and $\ol{\bm q}$
was denoted by the bold letter $\bm l$.
By a {\em bitangent} of the branch quartic of $\phi$,
we mean a line which is tangent to the quartic
at two {\em smooth} points of the quartic.
So if the surface $S$ is of type 
$\Aone,\Atwo$ or $\Athree$ 
we do not call the line $\bm l$ a bitangent because
in these cases $\bm q$ and $\ol{\bm q}$
are singular points of the quartic.
If $\ell$ is a bitangent of a quartic
in this sense and if $\ell\neq\bm l$, the inverse image $\phi\inv(\ell)$
splits into two components, and each component is a $(-1)$-curve on $S$
which is mapped isomorphically to the bitangent $\ell$.
Moreover these two components intersect each other  transversally at two points
over the tangent points.
Conversely, any $(-1)$-curve on $S_3$
 is mapped to a bitangent of the branch quartic,
and two $(-1)$-curves are mapped to the same bitangent
if and only if they intersect transversally at two points or they are exactly the pair of line components
in the real anti-canonical cycle.


In this subsection we investigate real bitangents of a real plane quartic
and determine its number.
We begin with the case  of type $\Azero$.

\begin{proposition}\label{prop:btan1}
If a twistor space $Z$ on $n\CP^2$ is of type 
$\Azero$ and $S$ is a real irreducible member of the pencil $|F|$, the branch quartic of 
the degree-two morphism $\phi:S\to\CP^2$ induced 
by the linear system $|D|$ has
exactly four real bitangents.
\end{proposition}

\proof
We find all $(-1)$-curves on the surface $S_3$ in the diagram \eqref{diagram1'},
and next find the partner for each one to detect  real bitangents.
First, the 6 exceptional curves of the blowing up $S_3\to S_0=\qdr$ 
are of course $(-1)$-curves.
Besides these, the strict transforms of
the following curves on $S_0=\qdr$ give $(-1)$-curves on $S_3$:
\begin{itemize}
\item the $(1,0)$-curves and the $(0,1)$-curves passing one of 
the 6 points $p_i$ and $\ol p_i$,
\item the $(1,1)$-curves  passing 3 of 
the 6 points,
\item the $(1,2)$-curves and the $(2,1)$-curves passing 5 of 
the 6 points,
\item the $(2,2)$-curves which have a
node at one of the 6 points
and which passes the remaining 5 points.
\end{itemize}
\noindent
The number of the $(-1)$-curves obtained from each item is, respectively,
$$6+6,\quad \begin{pmatrix}6\\3\end{pmatrix} = 20,\quad 6+6\qandq 6.$$
The images of all these $(-1)$-curves are bitangents.
Collecting these, we obtain fifty-six $(-1)$-curves.
On the other hand, as is well-known, 
the number of bitangents of a plane quartic is at most 28,
and the maximal number $28$ is attained by a smooth quartic.
These mean that each of the fifty-six $(-1)$-curves has the partner among these curves,
and the images of these 56 curves exhaust all bitangents.
In particular, the branch quartic is smooth.
Thus we have obtained that the branch quartic in the case of type $\Azero$
is smooth, and that any bitangent is the image of some of the above fifty-six
$(-1)$-curves.

Among these $(-1)$-curves,
each of the 6 exceptional curves of the blowing up $S_3\to S_0=\qdr$
intersects exactly one of the
strict transforms of the $(2,2)$-curves at two points from the nodal condition.
Therefore the partner of each of the 6 exceptional curves is 
among the strict transforms
of the nodal $(2,2)$-curves.
So from these we obtain 6 bitangents.
But since the real structure was a lift from 
that on $S_0$,
obviously none of these bitangents is real.
Similarly, the partners of the strict transforms
of the $(1,0)$-curves are
the strict transforms of the $(1,2)$-curves.
From these we obtain 6 bitangents.
Again none of these are real.
Also, the strict transforms of the $(0,1)$-curves and those of the
$(2,1)$-curves form 6 pairs, and they give 6 bitangents, which are again non-real.
So far we have obtained $6\times 3 = 18$ bitangents, and all of them are 
non-real.

The remainder is $(56-2\cdot18=)$ twenty $(-1)$-curves, and all of them are the strict transforms of the $(1,1)$-curves passing three of the six points.
To each of these $(1,1)$-curves, one can associate another $(1,1)$-curve
by taking the complement of the three points in the set of six points.
Hence we obtain $(20/2=)\,10$ pairs of $(-1)$-curves.
Each of these pairs give the same bitangent, and we get 10 bitangents.
Among these, the real one is obtained from the pairs
whose three points are mutually conjugate.
Concretely, the two $(1,1)$-curves passing the three points
$\{p_1,p_2,p_3\}$ and $\{\ol p_1,\ol p_2,\ol p_3\}$ are mutually conjugate,
and these give the same real bitangent.
The same is true for the $(1,1)$-curves passing the following three pairs of three points
\begin{gather}
\{\ol p_1, p_2, p_3\} \qandq \{p_1, \ol p_2, \ol p_3\},\notag\\
\{p_1,\ol p_2, p_3\} \qandq \{\ol p_1, p_2, \ol p_3\},\label{pa0}\\
\{p_1,p_2,\ol p_3\}  \qandq \{\ol p_1,\ol p_2, p_3\}\notag.
\end{gather}
In contrast, the $(1,1)$-curve through the three points $\{p_1,\ol p_1,p_2\}$
for example,
gives a non-real bitangent,
and the same thing holds for all remaining $(1,1)$-curves.
Thus we can conclude that the number of real bitangents is $1+3=4$.
\proofend

\medskip
Among these four real bitangents, 
 there is a special one
characterized by the property that 
$\phi\inv (l) = D$.
In the above proof this real bitangent is 
nothing but the first one which is determined by 
the $(1,1)$-curve through the three points
$p_1,p_2$ and $p_3$.
This real bitangent will be excluded in the sequel
since it is mapped to the ridge $\bm l$ of the scroll of planes,
and does not form a trope in the generalized sense.

Next we discuss the cases of type $\Aone$ and type $\Atwo$.
We recall from Section \ref{ss:acsfd} that in these cases the branch quartic has
$\Aone$-singularities
and $\Atwo$-singularities respectively
at the points $\bm q$ and $\ol{\bm q}$.
By our convention the line $\bm l$
is not a bitangent.

\begin{proposition}
\label{prop:btan2}
If a twistor space $Z$ on $n\CP^2$ is of type $\Aone$
(resp.\,type $\Atwo$) and $S$ is a 
real irreducible member of the pencil $|F|$, the branch quartic of 
the degree-two morphism $\phi:S\to\CP^2$ induced by the linear system $|D|$
has 
exactly 2 (resp.\,1) real bitangents.
\end{proposition}

\proof
Again we find all bitangents in the present sense  and  pick up real ones.
First we consider the case of type $\Aone$.
Let $\{p_1,p_2,p_3,\ol p_1,\ol p_2, \ol p_3\}$
be the six points on $S_0=\qdr$ as in Figure \ref{fig:bu_cfg2}.
The intermediate surface $S_3$ is obtained from $S_0$
by blowing up  these six points.
The anti-canonical cycle $C_0 + \ol C_0$ in Figure \ref{fig:bu_cfg2} determines a 
real anti-canonical cycle on $S_3$
consisting of four components, two of which are $(-1)$-curves and the remaining two are
$(-2)$-curves.
The $(-1)$-curves are the exceptional curves over
the point $p_3$ and $\ol p_3$,
and the $(-2)$-curves are the strict transforms
of the curves $C_0$ and $\ol C_0$.
The points $\bm q$ and $\ol{\bm q}$ are exactly 
the images of these $(-2)$-curves.

As we are not allowing bitangents to pass singularities of the quartic,
a $(-1)$-curve which corresponds to a bitangent has to be disjoint from the two $(-2)$-curves
in the real anti-canonical cycle.
Conversely, the image of a $(-1)$-curve on $S_3$
which does not intersect these $(-2)$-curves is a bitangent.
From the Pl\"ucker's formula \cite[Chapter 3,\,Section 4]{GH},
the presence of the two $\Aone$-singularities
of the quartic means that 
the number of bitangents in the present sense is at most 8.
In the following we concretely present sixteen $(-1)$-curves on the surface $S_3$ which do not
intersect the two $(-2)$-curves.

First we consider the pair of $(-1)$-curves
which are the strict transforms of the $(1,1)$-curves
passing respectively the points
\begin{align}\label{pa1}
\{\ol p_1, p_2, p_3\}
\qandq
\{p_1,\ol p_2, \ol p_3\}.
\end{align}
As an effect of the blowing up $S_3\to S_0=\qdr$, 
these two $(-1)$-curves do not intersect the two $(-2)$-curves on $S_3$,
and intersect each other at two points.
Hence these $(-1)$-curves are partners to each other,
and the image to $\CP^2$ is the same bitangent.
Moreover, since \eqref{pa1} is  a conjugate pair,
this bitangent is real.
In the same way, the pair of $(-1)$-curves
which are the strict transforms of the $(1,1)$-curves
passing respectively the points
\begin{align}\label{pa2}
\{p_1,\ol p_2, p_3\}
\qandq
\{\ol p_1, p_2, \ol p_3\}
\end{align}
give the same real bitangent.
Thus we have found two real bitangents.

Next the strict transform of the $(1,1)$-curve
passing the points $p_1,\ol p_1$ and $p_3$,
and that of the $(1,1)$-curve
passing the points $p_2,\ol p_2$ and $\ol p_3$
define the same bitangent.
Also, the strict transform of the $(1,1)$-curve
passing the points $p_1,\ol p_1$ and $\ol p_3$,
and that of the $(1,1)$-curve
passing the points $p_2,\ol p_2$ and $ p_3$
define the same bitangent.
These two bitangents are mutually conjugate.
Next the strict transform of the $(1,2)$-curve
passing all six but the point $p_3$ determines
a bitangent, and this bitangent is also obtained as the image 
of the strict transform of the $(1,0)$-curve passing the point $p_3$.
This bitangent is clearly non-real,
and we obtain another bitangent as its conjugation.
Similarly, the strict transform of the $(2,1)$-curve
passing all six but the point $p_3$ determines
a bitangent, and this is also obtained as the image 
of the strict transform of the $(0,1)$-curve passing the point $p_3$.
Again this bitangent is non-real
and we obtain another bitangent as its conjugation.

Thus we have obtained two real bitangents and six non-real ones.
By Pl\"ucker formula as above, these have to be all bitangents
of the branch quartic, and the quartic has no singularity
other than the two ordinary double points.
Thus the number of real bitangents is two in the case of type $\Aone$.

Next we consider the case of type $\Atwo$,
which is much easier than the case of type $\Aone$.
The points $\bm q$ and $\ol{\bm q}$ are $\Atwo$-singularities of the branch quartic.
From Pl\"ucker's formula, this means that 
the branch quartic has at most one bitangent
in the present sense.
The anti-canonical cycle $C_0+\ol C_0$
in Figure \ref{fig:bu_cfg2} determines a real anti-canonical cycle on $S_3$ which consists of 6 components.
Among the 6 components of the real anti-canonical cycle,
exactly two are $(-1)$-curves,
and the rest are $(-2)$-curves.
The image of the four $(-2)$-curves are the conjugate pair
of $\Atwo$-singularities, $\bm q$ and $\ol{\bm q}$.
The line $\bm l$ passing these two points is the image of the exceptional
curves over the points $p_3$ and $\ol p_3$.

We consider 
the pair of $(-1)$-curves on $S_3$
which are obtained as the strict transforms of the two $(1,1)$-curves
passing the following sets of three points respectively:
\begin{align}\label{pa3}
\{\ol p_1,  p_2, p_3\}
\qandq
\{p_1, \ol p_2,\ol p_3\}.
\end{align}
Here, passing the points $p_2$ and $p_3$ means
that it is tangent to the $(1,1)$-curve
through the points $p_1,p_2$ and $p_3$
(the curve $C_0$ in the figure)
at the point $p_2$.
These $(-1)$-curves do not intersect the four $(-2)$-curves,
and mutually intersect at two points.
Moreover, these are mutually conjugate.
Therefore these determine the same real bitangent.
Then as above, Pl\"ucker formula means that the quartic has
no singularities besides the two $\Atwo$-singularities,
and the above real bitangent is the unique bitangent.
In particular, the branch quartic has precisely one real bitangent.
\proofend
\medskip

We recall that for the case of type $\Athree$,
the branch quartic splits into 
two  conics which are tangent to each other at two points. 
These conics do not have a common tangent
other than the two tangent  at each of the two points.
So there is no bitangent in the present sense.

We now summarize what we have obtained
about real bitangents through
the proof of Propositions \ref{prop:btan1} and \ref{prop:btan2} and express the cohomology 
classes of the $(-1)$-curves on $S$
which are mapped to real bitangents.
Let $S$ be a real irreducible member of 
the pencil $|F|$, and
$\phi:S\to\CP^2$ the double covering map 
induced by the system $|D|$.
The surface $S$ can be obtained as the composition
\begin{align}\label{ee}
S\stackrel{\psi_1}\lras S_3\stackrel{\psi_0}\lras S_0=\qdr,
\end{align}
where $\psi_0$ is the blowing up at the six points
$p_1,p_2,p_3,\ol p_1,\ol p_2$ and $\ol p_3$,
and $\psi_1$ is a succession of pairwise blowing up 
at double points of the real anti-canonical cycle.
The configuration of the six points are as in Figure \ref{fig:bu_cfg2}.
We let 
$$
e_{\aaa}\qandq\ol e_{\aaa},\quad1\le \aaa\le 3
$$
 be the cohomology classes
on $S$ represented by the pull-back 
of the exceptional curves over 
the points $p_{\aaa}$ and $\ol p_{\aaa}$ respectively, 
under $\psi_1:S\to S_3$.
The classes of the exceptional curves
of  $\psi_1$ are denoted by $e_4,\dots,e_n$
and $\ol e_4,\dots, \ol e_n$.
Further we write $\psi:=\psi_0\circ\psi_1:S\to S_0=\qdr$.
From the concrete description
of real bitangents given as the strict transforms of the $(1,1)$-curves passing the three points \eqref{pa0}--\eqref{pa3},
we have the following
\begin{proposition}\label{prop:cc01}
Under the above notation,
if $\ell$ is a real bitangent of the branch quartic
of $\phi$,
the cohomology classes of the two irreducible
components 
(which are $(-1)$-curves) of the inverse image  $\phi\inv(\ell)$ are given by
\begin{align}\label{cc01}
\psi^*\ms O(1,1) - 
\ol e_{\aaa} -  e_{\beta} -   e_{\ccc}
\quad{\text{and}}\quad
\psi^*\ms O(1,1) - 
 e_{\aaa} - \ol e_{\beta} - \ol e_{\ccc},
\end{align}
where $({\aaa},{\beta},{\ccc})$ takes any one of the following
values.
\begin{itemize}
\item $({\aaa},{\beta},{\ccc}) = (1,2,3),\,(2,3,1),\,(3,1,2)$
in the case of type $\Azero$,
\item $({\aaa},{\beta},{\ccc}) = (1,2,3),\,(2,3,1)$
in the case of type $\Aone$,
\item $({\aaa},{\beta},{\ccc}) = (1,2,3)$
in the case of type $\Atwo$.
\end{itemize}
\end{proposition}

Simply speaking, the set of real bitangents
is in one-to-one correspondence with
the set of points among $p_1,p_2$ and $p_3$
which are not on the intersection of the
two $(1,1)$-curves $C_0$ and $\ol C_0$ 
on which the six points lie (see Figure \ref{fig:bu_cfg2}).
For example, in the case of type $\Azero$,
all three points $p_1,p_2$ and $p_3$
are not on the intersection, and
the branch quartic has three real bitangents.
Of course, the sum of the two classes
\eqref{cc01} is linearly equivalent to
the class of the divisor $D$,
and in this sense each of the classes
\eqref{cc01} is a `half' of
the effective divisor $D$.
Also we remark that,
for the case of type $\Azero$, as we have mentioned
right after the proof of Proposition \ref{prop:btan1},
the real bitangent whose inverse image
under $\phi$ is the divisor $D$ is excluded,
so that we have three classes as in the first item in
the present proposition.

From the geometric picture described in the last subsection,
we expect that
when a plane moves in the scroll,
 the tangent points of these bitangents constitute
a trope in the generalized sense.
In the next subsection we will show that this is really the case, again by working on
the blown-up space $\tilde Z$.

\subsection{Certain non-real line bundles on $\tilde Z$}
\label{ss:cc}
From the definition of the trope in the generalized sense,
the locus of tangent points of real bitangents form
a trope if the locus is lying on a hyperplane in $\CP^{m+2}$.
As discussed in  Section \ref{ss:rm}, this means that 
the member of the system $|mF|$ corresponding to
this hyperplane is  a real reducible member of $|mF|$,
and it will induce a constraint for the equation of the
branch divisor $B$ on the scroll of planes.
In the last subsection we have concretely specified
the cohomology classes on a real irreducible member $S$ of the pencil $|F|$,
in which one of the two $(-1)$-curves over
 a real bitangent belongs (Proposition \ref{prop:cc01}).
%
In this subsection
we find Chern classes on $\tilde Z$
such that their restriction to a real irreducible 
fiber of $\tilde f:\tilde Z\to\CP^1$ is one of 
the classes \eqref{cc01}.
This is possible by the following observation.

\begin{proposition}
\label{prop:cc02}
The two classes \eqref{cc01} on $S$ can be written in the form
\begin{align}\label{cc03}
D^{\rm h} + e_{\aaa} - \ol e_{\aaa}
\quad{\text{and}}\quad
\ol D^{\rm h} + \ol e_{\aaa} -  e_{\aaa}
\end{align}
respectively, where $D^{\rm h}$ is (the class of)
an effective curve on $S$
 which satisfies 
\begin{align}\label{half}
\Dh+\oDh = D.
\end{align}
Moreover the effective divisor $\Dh$ is 
independent on a choice of $\aaa$.
\end{proposition}

%


\proof
Evidently we have
\begin{align*}
\psi^*\ms O(1,1) - 
\ol e_{\aaa} -  e_{\beta} -  e_{\ccc}
&= \{\psi^*\ms O(1,1) -  e_{\aaa} - e_{\beta} -  e_{\ccc}\}
+ ( e_{\aaa} - \ol e_{\aaa})\\
&=
\{\psi^*\ms O(1,1) -  e_{1} - e_{2} -  e_{3}\}
+ ( e_{\aaa} - \ol e_{\aaa})
\end{align*}
where the last equality is from $\{\aaa,\bbb,\ccc\} = \{1,2,3\}$.
If we write $\psi$ as $\psi_0\circ\psi_1$
as in \eqref{ee},
the class $\psi_0^*\ms O(1,1) -  e_{\aaa} -  e_{\beta} -  e_{\ccc}$
is represented by a unique effective curve,
and this curve is a chain of smooth rational curves,
consisting of $1,2$ or $3$ components 
depending on types $\Azero,\Aone$ or $\Atwo$ respectively.
%
%
%
%
These curves are of course components of the real anti-canonical
cycle on $S_3$, and constitute exactly a half of the cycle.
Pulling back by $\psi_1:S\to S_3$
in \eqref{ee}, the class $\psi^*\ms O(1,1) -
 e_{\aaa} - e_{\beta} -  e_{\ccc}$ on $S$
is represented by an effective
divisor $\Dh$ which satisfies $\Dh+\oDh=D$.
Hence the assertion of the proposition follows simply by putting 
\begin{align}\label{flew}
D^{\rm h} = \psi^*\ms O(1,1) -
 e_{1} - e_{2} -  e_{3}.
\end{align}
This is obviously independent of $\aaa$.
\proofend

\medskip
From the relation \eqref{half}, the effective divisor $\Dh$
can be thought of as a half of the divisor $D$,
and the upper-script `h' stands for `half'.
Proposition \ref{prop:cc02} says that the class of each $(-1)$-curve on $S$ lying
 over a real bitangent
can be represented by a half of $D$
with a slight adjustment by the classes
$e_{\aaa}$ and $\ol e_{\aaa}$.

The formula \eqref{flew}  provides a practical
way to write down
the half $D^{\rm h}$ of the divisor $D$ in a concrete form; 
it is just the pull-back of the half of the 
real anti-canonical cycle on $S_3$, under the birational morphism
$\psi_1:S\to  S_3$.

Since the divisor $D$ on $S$ includes 
the line components $C_1$ and $\ol C_1$ by 
multiplicity one,
and since $D^{\rm h}$ and $\ol D^{\rm h}$ are effective
 satisfying $D^{\rm h}+\ol D^{\rm h}=D$,
the multiplicity of $C_1$ in $D^{\rm h}$ and $\ol D^{\rm h}$ is
either 1 or 0.
Now we make a distinction between the line components
$C_1$ and $\ol C_1$ as well as irreducible components $S\upi_+$ and
$S\upi_-$ of a reducible fiber $S\upi$ by the following rule.

\begin{convention}\label{conv}
{\em
We let $C_1$ to be the line component which is a component
of the above effective divisor $D^{\rm h}$.
(Therefore $\ol C_1$ is not a component of $D^{\rm h}$.)
Next,  for each $i\in\{1,\dots,k\}$,  we let the component $S\upi_+$ 
to be the one which includes  $C_1$, so that 
the component $S\upi_-$ 
is the one which includes $\ol C_1$.
\proofend
}
\end{convention}

In terms of the blowing up $\psi:S\to S_0=\qdr$,
this is equivalent to saying that 
\begin{itemize}
\item
In the case of type $\Azero$,
the component $C_1$ is the strict transform of the unique
$(1,1)$-curve through the points $ p_1, p_2$
and $p_3$,
\item
In the case of type $\Aone$, the component 
$C_1$ is the strict transform of the exceptional curve
of the point $\ol p_3$,
\item
In the case of type $\Atwo$, the component $C_1$ is the strict transform of the 
exceptional curve which arises when we make a blowing up
in the direction of the tangent vector $\ol p_3$.
\end{itemize}

Under this convention, from the formula \eqref{flew},
we immediately obtain the following information about 
the half $\Dh$.

\begin{proposition}\label{prop:half2}
The half $D^{\rm h}$ of $D$ on $S$ in 
Proposition \ref{prop:cc02} satisfies the following.
\begin{itemize}
\item If $S$ is of type $\Azero$ (and $n>3$), the divisor $D^{\rm h}$
always includes
each of the two neighboring 
components of $\ol C_1$ by multiplicity one.
Further, the divisor $\Dh$ always includes 
two neighboring components of $C_1$.
\item
If $S$ is of type $\Aone$ or $\Atwo$,
the divisor $D^{\rm h}$ includes one of the neighboring 
components of $\ol C_1$ by multiplicity one
and does not include the other neighboring component.
\end{itemize}
\end{proposition}

In view of Proposition \ref{prop:cc02}, 
for each $\aaa$, we define two non-effective divisors 
$\Dha$ and $\bDha$ on $S$ by 
\begin{align}\label{Dha}
\Dha = D^{\rm h} + e_{\aaa} - \ol e_{\aaa}
\qandq
\oDha = \ol D^{\rm h} - e_{\aaa} + \ol e_{\aaa}.
\end{align}
From Proposition \ref{prop:cc01},
each of these classes is represented by
a $(-1)$-curve on $S$
lying over a real bitangent, and 
the range of $\aaa$ is as follows.
\begin{align}\notag
{\text{type $\Azero$}}&\Ra \aaa\in \{1,2,3\},\\
{\text{type $\Aone$}}&\Ra \aaa\in \{1,2\}, \label{alpha}\\
{\text{type $\Atwo$}}&\Ra \aaa=1.\notag
\end{align}
The above two classes $\Dha$ and $\oDha$ 
satisfy the following 
nice numerical properties.
\begin{proposition}
\label{prop:Aint1}
Under Convention \ref{conv} for distinction
between $C_1$ and $\ol C_1$,
the intersection numbers $D^{\rm h}_{\aaa}.\,C_i$
and $D^{\rm h}_{\aaa}.\,\ol C_i$ on $S$
are all zero except 
$D^{\rm h}_{\aaa}.\,\ol C_1$,
which is one.

%
\end{proposition}

For this reason, Convention \ref{conv} is 
significant in the rest of this section.

\proof
The class $\Dha$ is represented by a $(-1)$-curve which is one of the two components
of the inverse image of a real bitangent
of the branch quartic.
These $(-1)$-curves are not a component of the real anti-canonical
cycle $C$ since we are excluding the ridge from bitangents.
(Recall the cycle $C$ is mapped to the ridge of the scroll.)
This means that $\Dha.\,C_i$ and $\Dha.\,\ol C_i$ are 
non-negative for any $i$.
Hence if $i\neq 1$ we have $\Dha.\,C_i = \Dha.\,\ol C_i=0$
since $D.\, C_i = D.\,\ol C_i = 0$ and $D = \Dha + \oDha$.

The remaining intersection numbers $\Dha.\,C_1$
and $\Dha.\,\ol C_1$ can be calculated from
the presentation \eqref{cc01} for $D^{\rm h}_{\aaa}$
and $\ol D^{\rm h}_{\aaa}$ as follows.
For the case of type $\Azero$, as above,
the line component
$C_1$ 
is the strict transforms of the  $(1,1)$-curve
through the three points $ p_1, p_2$ and $ p_3$.
This means 
$$
C_1 \simeq \psi^*\ms O(1,1) -   e_1 -  e_2 -  e_3 -
\eee,
$$
where the class $\eee$ is a (positive) linear combination
of the exceptional classes $e_4,\dots,e_n$
and $\ol e_4,\dots, \ol e_n$.
(Recall that these are classes of the exceptional
divisors of $\psi_1:S\to S_3$ respectively.)
On the other hand, from \eqref{cc01},
we have
\begin{align}\label{pb101}
D^{\rm h}_{\aaa} \simeq
\psi^*\ms O(1,1) -\ol e_{\aaa} -  e_{\beta} -  e_{\ccc},
\end{align}
where $(\aaa,\bbb,\ccc) = 
(1,2,3),\,(2,3,1)$ or $(3,1,2)$.
From these we get 
$D^{\rm h}_{\aaa}.\, C_1 = 2-2=0$.
Hence, we obtain 
$D^{\rm h}_{\aaa}.\, \ol C_1 = 
D.\,\ol C_1 - \oDha.\,\ol C_1 = 1 - \Dha.\,C_1=
1-0=1$.
%

For the case of type $\Aone$,
we have
\begin{align}\label{c1}
C_1 \simeq \ol e_3- \eee',
\end{align}
where $\eee'$ is again a (positive) linear combination
of the classes $e_4,\dots,e_n$
and $\ol e_4,\dots, \ol e_n$.
So noting that $\aaa$ is either 1 or 2,
the presentation \eqref{pb101} gives 
$\Dha.\,C_1 = 0.$
Hence we again obtain $\Dha.\,\ol C_1=1$.

For the case of type $\Atwo$,
we still have \eqref{c1} for some other $\eee'$.
So noting that $\aaa = 1$,
we again obtain $\Dha.\,C_1 = 0$
and $\Dha.\,\ol C_1=1$.
\proofend

\medskip
By using the effective divisor $\Dh$ on $S$
in Proposition \ref{prop:cc02},
we define
an effective divisor $\bDh$ on $\tilde Z$
as the formal extension of the divisor $\Dh$ on the whole of $\tilde Z$;
namely by replacing the curves $C_i$ and $\ol C_i$
appearing in the curve $D^{\rm h}$  by the exceptional 
divisors $E_i$ and $\ol E_i$ respectively.
So this is  analogous to define the divisor $\bm D$
on $\tilde Z$
from the divisor $D$ on $S$.
From the definition, we have,
$$
\ms O_{\tilde Z}(\bm D^{\rm h})|_S\simeq \ms O_S(D^{\rm h})
$$
for any real irreducible fiber $S$ of $\tilde f$.
The divisors $\bDh$ and $\obDh$
are `halves' of the divisor $\bm D$.
Of course we have the relation.
\begin{align}\label{division02}
\bm D^{\rm h} + \ol{\bm D}^{\rm h} =\bm D.
\end{align}


In order to obtain a Chern class on $\tilde Z$ whose restriction to $S$
is exactly $\Dha$,
for the moment we fix any total small resolution $\zeta:\tilde Z\to \hat Z$.
(We will choose a particular one later.)
We recall from Proposition \ref{prop:PP} that 
a real irreducible member $S$ in $|F|$ always
contains a family of twistor lines
parameterized by $S^1$, and by the twistor projection
$\pi:Z\to n\CP^2$,
each connected component of the complement of the union of these twistor lines is diffeomorphic 
to $n\CP^2\minus S^1$.
Moreover, the birational morphism $\psi:S\to S_0=\qdr$ does not 
blow up points on the image of the twistor lines.
Hence, via the degree-two map 
$\pi|_S:S\to n\CP^2$, each of the classes $e_1,\dots, e_n$ on $S$
defines a cohomology class on $n\CP^2$.
For each $\aaa$ in the range of \eqref{alpha},
let $\xi_{\aaa}\in H^2(n\CP^2,\ZZ)$ be the class defined by
$e_{\aaa}$.
(We can define this for any $\aaa\in\{1,\dots,k\}$, but 
we only use the above ones.)
Then from the choice we have 
$$(\pi^*\xi_{\aaa})|_S = e_{\aaa} - \ol e_{\aaa}$$
in $H^2(S,\ZZ)$.
As before let $\eta:\tilde Z\to Z$ be the birational morphism
which is the composition of the blowing up $\hat Z\to Z$
and the total small resolution $\zeta:\tilde Z\to \hat Z$.
Then we define a cohomology class
(or equivalently, a line bundle) 
$\msDha$ on $\tilde Z$ by
\begin{align}\label{lb002}
\msDha:=\bm D^{\rm h} + \eta^*\pi^*\xi_{\aaa}.
\end{align}
(This should read
``a half of $\bm D$, adjusted by $\xi_{\aaa}$''.)
Then if $S$ is any real irreducible fiber of $\tif$,
since $\eta$ is isomorphic on $S$,
we have
\begin{align*}
\msDha|_S &\simeq 
\bDh|_S + (\pi^*\xi_{\aaa})|_S\\
&\simeq \Dh + e_{\aaa} - \ol e_{\aaa}\\
&\simeq \Dha.
\end{align*}
Thus we have obtained the desired cohomology classes on $\tilde Z$.
Also from Proposition \ref{prop:Aint1},
if $S$ is a smooth fiber of $\tilde f:\tilde Z\to\CP^1$
and $C_i$ (resp.\,$\ol C_i$) means
the intersection $E_i\cap S$ (resp.\,$\ol E_i\cap S$),
we have
\begin{gather}\label{in2}
\msDha.\,C_i=0 \quad{\text{for any $i$}},\\
\msDha.\,\ol C_1 = 1,
\quad \msDha.\,\ol C_i = 0\quad{\text{for any $i\neq 1$}}.\label{in3}
\end{gather}
From these we immediately obtain
\begin{align}\label{in4}
h^0\big(\msDha|_C\big) = 1,
\end{align}
where $C$ means the cycle $S\cap E$.
Note that the discrepancy between \eqref{in2} and \eqref{in3}
reflects non-reality of the line bundle $\msDha$.
These intersection numbers will be soon used a number of times for investigating base curves of
the line bundle $\msDha$.

In the next subsection we discuss
how to find a real reducible member of the linear system
$\big|\ms D(m)\big|$ by making use of 
the line bundle $\msDha$. 

\subsection{Stable base curves of the non-real line bundle $\msDha$}\label{ss:dimh}
As we have discussed so far, our goal
in this section is to show that 
the system $\big|\ms D (m)\big|$ on $\tilde Z$
has a real reducible  member in the case of type $\Azero,\Aone$ 
and $\Atwo$.
From the result in the last subsection,
for this goal, we should investigate the 
line bundle $\msDha$ on $\tilde Z$.
But the line bundle $\msDha$ itself does
not admit a non-zero section since
$\msDha + \omsDha = \ms D$
holds 
and the line bundle $\ms D$ does not admit
a non-zero section other than defining sections of the divisor $\bm D$.
So as in the case for the real line bundle $\ms D$,
we put
\begin{align}\label{msDhal}
\msDha(l):=\msDha  \otimes\tilde f^*\ms O(l),\quad
l\in\ZZ.
\end{align}
We begin with the following easy property.

\begin{proposition}\label{prop:tau}
The direct image sheaf $\tif_*\msDha$ 
over $\CP^1$ is invertible
and of negative degree.
\end{proposition}

\proof
Since $\tif$ is a proper flat morphism to a smooth curve,
the direct image is invertible if $h^0\big(\msDha|_S\big)=1$ holds
for a generic fiber $S$ of $\tif$.
If $S$ is a real irreducible fiber,
we have $\msDha|_S\simeq \Dha$ and
this is represented by a $(-1)$-curve from the construction in
the last subsection.
Hence we have 
$h^0\big(\msDha|_S\big)=1$ for a generic real fiber $S$ of $\tif$.
This property automatically holds for generic $S$ which is not necessarily real,
and therefore the direct image $\tif_*\msDha$ is an invertible 
sheaf over $\CP^1$.
So we can write $\tif_*\msDha\simeq\ms O(-\tau)$
for some $\tau\in\ZZ$.
Since $H^0(\msDha) = 0$ as above, we have $\tau>0$.
\proofend

\begin{definition}{\em
In the following we always use the letter $\tau$ to mean 
the positive integer which satisfies $\tif_*\msDha\simeq\ms O_{\CP^1}(-\tau)$,
and we write $\Theta_{\aaa}$ for the unique member
of the system $\big|\msDha(\tau)\big|$ on $\tiZ$.
\proofend
}
\label{def:theta}
\end{definition}

Although we do not know if the divisor $\Theta_{\aaa}$ is irreducible,
we have the following.
\begin{proposition}\label{prop:tht}
The divisor $\Theta_{\aaa}$ does not contain any component of the divisor $E$.
\end{proposition}

\proof
Suppose that the divisor $\Theta_{\aaa}$ contains some component of $E$.
Let $S$ be any real irreducible fiber of $\tif$.
Since $\Theta_{\aaa}|_S\simeq \Dha$,
from Proposition \ref{prop:Aint1}, we obtain $S\cap E\subset
\Theta_{\aaa}\cap S$.
But the class $\Dha$ is represented by a $(-1)$-curve
over a real bitangent and cannot contain the real anti-canonical cycle $C$ on $S$.
\proofend

\medskip
Since $\Theta_{\aaa}\simeq \msDha(\tau)$, we have
$$\Theta_{\aaa} + \ol\Theta_{\aaa}\in \big|\ms D(2\tau)\big|.$$
Hence if $2\tau=m$ holds, this is a real reducible member 
of $\big|\ms D(m)\big|$ whose intersection with 
a real irreducible fiber $S$ of $\tif$ is a pair of
$(-1)$-curves over the same real bitangent.
Therefore $\Theta_{\aaa} + \ol\Theta_{\aaa}$ provides
the desired real reducible member of $\big|\ms D(m)\big|$.
However we have $2\tau\neq m$ in general.
For instance, the number $m$ is not necessarily even.
On the other hand, since $\big|\ms D(l)\big| = \bm D + \big|\tilde f^*\ms O(l)\big|$
when $l<m$ by Proposition \ref{prop:l<m},
from the property
$\msDha|_S \simeq \Dh_{\aaa}$,
we have $\Theta_{\aaa} + \ol\Theta_{\aaa}
\not\in\big|\bm D(l)\big|$ if $l<m$.
Hence we always have an inequality $2\tau\ge m$.

In the rest of this section, we shall show that 
{\em if the strict inequality $2\tau>m$ holds,
the divisor  $\Theta_{\aaa}$ always has
some degree-one divisors $S\upi_+$ and/or $S\upi_-$ as
irreducible components},
and that 
{\em by removing all such components from $\Theta_{\aaa}$,
we get a component of a real reducible member of $\big|\ms D(m)\big|$
whose restriction to any real irreducible fiber of $\tif$ is a 
$(-1)$-curve over a real bitangent.}
The proof will be completed at the end of the next subsection.

Our proof of this assertion relies on Propositions
\ref{prop:fvr} and
\ref{prop:szero},
which were about special sections of the line bundle $\ms D(m)$,
and the structure of the space $H^0\big(\ms D(m+l)\big)$ for
arbitrary $l$ respectively.
But we need several steps to complete the proof.
In this subsection, for the proof,
we investigate base curves of the line bundle 
of the form $\msDha(l)$, which lie on the exceptional divisor $E$.
As in the case for the line bundle $\ms D$,
this is done by computing intersection numbers of
$\msDha$ with components of the cycles $C\upi= S\upi\cap E$, $1\le i\le k$.

%
%
%
%
%
%
For this purpose, we fix a small resolution
$\zeta\upi$ of the nodal space $\hat Z$
for any $1\le i\le k$.
This time from the start we choose the following small resolutions,
for which we know that the number $e$ becomes the largest:
\begin{itemize}
\item
For $\zeta\upone$ and $\zeta\upk$,
we choose the small resolutions which
blow up the ridge
components $E_1$ and $\ol E_1$.
\item
For the remaining resolutions $\zeta\upi$,
$i\neq 1,k$,
we choose the ones which blow up the components
$E_i$ and $\ol E_i$.
\end{itemize}
We recall that we have denoted $\DDD\upi$ and $\ol\DDD\upi$
for the exceptional curves over the points
$z\upi$ and $\ol z\upi$ respectively.

In order to express base curves of $\msDha(l)$ in a concrete form,
we write the effective divisor $\bDh$ as
\begin{align}\label{Dh}
\bm D^{\rm h} = \sum_{i=1}^{k} d'_i E_i +
\sum_{i=1}^{k} d''_i \ol E_i.
\end{align}
Because the divisor $\Dh$ on $S$ is effective,
we have $d'_i\ge 0$ and $d''_i\ge 0$ for any $i$.
We have $d'_i\neq d''_i$ in general
because $\bm D^{\rm h}$ is  not a real divisor.
Also, some of these coefficients can be zero.
Since $\bDh+\obDh = \bm D$ as in \eqref{division02},
for any $i\in\{1,\dots,k\}$, we have
\begin{align}\label{}
d'_i + d''_i = d_i.
\end{align}

Since the divisor  $\bm D$ was real,
the difference of the multiplicities between
adjacent components $C_i$ and $C_{i+1}$,
and that of the components
 $\ol C_i$ and $\ol C_{i+1}$ were equal.
This property does not carry over to
the divisors $\bDh$ and $\obDh$,
since these divisors are not real.
But we have the following remark, which will be used later.

\begin{proposition}
\label{prop:cnt0}
For any index $i\neq 1,k$,
the differences of adjacent
multiplicities $d'_{i+1} - d'_i$ 
and $d''_{i+1} - d''_i$ cannot have a different
sign i.e.\,$(d'_{i+1} - d'_i)(d''_{i+1} - d''_i)\ge 0$.
\end{proposition}

We note that in the case $i=1,k$, the proposition does not necessarily
hold, where the indices are given cyclically in the following sense: $d'_0:=d''_k,\, d'_{k+1}:=d''_1,\,
d''_0:=d'_k$ and $d''_{k+1}:=d'_1$.

\medskip
\noindent{\em Proof of Proposition \ref{prop:cnt0}.}
This is an easy consequence of the formula \eqref{flew}
for the effective divisor $\Dh$ and induction on $n$.
Indeed, if $n=3$ and $Z$ is of type $\Azero$ or $\Aone$,
the assertion is vacuous because $k=1,2$ holds.
If $n=3$ and $Z$ is of type $\Atwo$, we have $k=3$
and the assertion makes sense only for $i=2$.
In this case, among the six components of the cycle $C$,
exactly three are included as multiplicity one
component in the divisor $\Dh$, and they constitute 
a chain (of length three).
Further the line component $C_1$ is an end of this chain,
and the remaining three components have multiplicity zero.
From this  $d'_2 = d'_3$ and $d''_2 = d''_3$ 
always hold,
and this means the assertion.

Next let $n>3$ and fix any index $i\neq 1,k$.
If the component $C_i$ is among 
the conjugate pair of the exceptional curves
of the final blowing up in the composition $\psi:S\to S_0=\qdr$
(see \eqref{ee}),
then from the formula \eqref{flew} 
we have $d'_{i}= d'_{i-1}+d'_{i+1}$ 
and  $d''_{i}= d''_{i-1}+d''_{i+1}$.
These mean
$d'_i\ge d'_{i+1}$ and $d''_i\ge d''_{i+1}$
since $d'_{i-1}\ge 0$ and $d''_{i-1}\ge 0$.
Therefore the assertion holds.
The same argument applies if the component $C_{i+1}$
is among 
the conjugate pair of the exceptional curves
of the final blowing up in $\psi$.
If both of the components
$C_i$ and $C_{i+1}$ are not among
the conjugate pair of the exceptional curves
of the final blowing up in $\psi:S\to S_0$,
then the assertion holds from the induction hypothesis.
\proofend

\medskip

Some of the multiplicities $d'_i$
and $d''_i$ can be concretely determined.
All these will be used later.

\begin{proposition}
\label{prop:cnt1}
The coefficients $d'_i$ and $d''_i$ for the divisor 
$\bDh$ satisfy the following.
\begin{itemize}
\item
Regardless of the type, $d'_1 = 1$ and $d''_1 = 0.$
\item
If $Z$ is of type $\Azero$ (and $n>3$), we always have $d'_k=d''_2=1$,
$d'_2>0$ and $d''_k>0$.
\item
If $Z$ is of type $\Aone$ or $\Atwo$, either 
$(d'_k,d''_2) = (1,0)$ or $(d'_k,d''_2) = (0,1)$ holds.
\end{itemize}
\end{proposition}

\proof
The first assertion immediately follows from $d_1=1$
and the promise that
the divisor $\Dh$ contains $C_1$ as in Convention \ref{conv}.
The second and the third assertions are just  
paraphrases of Proposition \ref{prop:half2}.
\proofend

\medskip
The intersection numbers  $\bDh.\,\DDD\upi$ and $\bDh.\,\ol\DDD\upi$
are expressed in terms of the coefficients $d'_j$ and $d''_j$.
For calculating $\msDha.\DDD\upi$ and $\msDha.\,\ol\DDD\upi$,
we also need the intersection numbers
$\big(\eta^*\pi^*\xi_{\aaa},\DDD_{i}\big)$ 
and
$\big(\eta^*\pi^*\xi_{\aaa},\ol\DDD_{i}\big)$.
For these, since the blown up points $p_{\aaa}$
on $S_0=\qdr$ (where $\aaa$ is in the range \eqref{alpha})
is not an intersection point of the two $(1,1)$-curves
$C_0$ and $\ol C_0$
as remarked right after Proposition \ref{prop:cc01},
the class $\xi_{\aaa}\in H^2(n\CP^2,\ZZ)$ is represented by 
a smooth real surface (isomorphic to $S^2$) which does not pass the double points of 
the cycle $\pi(C)$.
This means that the class $\eta^*\pi^*\xi_{\aaa}$ 
is represented by a 4-cycle which is disjoint from the curves
$\DDD\upi$ and $\ol\DDD\upi$.
Hence  
we have
\begin{align}\label{int259}
\big(\eta^*\pi^*\xi_{\aaa},\DDD\upi\big)=
\big(\eta^*\pi^*\xi_{\aaa},\ol\DDD\upi\big)=0
\quad{\text{for any $\aaa$ and $i$}}.
\end{align}
This means that, 
we simply have
\begin{align}\label{int260}
\big(\msDha,\DDD\upi\big) = \big(\bDh,\DDD\upi\big)
\qandq
\big(\msDha,\ol\DDD\upi\big) = \big(\bDh,\ol\DDD\upi\big).
\end{align}
Therefore, under the present choice of the small
resolution $\tiZ\to\hat Z$ as above,
by using the first item in  Proposition \ref{prop:cnt1}, we have
\begin{gather}\label{Dha30}
\msDha.\,\DDD\upi = d'_{i+1} - d'_i, \quad 1\le i<k,
\quad
\msDha.\,\DDD\upk = d'_k,\\
\msDha.\,\ol\DDD\upi = d''_{i+1} - d''_i, \quad 1\le i<k,
\quad
\msDha.\,\ol\DDD\upk = d''_k-1.
\label{Dha31}
\end{gather}


For later purpose, we also compute the intersection numbers 
with the curves $\Xi_i=E_i\cap E_{i+1}$ and $\ol\Xi_i=\ol E_i\cap \ol E_{i+1}$.
By the same reason to \eqref{int259}, 
we have
\begin{align}\label{int261}
\big(\eta^*\pi^*\xi_{\aaa},\Xi_i\big)=
\big(\eta^*\pi^*\xi_{\aaa},\ol\Xi_i\big)=0
\quad{\text{for any $\aaa$ and $i$}},
\end{align}
and hence
\begin{align}\label{int262}
\big(\msDha,\Xi_i\big) = \big(\bDh,\Xi_i\big)
\qandq
\big(\msDha,\ol\Xi_i\big) = \big(\bDh,\ol\Xi_i\big)
\quad{\text{for any $\aaa$ and $i$}}.
\end{align}
By using the above calculations, we investigate
base curves of the line bundle of the form $\msDha(l)$,
$l\in\ZZ$.
As in the case for the line bundle
of the form $\ms D(l)$, 
these will be contained in the cycles
$C\upi = S\upi\cap E$, $1\le i\le k$,
and will be again found through intersection numbers
as we mentioned.
But there is  considerable 
and interesting
difference coming from the following two factors.
\begin{itemize}
\item
As in Proposition \ref{prop:Aint1},
we have
$\Dha.\,C_1 = 0$ and 
$\Dha.\,\ol C_1 = 1$.
The former implies that, 
unlike for the the line bundle $\ms D$,
base curves in the cycle $C\upi$ can pass across the ridge component $E_1$
if $i\neq 1,k$.
\item
The line bundle $\msDha$ is not real.
This means $\msDha.\,\DDD\upi\neq\msDha.\,\ol\DDD\upi$
in general. 
\end{itemize}

In order to express base curves concretely,
as before, for any $i,j\in\{1,\dots,k\}$,
let $C\upi_j$ and $\ol C\upi_j$ be the components
of the cycle $C\upi\subset \tiZ$
which are mapped isomorphically to the components
$C_j$ and $\ol C_j$ of the cycle $C$ respectively by
the birational morphism $\eta:\tilde Z\to Z$.

As in the case for the line bundle $\ms D$,
the case $i=1,k$ requires another treatment due to relevance with the ridge components, and 
we first investigate base curves of $\msDha(l)$
which are on the fiber $S\upi$ when $i\neq 1,k$.
Recall that we denote the coefficients of the divisors $E_i$ 
and $\ol E_i$ in $\bDh$ by  $d'_i$ and $d''_i$ respectively.
In the following argument,
finding  base curves of $\msDha$
requires some care that was not needed
in the case of the line bundle $\ms D$.

\begin{figure}
\includegraphics{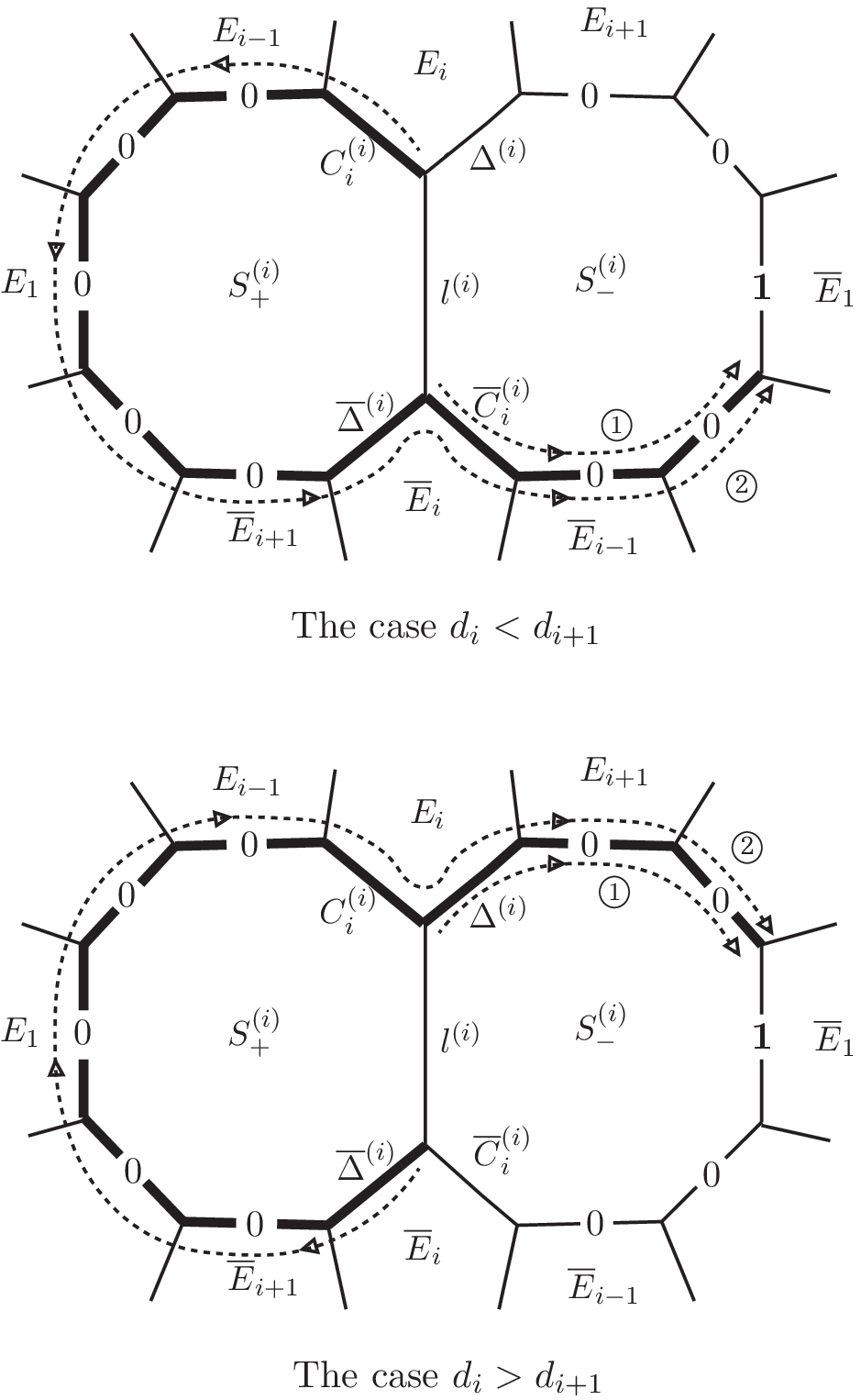}
\caption{stable base curves of $\msDha$ on $S\upi$}
\label{fig:rf8}
\end{figure}

Assume $i\neq 1,k$, and
suppose that $d'_i\le d'_{i+1}$ holds.
By Proposition \ref{prop:cnt0},
this is equivalent to $d''_i\le d''_{i+1}$,
or just $d_i\le d_{i+1}$.
Hence by \eqref{Dha30} and \eqref{Dha31},
we have
$\msDha.\,\DDD\upi\ge 0$ and
$\msDha.\,\ol\DDD\upi \ge 0$.
Since  $\msDha.\,C_i=0$ holds on a general fiber $S$ by Proposition \ref{prop:Aint1},
we have $\msDha.\,C\upi_i = - \msDha.\,\DDD\upi$
and $\msDha.\,\ol C\upi _i = - \msDha.\,\ol\DDD\upi$.
Therefore in the present assumption we obtain 
 $\msDha.\,C\upi_i \le 0$
and $\msDha.\,\ol C\upi_i\le 0$.
If the latter inequality is strict,
the chain
\begin{align}\label{sbcnr1}
\ol C\upi_i + \ol C\upi_{i-1} + \dots +
\ol C\upi_2
\end{align}
is a base curve of $\msDha(l)$ for any $l\in\ZZ$.
Moreover, by the same reason to the case of the line bundle $\ms D$,
the multiplicity on the divisor $E$ of this base curve equals $-\msDha.\,\ol C\upi_i = d''_{i+1}-d''_i$.
These base curves stop at the tail $\ol C\upi_2$ because
it hits the component $\ol E_1$ and we have
$\msDha.\,\ol C\upi_1 = 1$ by Proposition \ref{prop:Aint1}.
In the upper picture of Figure \ref{fig:rf8}, this chain is followed by 
the dotted curve with the letter \textcircled{\scriptsize 1}.
Let $L$ be the line bundle over $E$
which is obtained from the restriction
$\msDha|_{E}$ by subtracting the chain
\eqref{sbcnr1} with the above multiplicity.
(If $d'_i= d'_{i+1}$, there is no subtraction.)
Here, since the chain \eqref{sbcnr1} is not a Cartier divisor on $E$ at the point
$\ol C\upi_2\cap \ol C\upi_1$,
$L$ is not a line bundle in the usual sense.
But this is not a problem since we will not consider 
intersection number with the curve $\ol C\upi_1$.
Then as an effect of the subtraction, we have 
\begin{align}\label{int800}
L.\,\ol\DDD\upi = \msDha.\,\ol\DDD\upi - (d''_{i+1}-d''_i) =0.
\end{align}

Next, since the subtracted chain \eqref{sbcnr1} does not intersect the curve $C\upi_i$,
we have
\begin{align}\label{int594}
L.\,C\upi_i = 
\msDha.\,C\upi_i = d'_i - d'_{i+1}.
\end{align}
This is non-positive from the assumption.
If this is negative, noting  $\msDha.\,C\upi_1=0$,
the chain
$C\upi_+-\ol\DDD\upi$,
where $C\upi_+= S\upi_+\cap E$ as before,
is a base curve
of $\msDha(l)|_E$  with multiplicity $(d'_{i+1}-d'_i)$ for any $l\in\ZZ$.
{\em Moreover, from \eqref{int800}, this base curve does not stop at the tail  $\ol C\upi_{i+1}$,
and the adjacent chain 
\begin{align}\label{sbcnr2}
\ol\DDD\upi + \ol C\upi_i + \ol C\upi_{i-1} + \dots + \ol C\upi_2
\end{align}
is also included
by the same multiplicity $(d'_{i+1}-d'_i)$.}
Summing  these up,
the chain $C\upi_+$
is  fixed component of the line bundle $\msDha(l)|_E$  with multiplicity $(d'_{i+1}-d'_i)$ for any $l\in\ZZ$,
and the adjacent chain \eqref{sbcnr1} is fixed component
of the same line bundle with multiplicity
$$
(d''_{i+1} - d''_i) + (d'_{i+1}-d'_i)
= d_{i+1} - d_i.
$$
In the upper picture of Figure \ref{fig:rf8}, the longer chain of base curves
with multiplicity $(d'_{i+1} - d'_i)$ is followed by
the dotted curve with the letter \textcircled{\scriptsize 2}.
We write the sum of these base curves on $E$
(namely \textcircled{\scriptsize 1}$+$\textcircled{\scriptsize 2}
with the above individual multiplicities)
 by $\big(\Cha\big)\upi_E$,
 and  call it
 the {\em stable base curves} of $\msDha$.
By reality, the conjugate divisor $\big(\oCha\big)\upi_E$ is fixed component
of $\omsDha(l)\big|_E$ for any $l\in\ZZ$.
We mention that in this argument for obtaining
the base curves of $\msDha$ on $S\upi$, $i\neq 1,k$,
a trick is to first subtract the chain of base curves
which hits the component $\ol E_1$.

In a similar way, if the opposite inequality $d_i\ge d_{i+1}$ holds, 
the multiple chain 
$$
(d_i - d_{i+1})
\big(\DDD\upi + C\upi_{i+1} + \dots + C\upi_k\big)
$$
is fixed component of $\msDha(l)|_E$
for any $l$,
and the half chain $C\upi_+$
 is fixed component of $\msDha(l)|_E$ 
with multiplicity $(d''_i - d''_{i+1})$ for any $l$.
(Thus the situation is upside down to the previous case.
See the lower picture in Figure \ref{fig:rf8}.)
If we again write the sum of these chains by $\big(\Cha\big)\upi_E$,
then the conjugate curve $\big(\oCha\big)\upi_E$ is fixed component  of 
$\omsDha|_E$.

Next we compare the sum 
$\big(\Cha\big)\upi_E 
+
\big(\oCha\big)\upi_E $ with the stable base curve
$\CCC\upi_E$ of the line bundle $\ms D$.
Recall from Proposition \ref{prop:bc1} that if $d_i\le d_{i+1}$, 
we have
$$
\CCC\upi_E=(d_{i+1} - d_i)\Big\{
\big(C\upi_i + \dots + C\upi_2\big) +
\big(\ol C\upi_i + \dots + \ol C\upi_2
\big)\Big\}.
$$
From this
we obtain 
\begin{align}\notag
\big(\Cha\big)\upi_E 
+
\big(\oCha\big)\upi_E 
 &= \CCC\upi_E + (d'_{i+1} - d'_i)(C\upi_++C\upi_-)\\
 &=\CCC\upi_E + (d'_{i+1} - d'_i)C\upi.
 \label{krel1}
\end{align}
Thus, while the sum of the stable base curves of $\msDha$ and $\omsDha$
can be strictly greater than the stable base curves of $\ms D$,
the difference is a multiple of 
the cycle $C\upi$.
This observation will be significant later.

%

If the opposite inequality $d_i\ge d_{i+1}$ holds, 
in a similar way, we obtain the relation
\begin{align}\label{krel2}
\big(\Cha\big)\,\upi_E 
+ \big(\oCha\big)\upi_E 
= \CCC\upi_E + (d''_{i} - d''_{i+1})C\upi.
\end{align}
Thus we have obtained fixed components of the line bundles $\msDha(l)|_E$
and $\omsDha(l)\big|_E$, and also the comparison formulae
\eqref{krel1} and \eqref{krel2}, in the case $i\neq 1,k$.
We note that if the equality $d_i=d_{i+1}$ holds,
we have $d'_i=d'_{i+1}$ and $d''_i = d''_{i+1}$
by Proposition \ref{prop:cnt0}.
Hence the stable base curves of $\msDha$ and $\omsDha$ do not appear.
But stable base curves of $\ms D$ also do not appear.
In particular, the relation
$\ChaiE + \oChaiE = \CCC_E\upi$  holds trivially.

Next we show that similar relations hold on the remaining two 
reducible fibers $S\upone$ and $S\upk$.
We need to distinguish the case of type $\Azero$, 
and the cases of type $\Aone$ and $\Atwo$,
because of the differences for
possibility of the two numbers $d'_k$ and $d''_2$
appearing in Proposition \ref{prop:cnt1}.

First we consider the case of type $\Azero$.
For base curves on the fiber $S\upone$, we first have
$\msDha.\,\DDD\upone = d'_2 -1$ from \eqref{Dha30} as $d'_1=1$
and hence 
$\msDha.\,C\upone_1 = 1-d'_2$ from Proposition \ref{prop:Aint1}.
We note $d'_2>0$ holds in the case of type $\Azero$ as in Proposition \ref{prop:cnt1}.
Since $\msDha.\,C\upone_j = 0$ if $j>1$ by \eqref{in2},
if $d'_2>1$, the chain 
$C\upone_+ - \ol\DDD\upone$
is a base curve of $\msDha(l)$ for any $l\in\ZZ$.
(See Figure \ref{fig:rf8} after letting $i=1$.)
We have, from \eqref{int260},
\begin{align*}
\msDha.\,\ol\DDD\upone &= \big( d''_1 \ol E_1 + d''_2\ol E_2 \big).\,\ol\DDD\upone \\
&= - d''_1 + d''_2 
\end{align*}
and this is one since $d''_1 = 0$ and $d''_2=1$ by 
Proposition \ref{prop:cnt1}.
Therefore the above chain stops at the tail $\ol C\upone_2$.
The multiplicity of the above chain $C\upone_+ - \ol\DDD_1$,
when regarded as fixed component of the restriction $\msDha|_E$,
is $(d'_2-1)$.
Next, since $\msDha.\,\ol C\upone_1 = 1- \msDha.\,\ol\DDD\upone$
from Proposition \ref{prop:Aint1},
we have $\msDha.\,\ol C\upone_1 =0$.
Therefore no base curve appears on the other component $S\upone_-$.
From these, the curve 
\begin{align}\label{chain1}
\ChaoneE := (d'_2-1) \big(C\upone_+ - \ol\DDD\upone\big)
\end{align}
is fixed component of $\msDha(l)|_E$ for any $l\in\ZZ$.
If $d'_2=1$, this is the zero divisor.
By the real structure, 
the conjugate curve $\oChaoneE$ is fixed component
of $\omsDha(l)\big|_E$ for any $l\in\ZZ$.
%
On the other hand, by Proposition \ref{prop:bc2},
the stable base curve $\CCC\upone_E$ of the line bundle $\ms D$
was 
\begin{align}\label{chain3}
\CCC\upone_E = (d_2 - 2)\big(
C\upone - \DDD\upone - \ol\DDD\upone
\big).
\end{align}
Comparing this with the sum $\ChaoneE+\oChaoneE$,
we obtain the relation
\begin{align}\label{cmpr}
\ChaoneE + \oChaoneE &= \CCC_E\upone + (1-d''_2)C\upone\notag\\
&=\CCC_E\upone 
\end{align}
where the last equality is again from $d''_2 = 1$.
In the same way, 
 on another reducible fiber $S\upk$, 
 we obtain the stable base curve
$\ChakE$ of $\msDha|_E$ and $\oChakE$ of $\omsDha\big|_E$,
and they satisfy the relation
\begin{align}\label{cmpr1}
\ChakE + \oChakE = \CCC_E\upk. 
\end{align}

Next we consider the cases of type $\Aone$ and $\Atwo$.
In these cases, by Proposition \ref{prop:cnt1}, we always have 
$(d'_k,d''_2) = (0,1)$ or $(d'_k,d''_2)=(1,0)$.

First suppose that $(d'_k,d''_2) = (0,1)$ holds.
Then we readily obtain $\msDha.\,\ol\DDD\upone = 1$
and $\msDha.\ol C\upone_1 = 0$ from \eqref{int260},
\eqref{in3} and Proposition \ref{prop:cnt1}.
If $d'_2>0$,
from $\msDha.\,C\upone_1 = 1-d'_2$,
we obtain that the chain 
\eqref{chain1} is fixed component of $\msDha(l)|_E$ with 
multiplicity $(d'_2-1)$ for any $l\in\ZZ$.
Note that this chain stops at the tail $\ol C\upone_2$
since $\msDha.\,\ol\DDD\upone = 1>0$ as above for the adjacent component
$\ol\DDD\upone$.
No base curve appears on another component $S\upone_-$
since we have $\msDha.\,\ol C\upone_1 =0$ as above.
So we again define the stable base curve of $\msDha$,
 $\ChaoneE$, as
$(d'_2-1) \big(C\upone_+ - \ol\DDD\upone\big)$.
By the real structure,
the conjugate divisor $\oChaoneE$
is fixed component of $\omsDha(l)|_E$ for any $l\in\ZZ$.
On the other hand, as in Proposition \ref{prop:bc2},
the stable base curve $\CCC\upone_E$ of the line bundle $\ms D$
was the same as \eqref{chain3}.
From these, using $d''_2=1$, we obtain the same relation 
as \eqref{cmpr}.
%
%

Next, if ($(d'_k,d''_2)=(0,1)$ and) $d'_2=0$, we readily see that
this time just the half chain
$C\upone_- = S\upone_-\cap E$
is a base curve of $\msDha(l)|_E$ for any $l\in\ZZ$,
and its multiplicity is one.
So we put $\ChaoneE = C\upone_-$.
On the other hand,
by Proposition \ref{prop:bc2}
we have $\CCC\upone_E = 0$ in this situation
since $d_2 =d'_2+d''_2=0+1=1$.
From these we obtain the relation
\begin{align}\label{cmpr2}
\ChaoneE + \oChaoneE = \CCC_E\upone + C\upone.
\end{align}

Next suppose alternatively that $(d'_k,d''_2)=(1,0)$.
This implies
$\msDha.\,\ol\DDD\upone = 0$
and $\msDha.\ol C\upone_1 = 1$ from \eqref{int260},
\eqref{in3} and 
Proposition \ref{prop:cnt1}.
Since $d_2>0$ and $d_2 = d'_2 + d''_2$,
we obtain $d'_2>0$.
We still have $\msDha.\,C\upone_1 = -\msDha.\,\DDD\upone = 1-d'_2$.
 If $d'_2 = 1$, this is zero,
 and no base curve appears for the line bundle
$\msDha(l)|_E$.
This is the same  for $\ms D|_E$.
Hence we trivially have the relation
\begin{align}\label{cmpr3}
\ChaoneE + \oChaoneE = \CCC_E\upone. 
\end{align}
If $d'_2>1$, the chain $C\upone_+ - \ol\DDD\upone$ is fixed component of 
$\msDha(l)|_E$ with multiplicity $(d'_2-1)$ for any $l\in\ZZ$.
Further, as $\msDha.\,\ol\DDD\upone=0$ as above,
$\ol\DDD\upone$ is also a base curve.
But the multiplicity on $E$ of this component is only one.
Furthermore, no base curve appears on the other component $S\upone_-$.
So we put
$$
\ChaoneE = (d'_2 -1) \big( C\upone_+ - \ol\DDD\upone\big) + \ol\DDD\upone.
$$
On the other hand for the stable base curve of $\ms D$,
$\CCC\upone_E$ is as in \eqref{chain3}.
From these we obtain the relation
\begin{align}
\ChaoneE + \oChaoneE 
&= \CCC_E\upone + C\upone.\label{cmpr4}
\end{align}
In the same way, on the reducible fiber $S\upk$,
we also have the relation
\begin{align}\label{cmpr5}
\ChakE + \oChakE = \CCC_E\upk + n_kC\upk,
\end{align}
where $n_k$ is either 1 or 0 as in the case for $S\upone$.

Collecting all these base curves, in the same way to the line bundle $\ms D$,
we prepare the following notation.
\begin{definition}
\label{def:sbch}
{\em
We put
$$
\ChaE=\sum_{i=1}^k\big(\Cha\big)\upi_E
\qandq
\oChaE=\sum_{i=1}^k\big(\oCha\big)\upi_E,
$$
and call these
the {\em stable base curves} of the line bundles
$\msDha$ and $\omsDha$ respectively.
\proofend}
\end{definition}

Then from \eqref{krel1}, \eqref{krel2}
and \eqref{cmpr}--\eqref{cmpr5}
which compare the stable base curves
of $\msDha,\omsDha$ and $\ms D$,
we have

\begin{proposition}\label{prop:comparison1}
As before let $\CCC_E$ and $\ChaE$ be the stable base curves
of $\ms D$ and $\msDha$ respectively, with the multiplicities
on the divisor $E$ being taken into account as before.
Then we have the relation
\begin{align}\label{compa1}
\ChaE + \oChaE = \CCC_E + \sum_{i=1}^k n_i C\upi
\end{align}
for some non-negative integers $n_1,\dots, n_k$.
\end{proposition}

Of course we know the numbers $n_i$ concretely
from the above calculations,
but in the following only the relationship
\eqref{compa1} will be needed.

Though all components of $\ChaE$
are on $E$,
by the same reason to the stable base curve $\CCC_E$ of $\ms D$,
these base curves are not Cartier divisors on $E$.
But the part $\Ch_{\aaa,\,E-\ol E_1}:=
\ChaE|_{E-\ol E_1}
$ 
(resp.\,$\oCh_{\aaa,\,E-E_1}:=\oChaE|_{E-E_1}$) is a Cartier divisor on $E-\ol E_1$
(resp.\,$E- E_1$).
Moreover, from the above argument, the intersection numbers
of the line bundle $\msDha|_{E-\ol E_1}- \Ch_{\aaa,\,E-\ol E_1}$
 with 
any component of the chain $C\upi - \ol C\upi_1$
is zero for any index $i$.
This implies that there exists an integer $\tau'$
which satisfies
\begin{align}\label{plb1}
\msDha|_{E-\ol E_1}- \Ch_{\aaa,\,E-\ol E_1}
\simeq
\tif^*\ms O_{\CP^1}(-\tau').
\end{align}
This can be rewritten as
\begin{align}\label{plb2}
\msDha(\tau')|_{E-\ol E_1}\simeq \Ch_{\aaa,\,E-\ol E_1}.
\end{align}
By the real structure, this means
\begin{align}\label{plb3}
\omsDha(\tau')|_{E- E_1}\simeq \oCh_{\aaa,\,E- E_1}.
\end{align}
From these, as $\msDha + \omsDha \simeq\ms D$,  we obtain 
the relation
\begin{align}\label{plb4}
\ms D(2\tau')|_{E-E_1- \ol E_1}\simeq 
\Ch_{\aaa,\,E-E_1-\ol E_1} + \oCh_{\aaa,\,E-E_1-\ol E_1}.
\end{align}
From the property of the stable base curves of $\ms D$,
one might think that 
the relation \eqref{plb4} would imply a relation
\begin{align}\label{compa2}
\Ch_{\aaa,\,E-E_1-\ol E_1} + \oCh_{\aaa,\,E-E_1-\ol E_1} = \CCC_E + \sum_{i=1}^k n_i C\upi
\,\big|_{E-E_1-\ol E_1} 
\end{align}
for some non-negative integers $n_1,\dots, n_k$,
which would make the above calculations for the explicit forms
of the stable base curves of $\msDha$ unnecessary.
However, obtaining the inclusions 
$\CCC_E\upone\le \ChaoneE$ 
and
$\CCC_E\upk\le \ChakE$ 
requires concrete forms of $\ChaoneE$ and $\ChakE$
as we did as above.

\medskip
From \eqref{plb2}, it is easy to obtain the following result,
which is an analogue of Proposition \ref{prop:e2} for the
concrete expression for the
direct image $\tif_*\ms D|_E$.

\begin{proposition}\label{prop:l3}
Let $\tau'$ be as in \eqref{plb2}.
Then we have 
$$
\tilde f_* \big(\msDha|_E\big)\simeq\ms O_{\CP^1}(-\tau').
$$
\end{proposition}

\proof
If $S$ is a real irreducible fiber of $\tif:\tiZ\to\CP^1$,
we have $h^0\big(\msDha|_{E\cap S}\big) = 1$ by \eqref{in4}.
This implies that the direct image 
$\tif_*\msDha|_E$ is an invertible sheaf.
Hence it is enough to show 
\begin{align}\label{iso46}
H^0\big(\msDha(\tau')|_E\big)\simeq \CC.
\end{align}
Since the part $\Ch_{\aaa,\,E-\ol E_1}$ is fixed component of $\msDha(l)|_{E-\ol E_1}$
for any $l$, 
from \eqref{plb2}, we obtain
$$H^0\big(\msDha(\tau')|_{E-\ol E_1}\big)\simeq \CC.$$
Therefore in order to show \eqref{iso46} 
it suffices to show that 
the line bundle $\msDha|_{\ol E_1}$ satisfies the 
following property:
any section defined over the curves $\Xi_k\cup\ol \Xi_1\subset \ol E_1$ extends
in a unique way to the whole of $\ol E_1$.

As before let $(1,0)$ be the class on $\ol E_1$
represented by a fiber of the projection to $\ol C_1$,
and $(0,1)$ the class represented by a fiber of 
$\tif|_{\ol E_1}$.
Then $H^2(\ol E_1,\ZZ)$ is generated by $(1,0)$,
$(0,1)$, $\DDD\upk$ and $\ol\DDD\upone$.
In this notation, from $d''_1=0$, \eqref{Dha30}, \eqref{Dha31} and Proposition \ref{prop:Aint1},
we have
\begin{align}\label{dE1}
\msDha(\tau')|_{\ol E_1}\simeq
(1,\tau') - d'_k\DDD\upk - d''_2\ol\DDD\upone.
\end{align} 
On the other hand, 
we have $\Xi_k = (1,0) - \DDD\upk$
and $\ol\Xi_1 = (1,0)  - \ol\DDD\upone$.
Subtracting these from \eqref{dE1},
the kernel sheaf of the restriction homomorphism
from $\msDha(\tau')|_{\ol E_1}$ to the curves $\Xi_k\cup\ol\Xi_1$
is isomorphic to  $(-1,\tau') + (1-d'_k)\DDD\upk +	(1-d''_2)\ol\DDD\upone$.
By using that
$d_k'$ and $d_2''$ are either 1 or 0
from Proposition \ref{prop:cnt1},
it is easy to see that all cohomology groups of this class vanish.
Hence we obtain the unique extension property.
\proofend

\medskip
We note that the proof is a bit easier than Proposition \ref{prop:e2},
in that this time we do not need to subtract the stable base curves
lying on the component $\ol E_1$ when we prove the unique extension property.

The following property will soon be used in the next subsection.
Recall that any irreducible fiber $S$ of $\tif$ is smooth on 
the cycle $C=S\cap E$ by Proposition \ref{prop:sm1}.

\begin{proposition}\label{prop:restiso}
If $S$ is an irreducible fiber of $\tif$ which is not necessarily real
and $C=S\cap E$,
we have
\begin{align}\label{in5}
H^0\big(\msDha|_S\big)\simeq\CC
\qandq
H^0\big(\msDha|_C\big)\simeq\CC.
\end{align}
Further,
the restriction homomorphism
\begin{align}\label{res4}
H^0\big(\msDha|_S\big)
\lras
H^0\big(\msDha|_C\big)
\end{align}
is isomorphic.
\end{proposition}

\proof
Again the latter of \eqref{in5} was already shown in \eqref{in4}.
If $S$ is a real irreducible fiber, the class $\msDha|_S$
is represented by a $(-1)$-curve which is over a real bitangent.
In particular we have $H^0\big(\msDha|_S\big)\neq 0$
for such an $S$.
By upper semi-continuity, this means $H^0\big(\msDha|_S\big)\neq 0$ for any fiber $S$ of $\tif$.
Therefore, in order to show the former in \eqref{in5} and 
the isomorphicity of  the homomorphism \eqref{res4}, it is enough to 
prove that \eqref{res4} is injective
if $S$ is an irreducible fiber of $\tif$.
We note that from the real structure we also have
$H^0\big(\omsDha|_S\big)\neq 0$ for any fiber $S$ of $\tif$.

Let $S$ be an irreducible fiber of $\tif$ which is not necessarily real.
Suppose that there exists a non-zero element $t\in H^0\big(\msDha|_S\big)$  which satisfies $t|_C=0$.
Write $(t) = C + A$, where $A$ is an effective curve
(which might be zero) on $S$.
Let $t'$ be any non-zero element of $H^0\big(\omsDha|_S\big)$,
and write $(t')= A'$ by an effective curve $A'\neq 0$.
Since $\msDha + \omsDha \simeq\ms D$,
we have $(C+A)+A'\simeq \ms D|_S$.
Hence $A + A'\simeq \ms D|_S - C$.
By Proposition \ref{prop:van8}, the linear system 
$\big|\ms D|_S - C\big|$ on $S$ consists of a single member
$D-C$, where $D = \bm D|_S$.
Therefore $A+A' = D-C$.
Hence the curve $A'$ is a sub-divisor of $D-C$ and 
does not include the component $ C_1$.
On the other hand, since $A'\simeq \omsDha|_S$,
from Proposition \ref{prop:Aint1},
we have
\begin{align}\label{int823}
A'.\, \ol C_i = 0 \quad{\text{for any $i$,}}\qandq
A'.\, C_i = 
\begin{cases} 
1 & i=1,\\
0 & i\neq 1.
\end{cases}
\end{align}
(Note that the roles of $C_1$ and $\ol C_1$ are exchanged
since we are considering the conjugate bundle.)
In particular, the curve $A'$ intersects $C_1$.
But since $A'$ does not include $C_1$ and $\Supp A'\subset (C- C_1)$,
at least one of $ C_2$ and $\ol C_k$ has to be a component
of $A'$.
From the  intersection numbers \eqref{int823},
this readily means that $A'$ includes whole of the cycle $C$.
This  contradicts $\Supp A'\subset (C- C_1)$.
Hence $t$ has to be the zero section.
Therefore the homomorphism \eqref{res4} is injective.
\proofend

\subsection{Proof of the existence of real reducible members of $|mF|$}
\label{ss:redpf}
In this subsection, by using 
the results in the previous subsection,
we show the following result.
Recall that the class $\Dha$ (defined in \eqref{Dha}) on a real irreducible fiber $S$ of $\tif:\tiZ\to\CP^1$
is represented by a $(-1)$-curve 
over a real bitangent.

\begin{proposition}\label{prop:exred}
As before let $\aaa$ be as in the ranges \eqref{alpha} 
according to the type.
Then the linear system $\big|\ms D(m)\big|$ on the space $\tilde Z$
has a real reducible member of the form $\tilde T_{\aaa} + \sigma(\tilde  T_{\aaa})$,
which satisfies $\tilde T_{\aaa}|_S \simeq \Dh_{\aaa}$ for any real irreducible 
fiber $S$ of the morphism $\tilde f:\tilde Z\to\CP^1$.
\end{proposition}

Let $T_{\aaa}$ be the image of $\tilde T_{\aaa}$ to the twistor space under the birational morphism $\eta:\tilde Z\to Z$.
Then from the proposition we obtain the following expected result.
 
\begin{corollary}
\label{cor:rred}
Let $\aaa$ be as in Proposition \ref{prop:exred}.
Then the linear system $|mF|$ on the twistor space $Z$ admits
a real reducible member of the form
 $T_{\aaa} + \ol T_{\aaa}$
which satisfies $ T_{\aaa}|_S \simeq \Dh_{\aaa}$ for any real irreducible 
fundamental divisor $S$ on $Z$.\end{corollary}

In the rest of this subsection we prove Proposition \ref{prop:exred}.
To this end  we first consider the short exact sequence
\begin{align}\label{ses:45}
0\lras 
\msDha( - E)
\lras 
\msDha  
\lras
\msDha|_E
\lras 0.
\end{align}

\begin{proposition}\label{prop:van3}
We have
$\tilde f_*\msDha ( - E)
=0$.
\end{proposition}

\proof
It is enough to show 
$H^0\big(\msDha( - E)|_S\big)=0$ for 
a generic real fiber $S$ of $\tif$.
As above, for such an $S$, the class $\msDha|_S=\Dha$ is 
represented by a $(-1)$-curve 
over a  real bitangent and the $(-1)$-curve
cannot contain the real anti-canonical cycle $C$.
Hence we have $H^0\big(\msDha( - E)|_S\big)=
H^0(\Dha - C) =0$.
\proofend

\medskip
From this proposition,
by taking the direct image of \eqref{ses:45}, we obtain a long exact sequence
\begin{align}\label{ses:46}
0\lras 
\tilde f_*\msDha   
\lras
\tilde f_*(\msDha |_E )
\lras 
R^1\tilde f_*\msDha( - E)
\lras \dots
\end{align}
Recall that $\tif_*\msDha$ is also invertible by Proposition \ref{prop:tau}
and we wrote it as $\ms O(-\tau)$
where $\tau>0$ (Definition \ref{def:theta}).
Further we have $\tilde f_*(\msDha |_E )\simeq\ms O(-\tau')$ by Proposition \ref{prop:l3},
where  $\tau'$ is the number given in \eqref{plb1} or \eqref{plb2}.
Hence  from the exact sequence \eqref{ses:46} we have an inequality $\tau\ge\tau'$.
%
%


As before let $s_{\sst{\bm D}}$ be a section
of the line bundle $\ms D$ 
which satisfies $(s_{\sst{\bm D}})=\bm D$,
and $s_{\bm q}$ a non-real section of the line bundle $\ms D(m)$
as in Proposition \ref{prop:szero},
which corresponds to a non-real hyperplane in $\CP^{m+2}$ 
that passes the point $\bm q$ on the ridge $\bm l$.
The latter satisfies the vanishing property
\begin{align}\label{sqr0}
s_{\bm q}|_{ E_2 \cup \dots \cup  E_k} = 0.
\end{align}
Then by Proposition \ref{prop:fvr}, we have
\begin{align}\label{basic7}
H^0\big(\ms D(2\tau)\big)
=\tif^* H^0
\big(\ms O_{\CP^1}(2\tau)\big)\otimes s_{\sst{\bm D}} \, \bigoplus \,
\tif^*H^0\big(
\ms O_{\CP^1}(2\tau-m)\big)\otimes\langle s_{\bm q},\ol s_{\bm q}\rangle.
\end{align}
Now let $\ta$  be a generator
of the 1-dimensional vector space $ H^0\big(\msDha(\tau)\big)$.
We have $(\theta_{\aaa}) =\Theta_{\aaa}$,
where $\Theta_{\aaa}$ is the unique element
of the system $\big|\msDha(\tau)\big|$ as in 
Definition \ref{def:theta}.
Obviously we have $\ta\ota \in H^0\big(\ms D(2\tau)\big)$.
Hence from \eqref{basic7}, omitting the pull-back notation $\tif^*$
for simplicity, we can write
\begin{align}\label{red61}
\ta\ota = g s_{\sst{\bm D}} + hs_{\bm q} + \ol {h s_{\bm q}},
\end{align}
for some $g\in H^0\big(\ms O_{\CP^1}(2\tau)\big)$ and 
$h \in H^0\big(\ms O_{\CP^1}(2\tau-m)\big)$,
where $g$ is real.
In the sequel this relation plays a principal role.

We first determine the latter polynomial $h$ in a concrete form.
For this, we show the following proposition.
Recall that from Proposition \ref{prop:l3} 
we have $H^0\big(\msDha(\tau')|_E\big)\simeq\CC$.
For each index $i$,
we write $u\upi\in H^0\big(\tif^*\ms O_{\CP^1}(1)\big)$
for a real element which satisfies $\big(u\upi\big) = S\upi$.

\begin{proposition}\label{prop:comparison2}
Let $t_{\aaa}$ be a generator of $H^0\big(\msDha(\tau')|_E\big)\simeq\CC$.
Then on the divisor $E-\ol E_1$, we have the relation
\begin{align}\label{tae1}
\ta|_{E-\ol E_1}
= 
t_{\aaa}\prod_{i=1}^k \big(u\upi\big)^{n'_i}\,\Big|_{E-\ol E_1}
\end{align}
for some integers $n'_1,\dots,n'_k\ge 0$.
\end{proposition}

\proof
For simplicity of notation we write
$\theta'_{\aaa}:= \ta|_E$ in this proof.
By Proposition \ref{prop:tht}, this does not vanish
identically on any component of $E$.
By  \eqref{plb2},
over the divisor $E-\ol E_1$,
the section $t_{\aaa}$  defines the 
stable base curve $\ChaE$ of $\msDha$.
This means that on $E-\ol E_1$
an inequality $(t_{\aaa})|_{E-\ol E_1}\le (\theta'_{\aaa})|_{E-\ol E_1}$
as effective divisors holds.
Therefore, on $E-\ol E_1$,
the quotient $\theta'_{\aaa}/t_{\aaa}$ 
is a holomorphic section of the line bundle
$$
\msDha(\tau) - \msDha(\tau') \,\big| _{E-\ol E_1}.
$$
This is  isomorphic to $\tif^*\ms O_{\CP^1}(\tau-\tau')|_{E-\ol E_1}$.
Hence the support of the divisor
$(\theta'_{\aaa}/t_{\aaa})$ 
consists of fibers of $\tif |_{E-\ol E_1}:
E-\ol E_1\to\CP^1$.
Clearly such a fiber is contained in 
the zeroes of $\ta'=\ta|_{E}$.
Hence in order to show the assertion, it is enough to 
show that the restriction $\ta|_{E-\ol E_1}$ 
never vanishes identically on any {\em irreducible} fiber
of $\tif:\tiZ\to\CP^1$.

For this, 
let $S$ be any irreducible fiber of $\tif$, and put $C=S\cap E$ and $\ol C_1= S\cap \ol E_1$.
Suppose $\ta|_{C-\ol C_1}= 0$.
Then since $\msDha.\,\ol C_1=1$ as in \eqref{in3},
we have $\ta|_{\ol C_1}=0$ and hence $\ta|_C=0$.
By Proposition \ref{prop:restiso}, this implies $\ta|_S=0$.
Therefore dividing $\ta$ by a defining section 
of the fiber $S$, we obtain  a non-zero element 
of $H^0\big(\msDha(\tau-1)\big)$.
But by taking the direct image under $\tif$, 
recalling that $\tif_*\msDha\simeq\ms O(-\tau)$,
we have 
$$H^0\big(\msDha(\tau-1)\big)
\simeq H^0\big(\ms O_{\CP^1}(-\tau)\otimes
\ms O_{\CP^1}(\tau-1)\big)=
H^0\big(\ms O_{\CP^1}(-1)\big)
=0.
$$
This implies $\ta=0$, which contradicts
the choice of $\ta$.
\proofend


\medskip
The polynomial $h$ in the 
relation \eqref{red61} can be concretely written down
in terms of the integers $n_i$ and $n'_i$ as follows.

\begin{proposition}\label{prop:h}
Let $n_i$ and $n'_i$ be the non-negative
integers in Propositions \ref{prop:comparison1} and \ref{prop:comparison2} respectively.
Let $\mu\upi$ be the multiplicity of the cyclic base curves
as before.
Then for the polynomial $h$ in \eqref{red61},
we have 
\begin{align}\label{divide1}
h=
\prod_{i=1}^k \left(
u\upi\right)^{n_i+2n'_i-\mu\upi}.
\end{align}
Here we are regarding the defining section $u\upi$ of the reducible 
fiber $S\upi$ as a homogeneous polynomial of degree one
through the morphism $\tif:\tiZ\to\CP^1$.
\end{proposition}

Since the reducible fibers $S\upi$ are real,
in terms of homogeneous coordinates $(z_0,z_1)$
on $\CP^1$ with respect to which the real structure is 
given by $(z_0,z_1)\mapsto (\ol z_0,\ol z_1)$,
the conclusion \eqref{divide1} can be written more concretely as
\begin{align}\label{}
h(z_0,z_1) = \prod_{i=1}^k(z_0 - a_iz_1)^{n_i+2n'_i-\mu\upi}
\end{align}
for some distinct real numbers $a_1,\dots,a_k$.

\medskip
\noindent {\em Proof of Proposition \ref{prop:h}.}
We restrict the relation \eqref{red61} to 
the divisors $\ol E_2\cup \ol E_3\cup \dots \cup \ol E_k$.
Then since not only $s_{\sst{\bm D}}$ but also $\ol s_{\bm q}$ vanish identically there from \eqref{sqr0}, we obtain
\begin{align}\label{red62}
\ta\ota|_{\ol E_2\cup\dots \cup \ol E_k} = h s_{\bm q}|_{\ol E_2\cup\dots \cup\ol  E_k}.
\end{align}
As before let $\CCC_E$ be the stable base curve of $\ms D$.
Then for the right-hand side,
by  \eqref{szr0},
the divisor $\big(s_{\bm q}|_{\ol E_2\cup\dots \cup \ol E_k}\big)$ is exactly
the restriction of the sum
$$
\CCC_E + \sum_{i=1}^k \mu\upi C\upi
$$
to $\ol E_2\cup\dots \cup \ol E_k$.
On the other hand, from \eqref{tae1}, we have
\begin{align}\label{tae2}
\Big(\ta|_{E-\ol E_1}\Big)
= 
\ChaE + \sum_{i=1}^k n'_i C\upi\,\,\Big|_{E-\ol E_1}.
\end{align}
From the real structure, recalling that any $S\upi$ is real,
this means 
\begin{align}\label{tae3}
\Big(\ota|_{E- E_1}\Big)
= 
\oChaE + \sum_{i=1}^k n'_i C\upi\,\,\Big|_{E-E_1}.
\end{align}
%
%
%
%
From \eqref{tae2} and \eqref{tae3} we obtain
\begin{align}\label{tae4}
\Big(\ta\ota|_{\ol E_2\cup\dots \cup \ol E_k}\Big) 
= 
\ChaE + \oChaE + 2\sum_{i=1}^k n'_i C\upi
\,\,\Big|_{\ol E_2\cup\dots\cup \ol E_k}
\end{align}
Therefore, from \eqref{red62}, we obtain
\begin{align*}
\big(h|_{\ol E_2\cup\dots \cup \ol E_k}\big) &= 
\Big(
\ta\ota\big|_{\ol E_2\cup\dots \cup \ol E_k}
\Big)
- \big(s_{\bm q}|_{\ol E_2\cup\dots \cup \ol E_k}\big)\\
&=
\Big(
\ChaE + \oChaE + 2\sum_{i=1}^k n'_i C\upi
\Big)\,\Big|_{\ol E_2\cup\dots \cup \ol E_k}
-
\Big(\CCC_E + \sum_{i=1}^k \mu\upi C\upi
\Big)
\,\Big|_{\ol E_2\cup\dots \cup \ol E_k}\\
&= \sum_{i=1}^k (n_i+ 2n'_i-\mu\upi) C\upi
\,\Big|_{\ol E_2\cup\dots \cup \ol E_k},
\end{align*}
where the last equality is from 
Proposition \ref{prop:comparison1}.
Since $\big(u\upi\big) = S\upi$
and $S\upi\cap E=C\upi$, this implies  \eqref{divide1}.
\proofend

\medskip
Our next aim is to show that the polynomial $g$ in \eqref{red61}
is divisible by the polynomial $h$. For this purpose,  we need the following lemma.

\begin{lemma}\label{lemma:lb1}
For any index $i\in\{1,\dots,k\}$ and any non-negative integers $m_+$
and $m_-$, any section $t$ of the line bundle
\begin{align}\label{lb1}
\msDha - m_+ S\upi_+ - m_- S\upi_-\Big|_{S\upi}
\end{align}
over the reducible fiber $S\upi$
satisfies at least one of the following two properties.
\begin{itemize}
\item
the section $t$ vanishes at some point which does not belong
to the cycle $C\upi$,
\item
the section $t$ vanishes identically on the cycle $C\upi$.
\end{itemize}
\end{lemma}

\proof
In this proof for simplicity we write $L\upi$ for the line bundle  \eqref{lb1}
over $S\upi$.
As before,
let $l\upi=S\upi_+\cap S\upi_-$.
We have
\begin{align*}
\big(S\upi_+,l\upi\big)_{\tiZ} = 
\big(S\upi_+|_{S\upi_-},
l\upi\big)_{S\upi_-} = \big(l\upi,l\upi\big)_{S\upi_-} = 0.
\end{align*}
Hence $\big(S\upi_-,l\upi\big)_{\tiZ}=0$ by the real structure.
These imply
$
L\upi.\,l\upi = \msDha.\,l\upi
$
for any integers $m_+$ and $m_-$.
Further we have $\msDha.\,l\upi. = \bDh.\,l\upi$,
which is always positive because
the effective divisor $\bDh$ contains at least one
component of $E$ which intersects $l\upi$.
Therefore any section $t$ of the line bundle $L\upi$
vanishes at some point of $l\upi$.
If this point is not among the two points $l\upi\cap E$,
then $t$ satisfies the first item in the lemma.

\begin{figure}
\includegraphics{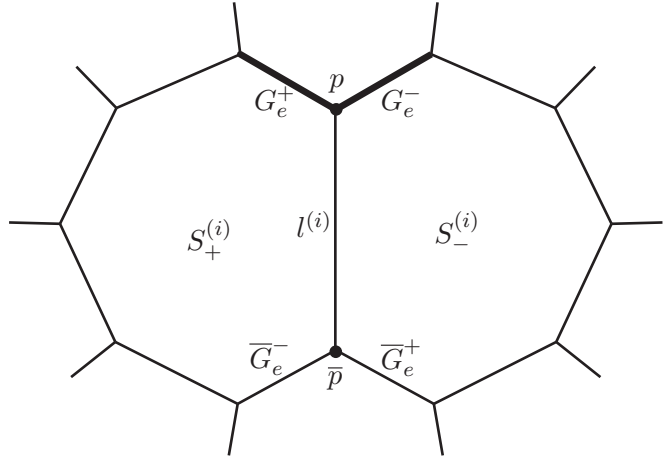}
\caption{Notations in the proof of Lemma \ref{lemma:lb1}.}
\label{fig:rf9}
\end{figure}

So suppose that 
the zero locus of $t|_{l\upi}$ is included in the two points $l\upi\cap E$.
In the following we write $p$ and $\ol p$ for these two points
(see Figure \ref{fig:rf9}),
and suppose that $t(p) = 0$.
We write $t_{\pm}:=t|_{S\upi\dpm}$ respectively.
Then since the zero locus of $t_+$ is not isolated on $S\upi_+$,
the section  $t_+$ vanishes on a curve on $S\upi_+$ that passes the point $p$.
If this curve is not the end component of the chain $C\upi_+=S\upi_+\cap E$
which passes the point $p$,
then the section $t$ again satisfies the first item in the lemma.
So suppose that the curve is the end component of $C\upi_+$.
By the same reason, we may suppose that $t_-$ vanishes
on the end component of $C\upi_-=S\upi_-\cap E$ that passes the point $p$.
We write $G_e^+$ and $G_e^-$ for these end components
of $C\upi_+$ and $C\upi_-$ respectively (see Figure \ref{fig:rf9}).
We show that under this situation the section $t$ satisfies
at least one of the two properties in the lemma.

For this, let $G$ be a component of the cycle $C\upi=S\upi\cap E$
which is not an end of the chains $C\upi_+$ nor $C\upi_-$.
Then we have $S\upi_+.\, G = S\upi_-.\,G=0$ since
when $G\subset C\upi_+$, we have
 $S\upi_+.\,G = (S\upi_+ + S\upi_-).\,G = 0$
and similarly 
$S\upi_-.\, G =0$
when $G\subset C\upi_-$.
Therefore
$$
L\upi.\, G = \msDha.\,G
$$
holds for any integers $m_+$ and $m_-$.
So by the intersection numbers
\eqref{in2} and \eqref{in3}, we have
\begin{align}\label{int200}
L\upi.\, G
= 
\begin{cases}
1 & {\text{if $G = \ol C\upi_1\,(= \ol E_1\cap S\upi)$,}}\\
0 & {\text{otherwise.}}
\end{cases}
\end{align}		
We also have
\begin{align*}
S\upi_+.\,\ol G_e^+ = 1
\end{align*}
and
\begin{align*}
S\upi_+.\,\ol G_e^- &= \big(S\upi - S\upi_-\big).\,\ol G_e^- = - S\upi_-.\,\ol G_e^-\\
&=- S\upi_-|_{S\upi_+}.\,\ol G_e^- = - l\upi .\, \ol G^-_e = -1.
\end{align*}
Therefore we have $S\upi_+.\,\big(\ol G_e^+ + \ol G_e^-\big) = 0$.
In a similar way we also obtain 
$S\upi_-.\,\big(\ol G_e^+ + \ol G_e^-\big) = 0$.
Hence we have
\begin{align*}
L\upi.\,
\big(\ol G_e^+ + \ol G_e^-\big) = \msDha.\,\big(\ol G_e^+ + \ol G_e^-\big)
\end{align*}
for any  $m_+$ and $m_-$.
Hence again from the intersection numbers \eqref{in2} and \eqref{in3}, it follows that
\begin{align}\label{int201}
L\upi.\,
\big(\ol G_e^+ + \ol G_e^-\big) = 
\begin{cases}
{\text{$1$ or $0$}} & i=1,k,\\
0 & i\neq 1,k.
\end{cases}
\end{align}

By using \eqref{int200} and \eqref{int201} we show that the section $t$ satisfies one of the properties
in the lemma when $i\neq 1,k$.
From the present choice of the small resolutions,
the line component $\ol C\upi_1$ is not an end of 
the chains $C\upi_+$ and $C\upi_-$.
From Convention \ref{conv},  $\ol C\upi_1\subset S\upi_-$.
Then by \eqref{int200}, since $t$ is supposed to 
vanish identically on $G_e^+$,
the section $t$ vanishes identically
on any component of $C\upi_+$ except possibly on the end $\ol G_e^-$.
If $L\upi.\,\ol G_e^-\le 0$, the section $t$ vanishes further on
the component $\ol G_e^-$.
This implies $t(\ol p) = 0$.
Hence the section $t_-=t|_{S\upi_-}$ vanishes identically on a curve on $S\upi_-$ through 
the point $\ol p$. 
If this curve is different from the end component $\ol G_e^+$, then the section
$t$ satisfies the property in the first item of the lemma.
Hence suppose that $t_-$ vanishes identically on $\ol G_e^+$.
Then because we have assumed that $t_-$ vanishes also on the component
$G_e^-$, 
the section $t_-$ vanishes whole on the chain $C\upi_-$ by \eqref{int200}, 
except possibly on $\ol C\upi_1$.
But since $L\upi.\,\ol C\upi_1 =1$ by \eqref{int200},
it follows that $t_-$ vanishes identically on the component
$\ol C\upi_1$ as well.
This means that $t$ vanishes identically on the cycle $C\upi$;
namely $t$ satisfies the property in the second item
in the lemma.
Hence we are done if 
 $L\upi.\,\ol G_e^-\le 0$.
If $L\upi.\,\ol G_e^->0$, we have $L\upi.\,\ol G_e^+<0$ by \eqref{int201}.
Hence in the same way to above,
we obtain $t|_{C\upi_-}= 0$,
and also  $t(\ol p)=0$.
Then we may suppose $t_+|_{\ol G_e^-}= 0$ since
otherwise $t$ satisfies the first item in the lemma.
From this, by using \eqref{int200}, we again obtain that $t$ vanishes identically on  the cycle $C\upi$.
Therefore the proof is over when $i\neq 1,k$.

Next we show that the section $t$ satisfies one of the two
properties in the lemma when $i= 1,k$.
In these cases, the line component $\ol C\upi_1$ is one of 
the four ends $G_e^+, G_e^-, \ol G_e^+$ and $\ol G_e^-$.
From Convention \ref{conv}, either 
$\ol C\upi_1 = G_e^-$ or $\ol C\upi_1 = \ol G_e^+$ holds.
In both cases, it follows from \eqref{int200}
that $t$ vanishes identically on $C\upi$
except possibly on the two ends $\ol G_e^+$ and $\ol G_e^-$.
If $L\upi.\,\ol G_e^- \le 0$, the section $t$ vanishes
 further on $\ol G_e^-$ and
from this it follows that $t$ vanishes on a curve on $S\upi_-$
which passes the point $\ol p$.
This again means that $t$ satisfies one of the two properties
in the lemma.
If $L\upi.\,\ol G_e^- > 0$, from \eqref{int201}, we have
$L\upi.\,\ol G_e^+< 1$.
This means that 
$t$ vanishes identically on $\ol G_e^+$, and from this
we again obtain that $t$ satisfies
one of the two properties in the lemma.
\proofend

\medskip
In order to complete a proof of 
Proposition \ref{prop:exred},
we next derive a constraint for the polynomial $g$ in 
\eqref{red61}.

\begin{proposition}\label{prop:g}
The polynomial $g$ in \eqref{red61} is divisible by 
the polynomial $h$.
\end{proposition}

\proof
We take any index $i$ for which 
$n_i + 2n'_i -\mu\upi>0$ holds, and 
restrict the relation \eqref{red61}
to the reducible fiber $S\upi$ of $\tif$.
Since $g$ is a pullback of a polynomial on $\CP^1$,
the restriction $g|_{S\upi}$ is  constant.
We write $c_i = g|_{S\upi}\in\CC$.
Then since $h|_{S\upi}= 0$ by Proposition \ref{prop:h}, 
from the relation \eqref{red61}, we obtain 
\begin{align}\label{red63}
\ta\ota|_{S\upi} = c_i s_{\sst{\bm D}}|_{S\upi}.
\end{align}
We shall show $c_i = 0$. 
By Lemma \ref{lemma:lb1}, 
the restriction $\ta|_{S\upi}$ vanishes
at some point outside the cycle $C\upi$, or
vanishes identically on $C\upi$.
If the former is the case, we immediately obtain $c_i=0$
since $s_{\sst{\bm D}}|_{S\upi}$ vanishes only on $C\upi$.
If the latter is the case, another section
$\ota$ also vanishes identically on $C\upi$ as $C\upi$ is real,
and hence the left-hand side of \eqref{red63} vanishes
on any component of $C\upi$ by the order at least two.
On the other hand, the section $s_{\sst{\bm D}}|_{S\upi}$ vanishes 
on the components $C\upi_1$ and $\ol C\upi_1$ by the order exactly one
because the coefficients of the ridge components $E_1$ and $\ol E_1$ for $\bm D$ are 
always one.
Therefore from \eqref{red63} we again obtain $c_i=0$.

From this, by the relation \eqref{red61}, if we put
\begin{align}\label{fact1}
U_1:=\prod_{n_i+2n'_i >\mu\upi} u\upi,
\end{align}
then $U_1 \mid g$.
(Namely $g$ is divisible by the product $U_1$.)
Hence by dividing both-sides of the relation \eqref{red61}
by $U_1$, noting $\ol U_1=U_1$ from reality of each $u\upi$, we obtain
\begin{align}\label{fact2}
\frac{\ta\ota}{U_1} 
= \frac g{U_1}  s_{\sst{\bm D}} + 
\frac h{U_1} s_{\bm q} + 
\frac {\ol h}{\ol U_1} \ol s_{\bm q}.
\end{align}
Of course, $g/U_1,h/U_1$ and $\ol h/\ol U_1$ are polynomials.
So if there is no index $i$ for which $n_i+2n'_i-\mu\upi>1$ holds,
we have $h=U_1$ and $g$ is divisible by $h$, and the proof is over.
But for later purpose we rewrite the relation \eqref{fact2} in more concrete form.

Each $u\upi$ can be written as
the product $u\upi_+ u\upi_-$, where $u\upi\dpm$ are
defining sections of the divisors $S\upi\dpm$ respectively.
Thus we can write
$$
U_1 = \prod_{n_i+2n'_i>\mu\upi} u\upi_+ u\upi_-.
$$
Then since the divisors $S\upi\dpm$ are irreducible,
from the relation \eqref{red61},
for each factor $u\upi$ in $U_1$, the section $\ta$ can be 
divisible by either $u\upi_+$ or $u\upi_-$.
Geometrically this means that the divisor $\Theta_{\aaa} = (\ta)$ includes
the component $S\upi_+$ or $S\upi_-$.
Let $V_1$ be the product of those which are factors of 
$\ta$.
Then $\ol V_1$ divides $\ota$, and 
\eqref{fact2} can be rewritten as
\begin{align}\label{fact3}
\frac{\ta}{V_1}
\frac{\ota}{\ol V_1}
= 
\frac g{U_1}  s_{\sst{\bm D}} + 
\frac h{U_1} s_{\bm q} + 
\frac {\ol h}{\ol U_1} \ol s_{\bm q}.
\end{align}
Here, since $V_1$ is the product of defining sections of
the divisor $S\upi_+$ and/or $S\upi_-$,
the quotient $\ta/V_1$ is a holomorphic section of the line bundle
of the form
\begin{align}\label{lb2}
\msDha(\tau') - \sum_{i=1}^k m_i S\upi_+ - \sum_{i=1}^k m'_i S\upi_-
\end{align}
where $m_i$ and $m'_i$ are either 0 or 1.
Now suppose that $n_i+2n'_i-\mu\upi>1$ for some $i$.
By restricting the relation \eqref{fact3} to $S\upi$ 
for such an index $i$,
since $h/U_1$ and $\ol h/\ol U_1$
can be further divided by $u\upi$ from  \eqref{divide1},
we obtain
\begin{align}\label{}
\frac{\ta}{V_1}
\frac{\ota}{\ol V_1}
\,\Big|_{S\upi}
= 
\frac g{U_1}\,\Big|_{S\upi} \cdot s_{\sst{\bm D}}\,|_{S\upi}.
\end{align}
The restriction $(g/U_1)|_{S\upi}$ is again
a constant.
Further, by Lemma \ref{lemma:lb1},
since the quotient $\ta/V_1$ is a section of the line bundle
of the form \eqref{lb2},
the left-hand side has a zero point
which does not belong to the cycle $C\upi$,
or vanishes along the components $C\upi_1$ and $\ol C\upi_1$
by the order at least two.
Comparing with the right-hand side,
this again means that the last constant is zero.
Therefore we can divide the relation \eqref{fact3}
by the product $U_2$ of all $u\upi$ for which $n_i + 2n'_i-\mu\upi>1$.
Furthermore,
as in the case of $U_1$, the product $U_2$ can be naturally factored as $U_2 = V_2\ol V_2$,
where $V_2$ is the product of all factors of 
the form $u\upi_+$ or $u\upi_-$ 
which divide $\ta/V_1$.
Thus dividing the relation \eqref{fact3} by $U_2 = \ol U_2$, we obtain the relation
\begin{align*}
\frac{\ta}{V_1V_2}
\frac{\ota}{\ol V_1 \ol V_2}
= 
\frac g{U_1U_2}  s_{\sst{\bm D}} + 
\frac h{U_1U_2} s_{\bm q} + 
\frac {\ol h}{\ol U_1\ol U_2} \ol s_{\bm q}.
\end{align*}
Again all quotients in the right-hand side are polynomials,
and in particular $U_1U_2$ divides $g$.
If there is no index $i$ for which $n_i + 2n'_i-\mu\upi>2$ holds,
we have $h=U_1U_2$ and $g$ is divisible by $h$, as desired.

By using Lemma \ref{lemma:lb1}, this process can be repeated until the coefficient
polynomial of $s_{\bm q}$ becomes a constant
(namely until the product of $U_1$, $U_2,\dots$ becomes $h$),
and finally we obtain a relation of the form
\begin{align}\label{fact4}
\frac{\ta}{V}
\frac{\ota}{\ol V}
= 
\frac gh  s_{\sst{\bm D}} + 
 s_{\bm q} + \ol s_{\bm q},
 \quad h = V\ol V,
\end{align}
where $V$ is the product of $V_1,V_2,\dots$.
In particular, the polynomial $g$ is divisible by $h$.
\proofend

\medskip
Now we are about to finish the proof of 
the existence of a real reducible member of
the linear system $|\ms D(m)|$.

\medskip\noindent
{\em Completion of the proof of Proposition \ref{prop:exred}.}
The last equation \eqref{fact4} is a relation in $H^0\big(\ms D(m)\big)$.
Therefore if we put 
\begin{align*}
\tilde T_{\aaa}:=\Big(
\frac{\ta}{V}\Big),
\end{align*}
the divisor $\tilde T_{\aaa} + \sigma(\tilde T_{\aaa})$
gives a real reducible member of the linear system $
\big|\ms D(m)\big|$.
Since the denominator $V$ is a product of
factors of the form $u\upi_+$ and $u\upi_-$,
we have $\tilde T_{\aaa}|_S\simeq (\ta)|_S
= \Theta_{\aaa}|_S$ for any fiber $S$ of $\tif$.
Since $\Theta_{\aaa}|_S\simeq\msDha|_S\simeq \Dha$
if $S$ is real and irreducible, we obtain $\tilde T_{\aaa}|_S\simeq\Dha$,
as desired.
\proofend

\medskip
The proof shows that, 
if $2\tau>m$ holds,
the divisors $(\ta)=\Theta_{\aaa}$ 
and $(\ota)=\ol\Theta_{\aaa}$ necessarily 
have components of the form $S\upi_+$ and/or $S\upi_-$,
and by removing all such components
we are able to obtain components of a real 
reducible member of $\big|\ms D(m)\big|$.

\section{Description of the 
twistor spaces by quartic polynomials}

In this section, assembling the results  obtained
so far, we derive a constraint for an equation of the quartic hypersurface which cuts out
the branch divisor of the double covering structure, from the scroll of planes.

We briefly recall what we have obtained so far.
Suppose $n>3$ and
let $Z$ be a Moishezon twistor space on $n\CP^2$ 
which is one of the four types $\Azero,\Aone, \Atwo$
and $\Athree$ (see Definition \ref{def:type}).
We have $\dim |F|=1$, and 
the base locus of the pencil $|F|$ is a curve $C$
which is a real anti-canonical cycle
of any real irreducible member $S$ of $|F|$.
The cycle $C$ consists of $2k$ 
irreducible components,
where $k$ is $n-2,n-1,n$ or $n+1$ according
as $Z$ is of type  $\Azero,\Aone, \Atwo$
or $\Athree$ respectively.
Let $\hat Z\to Z$ be the blowing up at the cycle $C$,
and take a total small resolution
$\tilde Z\to\hat Z$ which preserves the real structure.
We have the morphism $\tilde f:\tilde Z\to \CP^1$ induced by the pencil.
This has exactly $k$ reducible fibers,
and all of them are real.
We write $S\upone,\dots,S\upk$ for the reducible
fibers, and $\lmd\upi=\tif(S\upi)\in\CP^1$.
Let $\bm D$ be the effective divisor on $\tilde Z$ defined
as the formal extension of the divisor $D$ on the surface $S$ to $\tiZ$, and $\ms D=\ms O_{\tiZ}(\bm D)$
the associated invertible sheaf.
Then for the direct image sheaf, we have 
\begin{align}\label{d-im}
\tilde f_*\ms D\simeq \ms O(-m)^2\oplus \ms O
\end{align}
for some $m>0$,
and the number $m$ is independent of a choice of the
small resolution $\tilde Z\to\hat Z$.
We have the formula for this number $m$
in terms of the coefficients
of the effective divisor $D$ and 
local invariants  associated to neighborhoods in $\tiZ$ of
the reducible fibers $S\upi$.
Writing $\ms D(m) = \ms D\otimes\tif^*\ms O(m)$,
from \eqref{d-im}, we obtain
$h^0\big(\tilde Z,\ms D(m)\big) = m+3$, and 
the linear system $\big|\ms D (m)\big|$ on $\tilde Z$
induces a degree-two meromorphic map
$\tilde\Phi_m:\tilde Z\to Y_m\subset\CP^{m+2}$,
where $Y_m$ is a scroll of planes over a rational 
normal curve $\Lmd_m$ of degree $m$.
The branch divisor of $\tilde\Phi_m$ is a cut of the scroll
by a quartic hypersurface in $\CP^{m+2}$.
We denote by $\bm l$ for the ridge of the scroll $Y_m$.
The ridge $\bm l$ has two special points $\bm q$ and $\ol{\bm q}$
characterized by the property $B|_{\bm l} = 2\bm q + 2\ol{\bm q}$,
where $B$ is the branch divisor of the
degree-two meromorphic
map $\tilde\Phi_m:\tiZ\to Y_m$.
The target space of $\tif$ and the curve $\Lmd_m$ are naturally identified.
The branch divisor $B$ is a cut of
the scroll $Y_m$ by a quartic hypersurface in $\CP^{m+2}$.
Our goal is to derive an equation of 
this quartic hypersurface in a concrete form.

We denote $$\mathbb P\upone,\mathbb P\uptwo,\dots,\mathbb P\upk$$ for the planes of the scroll $Y_m$
which correspond to the reducible fibers
$S\upone,S\uptwo,\dots,S\upk$ respectively.
We have $\tilde\Phi_m\inv(\mathbb P\upi) = S\upi$ for
each $i$.
Therefore, since $S\upi$ is reducible,
the restriction $B|_{\mathbb P\upi}$ is of the form $2\ms C\upi$
for some curve $\ms C\upi\subset\mathbb P\upi$.
Since $B$ is real and of degree four on the plane, 
the curve $\ms C\upi$ is a real conic.
If $\tilde\Phi_m$ does not contract the divisors 
$S\upi_+$ and $S\upi_-$ to lower-dimensional subvarieties, the conic $\ms C\upi$ is irreducible
since it is the image of the twistor line $l\upi$.
But $\tilde\Phi_m$ can contract $S\upi_+$ and $S\upi_-$ to lower-dimensional subvarieties, and in that case the real conic $\ms C\upi$ can be two lines.
Such examples can be found in \cite[Prop.\,4.4]{Hon_Cre2}
and \cite[Section 5.2]{Hon_JAG2}.
In these examples two irreducible components of $S\upi$
are contracted to mutually distinct lines in $\mathbb P\upi$
and $\ms C\upi$ consists of these lines.
(In the latter article, such a conic $\ms C\upi$ is called a splitting double conic.)
We call these conics
\begin{align}\label{doubleconics}
\ms C\upone,\ms C\uptwo,\dots, \ms C\upk
\end{align}
the {\em double conics} on the branch divisor $B$.
All these satisfy $B|_{\mathbb P\upi} = 2\ms C\upi$,
and each passes the  points
$\bm q$ and $\ol{\bm q}$ on the ridge $\bm l$.
Also, any double conic does not include the ridge $\bm l$
since $\bm l\not\subset B$.

Let $\aaa$ be an index in the range \eqref{alpha}.
These indices were in one-to-one correspondence
with real bitangents to the branch quartic
of the degree-two morphism $\phi:S\to\CP^2$ induced by $|D|$;
namely the restriction $\tilde\Phi_m|_S$,
where $S$ is any real irreducible member of the pencil $|F|$.
The class of a $(-1)$-curve over the real bitangent
was denoted by $\Dha$.
By Proposition  \ref{prop:exred},
there exists a real member of the system $\big|\ms D(m)\big|$
which is of the form $\tilde T_{\aaa} + \sigma({\tilde T}_{\aaa})$,
and which satisfies $\tilde T_{\aaa}|_S\simeq \Dha$
for any real irreducible member $S\in |F|$.
We write $H_{\aaa}$ for the hyperplane in $\CP^{m+2}$
which corresponds to this member.
This is real.
Let $\ms T_{\aaa}$ be the associated trope in the generalized sense.
By definition, this is a curve on the branch divisor $B$ which satisfies
\begin{align}\label{trope3}
B|_{Y_m\cap H_{\aaa}} = 2\ms T_{\aaa}.
\end{align}

%
\begin{center}
\begin{table}
\begin{tabular}{c|c|c|c}
                 type $\Azero$ & type $\Aone$ & type $\Atwo$ & type $\Athree$\\ \hline
 $\ms T_1,\ms T_2, \ms T_3$      &    
$\ms T_1, \ms T_2$    &  $\ms T_1$      & none
\end{tabular} 
\medskip
\caption{Tropes in the generalized sense of the branch divisor}
\label{tbl:trope2}
\end{table}
\end{center}

In Table \ref{tbl:trope2},
we display all tropes in the generalized sense for each type.
As a curve in $\CP^{m+2}$,
the degree of each $\ms T_{\aaa}$ is $2m$.
Indeed, since $B$ is of degree four,
from \eqref{trope3}, we have $\ms T_{\aaa} \in| \ms O_{Y_m\cap H_{\aaa}}(2)|$.
Since $Y_m$ is a scroll over the rational normal curve of degree $m$,
the intersection $Y_m\cap H_{\aaa}$ can be identified with $ \ol{\mathbb F}_m$,
the ruled surface of degree $m$ with the $(-m)$-section contracted,
and the class $\ms O_{Y_m\cap H_{\aaa}}(1)$ is represented by an image of 
an $(+m)$-section of $\mathbb F_m\to\CP^1$
under the contraction map $\mathbb F_m\to\ol{\mathbb F}_m$.
Hence the intersection number of the trope $\ms T_{\aaa}$ with a 
generic hyperplane section is $1\cdot 2 \cdot m=2m$, as asserted.
We call the double conics $\ms C\upi$, $1\le i\le k$, and 
the tropes $\ms T_{\aaa}$ in the generalized sense
as the {\em double curves} of the branch divisor $B$.
The number of double curves on the branch divisor $B$ is 
$(n+1)$ in total, regardless of the type of $Z$.

Next we show the existence of a special quadric
in $\CP^{m+2}$.
This quadric will appear in a defining equation of the 
quartic hypersurface that gives the branch divisor $B$.

\begin{proposition}\label{prop:Qexist}
Let the notation be as above.
Then there exists a real quadratic hypersurface in $\CP^{m+2}$ 
whose intersection with the planes $\mathbb P\upi$ and the hyperplane sections
$Y_m\cap H_{\aaa}$ equals the double conic $\ms C\upi$ and the trope
$\ms T_{\aaa}$ in the generalized sense respectively for any $i$ and $\aaa$.
\end{proposition}

\proof
We first show that,  regardless of types, we can take four real hyperplanes $H_1,H_2,H_3$ and $H_4$ in $\CP^{m+2}$
which satisfy 
\begin{align}\label{dc0}
\left(\bigcup_{i=1}^k \ms C\upi 
\right)
\,\cup\,
\left(\bigcup_{\aaa} \ms T_{\aaa}
\right)
\subset
H_1\cup H_2\cup H_3\cup H_4,
\end{align}
where the index $\aaa$ runs in the range
as in Table \ref{tbl:trope2} according to the type.
%

For the hyperplane which includes the trope $\ms T_{\aaa}$ in the generalized sense, we just take the hyperplane $H_{\aaa}$.
When $Z$ is of type $\Azero$,
it remains to show that all double conics
$\ms C\upi$, $1\le i\le k= n-2$, are contained in 
a single real hyperplane, say $H_4$.
The rational normal curve $\Lmd_m$ is in $\CP^m$,
and we have an inequality $m\ge n-2$ as we obtained in \eqref{ineqc}.
Hence there exists a hyperplane $H'_4$
in $\CP^m$ which
satisfies 
\begin{align}\label{Azero'}
H'_4\cap \Lmd_m = \big\{\lmd\upone,\lmd\uptwo,\dots,\lmd\upk\big\}.
\end{align}
$H'_4$ can be taken as real since $\lmd\upone,\dots,\lmd\upk$ are real.
Then by using the linear projection $p:\CP^{m+2}\to\CP^m$ 
whose center is the ridge $\bm l$,
we let $H_4=p\inv(H'_4)$.
Obviously the hyperplane $H_4$ is real and contains
all double conics.

Similarly, if $Z$ is of type $\Aone$,
since $k=n-1\le m+1$ from the  inequality
$m\ge n-2$, 
we can take two real hyperplanes $H'_3$ and $H'_4$
in $\CP^m$ which satisfy
\begin{align}\label{Aone'}
(H'_3\cup H'_4)\cap \Lmd_m =\big \{\lmd\upone,\lmd\uptwo,\dots,\lmd\upk\big\}.
\end{align}
If we put $H_3=p\inv(H'_3)$ and $H_4=p\inv(H'_4)$,
these are real hyperplanes and 
the union $H_3\cup H_4$ contains all double conics.

If $Z$ is of type $\Atwo$, 
since $k=n\le m+2$ from the  inequality
$m\ge n-2$, 
we can take three real hyperplanes $H'_2,H'_3$
and $H'_4$ which satisfy
\begin{align}\label{Atwo'}
(H'_2\cup H'_3\cup H'_4)\cap \Lmd_m = 
\big\{\lmd\upone,\lmd\uptwo,\dots,\lmd\upk\big\}.
\end{align}
We put $H_2 = p\inv(H'_2),
H_3=p\inv(H'_3)$ and $H_4=p\inv(H'_4)$.
Then these are real hyperplanes and  the union  $H_2\cup H_3\cup H_4$
contains all double conics.

If $Z$ is of type $\Athree$,
since $k=n+1\le m+3$ from the  inequality
$m\ge n-2$, 
we can take four real hyperplanes $H'_1,H'_2,H'_3$
and $H'_4$ which satisfy
\begin{align}\label{Athree'}
(H'_1\cup H'_2\cup H'_3\cup H'_4)\cap \Lmd_m = 
\big\{\lmd\upone,\lmd\uptwo,\dots,\lmd\upk\big\}.
\end{align}
We put $H_{\aaa}=p\inv(H'_{\aaa})$ 
for each $\aaa$, $1\le \aaa\le 4$.
Then these are real hyperplanes and  the union  $H_1\cup H_2\cup H_3\cup H_4$
contains all double conics.
Thus we have shown the assertion at the beginning of this proof.

Since the branch divisor $B$ is cut of 
the scroll $Y_m$ by a quartic hypersurface
by Theorem \ref{thm:dc2},
the left-hand side of \eqref{dc0} is an element of the linear
system $\big|\ms O_{Y_m\cap (H_1\cup H_2\cup H_3\cup H_4)}(2)\big|$.
Therefore, in order to show the existence of 
a quadric hypersurface in $\CP^{m+2}$ which contains all the double curves, it is enough to show that
the restriction homomorphism 
$H^0\big(\ms O_{Y_m}(2)\big)
\to H^0\big(\ms O_{Y_m\cap (H_1\cup H_2\cup 
H_3\cup H_4)}(2) 
\big)$ is surjective.
For this, we have a standard exact sequence
\begin{align}\label{sexs}
0 \lras \ms O_{Y_m}(-2) \lras
\ms O_{Y_m}(2) \lras 
\ms O_{Y_m\cap (H_1\cup H_2\cup 
H_3\cup H_4)}(2) \lras 0,
\end{align}
and from this it suffices to show 
$H^1\big( \ms O_{Y_m}(-2)\big) = 0$.
But this follows easily from an exact sequence
$0 \lras \ms O_{Y_m}(-1)\lras \ms O_{Y_m}
\lras \ms O_{Y_m\cap H} \lras 0$,
where $H$ is a hyperplane in $\CP^{m+2}$ which does not contain the ridge $\bm l$,  vanishing of $H^1(\ms O_{Y_m})$, and also an exact sequence
$0 \lras \ms O_{Y_m}(-2)\lras \ms O_{Y_m}(-1)
\lras \ms O_{Y_m\cap H}(-1) \lras 0$.
The reality of the quadric immediately follows from the reality of 
the double curves \eqref{dc0} and the isomorphism
$H^0\big(\ms O_{Y_m}(2)\big)\simeq 
H^0\big(\ms O_{Y_m\cap (H_1\cup H_2\cup 
H_3\cup H_4)}(2)\big)$
implied by the exact sequence \eqref{sexs}.
\proofend

\medskip
We are ready to derive a concrete form of the 
equation of the quartic hypersurface
in $\CP^{m+2}$.
In order to state it precisely,
we take homogenous coordinates
$(z_0,z_1,\dots, z_{m+2})$  such that 
the ridge $\bm l$ is defined by 
$z_0 = z_1 = \dots = z_m=0$,
and such that the real structure induced from that on
the twistor space is given by the complex conjugation
$(z_0,z_1,\dots, z_{m+2})\mapsto
(\ol z_0,\ol z_1,\dots,\ol z_{m+2})$.

In these coordinates, the linear projection 
$p:\CP^{m+2}\to\CP^m$ from $\bm l$ is the map which drops
the last two coordinates $z_{m+1}$ and $z_{m+2}$.

\begin{theorem}\label{thm:main}
Let the situation be as stated in 
the beginning of this section.
Then in the above homogeneous coordinates,
a defining equation of the quartic hypersurface
which cuts out the branch divisor $B$ from the scroll $Y_m$ can be taken in the form
\begin{align}\label{BB}
 h_1h_2h_2h_3 = Q^2,
\end{align}
where  $h_i$  and $Q$ are  linear and
quadratic polynomials with real coefficients in $z_0,\dots, z_{m+2}$
respectively. Moreover,
exactly $1,2,3$ or $4$ among the four polynomials $h_1,\dots, h_4$ belong
to the ideal $(z_0,z_1,\dots, z_m)$ if $Z$ is of type $\Azero,\Aone,\Atwo$ or $\Athree$ respectively.
\end{theorem}

\proof
Let $R=R(z_0,\dots z_{m+2})$ be a homogeneous polynomial of degree four whose intersection with the scroll $Y_m$ is the branch divisor $B$.
We choose four real hyperplanes $H_1,H_2,H_3$ and $H_4$ in $\CP^{m+2}$
as in the proof of the previous proposition,
and let $h_{\aaa}$, $1\le \aaa\le 4,$ be 
a defining equation of $H_{\aaa}$.
Let $Q$ be a quadratic polynomial
which defines  the quadric in the previous proposition.
These polynomials can be taken as real ones because $H_{\aaa}$ and $Q$ are real.

For an algebraic subset $X\subset\CP^{m+2}$,
we denote by $I_X\subset\CC[z_0,\dots,z_{m+2}]$ for the defining ideal 
of $X$ in $\CP^{m+2}$.
On the planes $\mathbb P\upi$,
$1\le i\le k$,
the polynomial $R$ vanishes doubly on
the conic $\ms C\upi$.
Since the intersection of the quadric and the plane $\mathbb P\upi$ equals $\ms C\upi$,
there exists a constant $c_i\in\CC^*$
which satisfies $R-c_iQ^2\in I_{\mathbb P\upi}$.
Hence taking a difference, for any indices $i$ and $j$,
we obtain $(c_i - c_j)Q^2 \in I_{\mathbb P\upi}+ I_{\mathbb P\upj}$.
So if $c_i\neq c_j$, we have
$Q^2\in I_{\mathbb P\upi}+ I_{\mathbb P\upj}$.
Further we have $I_{\mathbb P\upi} + I_{\mathbb P\upj}= 
I_{\mathbb P\upi\cap\mathbb P\upj} = I_{\bm l}$.
So we have $Q^2\in I_{\bm l}$.
Hence $Q\in I_{\bm l}$.
Since $Q|_{\mathbb P\upi}$ defines the conic $\ms C\upi$,
this means $\bm l\subset \ms C\upi$.
As $\ms C\upi\subset B$, this means $\bm l\subset B$.
This contradicts $\bm l\not\subset B$.
Therefore $c_i= c_j$ for any indices $i$ and $j$.

Next, let $\aaa$ be an index such that 
 intersection
$H_{\aaa}\cap B$ is the trope $\ms T_{\aaa}$. 
Namely $\aaa$ is in the range presented in  \eqref{alpha}.
Then since the polynomial $R$ defines $\ms T_{\aaa}$ on 
the hyperplane section $Y_m\cap H_{\aaa}$ as a double curve
and $Q$ vanishes exactly on $\ms T_{\aaa}$ over the hyperplane section, there exists
a constant $c'_{\aaa}\in\CC^*$ which satisfies 
$R-c'_{\aaa}Q^2\in I_{Y_m\cap H_{\aaa}}$.
Suppose that $c'_{\aaa}\neq c_1$ for some $\aaa$.
Then taking a difference we have $Q^2\in I_{Y_m\cap H_{\aaa}}+ I_{\mathbb P\upone}$.
Further we have $I_{Y_m\cap H_{\aaa}}+ I_{\mathbb P\upone}
= I_{(Y_m\cap H_{\aaa})\cap\mathbb P\upone}
=I_{H_{\aaa}\cap\mathbb P\upone}$ as $\mathbb P\upone\subset Y_m$.
So $Q^2\in I_{H_{\aaa}\cap \mathbb P\upone}$.
Since $\bm l\not\subset H_{\aaa}$,
the intersection $H_{\aaa}\cap \mathbb P\upone$ is a line.
It follows that the conic $\ms C\upone$ contains this line
as an irreducible component.
This line is real since $H_{\aaa}$ and $\mathbb P\upone$ are real,
and is different from the ridge $\bm l$
since $\bm l\not\subset H_{\aaa}$.
Hence this line does not pass the points $\bm q$ nor $\ol{\bm q}$.
Also, even if $\ms C\upone$ has
another line as a component, it cannot be $\bm l$
since $\bm l\not\subset B$.
These imply $\{\bm q,\ol{\bm q}\}\cap\ms C\upone=\emptyset$,
which is a contradiction.
Therefore we obtain $c'_{\aaa}=c_1$ for any $\aaa$
in the range as in \eqref{alpha}.

Thus all the constants $c_i$ and $c'_{\aaa}$ are equal.
Hence by leaving the constant to
the polynomial $Q$, we may suppose that 
\begin{gather}\label{rq}
R - Q ^ 2 \in \bigcap_{i=1}^k I_{\mathbb P\upi}
= I_{\mathbb P\upone\cup\dots\cup\,\mathbb P\upk}\qandq\\
R - Q ^ 2 \in I_{H_{\aaa}\cap Y_m}
\quad{\text{for any $\aaa$ as in \eqref{alpha}}}\label{rq2}.
\end{gather}

If $Z$ is of type $\Azero$, 
from \eqref{Azero'} we have $ \mathbb P\upone\cup\dots\cup\,\mathbb P\upk
=H_4\cap Y_m$.
Therefore, from  \eqref{rq}, we have
$R-Q^2\in I_{H_4\cap Y_m}$.
Hence, together with \eqref{rq2}, we have $R-Q^2\in I_{H_{\aaa}\cap Y_m}$ for any $\aaa$, $1\le \aaa\le 4$.
Furthermore, for these $\aaa$, we have $I_{H_{\aaa}\cap Y_m} = I_{H_{\aaa}} + I_{Y_m}
= (h_{\aaa}) + I_{Y_m}$.
From these we can conclude $R-Q^2 \in (h_1h_2h_3h_4) + I_{Y_m}$.
This implies that we can take \eqref{BB} as an equation of $B$,
and that only $h_4$ belongs to the ideal $I_{\bm l}
=(z_1,\dots,z_m)$.

If $Z$ is of type $\Aone$, from \eqref{Aone'}, 
we have $\mathbb P\upone\cup\dots\cup\,\mathbb P\upk=  (H_3\cup H_4)\cap Y_m$.
Therefore, from  \eqref{rq}, we have
$R-Q^2\in I_{(H_3\cup H_4)\cap Y_m}=(h_3h_4) + I_{Y_m}$.
Also from \eqref{rq2} we have $R-Q^2 \in (h_{\aaa}) + I_{Y_m}$
if $\aaa=1,2$.
From these we again obtain $R-Q^2 \in (h_1h_2h_3h_4) + I_{Y_m}$.
and that exactly $h_3$ and $h_4$ belong to the ideal $I_{\bm l}$.

If $Z$ is of type $\Atwo$, from \eqref{Atwo'}, 
we have $\mathbb P\upone\cup\dots\cup\,\mathbb P\upk= (H_2\cup H_3\cup H_4)\cap Y_m$.
Therefore, from  \eqref{rq}, we have
$R-Q^2\in I_{(H_2\cup H_3\cup H_4)\cap Y_m}=(h_2h_3h_4) + I_{Y_m}$.
Also from \eqref{rq2} we have $R-Q^2 \in (h_{1}) + I_{Y_m}$.
From these we again obtain $R-Q^2 \in (h_1h_2h_3h_4) + I_{Y_m}$.
and that exactly $h_2,h_3$ and $h_4$ belong to the ideal $I_{\bm l}$.

Finally if $Z$ is of type $\Athree$, 
we have $\mathbb P\upone\cup\dots\cup\,\mathbb P\upk=  (H_1\cup H_2\cup H_3\cup H_4)\cap Y_m$ by \eqref{Athree'}.
Therefore, from  \eqref{rq}, we have
$R-Q^2\in (h_1h_2h_3h_4) + I_{Y_m}$.
Further
all $h_{\aaa}$ belong to the ideal $I_{\bm l}$.
\proofend

\medskip
One might wonder the meaning of the constraint for the linear polynomials
$h_1,\dots,h_4$ given in the final part of the theorem.
This is directly relevant to the  singularities
of the quartic on the planes of the scroll $Y_m$ as follows.
Let $(z_0,z_1,\dots,z_{m+2})$ be the homogeneous coordinates on $\CP^{m+2}$ 
as in Theorem \ref{thm:main}.
In particular the ridge $\bm l$ is defined by $z_0 = z_1 = \dots = z_m = 0$.
By a linear change of the two coordinates
$z_{m+1}$ and $z_{m+2}$, we can suppose that 
the points $\bm q$ and $\bm {\ol q}$ on $\bm l$ satisfy $z_{m+1}=0$ and $z_{m+2}=0$ respectively.
(So in these coordinates the real structure is not precisely
given by the complex conjugation.)
An equation of a generic plane of the scroll is of the form
\begin{align}\label{pln}
z_j = c_j\,z_0,\quad c_j\in\CC,\quad 1\le j\le m.
\end{align}
On these planes the triple $(z_0,z_{m+1}, z_{m+2})$ work
as homogeneous coordinates.
The ridge $\bm l$ is not contained in the quadric $Q=0$
because, on the plane $\mathbb P\upi$, $1\le i\le k$,
the equation $Q=0$ defines the conic $\ms C\upi$,
and this does not contain the ridge.
From this, 
the quadratic polynomial $Q$ includes at least one of the three monomials $z_{m+1}^2, z_{m+1}z_{m+2}$
and $z_{m+2}^2$.
But from the above choice of the coordinates
$z_{m+1}$ and $z_{m+2}$
and the fact $\bm q,\ol{\bm q}\in  \{Q=0\}$,
it follows that $Q$ does not include
the two monomials $z_{m+1}^2$ 
and $z_{m+2}^2$.
Therefore we may suppose
\begin{align}\label{}
Q = z_{m+1}z_{m+2} +  g_1(z_0,z_1,\dots, z_{m+2}),
\end{align}
for some quadratic polynomial $g_1$ that does not include
the monomials $z_{m+1}^2, z_{m+1}z_{m+2}$
and $z_{m+2}^2$.

As in the theorem, if $Z$ is of type ${\rm A}_{\nu}$
(where $0\le \nu\le 3$),
exactly $(\nu+1)$ among the four linear polynomials $h_1,h_2,h_3$ and $h_4$
belong to the ideal $(z_0,z_1,\dots,z_m)$.
Hence, {\em 
the restriction of the product $h_1h_2h_3h_4$
to the planes of the form \eqref{pln}
can be divided by $z_0^{\nu+1}$ but cannot be divided by $z_0^{\nu+2}$}.
From these, 
on the plane \eqref{pln} of the scroll, the equation of the quartic becomes
$$
\big\{
z_{m+1}z_{m+2} + g_2(z_0, z_{m+1}, z_{m+2})
\big\}^2 
= z_0^{\nu+1} g_3(z_0,z_{m+1},z_{m+2}),
$$
where $g_2$ is of degree two 
and $g_3$ is of degree $(3-\nu)$ and not divisible by $z_0$.
On these planes,
the equation of the ridge $\bm l$ is $z_0=0$, and 
the points $\bm q$ and $\ol{\bm q}$
respectively satisfy $(z_0,z_{m+1}, z_{m+2})=(0,0,1)$ and $(z_0,z_{m+1}, z_{m+2})=(0,1,0)$.
From this, we can easily read off that, on
these planes,
the branch quartic  has ${\rm A}_{\nu}$-singularities
at the points $\bm q$ and $\bm {\ol q}$.
Thus the constraint for the linear polynomials
$h_1,h_2,h_3$ and  $h_4$ in Theorem \ref{thm:main}
is directly related to the singularities 
of the quartics on the planes of the scroll $Y_m$.

\end{document}